\documentclass[10pt]{article}
\usepackage[all]{xy}
\usepackage{amsmath}
\usepackage{amsfonts}
\usepackage{amssymb}
\usepackage{amscd}
\usepackage{amsthm}
\usepackage{latexsym}
\usepackage{amsbsy}

\begin{document}

\title{\bf{Scalar Curvature for the  Noncommutative Two Torus}}
\author {\bf{Farzad Fathizadeh}\\
Department of Mathematics and Statistics\\
York University\\
Toronto, Ontario, Canada, M3J 1P3\\
ffathiza@mathstat.yorku.ca\\ 
\\
 \bf{Masoud Khalkhali}\\
 Department of Mathematics\\
 The University of Western Ontario\\
 London, Ontario, Canada, N6A 5B7\\
 masoud@uwo.ca
 }

\date{}
\maketitle
\begin{abstract}
We give  a local expression for the {\it scalar curvature}  of the noncommutative two torus $ A_{\theta} = C(\mathbb{T}_{\theta}^2)$ equipped
with an arbitrary translation invariant complex structure and Weyl factor.   This is achieved by evaluating the value
of the (analytic continuation of the)  {\it spectral zeta functional} $\zeta_a(s): = \text{Trace}(a \triangle^{-s})$ at $s=0$ as a linear functional in $a  \in  C^{\infty}(\mathbb{T}_{\theta}^2)$. A new,   purely noncommutative,  feature here
 is the appearance of the {\it  modular automorphism group} from the theory of type III factors  and  quantum statistical mechanics in the final formula for the  curvature.  This formula coincides with the   formula that was recently obtained independently by Connes and Moscovici  in their recent paper \cite{conmos2}.
 \end{abstract}

\tableofcontents

\section{Introduction}
In this paper we give  a local expression for the {\it scalar curvature}  of the noncommutative two torus $ A_{\theta} = C(\mathbb{T}_{\theta}^2)$ equipped
with an arbitrary translation invariant complex structure and Weyl factor.  More precisely,  for any complex number $\tau$  in the upper half
plane, representing the conformal class of a metric on  $\mathbb{T}_{\theta}^2$,
 and a Weyl factor given by a positive  invertible element $k \in C^{\infty}(\mathbb{T}_{\theta}^2)$, we give an explicit formula for an  element $R = R (\tau, k) \in    C^{\infty}(\mathbb{T}_{\theta}^2)$ that  is the scalar curvature of the underlying  noncommutative Riemannian manifold  $\mathbb{T}_{\theta}^2$.  This is achieved by evaluating the value
of the (analytic continuation of the)  {\it spectral zeta functional} $\zeta_a(s): = \text{Trace}(a \triangle^{-s})$ at $s=0$ as a linear functional in $a  \in  C^{\infty}(\mathbb{T}_{\theta}^2)$. A new,   purely noncommutative,  feature here
 is the appearance of the {\it  modular automorphism group} from the theory of type III factors  and  quantum statistical mechanics in the final formula for curvature.  This formula exactly reproduces the  formula that was recently obtained independently by Connes and Moscovici  in their recent paper \cite{conmos2}. It also reduces, for
$\tau = \sqrt{-1}$, to a formula that was earlier obtained by Alain Connes for the scalar curvature of the noncommutative two torus.

 Our main result (Theorem \ref{curfull} below) extends  and refines the recent work on {\it Gauss-Bonnet theorem}  for the noncommutative two torus that was initiated  in the pioneering work of Connes and Tretkoff in  \cite{contre} (cf. also \cite{cohcon, coh} for a preliminary version) and  its later generalization in  \cite{fatkha}. In fact after applying the standard trace
of the noncommutative torus to the scalar curvature $R$ one obtains, for all values of $\tau$ and $ k,$  the value 0. This is the Gauss-Bonnet theorem  for the noncommutative two torus
  and, in the commutative case, is equivalent to the classical Gauss-Bonnet theorem for a surface
of genus 1.

The backbone of the present  paper is  {\it  noncommutative differential geometry program}  \cite{con, con0, con1, conmar}. According to parts of this theory that is relevant here   the metric
information on a noncommutative space is fully encoded as a {\it spectral triple} on the noncommutative algebra of coordinates  on   that space. Various technical results corroborates, in fact fully justifies,  this vision. First of all,  {\it Connes' reconstruction theorem} \cite{con2}  guarantees that in the commutative case, the notion of spectral triple is strong enough to fully recover the Riemannian (spin)  manifold from  its natural spectral triple data defined using the Dirac operator acting on spinors. Secondly, as it is  shown in \cite{con0.5, con1, conmar},  ideas of spectral geometry, in particular formulation of several invariants of a Riemannian  manifold like  volume and scalar curvature in terms of asymptotics of the trace of the heat kernel of Laplacians and Dirac operators, have very natural extensions in the noncommutative setting and recover the classical results in the commutative case. Other relevant results are the Connes-Moscovici local index formula 
 \cite{conmos0} and Chamseddine-Connes spectral action principle \cite{chacon}.  In passing to the noncommutative case, sooner or later  one must face the prospect of type III algebras and the lack of trace on them. It was exactly for this reason that \emph{twisted spectral triples} were introduced by Connes and Moscovici in \cite{conmos1}. The spectral triple at the foundation of the present paper was defined in \cite{contre} and is in fact, via the right action corresponding to the Tomita anti-linear unitary map, a twisted spectral  triple.

One of the main technical tools  employed in this paper is Connes'
pseudodifferential operators and their symbol calculus on the noncommutative torus \cite{con}  and the use of the
asymptotic expansion of the heat kernel in computing  zeta values. This, however, by itself 
is not enough and, similar to \cite{cohcon, contre, fatkha}, one needs an extra and intricate argument to express $\zeta_a(0)$  in
terms of the modular operator defined by the Weyl factor. As a first step, the calculation of the
asymptotic expansion of the heat operator for arbitrary values of the conformal class is quite
involved and must be performed by a computer. We found it impossible to carry this step without
the use of symbolic calculations. Finally we should mention that, as is explained in \cite{contre, fatkha}, there is a
close relationship between the subject of this paper and scale invariance in spectral action \cite{chacon, chacon1}
on the one hand, and non-unimodular (or twisted) spectral triples \cite{conmos1} on the other hand.

This paper is organized as follows. In Section \ref{prel} we recall  a  twisted spectral triple
 on the noncommutative torus from \cite{contre} and the conformal structures of  this  noncommutative space. An important idea here is to
determine the conformal class of a metric by defining a complex structure on the noncommutative torus, and perturbing this metric by changing 
the tracial volume form to a KMS state by means of a Weyl factor given by an invertible positive smooth element \cite{contre}. In Section \ref{scacursec}
we give a spectral definition for the scalar curvature of the noncommutative torus equipped with a general metric.  We also recall the  pseudodifferential calculus \cite{con}  for the special case of the canonical dynamical system defining  the noncommutative torus and explain how this will
provide a method for computing a local expression for the scalar curvature of this noncommutative Riemannian manifold. In Section \ref{computation} we illustrate
the process of finding this local expression by means of   pseudodifferential calculus on the noncommutative torus  and heat kernel techniques. Another crucial technique
here, as in \cite{cohcon, contre, fatkha}, is to use the {\it modular automorphism} to permute elements of the noncommutative torus with the Weyl factor.  In fact this
prepares the ground for using functional calculus to write the final formula for the scalar curvature in a concise form. Considering the lengthy
computations and formulas in this section, the final concise formula shows some magical cancelations and simplifications after the necessary
rearrangements and permutations by means of the modular automorphism.
In Section \ref{scacurlogsec} we simplify our formula for the scalar curvature of the noncommutative torus in terms of  the logarithm of the Weyl
factor. Here again the modular automorphism is used  crucially to find some identities that relate the derivative of the Weyl factor and the derivative
of its logarithm with respect to the noncommutative coordinates of the noncommutative torus. At the end, for the convenience of the reader, we have 
recorded in Appendices, the lengthy formulas for the pseudodifferential symbols that appear in approximating  the resolvent of the Laplacian and contain 
the geometric information for computing the scalar curvature. 

The definition of the scalar curvature for spectral triples in terms of the second term of the heat expansion was given in \cite{conmar},  Definition 1.147 of Section  11.1. The refinement used here as well as in  \cite{conmos2}  is to introduce the chiral scalar curvature from which the scalar curvature using the Laplacian on functions is easily deduced,
(see also \cite{bhumar} for a variant).

We would like  to express our indebtedness to Alain Connes for motivating and enlightening discussions and for much help  during the various stages of the work on this paper.  At  several crucial stages he generously shared his insight and ideas with us and communicated  their  relevant  joint results   in \cite{conmos2}  with us.  This gave us a good chance of finding potential errors in the computations. In fact the idea of using the {\it full  Laplacian},  on functions and  $1$-forms,  as opposed to just  functions,  was suggested to us by him. While in the commutative case one can recover the curvature from zeta functionals from the Laplacian on functions, this is no more the case in the noncommutative case.  We would also like to heartily thank Henri Moscovici for a push in the right direction at an early stage. After the appearance of our Gauss-Bonnet paper  \cite{fatkha},  Henri and Alain   kindly pointed  out to us that the calculations
 in that 
  paper might be quite relevant for computing the  scalar curvature of the noncommutative two torus. Finally F. F. would like to thank  IHES for kind support and excellent environment during his visit  in Summer 2011 where  part of this work was carried out.

\section{Preliminaries} \label{prel}
Let $\Sigma$ be a closed, oriented,  2-dimensional smooth manifold equipped with a Riemannian metric $g$. 
The scalar curvature of $(\Sigma, g)$ can be expressed by a local formula in terms of the symbol of the Laplacian 
$\triangle_{g} = d^*d$, where $d$ is the de Rham differential operator acting on smooth functions on $\Sigma$. In fact 
using the Cauchy integral formula, for any $t >0$ one can write 
\begin{equation} \label{kerexp} 
e^{-t\triangle_g} =  \frac{1}{2\pi i} \int_C e^{-t \lambda} (\triangle_g - \lambda)^{-1} \, d \lambda, \nonumber
\end{equation}
where $C$ is a curve in the complex plane that goes around the non-negative
real axis in the clockwise direction without touching it. The operator $e^{-t \triangle_g}$ has a smooth kernel 
$K(t, x, y)$ and there is an asymptotic expansion of the form 
\[ 
K(t, x, x) \sim t^{-1} \sum_{n=0}^{\infty} e_{2n}(x, \triangle_g)t^{n} \qquad (t \to 0).
\]
The term $e_2(x, \triangle_g)$ turns out to be a constant multiple of the scalar curvature of $(\Sigma, g)$. 

As a first step towards computing the  \emph{scalar curvature} of the noncommutative two torus, we recall the 
notion of the perturbed spectral triple attached to $(\mathbb{T}_\theta^2, \tau, k)$, 
where $\tau \in \mathbb{C}\setminus \mathbb{R}$ represents the conformal class of a metric on the noncommutative two 
tours $\mathbb{T}_\theta^2$, and $k \in C^\infty(\mathbb{T}_\theta^2)$ is the Weyl factor by the aid of which one can 
vary inside the conformal class of the metric \cite{contre, fatkha}.   

\subsection{The irrational rotation algebra.}
Let  $\theta$ be an irrational number. Recall that the  irrational rotation $C^*$-algebra
 $A_{\theta}$ is, by definition, the universal unital $C^*$-algebra  generated by two unitaries $U, V$ satisfying
 $$VU=e^{2 \pi i \theta} UV.$$
One usually thinks of $A_{\theta}$ as the algebra of continuous functions on the noncommutative 2-torus
$\mathbb{T}_{\theta}^2$.
There is a continuous action of $\mathbb{T}^2$, $\mathbb{T}= \mathbb{R}/2\pi \mathbb{Z}$, on $A_{\theta}$ by $C^*$-algebra
automorphisms  $\{ \alpha_s\}$, $s\in \mathbb{R}^2$, defined by
\[\alpha_s(U^mV^n)=e^{is.(m,n)}U^mV^n.\]
The space of smooth elements for this action, that is those elements $a \in A_{\theta}$ for which the map $s \mapsto
\alpha_s (a)$  is $C^{\infty}$ will be denoted by
$A_{\theta}^{\infty}$. It is a dense subalgebra of $A_{\theta}$  which can be alternatively
described as the algebra of elements in $A_{\theta}$
whose (noncommutative) Fourier expansion has rapidly decreasing coefficients:
\[
A_{\theta}^{\infty}=\big \{\sum_{m,n\in \mathbb{Z}}a_{m,n}U^mV^n; \quad \sup_{m,n\in \mathbb{Z}} (|m|^k|n|^q|a_{m,n}|) < \infty, 
\forall k,q \in \mathbb{Z}  \big \}.
\]
There is a unique normalized trace   $\mathfrak{t}\,$  on $A_{\theta}$ whose restriction on smooth elements is given by
\[\mathfrak{t}\,(\sum_{m,n\in \mathbb{Z}}a_{m,n}U^mV^n)=a_{0,0}.\]

 The infinitesimal generators of the above action of $\mathbb{T}^2$ on $A_{\theta}$ are the
 derivations  $\delta_1, \, \delta_2: A_{\theta}^{\infty} \to A_{\theta}^{\infty}$ defined by
\[\delta_1(U)=U, \quad \delta_1(V)=0, \quad  \delta_2(U)=0, \quad  \delta_2(V)=V.\]
In fact, $\delta_1, \delta_2$ are analogues of the differential operators
$\frac{1}{i}\partial/\partial x, \frac{1}{i}\partial/\partial y$ acting on the smooth functions on
the ordinary two torus. We have $\delta_j(a^*)= -\delta_j(a)^* $ for $j=1, 2$ and all $a\in A_{\theta}^{\infty}$.
 Moreover, since $\mathfrak{t}\, \circ \delta_j =0$, for $j=1, 2$, we have the
 analogue of the
integration by parts formula:
\[ \mathfrak{t}\,(a\delta_j(b)) = -\mathfrak{t}\,(\delta_j(a)b), \qquad \forall a,b \in A_{\theta}^{\infty}. \]

We define an inner product on  $A_{\theta}$ by
\[ \langle a, b \rangle = \mathfrak{t}\,(b^*a), \qquad a,b \in A_{\theta}, \nonumber \]
and complete $A_{\theta}$ with respect to this inner product to obtain a Hilbert space denoted by
$\mathcal{H}_0$. The derivations $\delta_1, \delta_2$, as unbounded
operators on $\mathcal{H}_0$,  are formally selfadjoint and have unique extensions to selfadjoint operators.

\subsection{Conformal structures on $\mathbb{T}_{\theta}^2$.}
To any complex number $\tau = \tau_1+i\tau_2, \, \tau_1, \tau_2 \in \mathbb{R}$, with non-zero imaginary part, we can associate 
a complex structure on the noncommutative two torus by defining
\[ \partial = \delta_1 + \bar \tau \delta_2, \qquad \partial^*=  \delta_1 + \tau \delta_2. \]
To the conformal structure defined by $\tau$,  corresponds a  positive
Hochschild  two cocycle on  $A_{\theta}^{\infty}$ given by (\emph{cf.}
\cite{con1})
\[ \psi (a, b, c)=- \mathfrak{t}\, (a\partial b \partial^* c).\]

We note that   $\partial$ is an unbounded operator on  $\mathcal{H}_0$ and $\partial^*$ is
 its formal adjoint. The analogue of the space of $(1,0)-$forms on the ordinary two torus
is defined to be the Hilbert space completion of the space of finite sums $\sum a \partial b$, $a,b
\in A_{\theta}^{\infty}$, with respect to the inner product defined above, and it is denoted by
$\mathcal{H}^{(1,0)}$. 

Now we can vary inside the conformal class of the metric \cite{contre} by choosing a smooth
selfadjoint element $h=h^* \in A_{\theta}^{\infty}$, and define a linear functional
$\varphi$ on $A_{\theta}$ by
\[\varphi(a)= \mathfrak{t}\,(ae^{-h}), \qquad a \in A_{\theta}.\]
In fact, $\varphi$ is a positive linear functional which is not a trace, however,
it is a twisted trace, and satisfies the KMS condition at $\beta = 1$ for the
1-parameter group $\{\sigma_t \}$, $t \in \mathbb{R}$ of inner automorphisms
$\sigma_t= \Delta^{-it}$ where the modular operator for $\varphi$ is given by
(\emph{cf.} \cite{contre})
\[\Delta(x)=e^{-h}xe^{h};\]
moreover, the 1-parameter group of automorphisms $\sigma_t$ is generated by the derivation
$- \log \Delta$ where
\[
\log \Delta (x) = [-h,x], \qquad x \in A_{\theta}^{\infty}. 
\]

We  define an inner product $\langle \, , \, \rangle_{\varphi}$ on $A_{\theta}$ by
\[ \langle a,b  \rangle_{\varphi} = \varphi(b^*a), \qquad a,b \in A_{\theta}. \]
The Hilbert space obtained from completing $A_{\theta}$ with respect to this
inner product will be denoted by $\mathcal{H}_{\varphi}$.

\subsection{Spectral triple on $A_{\theta}$. } \label{spectrip}

In this subsection, we recall the Connes-Tretkoff ordinary and twisted spectral triple over $A_\theta$ and $A_\theta^{\textnormal{op}}$ respectively. 

Let us view the operator $\partial$ as an unbounded operator from $\mathcal{H}_\varphi$ to  $\mathcal{H}^{(1,0)}$ 
and denote it by $\partial_\varphi$. Similar to \cite{contre}, we construct an \emph{even spectral triple} by considering the 
left action of $A_{\theta}$ on the 
Hilbert space
\[   \mathcal{H} = \mathcal{H}_\varphi \oplus \mathcal{H}^{(1,0)},\]
and the operator
\[ D=
\left(\begin{array}{c c}
0 & \partial_\varphi^* \\
\partial_\varphi & 0
\end{array}\right)
: \mathcal{H} \to \mathcal{H}.\]
Then  the \emph{Laplacian} has the following form:
\begin{equation} 
\triangle := D^2 =  \left(\begin{array}{c c}
 \partial_\varphi^*\partial_\varphi & 0  \\
0 & \partial_\varphi \partial_\varphi^*
\end{array}\right). \nonumber
 \end{equation}
 We also note that the \emph{grading} is given by 
 \[\gamma =  \left(\begin{array}{c c}
 1 & 0  \\
0 & -1
\end{array}\right) : \mathcal{H} \to \mathcal{H}.
\]

It is shown in \cite{fatkha, contre} that the operator 
\[
\partial_\varphi^* \partial_\varphi: \mathcal{H}_\varphi \to \mathcal{H}_\varphi,
\] 
is anti-unitarily equivalent to
\[k \partial^* \partial k :\mathcal{H}_0 \to \mathcal{H}_0,\] 
where $k:= e^{h/2}$ acts on $\mathcal{H}_0$ by left multiplication. In a similar manner, we have the following equivalence for the 
other half of the Laplacian:

\newtheorem{main1}{Lemma}[section]
\begin{main1}
The operator $\partial_\varphi \partial_\varphi^*: \mathcal{H}^{(1, 0)} \to \mathcal{H}^{(1, 0)} $ is anti-unitarily equivalent to 
\[\partial^* k^2 \partial : \mathcal{H}^{(1,0)} \to \mathcal{H}^{(1,0)},\]
where $k^2$ acts by left multiplication.
\begin{proof}
One can easily see that the formal adjoint of $\partial_\varphi: \mathcal{H}_{\varphi} \to \mathcal{H}^{(1,0)}$ is given by 
$ R_{k^2} \partial^*,$
where $R_{k^2}$ denotes the right multiplication by $k^2$. Let $J$ be the involution on $\mathcal{H}^{(1,0)}$ given by $J(a)=a^*$. Then 
we have:
\[ J \partial_\varphi \partial_\varphi^* J =  J \partial R_{k^2}  \partial^* J =  J \partial JJ R_{k^2}  JJ \partial^* J = \partial^* k^2 \partial.\]
\end{proof}
\end{main1}

In \cite{contre}, a twisted spectral triple is also constructed over $A_\theta^{\textnormal{op}}$. In fact, considering the Tomita anti-linear unitary map $J_\varphi$ 
in $\mathcal{H}_\varphi$, and the corresponding unitary right action of $A_\theta$ in $\mathcal{H}_\varphi$ given by $a \mapsto J_\varphi a^* J_\varphi$, 
it is shown that $(A_\theta^{\textnormal{op}}, \mathcal{H}, D)$ is a twisted spectral triple in the sense that the twisted commutators 
\[ D a^{\textnormal{op}} - (k^{-1}ak)^{\textnormal{op}}D\]
are bounded operators for all $a \in A_\theta$.

\section{Scalar Curvature} \label{scacursec}

In their book \cite{conmar}, Connes and Marcolli give a definition for the {\it scalar curvature}  of {\it  spectral triples of metric dimension 4.} This uses residues of 
the zeta function at its poles and cannot be applied to spectral triples of metric dimension 2, as  is  the case for the noncommutative two torus. 
 For spectral triples of metric dimension 2,  it is the value of the zeta functional at $s=0$ that gives the scalar curvature.  The general definition of \emph{scalar curvature} for spectral triples of metric dimension 2,  reduces to the following definition in the case of 
the noncommutative two torus (cf. also \cite{conmos2} for further explanations,  motivations, and extensions, 
and \cite{bhumar} for a related proposal).    The scalar curvature of the spectral triple attached to $(\mathbb{T}_\theta^2, \tau, k)$ in Subsection \ref{spectrip} 
is the unique element $R \in A_\theta^\infty$ satisfying the equation 
\[  \textnormal{Trace}\, (a \triangle^{-s})_{|_{s=0}} + \textnormal{Trace}\,(a P) = \mathfrak{t}\,(a R), \qquad \forall a \in A_\theta^\infty,\]
where $P$ is the projection onto the kernel of $\triangle$.   
The first term on the left hand side of this formula denotes the value at the origin, $\zeta_a(0)$, of the zeta function
\[\zeta_a(s) := \textnormal{Trace} \,(a \triangle^{-s}), \qquad \textnormal{Re}(s) >> 0.\]
This function has a holomorphic continuation to $\mathbb{C} \setminus \{1\}$, in particular its value at the origin is defined (\emph{cf.} the 
proof of Proposition \ref{localexpcur}). 

In a similar manner, for the graded case, where the additional data of grading $\gamma$ is involved, the \emph{ chiral scalar 
curvature} $R^\gamma$ is the unique  element $R^\gamma \in A_\theta^\infty$ which satisfies the equation 
\[  \textnormal{Trace}\, (\gamma a \triangle^{-s})_{|_{s=0}} = \mathfrak{t}\,(a R^\gamma), \qquad \forall a \in A_\theta^\infty.\] 

Proposition \ref{localexpcur} will provide the means for finding a local expression for $R$ and $R^\gamma$.
First we recall the pseudodifferential calculus that we shall use for finding these elements.

\subsection{Connes' pseudodifferential operators on $\mathbb{T}_\theta^2$.}
For a non-negative integer $n$, the space of differential operators
on $A_{\theta}^{\infty}$ of order at most $n$ is defined to be the
vector space of operators of the form
\[ \sum_{j_1+j_2 \leq n } a_{j_1,j_2}
\delta_1^{j_1} \delta_2^{j_2}, \qquad j_1, \, j_2 \geq 0, \qquad a_{j_1,j_2}
\in A_{\theta}^{\infty}.\]

The notion of a differential operator on $A_{\theta}^{\infty}$ can be generalized
to the notion of a pseudodifferential operator using operator valued symbols \cite{con}. In fact 
this is achieved by considering the pseudodifferential calculus associated to $C^*$-dynamical systems 
\cite{con, baj}, for the canonical dynamical system $(A_\theta^\infty, \{\alpha_s \})$. In the
sequel, we shall use the notation $\partial_1=\frac{\partial}{\partial \xi_1}$, $\partial_2=\frac{\partial}{\partial \xi_2}$.

\newtheorem{symboldef}[main1]{Definition}
\begin{symboldef}
For an integer $n$, a smooth map $\rho: \mathbb{R}^2 \to A_{\theta}^{\infty}$ is  said
to be a symbol of order $n$, if for all non-negative integers $i_1, i_2, j_1,
j_2,$
\[ ||\delta_1^{i_1} \delta_2^{i_2} \partial_1^{j_1} \partial_2^{j_2} \rho(\xi) ||
\leq c (1+|\xi|)^{n-j_1-j_2},\]
where $c$ is a constant, and if there exists a smooth map $k: \mathbb{R}^2 \to
A_{\theta}^{\infty}$ such that
\[\lim_{\lambda \to \infty} \lambda^{-n} \rho(\lambda\xi_1, \lambda\xi_2) = k (\xi_1, \xi_2).\]
The space of symbols of order $n$ is denoted by $S_n$.
\end{symboldef}

To a symbol $\rho$ of order $n$, one can associate an operator on $A_{\theta}^{\infty}$,
denoted by $P_{\rho}$, given by
\[ P_{\rho}(a) = (2 \pi)^{-2} \int \int e^{-is \cdot \xi} \rho(\xi) \alpha_s(a) \,ds \,
d\xi. \]
The operator $P_{\rho}$ is said to be a pseudodifferential operator of order $n$. For
example, the differential operator $\sum_{j_1+j_2 \leq n } a_{j_1,j_2}
\delta_1^{j_1} \delta_2^{j_2}$ is associated with the symbol $\sum_{j_1+j_2 \leq n } a_{j_1,j_2}
\xi_1^{j_1} \xi_2^{j_2}$ via the above formula.

\newtheorem{symbolequivalence}[main1]{Definition}
\begin{symbolequivalence}
Two symbols $\rho$, $\rho'\in S_k$  are said to be equivalent if and only if $\rho-\rho'$ is in
$S_n$ for all integers $n$. The equivalence of the symbols will be denoted by  $\rho \sim \rho'$.
\end{symbolequivalence}

The following lemma shows that the space of pseudodifferential operators
is an algebra and one can find the symbol of the product of
pseudodifferential operators up to the above equivalence. Also,
the adjoint of a pseudodifferential operator, with respect to the inner
product defined on $\mathcal{H}_0$ in Section
\ref{prel}, is a pseudodifferential operator with the symbol given in the following proposition 
up to the equivalence (\emph{cf.} \cite{con}).

\newtheorem{symbolcalculus}[main1]{Proposition}
\begin{symbolcalculus} \label{symbolcalculus}
Let $P$ and $Q$ be pseudodifferential operators with symbols
$\rho$ and $\rho'$ respectively. Then the adjoint $P^*$ and
the product $PQ$ are pseudodifferential operators with the following
symbols
\[
\sigma(P^*) \sim \sum_{\ell_1, \ell_2 \geq 0} \frac{1}{\ell_1! \ell_2!}
\partial_1^{\ell_1} \partial_2^{\ell_2}\delta_1^{\ell_1}\delta_2^{\ell_2}
(\rho(\xi))^*,
\]

\[
\sigma (P Q) \sim \sum_{\ell_1, \ell_2 \geq 0} \frac{1}{\ell_1! \ell_2!}
\partial_1^{\ell_1} \partial_2^{\ell_2}(\rho (\xi))
\delta_1^{\ell_1}\delta_2^{\ell_2} (\rho'(\xi)).
\]

\end{symbolcalculus}

\newtheorem{elliptic}[main1]{Definition}
\begin{elliptic} \label{elliptic}
Let $\rho$ be a symbol of order $n$. It is said to be elliptic if $\rho(\xi)$ is invertible for $\xi \neq 0$, and if 
there exists a constant $c$ such that
\[ || \rho(\xi)^{-1} || \leq c (1+|\xi|)^{-n} \]
for sufficiently large $|\xi|.$
\end{elliptic}

For example, the flat Laplacian $ \delta_1^2 + 2 \tau_1 \delta_1\delta_2 +
|\tau|^2\delta_2^2$ is  an elliptic pseudodifferential operator (\emph{cf.} \cite{contre,fatkha}).

\subsection{Local expression for scalar curvature.} \label{recipe}

Here we explain how one can find a local expression for the scalar curvature of the noncommutative two torus. This will justify the lengthy 
computations in the following sections.  

Using the Cauchy integral formula, one has
\begin{equation} \label{cauchyint} 
 e^{-t\triangle} =  \frac{1}{2\pi i} \int_C e^{-t \lambda} (\triangle - \lambda)^{-1} \, d \lambda
\end{equation} 
where $C$ is a curve in the complex plane that goes around the non-negative
real axis in the {\it clockwise} direction without touching it. Similar to the formula in \cite{gil}, one can 
approximate the inverse of the operator
$(\triangle-\lambda)$ by a pseudodifferential
operator $B_{\lambda}$ whose symbol has an expansion of the form
\[ b_0(\xi, \lambda) +  b_1(\xi, \lambda) +  b_2(\xi, \lambda) + \cdots \]
where $b_j(\xi, \lambda)$ is a symbol of order $-2-j$, and
\[ \sigma(B_{\lambda} (\triangle-\lambda) ) \sim 1 .\]

\newtheorem{localexpcur}[main1]{Proposition}
\begin{localexpcur} \label{localexpcur}
The scalar curvature of the ungraded spectral triple attached to $(\mathbb{T}_\theta^2, \tau, k)$ is equal to  
\[  \frac{1}{2 \pi i} \int_{\mathbb{R}^2} \int_C e^{-\lambda} b_2(\xi, \lambda) \, d\lambda \, d\xi ,    \]
where $b_2$ is defined as above.
\begin{proof}
Using the Mellin transform we have
\[\triangle^{-s} = \frac{1}{\Gamma(s)} \int_0^\infty (e^{-t \triangle} - P)t^{s-1} \,dt,\]
where $P$ denotes the orthogonal projection on $\text{Ker} (\triangle)$.
Therefore for any $a \in A_\theta^\infty$, we have
\[ a \triangle^{-s} = \frac{1}{\Gamma(s)} \int_0^\infty a(e^{-t \triangle} - P)t^{s-1} \,dt,\]
which gives 
\begin{eqnarray}
 \text{Trace} (a \triangle^{-s}) &=& \frac{1}{\Gamma(s)} \int_0^\infty \big ( \text{Trace}(a e^{-t \triangle}) - \text{Trace}(aP) \big )t^{s-1} \,dt \nonumber \\
 &=& \frac{1}{\Gamma(s)} \int_0^\infty \big ( \text{Trace}(a e^{-t \triangle}) - 2 \, \mathfrak{t}\,(a) \big )t^{s-1} \,dt. \nonumber
\end{eqnarray}
Appealing to the Cauchy integral formula \eqref{cauchyint} and using similar arguments to those in \cite{gil}, one can derive an asymptotic expansion: 
\[ \text{Trace}(ae^{-t \triangle}) \sim t^{-1} \sum_{n=0}^{\infty} B_{2n} (a, \triangle) \, t^n \qquad (t \to 0).  \]
Using this asymptotic expansion and the fact that $\Gamma$ has a simple pole at $0$ with residue $1$, one can see that the zeta function 
\[\zeta_a(s) = \textnormal{Trace} \,(a \triangle^{-s}), \qquad \textnormal{Re}(s) >> 0,\]
has a meromorphic extension to the whole plane with no pole at $0$, and
\[  \zeta_a(0) = B_2(a, \triangle) - 2 \mathfrak{t}\,(a).\]
Also one can see that 
\begin{eqnarray}
B_2(a, \triangle) &=& \frac{1}{2 \pi i} \int \int_C e^{-\lambda} \, \mathfrak{t}\,(a b_2(\xi, \lambda))\, d \lambda \, d \xi \nonumber \\
&=& \frac{1}{2 \pi i} \mathfrak{t}\, \Big (a \int \int_C e^{-\lambda} b_2(\xi, \lambda)\, d \lambda \, d \xi \Big ), \nonumber
\end{eqnarray}
where, as above, $b_2$ is the symbol of order $-4$ in $\sigma(B_\lambda)$. Hence the scalar curvature is equal to 
\[  \frac{1}{2 \pi i} \int_{\mathbb{R}^2} \int_C e^{-\lambda} b_2(\xi, \lambda) \, d\lambda \, d\xi . \]
\end{proof}
\end{localexpcur}
Note that for the purpose of computing the scalar curvature, using a homogeneity argument, one can set $\lambda = -1$
and multiply the final answer by $-1$  (\emph{cf.} \cite{contre, fatkha}). In the sequel, we will write $b_2$ for $b_2(\xi, -1)$.

\section{Computation of the Scalar Curvature} \label{computation}

Following the recipe given in Subsection \ref{recipe} we compute the two components of the scalar curvature for the noncommutative two 
torus corresponding to the Laplacian of the perturbed spectral triple attached to $(\mathbb{T}_\theta^2, \tau, k)$.

\subsection{The computations for $k \partial^* \partial k$.} \label{firstcomputation}
In \cite{fatkha}, it is shown that the symbol of the operator 
$k \partial^* \partial k$, `\emph{the Laplacian on functions}', is equal to  $a_2(\xi)+a_1(\xi) + a_0(\xi)$
where
\[  a_2(\xi)= \xi_1^2k^2+|\tau|^2\xi_2^2k^2+2\tau_1\xi_1\xi_2k^2,  \]
\[a_1(\xi)=2\xi_1k\delta_1(k) + 2|\tau|^2\xi_2k\delta_2(k) +
2\tau_1\xi_1k\delta_2(k)+2\tau_1\xi_2k\delta_1(k),\]
\[ a_0(\xi)= k\delta_1^2(k)+ |\tau|^2k\delta_2^2(k) + 2\tau_1k\delta_1\delta_2(k).\]
It is also shown that one can find the terms $b_j$ inductively. In fact the equation   
\begin{eqnarray}
&& (b_0+b_1+b_2+\cdots)\,\sigma(k \partial^* \partial k+1) = \nonumber \\
&& \qquad  \qquad (b_0+b_1+b_2+\cdots) ((a_2+1)+a_1+a_0) \sim 1, \nonumber
\end{eqnarray} 
has a solution where each $b_j$ can be chosen to be a symbol of order $-2-j$. In fact, treating 1 as a symbol of 
order 2, we let 
$a_2'=a_2+1, a_1'=a_1, a_0'=a_0$.  Then, the above equation yields
\[ \sum_{\substack{j, \ell_1, \ell_2 \geq 0,\\ k=0, 1, 2}} \frac{1}{\ell_1! \ell_2!}
\partial_1^{\ell_1} \partial_2^{\ell_2}(b_j)
\delta_1^{\ell_1}\delta_2^{\ell_2} (a_k') \sim 1.\]
By decomposing the latter into terms of order $-n$, $n=0,1,2, \dots$, we find
\begin{equation}
b_0=a_2'^{-1}=(a_2+1)^{-1}=(\xi_1^2k^2+|\tau|^2\xi_2^2k^2+
2\tau_1\xi_1\xi_2k^2+1)^{-1},\nonumber
\end{equation}
\begin{eqnarray}
b_1&=&
-(b_0a_1b_0  + \partial_1(b_0)\delta_1(a_2)b_0 +
\partial_2(b_0)\delta_2(a_2)b_0), \nonumber
\end{eqnarray}
\begin{eqnarray}
b_2&=&
-(b_0a_0b_0 + b_1a_1b_0 + \partial_1(b_0)\delta_1(a_1)b_0 +
\partial_2(b_0)\delta_2(a_1)b_0 + \nonumber \\ 
&& \quad \partial_1(b_1)\delta_1(a_2)b_0 + 
\partial_2(b_1)\delta_2(a_2)b_0 + (1/2)\partial_{11}(b_0)\delta_1^2(a_2)b_0 + \nonumber \\
&& \quad (1/2)\partial_{22}(b_0)\delta_2^2(a_2)b_0  +
\partial_{12}(b_0)\delta_{12}(a_2)b_0). \nonumber
\end{eqnarray}
After a direct computation, we find a lengthy formula for $b_2$ which we record in Appendix \ref{firstb2}.
In order to integrate $b_2$  over the $\xi$-plane we pass to the following coordinates
\begin{equation} \label{changevar}
\xi_1= r \cos \theta - r \frac{\tau_1}{\tau_2} \sin \theta,
\qquad
\xi_2=\frac{r}{\tau_2}\sin \theta ,
\end{equation}
where $\theta$ ranges from 0 to $2\pi$ and $r$ ranges from 0 to $\infty$. After the integration with respect to $\theta$, 
up to a factor of $\frac{r}{\tau_2}$ which is the Jacobian of the change of variables, one gets \\

\noindent
$
4|\tau|^2\pi r^6b_0^2k^2\delta_2(k)b_0^2k^3\delta_2(k)b_0k                                                                                                                                     
+ 4\tau_1\pi r^6b_0^2k^2\delta_2(k)b_0^2k^3\delta_1(k)b_0k
\\
+4|\tau|^2\pi r^6b_0^2k^2\delta_2(k)b_0^2k^4\delta_2(k)b_0                                                                                                                                    
+4\tau_1\pi r^6b_0^2k^2\delta_2(k)b_0^2k^4\delta_1(k)b_0                                                                                                                                     
\\
+4\tau_1\pi r^6b_0^2k^2\delta_1(k)b_0^2k^3\delta_2(k)b_0k                                                                                                                                    
+4\pi r^6b_0^2k^2\delta_1(k)b_0^2k^3\delta_1(k)b_0k     
\\                                                                                                                                
+4\tau_1\pi r^6b_0^2k^2\delta_1(k)b_0^2k^4\delta_2(k)b_0                                                                                                                                     
+4\pi r^6b_0^2k^2\delta_1(k)b_0^2k^4\delta_1(k)b_0                                                                                                                                     
\\
+4|\tau|^2\pi r^6b_0^2k^3\delta_2(k)b_0^2k^2\delta_2(k)b_0k                                                                                                                                     
+4\tau_1\pi r^6b_0^2k^3\delta_2(k)b_0^2k^2\delta_1(k)b_0k                                                                                                                                     
\\
+4|\tau|^2\pi r^6b_0^2k^3\delta_2(k)b_0^2k^3\delta_2(k)b_0                                                                                                                                    
+4\tau_1\pi r^6b_0^2k^3\delta_2(k)b_0^2k^3\delta_1(k)b_0                                                                                                                                     
\\
+4\tau_1\pi r^6b_0^2k^3\delta_1(k)b_0^2k^2\delta_2(k)b_0k                                                                                                                                    
+4\pi r^6b_0^2k^3\delta_1(k)b_0^2k^2\delta_1(k)b_0k                                                                                                                                     
\\
+4\tau_1\pi r^6b_0^2k^3\delta_1(k)b_0^2k^3\delta_2(k)b_0                                                                                                                                    
+4\pi r^6b_0^2k^3\delta_1(k)b_0^2k^3\delta_1(k)b_0                                                                                                                                    
\\
+8|\tau|^2\pi r^6b_0^3k^4\delta_2(k)b_0k\delta_2(k)b_0k                                                                                                                                     
+8\tau_1\pi r^6b_0^3k^4\delta_2(k)b_0k\delta_1(k)b_0k                                                                                                                                     
\\
+8|\tau|^2\pi r^6b_0^3k^4\delta_2(k)b_0k^2\delta_2(k)b_0                                                                                                                                     
+8\tau_1\pi r^6b_0^3k^4\delta_2(k)b_0k^2\delta_1(k)b_0                                                                                                                                   
\\
+8\tau_1\pi r^6b_0^3k^4\delta_1(k)b_0k\delta_2(k)b_0k                                                                                                                                    
+8\pi r^6b_0^3k^4\delta_1(k)b_0k\delta_1(k)b_0k                                                                                                                                     
\\
+8\tau_1\pi r^6b_0^3k^4\delta_1(k)b_0k^2\delta_2(k)b_0                                                                                                                                    
+8\pi r^6b_0^3k^4\delta_1(k)b_0k^2\delta_1(k)b_0                                                                                                                                    
\\
+8|\tau|^2\pi r^6b_0^3k^5\delta_2(k)b_0\delta_2(k)b_0k                                                                                                                                     
+8\tau_1\pi r^6b_0^3k^5\delta_2(k)b_0\delta_1(k)b_0k                                                                                                                                     
\\
+8|\tau|^2\pi r^6b_0^3k^5\delta_2(k)b_0k\delta_2(k)b_0                                                                                                                                    
+8\tau_1\pi r^6b_0^3k^5\delta_2(k)b_0k\delta_1(k)b_0                                                                                                                                     
\\
+8\tau_1\pi r^6b_0^3k^5\delta_1(k)b_0\delta_2(k)b_0k                                                                                                                                     
+8\pi r^6b_0^3k^5\delta_1(k)b_0\delta_1(k)b_0k                                                                                                                                    
\\
+8\tau_1\pi r^6b_0^3k^5\delta_1(k)b_0k\delta_2(k)b_0                                                                                                                                    
+8\pi r^6b_0^3k^5\delta_1(k)b_0k\delta_1(k)b_0                                                                                                                                     
\\
-4|\tau|^2\pi r^4b_0k\delta_2(k)b_0^2k^2\delta_2(k)b_0k                                                                                                                                     
-4\tau_1\pi r^4b_0k\delta_2(k)b_0^2k^2\delta_1(k)b_0k                                                                                                                                    
\\
-4|\tau|^2\pi r^4b_0k\delta_2(k)b_0^2k^3\delta_2(k)b_0                                                                                                                                     
-4\tau_1\pi r^4b_0k\delta_2(k)b_0^2k^3\delta_1(k)b_0                                                                                                                                    
\\
-4\tau_1\pi r^4b_0k\delta_1(k)b_0^2k^2\delta_2(k)b_0k                                                                                                                                     
-4\pi r^4b_0k\delta_1(k)b_0^2k^2\delta_1(k)b_0k                                                                                                                                   
\\
-4\tau_1\pi r^4b_0k\delta_1(k)b_0^2k^3\delta_2(k)b_0                                                                                                                                    
-4\pi r^4b_0k\delta_1(k)b_0^2k^3\delta_1(k)b_0                                                                                                                                    
\\
-8|\tau|^2\pi r^4b_0^2k^2\delta_2(k)b_0k\delta_2(k)b_0k                                                                                                                                     
-8\tau_1\pi r^4b_0^2k^2\delta_2(k)b_0k\delta_1(k)b_0k                                                                                                                                     
\\
-12|\tau|^2\pi r^4b_0^2k^2\delta_2(k)b_0k^2\delta_2(k)b_0                                                                                                                                     
-12\tau_1\pi r^4b_0^2k^2\delta_2(k)b_0k^2\delta_1(k)b_0                                                                                                                                    
\\
-8\tau_1\pi r^4b_0^2k^2\delta_1(k)b_0k\delta_2(k)b_0k                                                                                                                                    
-8\pi r^4b_0^2k^2\delta_1(k)b_0k\delta_1(k)b_0k                                                                                                                                    
\\
-12\tau_1\pi r^4b_0^2k^2\delta_1(k)b_0k^2\delta_2(k)b_0                                                                                                                                     
-12\pi r^4b_0^2k^2\delta_1(k)b_0k^2\delta_1(k)b_0                                                                                                                                     
\\
-12|\tau|^2\pi r^4b_0^2k^3\delta_2(k)b_0\delta_2(k)b_0k                                                                                                                                     
-12\tau_1\pi r^4b_0^2k^3\delta_2(k)b_0\delta_1(k)b_0k                                                                                                                                     
\\
-16|\tau|^2\pi r^4b_0^2k^3\delta_2(k)b_0k\delta_2(k)b_0                                                                                                                                    
-16\tau_1\pi r^4b_0^2k^3\delta_2(k)b_0k\delta_1(k)b_0                                                                                                                                     
\\
-12\tau_1\pi r^4b_0^2k^3\delta_1(k)b_0\delta_2(k)b_0k                                                                                                                                    
-12\pi r^4b_0^2k^3\delta_1(k)b_0\delta_1(k)b_0k                                                                                                                                     
\\
-16\tau_1\pi r^4b_0^2k^3\delta_1(k)b_0k\delta_2(k)b_0                                                                                                                                     
-16\pi r^4b_0^2k^3\delta_1(k)b_0k\delta_1(k)b_0                                                                                                                                    
\\
-4\pi r^4b_0^3k^4\delta_1^2(k)b_0k                                                                                                                                    
-4|\tau|^2\pi r^4b_0^3k^4\delta_2^2(k)b_0k                                                                                                                                    
-8\tau_1\pi r^4b_0^3k^4\delta_1\delta_2(k)b_0k                                                                                                                                     
\\-8|\tau|^2\pi r^4b_0^3k^4\delta_2(k)^2b_0                                                                                                                                     
-8\tau_1\pi r^4b_0^3k^4\delta_2(k)\delta_1(k)b_0                                                                                                                                   
-8\tau_1\pi r^4b_0^3k^4\delta_1(k)\delta_2(k)b_0                                                                                                                                    
\\
-8\pi r^4b_0^3k^4\delta_1(k)^2b_0                                                                                                                                    
-4\pi r^4b_0^3k^5\delta_1^2(k)b_0                                                                                                                                     
-4|\tau|^2\pi r^4b_0^3k^5\delta_2^2(k)b_0                                                                                                                                     
\\
-8\tau_1\pi r^4b_0^3k^5\delta_1\delta_2(k)b_0                                                                                                                                     
+4|\tau|^2\pi r^2b_0k\delta_2(k)b_0\delta_2(k)b_0k                                                                                                                                   
\\
+4\tau_1\pi r^2b_0k\delta_2(k)b_0\delta_1(k)b_0k                                                                                                                                     
+8|\tau|^2\pi r^2b_0k\delta_2(k)b_0k\delta_2(k)b_0                                                                                                                                     
\\
+8\tau_1\pi r^2b_0k\delta_2(k)b_0k\delta_1(k)b_0                                                                                                                                    
+4\tau_1\pi r^2b_0k\delta_1(k)b_0\delta_2(k)b_0k                                                                                                                                     
\\
+4\pi r^2b_0k\delta_1(k)b_0\delta_1(k)b_0k                                                                                                                                    
+8\tau_1\pi r^2b_0k\delta_1(k)b_0k\delta_2(k)b_0                                                                                                                                    
\\
+8\pi r^2b_0k\delta_1(k)b_0k\delta_1(k)b_0                                                                                                                                     
+2\pi r^2b_0^2k^2\delta_1^2(k)b_0k                                                                                                                                    
+2|\tau|^2\pi r^2b_0^2k^2\delta_2^2(k)b_0k                                                                                                                                    
\\
+4\tau_1\pi r^2b_0^2k^2\delta_1\delta_2(k)b_0k                                                                                                                                    
+8|\tau|^2\pi r^2b_0^2k^2\delta_2(k)^2b_0                                                                                                                                     
+8\tau_1\pi r^2b_0^2k^2\delta_2(k)\delta_1(k)b_0                                                                                                                                    
\\
+8\tau_1\pi r^2b_0^2k^2\delta_1(k)\delta_2(k)b_0                                                                                                                                    
+8\pi r^2b_0^2k^2\delta_1(k)^2b_0                                                                                                                                     
+6\pi r^2b_0^2k^3\delta_1^2(k)b_0                                                                                                                                     
\\
+6|\tau|^2\pi r^2b_0^2k^3\delta_2^2(k)b_0                                                                                                                                     
+12\tau_1\pi r^2b_0^2k^3\delta_1\delta_2(k)b_0                                                                                                                                     
-2\pi b_0k\delta_1^2(k)b_0                                                                                                                                     
\\
-2|\tau|^2\pi b_0k\delta_2^2(k)b_0                                                                                                                                    
-4\tau_1\pi b_0k\delta_1\delta_2(k)b_0
 $\\
  
\noindent where
\[b_0=(r^2k^2+1)^{-1}.\]
Since we are in the noncommutative case, where 
$b_0=(r^2k^2+1)^{-1}$ and $\delta_j(k), j=1,2,$ do not necessarily commute, for the computation of 
$\int_{0}^{\infty} \bullet \, r dr$ of these terms, we need to use the modular  automorphism 
$\Delta$ to permute $k$ with elements of $A_{\theta}^\infty$. In the next two subsections we explain how this 
calculation is performed for the above types of terms.

\subsubsection{Terms with two factors of the form $b_0^i$, $i \geq 1$.}

These are the following terms:\\

\noindent
$  
-4\pi r^4b_0^3k^4\delta_1^2(k)b_0k                                                                                                                                    
-4|\tau|^2\pi r^4b_0^3k^4\delta_2^2(k)b_0k                                                                                                                                    
-8\tau_1\pi r^4b_0^3k^4\delta_1\delta_2(k)b_0k                                                                                                                                     
\\-8|\tau|^2\pi r^4b_0^3k^4\delta_2(k)^2b_0                                                                                                                                     
-8\tau_1\pi r^4b_0^3k^4\delta_2(k)\delta_1(k)b_0                                                                                                                                   
-8\tau_1\pi r^4b_0^3k^4\delta_1(k)\delta_2(k)b_0                                                                                                                                    
\\
-8\pi r^4b_0^3k^4\delta_1(k)^2b_0                                                                                                                                    
-4\pi r^4b_0^3k^5\delta_1^2(k)b_0                                                                                                                                     
-4|\tau|^2\pi r^4b_0^3k^5\delta_2^2(k)b_0                                                                                                                                     
\\
-8\tau_1\pi r^4b_0^3k^5\delta_1\delta_2(k)b_0    
+2\pi r^2b_0^2k^2\delta_1^2(k)b_0k                                                                                                                                    
+2|\tau|^2\pi r^2b_0^2k^2\delta_2^2(k)b_0k                                                                                                                                    
\\
+4\tau_1\pi r^2b_0^2k^2\delta_1\delta_2(k)b_0k                                                                                                                                    
+8|\tau|^2\pi r^2b_0^2k^2\delta_2(k)^2b_0                                                                                                                                     
+8\tau_1\pi r^2b_0^2k^2\delta_2(k)\delta_1(k)b_0                                                                                                                                    
\\
+8\tau_1\pi r^2b_0^2k^2\delta_1(k)\delta_2(k)b_0                                                                                                                                    
+8\pi r^2b_0^2k^2\delta_1(k)^2b_0                                                                                                                                     
+6\pi r^2b_0^2k^3\delta_1^2(k)b_0                                                                                                                                     
\\
+6|\tau|^2\pi r^2b_0^2k^3\delta_2^2(k)b_0                                                                                                                                     
+12\tau_1\pi r^2b_0^2k^3\delta_1\delta_2(k)b_0                                                                                                                                     
-2\pi b_0k\delta_1^2(k)b_0                                                                                                                                     
-2|\tau|^2\pi b_0k\delta_2^2(k)b_0                                                                                                                                    
\\
-4\tau_1\pi b_0k\delta_1\delta_2(k)b_0.
$\\

The computation of $\int_{0}^{\infty} \bullet \, r dr$ of these terms is
achieved by the following lemma of Connes and Tretkoff proved
in \cite{contre}. 

\newtheorem{integration}[main1]{Lemma}
\begin{integration} \label{integration}
For any $\rho \in A_{\theta}^{\infty}$ and every non-negative
integer $m$, one has
\[\int_0^{\infty}\frac{k^{2m+2}u^m}{(k^2u+1)^{m+1}}\rho
\frac{1}{(k^2u+1)}du=\mathcal{D}_m(\rho),\]
where $\mathcal{D}_m=\mathcal{L}_m(\Delta)$, $\Delta$ is
the modular automorphism introduced in
Section \ref{prel}, and $\mathcal{L}_m$ is the
modified logarithm:
\begin{eqnarray}
\mathcal{L}_m(u)
&=&\int_0^{\infty}\frac{x^m}{(x+1)^{m+1}}\frac{1}{(xu+1)}dx \nonumber \\
&=&(-1)^m(u-1)^{-(m+1)}\big( \log u - \sum_{j=1}^{m}(-1)^{j+1}
\frac{(u-1)^j}{j} \big). \nonumber
\end{eqnarray}
\end{integration}

Using this lemma, one can see that $\int_{0}^{\infty} \bullet \, r dr$ of the above terms, up to an overall factor of $\pi$, is equal to \\

\noindent
$
-2 \mathcal{D}_2 \Delta^{1/2} (k^{-1} \delta_1^2(k)) 
-2 |\tau|^2 \mathcal{D}_2 \Delta^{1/2} (k^{-1} \delta_2^2(k))
-4 \tau_1 \mathcal{D}_2\Delta^{1/2} (k^{-1} \delta_1\delta_2(k)) \\
-4 |\tau|^2 \mathcal{D}_2 (k^{-2} \delta_2(k)^2) 
-4\tau_1 \mathcal{D}_2 (k^{-2} \delta_2(k) \delta_1(k)) 
-4\tau_1 \mathcal{D}_2 (k^{-2} \delta_1(k) \delta_2(k)) \\
-4  \mathcal{D}_2 (k^{-2} \delta_1(k)^2)
-2 \mathcal{D}_2(k^{-1} \delta_1^2(k))
-2 |\tau|^2 \mathcal{D}_2(k^{-1} \delta_2^2(k))
-4 \tau_1 \mathcal{D}_2(k^{-1} \delta_1 \delta_2(k))\\
+\mathcal{D}_1\Delta^{1/2} (k^{-1} \delta_1^2(k))
+|\tau|^2 \mathcal{D}_1\Delta^{1/2} (k^{-1} \delta_2^2(k))
+2 \tau_1 \mathcal{D}_1 \Delta^{1/2} (k^{-1} \delta_1\delta_2(k))\\
+4 |\tau|^2 \mathcal{D}_1(k^{-2} \delta_2(k)^2)
+4 \tau_1 \mathcal{D}_1 (k^{-2} \delta_2(k) \delta_1(k))
+4 \tau_1 \mathcal{D}_1 (k^{-2} \delta_1(k) \delta_2(k))\\
+4  \mathcal{D}_1 (k^{-2} \delta_1(k)^2)
+3 \mathcal{D}_1 (k^{-1} \delta_1^2(k))
+3 |\tau|^2 \mathcal{D}_1 (k^{-1} \delta_2^2(k))
+6 \tau_1 \mathcal{D}_1(k^{-1} \delta_1 \delta_2(k)) \\
- \mathcal{D}_0 (k^{-1} \delta_1^2(k))
- |\tau|^2 \mathcal{D}_0 (k^{-1} \delta_2^2(k))
-2 \tau_1  \mathcal{D}_0  (k^{-1} \delta_1\delta_2(k)).
$\\

Hence, up to an overall factor of $\pi$, $\int_0^\infty \bullet \,rdr$ of the terms with two positive powers of $b_0$ is equal to 
\begin{eqnarray}\label{twob0int}
&&f_1(\Delta)(k^{-1} \delta_1^2(k)) + f_2(\Delta)(k^{-2} \delta_1(k)^2) \nonumber \\ 
&+&   |\tau|^2 f_1(\Delta)(k^{-1} \delta_2^2(k)) + |\tau|^2 f_2(\Delta)(k^{-2} \delta_2(k)^2) \nonumber \\
&+& \tau_1 f_1(\Delta) (k^{-1} \delta_1\delta_2(k)) + \tau_1 f_2(\Delta)(k^{-2}\delta_1(k)\delta_2(k)) \nonumber \\
&+& \tau_1 f_1(\Delta) (k^{-1} \delta_2\delta_1(k)) + \tau_1 f_2(\Delta)(k^{-2}\delta_2(k)\delta_1(k)),
\end{eqnarray}
where 
\begin{eqnarray} \label{fone}
f_1(u) &:=& -2 \mathcal{L}_2(u) u^{1/2} - 2 \mathcal{L}_2(u) + \mathcal{L}_1(u)u^{1/2} + 3 \mathcal{L}_1(u) - \mathcal{L}_0(u) \nonumber \\
&=& -\frac{u^{1/2} (2-2u+(1+u) \log{u})}{(-1+u^{1/2})^3(1+u^{1/2})^2},
\end{eqnarray}
and 
\begin{equation} \label{ftwo}
f_2(u):= -4 \mathcal{L}_2(u) + 4 \mathcal{L}_1(u) = 2\frac{-1+u^2-2u \log{u}}{(-1+u)^3}.
\end{equation}

\subsubsection{Terms with three factors of the form $b_0^i$, $i \geq 1$.} 

These terms are the following:\\

\noindent
$
4|\tau|^2\pi r^6b_0^2k^2\delta_2(k)b_0^2k^3\delta_2(k)b_0k                                                                                                                                     
+ 4\tau_1\pi r^6b_0^2k^2\delta_2(k)b_0^2k^3\delta_1(k)b_0k
\\
+4|\tau|^2\pi r^6b_0^2k^2\delta_2(k)b_0^2k^4\delta_2(k)b_0                                                                                                                                    
+4\tau_1\pi r^6b_0^2k^2\delta_2(k)b_0^2k^4\delta_1(k)b_0                                                                                                                                     
\\
+4\tau_1\pi r^6b_0^2k^2\delta_1(k)b_0^2k^3\delta_2(k)b_0k                                                                                                                                    
+4\pi r^6b_0^2k^2\delta_1(k)b_0^2k^3\delta_1(k)b_0k     
\\                                                                                                                                
+4\tau_1\pi r^6b_0^2k^2\delta_1(k)b_0^2k^4\delta_2(k)b_0                                                                                                                                     
+4\pi r^6b_0^2k^2\delta_1(k)b_0^2k^4\delta_1(k)b_0                                                                                                                                     
\\
+4|\tau|^2\pi r^6b_0^2k^3\delta_2(k)b_0^2k^2\delta_2(k)b_0k                                                                                                                                     
+4\tau_1\pi r^6b_0^2k^3\delta_2(k)b_0^2k^2\delta_1(k)b_0k                                                                                                                                     
\\
+4|\tau|^2\pi r^6b_0^2k^3\delta_2(k)b_0^2k^3\delta_2(k)b_0                                                                                                                                    
+4\tau_1\pi r^6b_0^2k^3\delta_2(k)b_0^2k^3\delta_1(k)b_0                                                                                                                                     
\\
+4\tau_1\pi r^6b_0^2k^3\delta_1(k)b_0^2k^2\delta_2(k)b_0k                                                                                                                                    
+4\pi r^6b_0^2k^3\delta_1(k)b_0^2k^2\delta_1(k)b_0k                                                                                                                                     
\\
+4\tau_1\pi r^6b_0^2k^3\delta_1(k)b_0^2k^3\delta_2(k)b_0                                                                                                                                    
+4\pi r^6b_0^2k^3\delta_1(k)b_0^2k^3\delta_1(k)b_0                                                                                                                                    
\\
+8|\tau|^2\pi r^6b_0^3k^4\delta_2(k)b_0k\delta_2(k)b_0k                                                                                                                                     
+8\tau_1\pi r^6b_0^3k^4\delta_2(k)b_0k\delta_1(k)b_0k                                                                                                                                     
\\
+8|\tau|^2\pi r^6b_0^3k^4\delta_2(k)b_0k^2\delta_2(k)b_0                                                                                                                                     
+8\tau_1\pi r^6b_0^3k^4\delta_2(k)b_0k^2\delta_1(k)b_0                                                                                                                                   
\\
+8\tau_1\pi r^6b_0^3k^4\delta_1(k)b_0k\delta_2(k)b_0k                                                                                                                                    
+8\pi r^6b_0^3k^4\delta_1(k)b_0k\delta_1(k)b_0k                                                                                                                                     
\\
+8\tau_1\pi r^6b_0^3k^4\delta_1(k)b_0k^2\delta_2(k)b_0                                                                                                                                    
+8\pi r^6b_0^3k^4\delta_1(k)b_0k^2\delta_1(k)b_0                                                                                                                                    
\\
+8|\tau|^2\pi r^6b_0^3k^5\delta_2(k)b_0\delta_2(k)b_0k                                                                                                                                     
+8\tau_1\pi r^6b_0^3k^5\delta_2(k)b_0\delta_1(k)b_0k                                                                                                                                     
\\
+8|\tau|^2\pi r^6b_0^3k^5\delta_2(k)b_0k\delta_2(k)b_0                                                                                                                                    
+8\tau_1\pi r^6b_0^3k^5\delta_2(k)b_0k\delta_1(k)b_0                                                                                                                                     
\\
+8\tau_1\pi r^6b_0^3k^5\delta_1(k)b_0\delta_2(k)b_0k                                                                                                                                     
+8\pi r^6b_0^3k^5\delta_1(k)b_0\delta_1(k)b_0k                                                                                                                                    
\\
+8\tau_1\pi r^6b_0^3k^5\delta_1(k)b_0k\delta_2(k)b_0                                                                                                                                    
+8\pi r^6b_0^3k^5\delta_1(k)b_0k\delta_1(k)b_0                                                                                                                                     
\\
-4|\tau|^2\pi r^4b_0k\delta_2(k)b_0^2k^2\delta_2(k)b_0k                                                                                                                                     
-4\tau_1\pi r^4b_0k\delta_2(k)b_0^2k^2\delta_1(k)b_0k                                                                                                                                    
\\
-4|\tau|^2\pi r^4b_0k\delta_2(k)b_0^2k^3\delta_2(k)b_0                                                                                                                                     
-4\tau_1\pi r^4b_0k\delta_2(k)b_0^2k^3\delta_1(k)b_0                                                                                                                                    
\\
-4\tau_1\pi r^4b_0k\delta_1(k)b_0^2k^2\delta_2(k)b_0k                                                                                                                                     
-4\pi r^4b_0k\delta_1(k)b_0^2k^2\delta_1(k)b_0k                                                                                                                                   
\\
-4\tau_1\pi r^4b_0k\delta_1(k)b_0^2k^3\delta_2(k)b_0                                                                                                                                    
-4\pi r^4b_0k\delta_1(k)b_0^2k^3\delta_1(k)b_0                                                                                                                                    
\\
-8|\tau|^2\pi r^4b_0^2k^2\delta_2(k)b_0k\delta_2(k)b_0k                                                                                                                                     
-8\tau_1\pi r^4b_0^2k^2\delta_2(k)b_0k\delta_1(k)b_0k                                                                                                                                     
\\
-12|\tau|^2\pi r^4b_0^2k^2\delta_2(k)b_0k^2\delta_2(k)b_0                                                                                                                                     
-12\tau_1\pi r^4b_0^2k^2\delta_2(k)b_0k^2\delta_1(k)b_0                                                                                                                                    
\\
-8\tau_1\pi r^4b_0^2k^2\delta_1(k)b_0k\delta_2(k)b_0k                                                                                                                                    
-8\pi r^4b_0^2k^2\delta_1(k)b_0k\delta_1(k)b_0k                                                                                                                                    
\\
-12\tau_1\pi r^4b_0^2k^2\delta_1(k)b_0k^2\delta_2(k)b_0                                                                                                                                     
-12\pi r^4b_0^2k^2\delta_1(k)b_0k^2\delta_1(k)b_0                                                                                                                                     
\\
-12|\tau|^2\pi r^4b_0^2k^3\delta_2(k)b_0\delta_2(k)b_0k                                                                                                                                     
-12\tau_1\pi r^4b_0^2k^3\delta_2(k)b_0\delta_1(k)b_0k                                                                                                                                     
\\
-16|\tau|^2\pi r^4b_0^2k^3\delta_2(k)b_0k\delta_2(k)b_0                                                                                                                                    
-16\tau_1\pi r^4b_0^2k^3\delta_2(k)b_0k\delta_1(k)b_0                                                                                                                                     
\\
-12\tau_1\pi r^4b_0^2k^3\delta_1(k)b_0\delta_2(k)b_0k                                                                                                                                    
-12\pi r^4b_0^2k^3\delta_1(k)b_0\delta_1(k)b_0k                                                                                                                                     
\\
-16\tau_1\pi r^4b_0^2k^3\delta_1(k)b_0k\delta_2(k)b_0                                                                                                                                     
-16\pi r^4b_0^2k^3\delta_1(k)b_0k\delta_1(k)b_0        
\\
+4|\tau|^2\pi r^2b_0k\delta_2(k)b_0\delta_2(k)b_0k                                                                                                                                   
+4\tau_1\pi r^2b_0k\delta_2(k)b_0\delta_1(k)b_0k                                                                                                                                     
\\
+8|\tau|^2\pi r^2b_0k\delta_2(k)b_0k\delta_2(k)b_0                                                                                                                                     
+8\tau_1\pi r^2b_0k\delta_2(k)b_0k\delta_1(k)b_0                                                                                                                                    
\\
+4\tau_1\pi r^2b_0k\delta_1(k)b_0\delta_2(k)b_0k                                                                                                                                     
+4\pi r^2b_0k\delta_1(k)b_0\delta_1(k)b_0k                                                                                                                                    
\\+8\tau_1\pi r^2b_0k\delta_1(k)b_0k\delta_2(k)b_0                                                                                                                                    
+8\pi r^2b_0k\delta_1(k)b_0k\delta_1(k)b_0.  
$\\

For computing $\int_{0}^{\infty} \bullet \, rdr$ of these terms we will use the following lemma in which two variable modified 
logarithm functions appear. This lemma can be proved along the same lines as in \cite{contre}.

\newtheorem{integrationtwo}[main1]{Lemma}
\begin{integrationtwo} \label{integrationtwo}
For any $\rho, \rho' \in A_\theta^\infty$ and positive integers $m, m'$, we have
\begin{eqnarray}
\int_0^\infty \frac{1}{(k^2u+1)^m} \rho \frac{k^{2(m+m')} u^{m+m'-1}}{(k^2u+1)^{m'}}\rho' 
\frac{1}{k^2u+1} \, du 
=\mathcal{D}_{m,m'}(\Delta_{(1)}, \Delta_{(2)}) (\rho \rho'),  \nonumber
\end{eqnarray}
where the function $\mathcal{D}_{m, m'}$ is defined by 
\[\mathcal{D}_{m, m'}(u, v) = \int_0^\infty \frac{1}{(xu^{-1}+1)^m} \frac{x^{m+m'-1}}{(x+1)^{m'}}
 \frac{1}{xv+1} \,dx,\]
and $\Delta_{(1)}$ and $\Delta_{(2)}$ respectively denote the action of $\Delta$ on the first and second factor of the product.

\begin{proof}
Using the change of variable $s = \log u + h,$ we have:
\begin{eqnarray}
&& \int_0^\infty \frac{1}{(k^2u+1)^m} \rho \frac{k^{2(m+m')} u^{m+m'-1}}{(k^2u+1)^{m'}}\rho' \frac{1}{k^2u+1} \, du \nonumber \\
&=& \int_{-\infty}^\infty \frac{1}{(e^s+1)^m} \rho \frac{e^{s(m+m'-1)}k^2}{(e^s+1)^{m'}} \rho' \frac{k^{-2}}{e^s+1} \,d(e^s) \nonumber \\
&=&  \int_{-\infty}^\infty \frac{1}{(e^s+1)^m} \rho \frac{e^{s(m+m'-1/2)}}{(e^s+1)^{m'}} \Delta^{-1/2}(\rho') \frac{e^{s/2}}{e^s+1} \,ds \nonumber \\
&=&  \int_{-\infty}^\infty \frac{1}{(e^s+1)^m} \rho \frac{e^{s(m+m'-1/2)}}{(e^s+1)^{m'}} \Delta^{-1/2}(\rho')
 \int_{-\infty}^\infty\frac{e^{its}}{e^{\pi t}+e^{-\pi t}} \,dt \,ds \nonumber 
\end{eqnarray}
\begin{eqnarray}
&=&  \int_{-\infty}^\infty \frac{1}{(e^s+1)^m} \rho \frac{e^{s(m+m'-1/2)}}{(e^s+1)^{m'}} 
 \int_{-\infty}^\infty\frac{e^{its}}{e^{\pi t}+e^{-\pi t}} \Delta^{it-1/2}(\rho')\,dt \,ds \nonumber \\
 &=& \int_{-\infty}^\infty \frac{1}{(e^s+1)^m} \rho \frac{e^{s(m+m'-1/2)}}{(e^s+1)^{m'}} 
 \frac{e^{(s+\log \Delta)/2}}{e^{s+\log \Delta}+1} \Delta^{-1/2}(\rho') \,ds \nonumber \\
&=& \int_{-\infty}^\infty \frac{e^{ms/2}}{(e^s+1)^m}  \Delta^{m/2} (\rho) \frac{e^{s(m/2+m'-1/2)}}{(e^s+1)^{m'}} 
 \frac{e^{(s+\log \Delta)/2}}{e^{s+\log \Delta}+1} \Delta^{-1/2}(\rho') \,ds \nonumber \\ 
&=& \int_{-\infty}^\infty \prod_{j=1}^m \int_{-\infty}^\infty \frac{e^{it_js}}{e^{\pi t_j}+e^{-\pi t_j}} \,dt_j \, \Delta^{m/2} (\rho) \frac{e^{s(m/2+m'-1/2)}}{(e^s+1)^{m'}} 
\times \nonumber \\
&& \qquad \qquad \qquad \qquad  \qquad \qquad \qquad \qquad \qquad \frac{e^{(s+\log \Delta)/2}}{e^{s+\log \Delta}+1} \Delta^{-1/2}(\rho') \,ds \nonumber 
\end{eqnarray}
\begin{eqnarray}
&=& \int_{-\infty}^\infty  \int_{-\infty}^\infty \cdots \int_{-\infty}^\infty  \frac{\Delta^{-is \sum_{j=1}^m t_j +m/2} (\rho)e^{is \sum_{j=1}^m t_j}}{ \prod_{j=1}^m (e^{\pi t_j}+e^{-\pi t_j})} \,dt_1 \cdots \,dt_m \times  \nonumber \\
&& \qquad \qquad \qquad \qquad  \qquad \qquad  \,  \frac{e^{s(m/2+m'-1/2)}}{(e^s+1)^{m'}} 
 \frac{e^{(s+\log \Delta)/2}}{e^{s+\log \Delta}+1} \Delta^{-1/2}(\rho') \,ds \nonumber \\   
 &=& \int_{-\infty}^\infty  \int_{-\infty}^\infty \cdots \int_{-\infty}^\infty  \frac{\Delta_{(1)}^{-is \sum_{j=1}^m t_j} e^{is \sum_{j=1}^m t_j}}{ \prod_{j=1}^m (e^{\pi t_j}+e^{-\pi t_j})} \,dt_1 \cdots \,dt_m \times  \nonumber \\
&& \qquad \qquad \qquad   \,  \frac{e^{s(m/2+m'-1/2)}}{(e^s+1)^{m'}} 
 \frac{e^{(s+\log \Delta_{(2)})/2}}{e^{s+\log \Delta_{(2)}}+1} \,ds \, (  \Delta^{m/2}(\rho) \Delta^{-1/2}(\rho')) \nonumber
\end{eqnarray}
\[
=\int_{-\infty}^\infty \frac{e^{m(s- \log \Delta_{(1)})/2}}{(e^{s-\log \Delta_{(1)}}+1)^m}\frac{e^{s(m/2+m'-1/2)}}{(e^s+1)^{m'}} 
 \frac{e^{(s+\log \Delta_{(2)})/2}}{e^{s+\log \Delta_{(2)}}+1} \,ds \,  (  \Delta^{m/2}(\rho) \Delta^{-1/2}(\rho')) \nonumber 
\]
\begin{eqnarray}
&=& \int_{-\infty}^\infty \frac{e^{s/2} \Delta_{(1)}^{-m/2}}{(e^{s-\log \Delta_{(1)}}+1)^m}\frac{e^{s(m+m'-1)}}{(e^s+1)^{m'}} 
 \frac{e^{(s+\log \Delta_{(2)})/2}}{e^{s+\log \Delta_{(2)}}+1} \,ds \, (  \Delta^{m/2}(\rho) \Delta^{-1/2}(\rho')) \nonumber \\ 
&=& \int_0^\infty \frac{1}{(x \Delta_{(1)}^{-1}+1)^m}\frac{x^{m+m'-1}}{(x+1)^{m'}} 
 \frac{1}{x\Delta_{(2)}+1} \,dx \, (  \rho \rho'). \nonumber 
\end{eqnarray} 
\end{proof}
\end{integrationtwo}
In this paper, we only need these cases for our computations:
\begin{eqnarray}
\mathcal{D}_{1,1}(u, v)&=& ((-1 + v) \log[1/u] - (-1 + 1/u) \log[v])/ \nonumber \\ 
&& \qquad \qquad \qquad \qquad((-1 + 1/u) (-1 + v) (-(1/u) +v)); \nonumber
\end{eqnarray}
\begin{eqnarray}
&&\mathcal{D}_{2,2}(u,v)=(u ((-1 + v) ((-1 + 1/u) (1/u - v) (1 + 1/ u^2 - (1 + 1/u) v) + \nonumber \\
&& \qquad ((-1 + 3/u - 2 v) (-1 + v) \log[1/u])/u) - ((-1 + 1/u)^3 \log[v])/u))/ \nonumber \\ 
&&\qquad ((-1 + 1/u)^3 (1/u - v)^2 (-1 + v)^2); \nonumber
\end{eqnarray}
\begin{eqnarray}
&&\mathcal{D}_{1,2}(u,v)=  ((-1 + v)^2 \log[1/u] + (-1 + 1/u) ((1/u - v) (-1 + v) - \nonumber \\
&& \qquad \qquad (-1 + 1/u) \log[v]))/ ((-1 +1/u)^2 (1/u - v) (-1 + v)^2); \nonumber
\end{eqnarray}
\begin{eqnarray}
\mathcal{D}_{2,1}(u,v) &=& (
  u ((-1 + v) ((-1 + 1/u) (1/u - v) + ((1 - 2/u + v) \log[1/u])/ \nonumber \\
        &&u) + ((-1 + 1/u)^2 \log[v])/u))/((-1 + 1/u)^2 (1/u - v)^2 (-1 +v)); \nonumber
\end{eqnarray}
\begin{eqnarray}
\mathcal{D}_{3,1}(u,v)&=&(u^2 ((-1 + v) ((-1 + 1/u) (1/u - v) (5/u^2 + v - (3 (1 + v))/u) - \nonumber \\ 
&&(2 (1 + 3/u^2 + v + v^2 - (3 (1 + v))/u) \log[1/u])/u^2) + \nonumber \\
&&( 2 (-1 + 1/u)^3 \log[v])/u^2))/(2 (-1 + 1/u)^3 (1/u - v)^3 (-1 + v)). \nonumber
\end{eqnarray}

Using this lemma, $\int_0^\infty \bullet \,r dr$ of the terms listed in the beginning of this subsection, up to an overall factor of $\pi$, is equal to\\

\noindent
$
2 \mathcal{D}_{2,2}(\Delta_{(1)}, \Delta_{(2)})( \Delta^{-1}(\delta_1(k)k^{-1}) \Delta^{1/2} (k^{-1}\delta_1(k)))\\
+ 2  \mathcal{D}_{2,2}(\Delta_{(1)}, \Delta_{(2)})( \Delta^{-1}( \delta_1(k) k^{-1}) (k^{-1}\delta_1(k)) )\\
+ 2  \mathcal{D}_{2,2}(\Delta_{(1)}, \Delta_{(2)})( \Delta^{-3/2}( \delta_1(k) k^{-1}) \Delta^{1/2} (k^{-1}\delta_1(k)) )\\
+ 2  \mathcal{D}_{2,2}(\Delta_{(1)}, \Delta_{(2)})( \Delta^{-3/2}( \delta_1(k) k^{-1}) (k^{-1}\delta_1(k)) )\\
+ 4 \mathcal{D}_{3,1}(\Delta_{(1)}, \Delta_{(2)})( \Delta^{-2}(\delta_1(k)k^{-1} ) \Delta^{1/2}(k^{-1} \delta_1(k)) )\\
+ 4 \mathcal{D}_{3,1}(\Delta_{(1)}, \Delta_{(2)})( \Delta^{-2}( \delta_1(k)k^{-1}) (k^{-1}\delta_1(k)) )\\
+ 4 \mathcal{D}_{3,1}(\Delta_{(1)}, \Delta_{(2)})(\Delta^{-5/2}(\delta_1(k)k^{-1} )\Delta^{1/2}(k^{-1} \delta_1(k)) )\\
+ 4 \mathcal{D}_{3,1}(\Delta_{(1)}, \Delta_{(2)})(\Delta^{-5/2}( \delta_1(k)k^{-1})(k^{-1} \delta_1(k)) )\\
- 2 \mathcal{D}_{1,2}(\Delta_{(1)}, \Delta_{(2)})(\Delta^{-1/2}( \delta_1(k)k^{-1})\Delta^{1/2}(k^{-1}\delta_1(k)) )\\
- 2 \mathcal{D}_{1,2}(\Delta_{(1)}, \Delta_{(2)})(\Delta^{-1/2}( \delta_1(k)k^{-1})(k^{-1}\delta_1(k)) )\\
- 4 \mathcal{D}_{2,1}(\Delta_{(1)}, \Delta_{(2)})(\Delta^{-1}(\delta_1(k)k^{-1} )\Delta^{1/2}(k^{-1} \delta_1(k)) )\\
- 6 \mathcal{D}_{2,1}(\Delta_{(1)}, \Delta_{(2)})(\Delta^{-1}( \delta_1(k)k^{-1}) (k^{-1}\delta_1(k)) )\\
- 6 \mathcal{D}_{2,1}(\Delta_{(1)}, \Delta_{(2)})(\Delta^{-3/2}(\delta_1(k)k^{-1} ) \Delta^{1/2}(k^{-1} \delta_1(k)) )\\
- 8 \mathcal{D}_{2,1}(\Delta_{(1)}, \Delta_{(2)})(\Delta^{-3/2}( \delta_1(k)k^{-1}) (k^{-1}\delta_1(k)) )\\
+2 \mathcal{D}_{1,1}(\Delta_{(1)}, \Delta_{(2)})(\Delta^{-1/2}(\delta_1(k)k^{-1} ) \Delta^{1/2}(k^{-1} \delta_1(k)) )\\
+4 \mathcal{D}_{1,1}(\Delta_{(1)}, \Delta_{(2)})(\Delta^{-1/2}( \delta_1(k) k^{-1}) (k^{-1}\delta_1(k)) )\\
+|\tau|^2 \Big ( 2 \mathcal{D}_{2,2}(\Delta_{(1)}, \Delta_{(2)})( \Delta^{-1}(\delta_2(k)k^{-1}) \Delta^{1/2} (k^{-1}\delta_2(k)))\\
+ 2  \mathcal{D}_{2,2}(\Delta_{(1)}, \Delta_{(2)})( \Delta^{-1}( \delta_2(k) k^{-1}) (k^{-1}\delta_2(k)) )\\
+ 2  \mathcal{D}_{2,2}(\Delta_{(1)}, \Delta_{(2)})( \Delta^{-3/2}( \delta_2(k) k^{-1}) \Delta^{1/2} (k^{-1}\delta_2(k)) )\\
+ 2  \mathcal{D}_{2,2}(\Delta_{(1)}, \Delta_{(2)})( \Delta^{-3/2}( \delta_2(k) k^{-1}) (k^{-1}\delta_2(k)) )\\
+ 4 \mathcal{D}_{3,1}(\Delta_{(1)}, \Delta_{(2)})( \Delta^{-2}(\delta_2(k)k^{-1} ) \Delta^{1/2}(k^{-1} \delta_2(k)) )\\
+ 4 \mathcal{D}_{3,1}(\Delta_{(1)}, \Delta_{(2)})( \Delta^{-2}( \delta_2(k)k^{-1}) (k^{-1}\delta_2(k)) )\\
+ 4 \mathcal{D}_{3,1}(\Delta_{(1)}, \Delta_{(2)})(\Delta^{-5/2}(\delta_2(k)k^{-1} )\Delta^{1/2}(k^{-1} \delta_2(k)) )\\
+ 4 \mathcal{D}_{3,1}(\Delta_{(1)}, \Delta_{(2)})(\Delta^{-5/2}( \delta_2(k)k^{-1})(k^{-1} \delta_2(k)) )\\
- 2 \mathcal{D}_{1,2}(\Delta_{(1)}, \Delta_{(2)})(\Delta^{-1/2}( \delta_2(k)k^{-1})\Delta^{1/2}(k^{-1}\delta_2(k)) )\\
- 2 \mathcal{D}_{1,2}(\Delta_{(1)}, \Delta_{(2)})(\Delta^{-1/2}( \delta_2(k)k^{-1})(k^{-1}\delta_2(k)) )\\
- 4 \mathcal{D}_{2,1}(\Delta_{(1)}, \Delta_{(2)})(\Delta^{-1}(\delta_2(k)k^{-1} )\Delta^{1/2}(k^{-1} \delta_2(k)) )\\
- 6 \mathcal{D}_{2,1}(\Delta_{(1)}, \Delta_{(2)})(\Delta^{-1}( \delta_2(k)k^{-1}) (k^{-1}\delta_2(k)) )\\
- 6 \mathcal{D}_{2,1}(\Delta_{(1)}, \Delta_{(2)})(\Delta^{-3/2}(\delta_2(k)k^{-1} ) \Delta^{1/2}(k^{-1} \delta_2(k)) )\\
- 8 \mathcal{D}_{2,1}(\Delta_{(1)}, \Delta_{(2)})(\Delta^{-3/2}( \delta_2(k)k^{-1}) (k^{-1}\delta_2(k)) )\\
+2 \mathcal{D}_{1,1}(\Delta_{(1)}, \Delta_{(2)})(\Delta^{-1/2}(\delta_2(k)k^{-1} ) \Delta^{1/2}(k^{-1} \delta_2(k)) )\\
+4 \mathcal{D}_{1,1}(\Delta_{(1)}, \Delta_{(2)})(\Delta^{-1/2}( \delta_2(k) k^{-1}) (k^{-1}\delta_2(k)) ) \Big )\\
+\tau_1 \Big ( 2 \mathcal{D}_{2,2}(\Delta_{(1)}, \Delta_{(2)})( \Delta^{-1}(\delta_1(k)k^{-1}) \Delta^{1/2} (k^{-1}\delta_2(k)))\\
+ 2  \mathcal{D}_{2,2}(\Delta_{(1)}, \Delta_{(2)})( \Delta^{-1}( \delta_1(k) k^{-1}) (k^{-1}\delta_2(k)) )\\
+ 2  \mathcal{D}_{2,2}(\Delta_{(1)}, \Delta_{(2)})( \Delta^{-3/2}( \delta_1(k) k^{-1}) \Delta^{1/2} (k^{-1}\delta_2(k)) )\\
+ 2  \mathcal{D}_{2,2}(\Delta_{(1)}, \Delta_{(2)})( \Delta^{-3/2}( \delta_1(k) k^{-1}) (k^{-1}\delta_2(k)) )\\
+ 4 \mathcal{D}_{3,1}(\Delta_{(1)}, \Delta_{(2)})( \Delta^{-2}(\delta_1(k)k^{-1} ) \Delta^{1/2}(k^{-1} \delta_2(k)) )\\
+ 4 \mathcal{D}_{3,1}(\Delta_{(1)}, \Delta_{(2)})( \Delta^{-2}( \delta_1(k)k^{-1}) (k^{-1}\delta_2(k)) )\\
+ 4 \mathcal{D}_{3,1}(\Delta_{(1)}, \Delta_{(2)})(\Delta^{-5/2}(\delta_1(k)k^{-1} )\Delta^{1/2}(k^{-1} \delta_2(k)) )\\
+ 4 \mathcal{D}_{3,1}(\Delta_{(1)}, \Delta_{(2)})(\Delta^{-5/2}( \delta_1(k)k^{-1})(k^{-1} \delta_2(k)) )\\
- 2 \mathcal{D}_{1,2}(\Delta_{(1)}, \Delta_{(2)})(\Delta^{-1/2}( \delta_1(k)k^{-1})\Delta^{1/2}(k^{-1}\delta_2(k)) )\\
- 2 \mathcal{D}_{1,2}(\Delta_{(1)}, \Delta_{(2)})(\Delta^{-1/2}( \delta_1(k)k^{-1})(k^{-1}\delta_2(k)) )\\
- 4 \mathcal{D}_{2,1}(\Delta_{(1)}, \Delta_{(2)})(\Delta^{-1}(\delta_1(k)k^{-1} )\Delta^{1/2}(k^{-1} \delta_2(k)) )\\
- 6 \mathcal{D}_{2,1}(\Delta_{(1)}, \Delta_{(2)})(\Delta^{-1}( \delta_1(k)k^{-1}) (k^{-1}\delta_2(k)) )\\
- 6 \mathcal{D}_{2,1}(\Delta_{(1)}, \Delta_{(2)})(\Delta^{-3/2}(\delta_1(k)k^{-1} ) \Delta^{1/2}(k^{-1} \delta_2(k)) )\\
- 8 \mathcal{D}_{2,1}(\Delta_{(1)}, \Delta_{(2)})(\Delta^{-3/2}( \delta_1(k)k^{-1}) (k^{-1}\delta_2(k)) )\\
+2 \mathcal{D}_{1,1}(\Delta_{(1)}, \Delta_{(2)})(\Delta^{-1/2}(\delta_1(k)k^{-1} ) \Delta^{1/2}(k^{-1} \delta_2(k)) )\\
+4 \mathcal{D}_{1,1}(\Delta_{(1)}, \Delta_{(2)})(\Delta^{-1/2}( \delta_1(k) k^{-1}) (k^{-1}\delta_2(k)) )\\
+2 \mathcal{D}_{2,2}(\Delta_{(1)}, \Delta_{(2)})( \Delta^{-1}(\delta_2(k)k^{-1}) \Delta^{1/2} (k^{-1}\delta_1(k)))\\
+ 2  \mathcal{D}_{2,2}(\Delta_{(1)}, \Delta_{(2)})( \Delta^{-1}( \delta_2(k) k^{-1}) (k^{-1}\delta_1(k)) )\\
+ 2  \mathcal{D}_{2,2}(\Delta_{(1)}, \Delta_{(2)})( \Delta^{-3/2}( \delta_2(k) k^{-1}) \Delta^{1/2} (k^{-1}\delta_1(k)) )\\
+ 2  \mathcal{D}_{2,2}(\Delta_{(1)}, \Delta_{(2)})( \Delta^{-3/2}( \delta_2(k) k^{-1}) (k^{-1}\delta_1(k)) )\\
+ 4 \mathcal{D}_{3,1}(\Delta_{(1)}, \Delta_{(2)})( \Delta^{-2}(\delta_2(k)k^{-1} ) \Delta^{1/2}(k^{-1} \delta_1(k)) )\\
+ 4 \mathcal{D}_{3,1}(\Delta_{(1)}, \Delta_{(2)})( \Delta^{-2}( \delta_2(k)k^{-1}) (k^{-1}\delta_1(k)) )\\
+ 4 \mathcal{D}_{3,1}(\Delta_{(1)}, \Delta_{(2)})(\Delta^{-5/2}(\delta_2(k)k^{-1} )\Delta^{1/2}(k^{-1} \delta_1(k)) )\\
+ 4 \mathcal{D}_{3,1}(\Delta_{(1)}, \Delta_{(2)})(\Delta^{-5/2}( \delta_2(k)k^{-1})(k^{-1} \delta_1(k)) )\\
- 2 \mathcal{D}_{1,2}(\Delta_{(1)}, \Delta_{(2)})(\Delta^{-1/2}( \delta_2(k)k^{-1})\Delta^{1/2}(k^{-1}\delta_1(k)) )\\
- 2 \mathcal{D}_{1,2}(\Delta_{(1)}, \Delta_{(2)})(\Delta^{-1/2}( \delta_2(k)k^{-1})(k^{-1}\delta_1(k)) )\\
- 4 \mathcal{D}_{2,1}(\Delta_{(1)}, \Delta_{(2)})(\Delta^{-1}(\delta_2(k)k^{-1} )\Delta^{1/2}(k^{-1} \delta_1(k)) )\\
- 6 \mathcal{D}_{2,1}(\Delta_{(1)}, \Delta_{(2)})(\Delta^{-1}( \delta_2(k)k^{-1}) (k^{-1}\delta_1(k)) )\\
- 6 \mathcal{D}_{2,1}(\Delta_{(1)}, \Delta_{(2)})(\Delta^{-3/2}(\delta_2(k)k^{-1} ) \Delta^{1/2}(k^{-1} \delta_1(k)) )\\
- 8 \mathcal{D}_{2,1}(\Delta_{(1)}, \Delta_{(2)})(\Delta^{-3/2}( \delta_2(k)k^{-1}) (k^{-1}\delta_1(k)) )\\
+2 \mathcal{D}_{1,1}(\Delta_{(1)}, \Delta_{(2)})(\Delta^{-1/2}(\delta_2(k)k^{-1} ) \Delta^{1/2}(k^{-1} \delta_1(k)) )\\
+4 \mathcal{D}_{1,1}(\Delta_{(1)}, \Delta_{(2)})(\Delta^{-1/2}( \delta_2(k) k^{-1}) (k^{-1}\delta_1(k)) )  \Big ). \\
$

Putting together the latter terms with the ones obtained in \eqref{twob0int}, up to an overall factor of $\pi$, 
we find the following expression

\begin{eqnarray}\label{threeb0int}
&&f_1(\Delta)(k^{-1} \delta_1^2(k)) + f_2(\Delta)(k^{-2} \delta_1(k)^2) \nonumber \\
&+& F(\Delta_{(1)}, \Delta_{(2)})((\delta_1(k)k^{-1})(k^{-1}\delta_1(k))) \nonumber \\ 
&+&   |\tau|^2 f_1(\Delta)(k^{-1} \delta_2^2(k)) + |\tau|^2 f_2(\Delta)(k^{-2} \delta_2(k)^2) \nonumber \\
&+& |\tau|^2 F(\Delta_{(1)}, \Delta_{(2)})((\delta_2(k)k^{-1})(k^{-1}\delta_2(k))) \nonumber \\
&+& \tau_1 f_1(\Delta) (k^{-1} \delta_1\delta_2(k)) + \tau_1 f_2(\Delta)(k^{-2}\delta_1(k)\delta_2(k)) \nonumber \\
&+&  \tau_1 F(\Delta_{(1)}, \Delta_{(2)})((\delta_1(k)k^{-1})(k^{-1}\delta_2(k))) \nonumber \\
&+& \tau_1 f_1(\Delta) (k^{-1} \delta_2\delta_1(k)) + \tau_1 f_2(\Delta)(k^{-2}\delta_2(k)\delta_1(k)) \nonumber \\
&+& \tau_1 F(\Delta_{(1)}, \Delta_{(2)})((\delta_2(k)k^{-1})(k^{-1}\delta_1(k))), 
\end{eqnarray}
where as given by the formulas \eqref{fone} and \eqref{ftwo} we have
\begin{eqnarray}
f_1(u) = -\frac{u^{1/2} (2-2u+(1+u) \log{u})}{(-1+u^{1/2})^3(1+u^{1/2})^2}, \nonumber
\end{eqnarray}
\begin{equation}
f_2(u)= 2\frac{-1+u^2-2u \log{u}}{(-1+u)^3}, \nonumber
\end{equation}
and
\begin{eqnarray} \label{F(u,v)}
F(u,v) &:=& 2 \mathcal{D}_{2,2}(u, v)u^{-1}v^{1/2} 
+ 2  \mathcal{D}_{2,2}(u, v)u^{-1} 
+ 2  \mathcal{D}_{2,2}(u, v)u^{-3/2} v^{1/2} \nonumber \\
&&+ 2  \mathcal{D}_{2,2}(u,v) u^{-3/2}
+ 4 \mathcal{D}_{3,1}(u,v) u^{-2}v^{1/2}
+ 4 \mathcal{D}_{3,1}(u,v)u^{-2}  \nonumber \\
&&+ 4 \mathcal{D}_{3,1}(u,v) u^{-5/2} v^{1/2}
+ 4 \mathcal{D}_{3,1}(u,v)u^{-5/2}
- 2 \mathcal{D}_{1,2}(u,v)u^{-1/2} v^{1/2} \nonumber \\
&&- 2 \mathcal{D}_{1,2}(u,v)u^{-1/2}
- 4 \mathcal{D}_{2,1}(u,v)u^{-1}v^{1/2}
- 6 \mathcal{D}_{2,1}(u,v)u^{-1} \nonumber \\
&&- 6 \mathcal{D}_{2,1}(u,v) u^{-3/2} v^{1/2}
- 8 \mathcal{D}_{2,1}(u,v) u^{-3/2} 
+2 \mathcal{D}_{1,1}(u,v) u^{-1/2} v^{1/2} \nonumber \\
&&+ 4 \mathcal{D}_{1,1}(u,v) u^{-1/2} \nonumber
\end{eqnarray}
\begin{eqnarray}
&=& (2 u (-(((-1 + u v) (1 + \sqrt{u} (-1 - \sqrt{v} - (-2 + \sqrt{u} + u) v + 
          u v^{3/2})))/ \nonumber \\
          &&((-1 + \sqrt{u}) (-1 + \sqrt{v}))) + (
   \sqrt{u} \sqrt{
    v} (-1 - \sqrt{u} + u + u (-2 - \sqrt{u} + 2 u) \sqrt{v} \nonumber \\
    &&+ 
      u (-1 + \sqrt{u} + u) v + u^{5/2} v^{3/2}) \log{
     u})/((-1 + \sqrt{u})^2 (1 + \sqrt{u})) + (
   \sqrt{v} (1 -  \nonumber \\
   &&\sqrt{u} \sqrt{
       v} (-1 - \sqrt{v} + v + u v (-1 + \sqrt{v} + v) + 
         \sqrt{u} (-2 + \sqrt{v} + 2 v))) \log{
     v})/ \nonumber \\
     &&((-1 + \sqrt{v})^2 (1 + \sqrt{v}))))/
     (-1 + u v)^3. \nonumber
\end{eqnarray}

Note that in \eqref{threeb0int}, $\Delta_{(i)}$, i=1,2, signifies the action of $\Delta$ on the $i$-th factor of the product.

\subsection{The computations for $\partial^* k^2 \partial$.}

In order to compute the second component of the scalar curvature corresponding to $\partial^* k^2 \partial$, `\emph{the Laplacian on $(1,0)$-forms}',
 we first find the symbol 
of this operator:

\newtheorem{symbolsec}[main1]{Lemma}
\begin{symbolsec}
The symbol of $\partial^* k^2 \partial$ is equal to $c_2(\xi) + c_1(\xi)$ where 
\[c_2(\xi) =  \xi_1^2 k^2 + 2 \tau_1 \xi_1 \xi_2 k^2 + |\tau|^2 \xi_2^2 k^2, \]
\[c_1(\xi) =  (\delta_1(k^2)  + \tau \delta_2(k^2) )\xi_1 + (\bar \tau \delta_1(k^2)  + |\tau|^2 \delta_2(k^2) )\xi_2. \]
\begin{proof}
It follows easily from the composition rule explained in Proposition \ref{symbolcalculus} and the fact that the symbols of 
$\partial^*$, left multiplication by $k^2$, and $\partial$ are equal to $\xi_1+ \tau \xi_2$, $k^2$, and $\xi_1 + \bar \tau \xi_2.$ 
\end{proof}
\end{symbolsec}

Following the same method as in Subsection \ref{firstcomputation}, after a direct computation the corresponding $b_2$ term for 
the second half of the Laplacian $\triangle$, namely $\partial^* k^2 \partial$, is also given by a lengthy formula which is recorded in Appendix 
\ref{secondb2}. It is interesting to note that unlike the corresponding term for the first half, we have now terms with complex coefficient $i$ 
in front.

To integrate the second $b_2$  over the $\xi$-plane we use the change of variables \eqref{changevar}, namely we let
\begin{equation} \label{changevar}
\xi_1= r \cos \theta - r \frac{\tau_1}{\tau_2} \sin \theta,
\qquad
\xi_2=\frac{r}{\tau_2}\sin \theta , \nonumber
\end{equation}
where $\theta$ ranges from 0 to $2\pi$ and $r$ ranges from 0 to $\infty$.

After the integration with respect to $\theta$, up to a factor of $\frac{r}{\tau_2}$ which is the Jacobian of the change of variables, 
one gets \\

\noindent
$
4|\tau|^2\pi r^6b_0^2k^2\delta_2(k)b_0^2k^3\delta_2(k)b_0k+                                                        
4\tau_1\pi r^6b_0^2k^2\delta_2(k)b_0^2k^3\delta_1(k)b_0k+             \\                                         
4|\tau|^2\pi r^6b_0^2k^2\delta_2(k)b_0^2k^4\delta_2(k)b_0+                                                    
4\tau_1\pi r^6b_0^2k^2\delta_2(k)b_0^2k^4\delta_1(k)b_0+                 \\                               
4\tau_1\pi r^6b_0^2k^2\delta_1(k)b_0^2k^3\delta_2(k)b_0k+                                                 
4\pi r^6b_0^2k^2\delta_1(k)b_0^2k^3\delta_1(k)b_0k+                 \\                               
4\tau_1\pi r^6b_0^2k^2\delta_1(k)b_0^2k^4\delta_2(k)b_0+                                                
4\pi r^6b_0^2k^2\delta_1(k)b_0^2k^4\delta_1(k)b_0+                     \\                         
4|\tau|^2\pi r^6b_0^2k^3\delta_2(k)b_0^2k^2\delta_2(k)b_0k+                                           
4\tau_1\pi r^6b_0^2k^3\delta_2(k)b_0^2k^2\delta_1(k)b_0k+             \\                             
4|\tau|^2\pi r^6b_0^2k^3\delta_2(k)b_0^2k^3\delta_2(k)b_0+                                        
4\tau_1\pi r^6b_0^2k^3\delta_2(k)b_0^2k^3\delta_1(k)b_0+                  \\                                 
4\tau_1\pi r^6b_0^2k^3\delta_1(k)b_0^2k^2\delta_2(k)b_0k+                                                   
4\pi r^6b_0^2k^3\delta_1(k)b_0^2k^2\delta_1(k)b_0k+             \\                                      
4\tau_1\pi r^6b_0^2k^3\delta_1(k)b_0^2k^3\delta_2(k)b_0+                                                   
4\pi r^6b_0^2k^3\delta_1(k)b_0^2k^3\delta_1(k)b_0+               \\                                    
8|\tau|^2\pi r^6b_0^3k^4\delta_2(k)b_0k\delta_2(k)b_0k+                                                  
8\tau_1\pi r^6b_0^3k^4\delta_2(k)b_0k\delta_1(k)b_0k+              \\                                     
8|\tau|^2\pi r^6b_0^3k^4\delta_2(k)b_0k^2\delta_2(k)b_0+                                                  
8\tau_1\pi r^6b_0^3k^4\delta_2(k)b_0k^2\delta_1(k)b_0+             \\                                      
8\tau_1\pi r^6b_0^3k^4\delta_1(k)b_0k\delta_2(k)b_0k+                                            
8\pi r^6b_0^3k^4\delta_1(k)b_0k\delta_1(k)b_0k+                 \\                                  
8\tau_1\pi r^6b_0^3k^4\delta_1(k)b_0k^2\delta_2(k)b_0+                                                  
8\pi r^6b_0^3k^4\delta_1(k)b_0k^2\delta_1(k)b_0+            \\                                       
8|\tau|^2\pi r^6b_0^3k^5\delta_2(k)b_0\delta_2(k)b_0k+                                                  
8\tau_1\pi r^6b_0^3k^5\delta_2(k)b_0\delta_1(k)b_0k+           \\                                        
8|\tau|^2\pi r^6b_0^3k^5\delta_2(k)b_0k\delta_2(k)b_0+                                                 
8\tau_1\pi r^6b_0^3k^5\delta_2(k)b_0k\delta_1(k)b_0+           \\                                        
8\tau_1\pi r^6b_0^3k^5\delta_1(k)b_0\delta_2(k)b_0k+                                                   
8\pi r^6b_0^3k^5\delta_1(k)b_0\delta_1(k)b_0k+            \\                                       
8\tau_1\pi r^6b_0^3k^5\delta_1(k)b_0k\delta_2(k)b_0+                                                  
8\pi r^6b_0^3k^5\delta_1(k)b_0k\delta_1(k)b_0           \\                                        
-8|\tau|^2\pi r^4b_0^2k^2\delta_2(k)b_0k\delta_2(k)b_0k                                                   
-8\tau_1\pi r^4b_0^2k^2\delta_2(k)b_0k\delta_1(k)b_0k    \\                                              
-8|\tau|^2\pi r^4b_0^2k^2\delta_2(k)b_0k^2\delta_2(k)b_0                                                
-8\tau_1\pi r^4b_0^2k^2\delta_2(k)b_0k^2\delta_1(k)b_0     \\                                              
-8\tau_1\pi r^4b_0^2k^2\delta_1(k)b_0k\delta_2(k)b_0k                                                 
-8\pi r^4b_0^2k^2\delta_1(k)b_0k\delta_1(k)b_0k     \\                                             
-8\tau_1\pi r^4b_0^2k^2\delta_1(k)b_0k^2\delta_2(k)b_0                                                  
-8\pi r^4b_0^2k^2\delta_1(k)b_0k^2\delta_1(k)b_0   \\                                               
-8|\tau|^2\pi r^4b_0^2k^3\delta_2(k)b_0\delta_2(k)b_0k                                                 
-8\tau_1\pi r^4b_0^2k^3\delta_2(k)b_0\delta_1(k)b_0k   \\                                               
-8|\tau|^2\pi r^4b_0^2k^3\delta_2(k)b_0k\delta_2(k)b_0                                               
-8\tau_1\pi r^4b_0^2k^3\delta_2(k)b_0k\delta_1(k)b_0  \\                                                 
-8\tau_1\pi r^4b_0^2k^3\delta_1(k)b_0\delta_2(k)b_0k                                                 
-8\pi r^4b_0^2k^3\delta_1(k)b_0\delta_1(k)b_0k \\                                                  
-8\tau_1\pi r^4b_0^2k^3\delta_1(k)b_0k\delta_2(k)b_0                                                 
-8\pi r^4b_0^2k^3\delta_1(k)b_0k\delta_1(k)b_0 \\   
-4\pi r^4b_0^3k^4\delta_1^2(k)b_0k                                                   
-4|\tau|^2\pi r^4b_0^3k^4\delta_2^2(k)b_0k   \\                                                
-8\tau_1\pi r^4b_0^3k^4\delta_1\delta_2(k)b_0k                                                 
-8|\tau|^2\pi r^4b_0^3k^4\delta_2(k)^2b_0       \\                                            
-8\tau_1\pi r^4b_0^3k^4\delta_2(k)\delta_1(k)b_0                                                 
-8\tau_1\pi r^4b_0^3k^4\delta_1(k)\delta_2(k)b_0   \\                                               
-8\pi r^4b_0^3k^4\delta_1(k)^2b_0                                                
-4\pi r^4b_0^3k^5\delta_1^2(k)b_0         \\                                          
-4|\tau|^2\pi r^4b_0^3k^5\delta_2^2(k)b_0                                               
-8\tau_1\pi r^4b_0^3k^5\delta_1\delta_2(k)b_0     \\                                              
-2|\tau|^2\pi r^4b_0\delta_2(k^2)b_0^2k^2\delta_2(k)b_0k                                                 
-2(\tau_1 +i \tau_2)\pi r^4b_0\delta_2(k^2)b_0^2k^2\delta_1(k)b_0k  \\                                                
-2|\tau|^2\pi r^4b_0\delta_2(k^2)b_0^2k^3\delta_2(k)b_0                                                
-2(\tau_1 + i \tau_2)\pi r^4b_0\delta_2(k^2)b_0^2k^3\delta_1(k)b_0    \\                                               
-2(\tau_1 -i \tau_2)\pi r^4b_0\delta_1(k^2)b_0^2k^2\delta_2(k)b_0k                                                 
-2\pi r^4b_0\delta_1(k^2)b_0^2k^2\delta_1(k)b_0k    \\                                               
-2(\tau_1 -i \tau_2)\pi r^4b_0\delta_1(k^2)b_0^2k^3\delta_2(k)b_0                                                
-2\pi r^4b_0\delta_1(k^2)b_0^2k^3\delta_1(k)b_0   \\                                                 
-2|\tau|^2\pi r^4b_0^2k^2\delta_2(k^2)b_0\delta_2(k)b_0k                                                  
-2(\tau_1 +i \tau_2)\pi r^4b_0^2k^2\delta_2(k^2)b_0\delta_1(k)b_0k \\                                                   
-2|\tau|^2\pi r^4b_0^2k^2\delta_2(k^2)b_0k\delta_2(k)b_0                                                
-2(\tau_1 +i \tau_2)\pi r^4b_0^2k^2\delta_2(k^2)b_0k\delta_1(k)b_0    \\                                                   
-2(\tau_1 -i \tau_2)\pi r^4b_0^2k^2\delta_1(k^2)b_0\delta_2(k)b_0k                                                  
-2\pi r^4b_0^2k^2\delta_1(k^2)b_0\delta_1(k)b_0k        \\                                           
-2(\tau_1 -i \tau_2)\pi r^4b_0^2k^2\delta_1(k^2)b_0k\delta_2(k)b_0                                                  
-2\pi r^4b_0^2k^2\delta_1(k^2)b_0k\delta_1(k)b_0     \\                                             
-2|\tau|^2\pi r^4b_0^2k^2\delta_2(k)b_0k\delta_2(k^2)b_0                                            
-2(\tau_1 -i \tau_2)\pi r^4b_0^2k^2\delta_2(k)b_0k\delta_1(k^2)b_0     \\                                             
-2(\tau_1 +i \tau_2)\pi r^4b_0^2k^2\delta_1(k)b_0k\delta_2(k^2)b_0                                             
-2\pi r^4b_0^2k^2\delta_1(k)b_0k\delta_1(k^2)b_0     \\                                              
-2|\tau|^2\pi r^4b_0^2k^3\delta_2(k)b_0\delta_2(k^2)b_0                                         
-2(\tau_1 -i \tau_2)\pi r^4b_0^2k^3\delta_2(k)b_0\delta_1(k^2)b_0        \\                                           
-2(\tau_1 +i \tau_2)\pi r^4b_0^2k^3\delta_1(k)b_0\delta_2(k^2)b_0                                                
-2\pi r^4b_0^2k^3\delta_1(k)b_0\delta_1(k^2)b_0 +    \\    
2\pi r^2b_0^2k^2\delta_1^2(k)b_0k+                                                 
2|\tau|^2\pi r^2b_0^2k^2\delta_2^2(k)b_0k+    \\                                               
4\tau_1\pi r^2b_0^2k^2\delta_1\delta_2(k)b_0k+                                                 
4|\tau|^2\pi r^2b_0^2k^2\delta_2(k)^2b_0+      \\                                             
4\tau_1\pi r^2b_0^2k^2\delta_2(k)\delta_1(k)b_0+                                                    
4\tau_1\pi r^2b_0^2k^2\delta_1(k)\delta_2(k)b_0+      \\                                             
4\pi r^2b_0^2k^2\delta_1(k)^2b_0+                                             
2\pi r^2b_0^2k^3\delta_1^2(k)b_0+    \\                                               
2|\tau|^2\pi r^2b_0^2k^3\delta_2^2(k)b_0+                                                 
4\tau_1\pi r^2b_0^2k^3\delta_1\delta_2(k)b_0+    \\                                               
2|\tau|^2\pi r2b_0\delta_2(k^2)b_0\delta_2(k)b_0k+                                                 
2(\tau_1 +i \tau_2)\pi r^2b_0\delta_2(k^2)b_0\delta_1(k)b_0k+      \\                                             
2|\tau|^2\pi r^2b_0\delta_2(k^2)b_0k\delta_2(k)b_0+                                                 
2(\tau_1 +i \tau_2)\pi r^2b_0\delta_2(k^2)b_0k\delta_1(k)b_0+     \\                                              
2(\tau_1 -i \tau_2)\pi r^2b_0\delta_1(k^2)b_0\delta_2(k)b_0k+                                               
2\pi r^2b_0\delta_1(k^2)b_0\delta_1(k)b_0k+    \\                                               
2(\tau_1 -i \tau_2)\pi r^2b_0\delta_1(k^2)b_0k\delta_2(k)b_0+                                                  
2\pi r^2b_0\delta_1(k^2)b_0k\delta_1(k)b_0+   \\                                                
2\pi r^2b_0^2k^2\delta_1^2(k^2)b_0+                                          
2|\tau|^2\pi r^2b_0^2k^2\delta_2^2(k^2)b_0+    \\                                               
4\tau_1\pi r^2b_0^2k^2\delta_1\delta_2(k^2)b_0+                                               
0b_0\delta_2(k^2)b_0\delta_2(k^2)b_0+       \\                                            
0b_0\delta_2(k^2)b_0\delta_1(k^2)b_0+                                              
0b_0\delta_1(k^2)b_0\delta_2(k^2)b_0+           \\                                        
0b_0\delta_1(k^2)b_0\delta_1(k^2)b_0,
$ \\

\noindent
 where
\[b_0=(r^2k^2+1)^{-1}.\]

\subsubsection{Terms with two factors of the form $b_0^i$, $i \geq 1$.}

These are the following terms:\\

\noindent
$
-4\pi r^4b_0^3k^4\delta_1^2(k)b_0k                                                  
-4|\tau|^2\pi r^4b_0^3k^4\delta_2^2(k)b_0k                                                 
-8\tau_1\pi r^4b_0^3k^4\delta_1\delta_2(k)b_0k \\                                             
-8|\tau|^2\pi r^4b_0^3k^4\delta_2(k)^2b_0                                                  
-8\tau_1\pi r^4b_0^3k^4\delta_2(k)\delta_1(k)b_0                                                   
-8\tau_1\pi r^4b_0^3k^4\delta_1(k)\delta_2(k)b_0   \\                                               
-8\pi r^4b_0^3k^4\delta_1(k)^2b_0                                                  
-4\pi r^4b_0^3k^5\delta_1^2(k)b_0                                                  
-4|\tau|^2\pi r^4b_0^3k^5\delta_2^2(k)b_0 \\                                                 
-8\tau_1\pi r^4b_0^3k^5\delta_1\delta_2(k)b_0 +     
2\pi r^2b_0^2k^2\delta_1^2(k)b_0k+                                                 
2|\tau|^2\pi r^2b_0^2k^2\delta_2^2(k)b_0k+    \\                                               
4\tau_1\pi r^2b_0^2k^2\delta_1\delta_2(k)b_0k+                                                  
4|\tau|^2\pi r^2b_0^2k^2\delta_2(k)^2b_0+                                                  
4\tau_1\pi r^2b_0^2k^2\delta_2(k)\delta_1(k)b_0+ \\                                                    
4\tau_1\pi r^2b_0^2k^2\delta_1(k)\delta_2(k)b_0+                                                   
4\pi r^2b_0^2k^2\delta_1(k)^2b_0+                                                
2\pi r^2b_0^2k^3\delta_1^2(k)b_0+    \\                                               
2|\tau|^2\pi r^2b_0^2k^3\delta_2^2(k)b_0+                                                  
4\tau_1\pi r^2b_0^2k^3\delta_1\delta_2(k)b_0+      
2\pi r^2b_0^2k^2\delta_1^2(k^2)b_0+                \\                                 
2|\tau|^2\pi r^2b_0^2k^2\delta_2^2(k^2)b_0+                                                   
4\tau_1\pi r^2b_0^2k^2\delta_1\delta_2(k^2)b_0.$ \\

Using Lemma \ref{integration}, $\int_0^\infty \bullet \, r dr$ of these terms, up to an overall factor of $\pi$, is equal to \\

\noindent
$
-2 \mathcal{D}_2(\Delta^{1/2}(k^{-1}\delta_1^2(k))) 
-2 |\tau|^2 \mathcal{D}_2(\Delta^{1/2}(k^{-1}\delta_2^2(k)))
-4 \tau_1 \mathcal{D}_2(\Delta^{1/2}(k^{-1}\delta_1 \delta_2(k))) \\
-4 |\tau|^2 \mathcal{D}_2(k^{-2}\delta_2(k)^2)
-4 \tau_1 \mathcal{D}_2(k^{-2}\delta_2(k) \delta_1(k))
-4 \tau_1 \mathcal{D}_2(k^{-2}\delta_1(k) \delta_2(k))\\
-4 \mathcal{D}_2(k^{-2}\delta_1(k)^2) 
-2 \mathcal{D}_2(k^{-1}\delta_1^2(k)) 
-2 |\tau|^2 \mathcal{D}_2(k^{-1}\delta_2^2(k)) \\
-4 \tau_1 \mathcal{D}_2(k^{-1} \delta_1\delta_2(k))
+\mathcal{D}_1(\Delta^{1/2}(k^{-1}\delta_1^2(k)))
+|\tau|^2 \mathcal{D}_1(\Delta^{1/2}(k^{-1}\delta_2^2(k)))\\
+2 \tau_1 \mathcal{D}_1 ( \Delta^{1/2}(k^{-1} \delta_1\delta_2(k)))
+2|\tau|^2 \mathcal{D}_1(k^{-2}\delta_2(k)^2)
+2 \tau_1 \mathcal{D}_1 (k^{-2} \delta_2(k)\delta_1(k))\\
+2 \tau_1 \mathcal{D}_1 (k^{-2} \delta_1(k)\delta_2(k))
+2 \mathcal{D}_1 (k^{-2} \delta_1(k)^2)
+ \mathcal{D}_1 (k^{-1} \delta_1^2(k))\\
+ |\tau|^2 \mathcal{D}_1 (k^{-1} \delta_2^2(k))
+2 \tau_1 \mathcal{D}_1 (k^{-1} \delta_1\delta_2(k))
+  \mathcal{D}_1 (k^{-2} \delta_1^2(k^2))\\
+  |\tau|^2 \mathcal{D}_1 (k^{-2} \delta_2^2(k^2))
+2 \tau_1 \mathcal{D}_1 (k^{-2} \delta_1\delta_2(k^2)).
$ \\

Therefore, up to an overall factor of $\pi$, the $\int_0^\infty \bullet \,r dr$ of the terms with two factors of powers of $b_0$ is 
equal to 

\begin{eqnarray} \label{twob0inttwo}
&&g_1(\Delta)(k^{-1} \delta_1^2(k)) + g_2(\Delta)(k^{-2} \delta_1(k)^2) \nonumber \\
&+& |\tau|^2 g_1(\Delta)(k^{-1} \delta_2^2(k)) +|\tau|^2 g_2(\Delta)(k^{-2} \delta_2(k)^2) \nonumber \\
&+& \tau_1 g_1(\Delta)(k^{-1} \delta_1\delta_2(k)) + \tau_1 g_2(\Delta)(k^{-2} \delta_1(k)\delta_2(k)) \nonumber \\
&+& \tau_1 g_1(\Delta)(k^{-1} \delta_2\delta_1(k)) + \tau_1 g_2(\Delta)(k^{-2} \delta_2(k)\delta_1(k)),
\end{eqnarray} 
where
\begin{eqnarray} \label{gone}
g_1(u)&:=& -2 \mathcal{L}_2(u) u^{1/2} -2 \mathcal{L}_2(u) + 2 \mathcal{L}_1(u) u^{1/2} +2 \mathcal{L}_1(u) \nonumber \\
&=& \frac{-1+u^2-2u \log u}{(-1+u^{1/2})^3(1+u^{1/2})^2},
\end{eqnarray}
and

\begin{eqnarray} \label{gtwo}
g_2(u)&:=& -4 \mathcal{L}_2(u) +4 \mathcal{L}_1(u) = 2\frac{-1+u^2-2u \log u}{(-1+u)^3}.
\end{eqnarray}

\subsubsection{Terms with three factors of the form $b_0^i$, $i \geq 1$.}

These are the following terms:\\

\noindent
$
4|\tau|^2\pi r^6b_0^2k^2\delta_2(k)b_0^2k^3\delta_2(k)b_0k+                                                       
4\tau_1\pi r^6b_0^2k^2\delta_2(k)b_0^2k^3\delta_1(k)b_0k+             \\                                         
4|\tau|^2\pi r^6b_0^2k^2\delta_2(k)b_0^2k^4\delta_2(k)b_0+                                                    
4\tau_1\pi r^6b_0^2k^2\delta_2(k)b_0^2k^4\delta_1(k)b_0+                 \\                               
4\tau_1\pi r^6b_0^2k^2\delta_1(k)b_0^2k^3\delta_2(k)b_0k+                                                  
4\pi r^6b_0^2k^2\delta_1(k)b_0^2k^3\delta_1(k)b_0k+                 \\                               
4\tau_1\pi r^6b_0^2k^2\delta_1(k)b_0^2k^4\delta_2(k)b_0+                                                
4\pi r^6b_0^2k^2\delta_1(k)b_0^2k^4\delta_1(k)b_0+                     \\                         
4|\tau|^2\pi r^6b_0^2k^3\delta_2(k)b_0^2k^2\delta_2(k)b_0k+                                           
4\tau_1\pi r^6b_0^2k^3\delta_2(k)b_0^2k^2\delta_1(k)b_0k+             \\                             
4|\tau|^2\pi r^6b_0^2k^3\delta_2(k)b_0^2k^3\delta_2(k)b_0+                                         
4\tau_1\pi r^6b_0^2k^3\delta_2(k)b_0^2k^3\delta_1(k)b_0+                  \\                                 
4\tau_1\pi r^6b_0^2k^3\delta_1(k)b_0^2k^2\delta_2(k)b_0k+                                                   
4\pi r^6b_0^2k^3\delta_1(k)b_0^2k^2\delta_1(k)b_0k+             \\                                      
4\tau_1\pi r^6b_0^2k^3\delta_1(k)b_0^2k^3\delta_2(k)b_0+                                                   
4\pi r^6b_0^2k^3\delta_1(k)b_0^2k^3\delta_1(k)b_0+               \\                                    
8|\tau|^2\pi r^6b_0^3k^4\delta_2(k)b_0k\delta_2(k)b_0k+                                                  
8\tau_1\pi r^6b_0^3k^4\delta_2(k)b_0k\delta_1(k)b_0k+              \\                                     
8|\tau|^2\pi r^6b_0^3k^4\delta_2(k)b_0k^2\delta_2(k)b_0+                                                  
8\tau_1\pi r^6b_0^3k^4\delta_2(k)b_0k^2\delta_1(k)b_0+             \\                                      
8\tau_1\pi r^6b_0^3k^4\delta_1(k)b_0k\delta_2(k)b_0k+                                                   
8\pi r^6b_0^3k^4\delta_1(k)b_0k\delta_1(k)b_0k+                 \\                                  
8\tau_1\pi r^6b_0^3k^4\delta_1(k)b_0k^2\delta_2(k)b_0+                                                   
8\pi r^6b_0^3k^4\delta_1(k)b_0k^2\delta_1(k)b_0+            \\                                       
8|\tau|^2\pi r^6b_0^3k^5\delta_2(k)b_0\delta_2(k)b_0k+                                                   
8\tau_1\pi r^6b_0^3k^5\delta_2(k)b_0\delta_1(k)b_0k+           \\                                        
8|\tau|^2\pi r^6b_0^3k^5\delta_2(k)b_0k\delta_2(k)b_0+                                                   
8\tau_1\pi r^6b_0^3k^5\delta_2(k)b_0k\delta_1(k)b_0+           \\                                        
8\tau_1\pi r^6b_0^3k^5\delta_1(k)b_0\delta_2(k)b_0k+                                                   
8\pi r^6b_0^3k^5\delta_1(k)b_0\delta_1(k)b_0k+            \\                                       
8\tau_1\pi r^6b_0^3k^5\delta_1(k)b_0k\delta_2(k)b_0+                                               
8\pi r^6b_0^3k^5\delta_1(k)b_0k\delta_1(k)b_0           \\                                       
-8|\tau|^2\pi r^4b_0^2k^2\delta_2(k)b_0k\delta_2(k)b_0k                                                  
-8\tau_1\pi r^4b_0^2k^2\delta_2(k)b_0k\delta_1(k)b_0k    \\                                              
-8|\tau|^2\pi r^4b_0^2k^2\delta_2(k)b_0k^2\delta_2(k)b_0                                                 
-8\tau_1\pi r^4b_0^2k^2\delta_2(k)b_0k^2\delta_1(k)b_0     \\                                              
-8\tau_1\pi r^4b_0^2k^2\delta_1(k)b_0k\delta_2(k)b_0k                                                 
-8\pi r^4b_0^2k^2\delta_1(k)b_0k\delta_1(k)b_0k     \\                                             
-8\tau_1\pi r^4b_0^2k^2\delta_1(k)b_0k^2\delta_2(k)b_0                                                
-8\pi r^4b_0^2k^2\delta_1(k)b_0k^2\delta_1(k)b_0   \\                                               
-8|\tau|^2\pi r^4b_0^2k^3\delta_2(k)b_0\delta_2(k)b_0k                                             
-8\tau_1\pi r^4b_0^2k^3\delta_2(k)b_0\delta_1(k)b_0k   \\                                               
-8|\tau|^2\pi r^4b_0^2k^3\delta_2(k)b_0k\delta_2(k)b_0                                               
-8\tau_1\pi r^4b_0^2k^3\delta_2(k)b_0k\delta_1(k)b_0  \\                                                 
-8\tau_1\pi r^4b_0^2k^3\delta_1(k)b_0\delta_2(k)b_0k                                             
-8\pi r^4b_0^2k^3\delta_1(k)b_0\delta_1(k)b_0k \\                                                  
-8\tau_1\pi r^4b_0^2k^3\delta_1(k)b_0k\delta_2(k)b_0                                                
-8\pi r^4b_0^2k^3\delta_1(k)b_0k\delta_1(k)b_0 \\  
-2|\tau|^2\pi r^4b_0\delta_2(k^2)b_0^2k^2\delta_2(k)b_0k                                           
-2(\tau_1 + i \tau_2)\pi r^4b_0\delta_2(k^2)b_0^2k^2\delta_1(k)b_0k  \\                                                
-2|\tau|^2\pi r^4b_0\delta_2(k^2)b_0^2k^3\delta_2(k)b_0                                                
-2(\tau_1 +i \tau_2)\pi r^4b_0\delta_2(k^2)b_0^2k^3\delta_1(k)b_0    \\                                               
-2(\tau_1 -i \tau_2)\pi r^4b_0\delta_1(k^2)b_0^2k^2\delta_2(k)b_0k                                                 
-2\pi r^4b_0\delta_1(k^2)b_0^2k^2\delta_1(k)b_0k    \\                                               
-2(\tau_1 -i \tau_2)\pi r^4b_0\delta_1(k^2)b_0^2k^3\delta_2(k)b_0                                                 
-2\pi r^4b_0\delta_1(k^2)b_0^2k^3\delta_1(k)b_0   \\                                                 
-2|\tau|^2\pi r^4b_0^2k^2\delta_2(k^2)b_0\delta_2(k)b_0k                                                 
-2(\tau_1 +i \tau_2)\pi r^4b_0^2k^2\delta_2(k^2)b_0\delta_1(k)b_0k \\                                                   
-2|\tau|^2\pi r^4b_0^2k^2\delta_2(k^2)b_0k\delta_2(k)b_0                                                
-2(\tau_1 +i \tau_2)\pi r^4b_0^2k^2\delta_2(k^2)b_0k\delta_1(k)b_0    \\                                                   
-2(\tau_1 -i \tau_2)\pi r^4b_0^2k^2\delta_1(k^2)b_0\delta_2(k)b_0k                                                   
-2\pi r^4b_0^2k^2\delta_1(k^2)b_0\delta_1(k)b_0k        \\                                           
-2(\tau_1 -i \tau_2)\pi r^4b_0^2k^2\delta_1(k^2)b_0k\delta_2(k)b_0                                                  
-2\pi r^4b_0^2k^2\delta_1(k^2)b_0k\delta_1(k)b_0     \\                                             
-2|\tau|^2\pi r^4b_0^2k^2\delta_2(k)b_0k\delta_2(k^2)b_0                                                 
-2(\tau_1 -i \tau_2)\pi r^4b_0^2k^2\delta_2(k)b_0k\delta_1(k^2)b_0     \\                                             
-2(\tau_1 +i \tau_2)\pi r^4b_0^2k^2\delta_1(k)b_0k\delta_2(k^2)b_0                                                 
-2\pi r^4b_0^2k^2\delta_1(k)b_0k\delta_1(k^2)b_0     \\                                              
-2|\tau|^2\pi r^4b_0^2k^3\delta_2(k)b_0\delta_2(k^2)b_0                                                 
-2(\tau_1 -i \tau_2)\pi r^4b_0^2k^3\delta_2(k)b_0\delta_1(k^2)b_0        \\                                           
-2(\tau_1 +i \tau_2)\pi r^4b_0^2k^3\delta_1(k)b_0\delta_2(k^2)b_0                                                 
-2\pi r^4b_0^2k^3\delta_1(k)b_0\delta_1(k^2)b_0 +    \\    
2|\tau|^2\pi r2b_0\delta_2(k^2)b_0\delta_2(k)b_0k+                                                 
2(\tau_1 +i \tau_2)\pi r^2b_0\delta_2(k^2)b_0\delta_1(k)b_0k+      \\                                             
2|\tau|^2\pi r^2b_0\delta_2(k^2)b_0k\delta_2(k)b_0+                                                 
2(\tau_1 +i \tau_2)\pi r^2b_0\delta_2(k^2)b_0k\delta_1(k)b_0+     \\                                              
2(\tau_1 -i \tau_2)\pi r^2b_0\delta_1(k^2)b_0\delta_2(k)b_0k+                                                   
2\pi r^2b_0\delta_1(k^2)b_0\delta_1(k)b_0k+    \\                                               
2(\tau_1 -i \tau_2)\pi r^2b_0\delta_1(k^2)b_0k\delta_2(k)b_0+                                                 
2\pi r^2b_0\delta_1(k^2)b_0k\delta_1(k)b_0  +\\         
0b_0\delta_2(k^2)b_0\delta_2(k^2)b_0+                                                   
0b_0\delta_2(k^2)b_0\delta_1(k^2)b_0+         \\                                          
0b_0\delta_1(k^2)b_0\delta_2(k^2)b_0+                                                   
0b_0\delta_1(k^2)b_0\delta_1(k^2)b_0.
$\\

Using Lemma \ref{integrationtwo} we compute $\int_0^\infty \bullet \, r  dr$ of these terms, and the result, 
up to an overall factor of $\pi$, is equal to:\\

\noindent  
$
2 \mathcal{D}_{2,2}(\Delta_{(1)}, \Delta_{(2)}) ( \Delta^{-1}(\delta_1(k)k^{-1}) \Delta^{1/2}(k^{-1}\delta_1(k))) \\
+ 2 \mathcal{D}_{2,2}(\Delta_{(1)}, \Delta_{(2)}) ( \Delta^{-1}(\delta_1(k)k^{-1}) (k^{-1}\delta_1(k))) \\
+ 2 \mathcal{D}_{2,2}(\Delta_{(1)}, \Delta_{(2)}) ( \Delta^{-3/2}(\delta_1(k)k^{-1}) \Delta^{1/2}(k^{-1}\delta_1(k))) \\
+ 2 \mathcal{D}_{2,2}(\Delta_{(1)}, \Delta_{(2)}) ( \Delta^{-3/2}(\delta_1(k)k^{-1}) (k^{-1}\delta_1(k))) \\
+ 4 \mathcal{D}_{3,1}(\Delta_{(1)}, \Delta_{(2)}) ( \Delta^{-2}(\delta_1(k)k^{-1}) \Delta^{1/2}(k^{-1}\delta_1(k))) \\
+ 4 \mathcal{D}_{3,1}(\Delta_{(1)}, \Delta_{(2)}) ( \Delta^{-2}(\delta_1(k)k^{-1}) (k^{-1}\delta_1(k))) \\
+ 4 \mathcal{D}_{3,1}(\Delta_{(1)}, \Delta_{(2)}) ( \Delta^{-5/2}(\delta_1(k)k^{-1}) \Delta^{1/2}(k^{-1}\delta_1(k))) \\
+ 4 \mathcal{D}_{3,1}(\Delta_{(1)}, \Delta_{(2)}) ( \Delta^{-5/2}(\delta_1(k)k^{-1}) (k^{-1}\delta_1(k))) \\
- 4 \mathcal{D}_{2,1}(\Delta_{(1)}, \Delta_{(2)}) ( \Delta^{-1}(\delta_1(k)k^{-1}) \Delta^{1/2} (k^{-1}\delta_1(k))) \\
- 4 \mathcal{D}_{2,1}(\Delta_{(1)}, \Delta_{(2)}) ( \Delta^{-1}(\delta_1(k)k^{-1})  (k^{-1}\delta_1(k))) \\
- 4 \mathcal{D}_{2,1}(\Delta_{(1)}, \Delta_{(2)}) ( \Delta^{-3/2}(\delta_1(k)k^{-1})  \Delta^{1/2}(k^{-1}\delta_1(k))) \\
- 4 \mathcal{D}_{2,1}(\Delta_{(1)}, \Delta_{(2)}) ( \Delta^{-3/2}(\delta_1(k)k^{-1})  (k^{-1}\delta_1(k))) \\
-  \mathcal{D}_{1,2}(\Delta_{(1)}, \Delta_{(2)}) ( (\delta_1(k)k^{-1}) \Delta^{1/2} (k^{-1}\delta_1(k))) \\
-  \mathcal{D}_{1,2}(\Delta_{(1)}, \Delta_{(2)}) (\Delta^{-1/2} (\delta_1(k)k^{-1}) \Delta^{1/2} (k^{-1}\delta_1(k))) \\
-  \mathcal{D}_{1,2}(\Delta_{(1)}, \Delta_{(2)}) ( \Delta^{-1/2}(\delta_1(k)k^{-1}) (k^{-1}\delta_1(k))) \\
-  \mathcal{D}_{1,2}(\Delta_{(1)}, \Delta_{(2)}) ( (\delta_1(k)k^{-1}) (k^{-1}\delta_1(k))) \\
-  \mathcal{D}_{2,1}(\Delta_{(1)}, \Delta_{(2)}) ( \Delta^{-3/2}(\delta_1(k)k^{-1}) \Delta^{1/2}(k^{-1}\delta_1(k))) \\
-  \mathcal{D}_{2,1}(\Delta_{(1)}, \Delta_{(2)}) ( \Delta^{-1}(\delta_1(k)k^{-1}) \Delta^{1/2}(k^{-1}\delta_1(k))) \\
-  \mathcal{D}_{2,1}(\Delta_{(1)}, \Delta_{(2)}) ( \Delta^{-3/2}(\delta_1(k)k^{-1}) (k^{-1}\delta_1(k))) \\
-  \mathcal{D}_{2,1}(\Delta_{(1)}, \Delta_{(2)}) ( \Delta^{-1}(\delta_1(k)k^{-1}) (k^{-1}\delta_1(k))) \\
-  \mathcal{D}_{2,1}(\Delta_{(1)}, \Delta_{(2)}) ( \Delta^{-1}(\delta_1(k)k^{-1}) (k^{-1}\delta_1(k))) \\
-  \mathcal{D}_{2,1}(\Delta_{(1)}, \Delta_{(2)}) ( \Delta^{-1}(\delta_1(k)k^{-1}) \Delta^{1/2}(k^{-1}\delta_1(k))) \\
-  \mathcal{D}_{2,1}(\Delta_{(1)}, \Delta_{(2)}) ( \Delta^{-3/2}(\delta_1(k)k^{-1}) (k^{-1}\delta_1(k))) \\
-  \mathcal{D}_{2,1}(\Delta_{(1)}, \Delta_{(2)}) ( \Delta^{-3/2}(\delta_1(k)k^{-1}) \Delta^{1/2}(k^{-1}\delta_1(k))) \\
+ \mathcal{D}_{1,1}(\Delta_{(1)}, \Delta_{(2)}) ( \Delta^{-1/2}(\delta_1(k)k^{-1}) \Delta^{1/2}(k^{-1}\delta_1(k))) \\
+ \mathcal{D}_{1,1}(\Delta_{(1)}, \Delta_{(2)}) ( (\delta_1(k)k^{-1}) \Delta^{1/2}(k^{-1}\delta_1(k))) \\
+ \mathcal{D}_{1,1}(\Delta_{(1)}, \Delta_{(2)}) ( \Delta^{-1/2}(\delta_1(k)k^{-1}) (k^{-1}\delta_1(k))) \\
+ \mathcal{D}_{1,1}(\Delta_{(1)}, \Delta_{(2)}) ( (\delta_1(k)k^{-1}) (k^{-1}\delta_1(k))) \\
+|\tau|^2 \Big (   
2 \mathcal{D}_{2,2}(\Delta_{(1)}, \Delta_{(2)}) ( \Delta^{-1}(\delta_2(k)k^{-1}) \Delta^{1/2}(k^{-1}\delta_2(k))) \\
+2 \mathcal{D}_{2,2}(\Delta_{(1)}, \Delta_{(2)}) ( \Delta^{-1}(\delta_2(k)k^{-1}) (k^{-1}\delta_2(k))) \\
+2 \mathcal{D}_{2,2}(\Delta_{(1)}, \Delta_{(2)}) ( \Delta^{-3/2}(\delta_2(k)k^{-1}) \Delta^{1/2}(k^{-1}\delta_2(k))) \\
+ 2 \mathcal{D}_{2,2}(\Delta_{(1)}, \Delta_{(2)}) ( \Delta^{-3/2}(\delta_2(k)k^{-1}) (k^{-1}\delta_2(k))) \\
+ 4 \mathcal{D}_{3,1}(\Delta_{(1)}, \Delta_{(2)}) ( \Delta^{-2}(\delta_2(k)k^{-1}) \Delta^{1/2}(k^{-1}\delta_2(k))) \\
+ 4 \mathcal{D}_{3,1}(\Delta_{(1)}, \Delta_{(2)}) ( \Delta^{-2}(\delta_2(k)k^{-1}) (k^{-1}\delta_2(k))) \\
+ 4 \mathcal{D}_{3,1}(\Delta_{(1)}, \Delta_{(2)}) ( \Delta^{-5/2}(\delta_2(k)k^{-1}) \Delta^{1/2}(k^{-1}\delta_2(k))) \\
+ 4 \mathcal{D}_{3,1}(\Delta_{(1)}, \Delta_{(2)}) ( \Delta^{-5/2}(\delta_2(k)k^{-1}) (k^{-1}\delta_2(k))) \\
- 4 \mathcal{D}_{2,1}(\Delta_{(1)}, \Delta_{(2)}) ( \Delta^{-1}(\delta_2(k)k^{-1}) \Delta^{1/2} (k^{-1}\delta_2(k))) \\
- 4 \mathcal{D}_{2,1}(\Delta_{(1)}, \Delta_{(2)}) ( \Delta^{-1}(\delta_2(k)k^{-1})  (k^{-1}\delta_2(k))) \\
- 4 \mathcal{D}_{2,1}(\Delta_{(1)}, \Delta_{(2)}) ( \Delta^{-3/2}(\delta_2(k)k^{-1})  \Delta^{1/2}(k^{-1}\delta_2(k))) \\
- 4 \mathcal{D}_{2,1}(\Delta_{(1)}, \Delta_{(2)}) ( \Delta^{-3/2}(\delta_2(k)k^{-1})  (k^{-1}\delta_2(k))) \\
-  \mathcal{D}_{1,2}(\Delta_{(1)}, \Delta_{(2)}) ( (\delta_2(k)k^{-1}) \Delta^{1/2} (k^{-1}\delta_2(k))) \\
-  \mathcal{D}_{1,2}(\Delta_{(1)}, \Delta_{(2)}) (\Delta^{-1/2} (\delta_2(k)k^{-1}) \Delta^{1/2} (k^{-1}\delta_2(k))) \\
-  \mathcal{D}_{1,2}(\Delta_{(1)}, \Delta_{(2)}) ( \Delta^{-1/2}(\delta_2(k)k^{-1}) (k^{-1}\delta_2(k))) \\
-  \mathcal{D}_{1,2}(\Delta_{(1)}, \Delta_{(2)}) ( (\delta_2(k)k^{-1}) (k^{-1}\delta_2(k))) \\
-  \mathcal{D}_{2,1}(\Delta_{(1)}, \Delta_{(2)}) ( \Delta^{-3/2}(\delta_2(k)k^{-1}) \Delta^{1/2}(k^{-1}\delta_2(k))) \\
-  \mathcal{D}_{2,1}(\Delta_{(1)}, \Delta_{(2)}) ( \Delta^{-1}(\delta_2(k)k^{-1}) \Delta^{1/2}(k^{-1}\delta_2(k))) \\
-  \mathcal{D}_{2,1}(\Delta_{(1)}, \Delta_{(2)}) ( \Delta^{-3/2}(\delta_2(k)k^{-1}) (k^{-1}\delta_2(k))) \\
-  \mathcal{D}_{2,1}(\Delta_{(1)}, \Delta_{(2)}) ( \Delta^{-1}(\delta_2(k)k^{-1}) (k^{-1}\delta_2(k))) \\
-  \mathcal{D}_{2,1}(\Delta_{(1)}, \Delta_{(2)}) ( \Delta^{-1}(\delta_2(k)k^{-1}) (k^{-1}\delta_2(k))) \\
-  \mathcal{D}_{2,1}(\Delta_{(1)}, \Delta_{(2)}) ( \Delta^{-1}(\delta_2(k)k^{-1}) \Delta^{1/2}(k^{-1}\delta_2(k))) \\
-  \mathcal{D}_{2,1}(\Delta_{(1)}, \Delta_{(2)}) ( \Delta^{-3/2}(\delta_2(k)k^{-1}) (k^{-1}\delta_2(k))) \\
-  \mathcal{D}_{2,1}(\Delta_{(1)}, \Delta_{(2)}) ( \Delta^{-3/2}(\delta_2(k)k^{-1}) \Delta^{1/2}(k^{-1}\delta_2(k))) \\
+ \mathcal{D}_{1,1}(\Delta_{(1)}, \Delta_{(2)}) ( \Delta^{-1/2}(\delta_2(k)k^{-1}) \Delta^{1/2}(k^{-1}\delta_2(k))) \\
+ \mathcal{D}_{1,1}(\Delta_{(1)}, \Delta_{(2)}) ( (\delta_2(k)k^{-1}) \Delta^{1/2}(k^{-1}\delta_2(k))) \\
+ \mathcal{D}_{1,1}(\Delta_{(1)}, \Delta_{(2)}) ( \Delta^{-1/2}(\delta_2(k)k^{-1}) (k^{-1}\delta_2(k))) \\
+ \mathcal{D}_{1,1}(\Delta_{(1)}, \Delta_{(2)}) ( (\delta_2(k)k^{-1}) (k^{-1}\delta_2(k))) \Big ) \\
+ \tau_1 \Big (  
2 \mathcal{D}_{2,2}(\Delta_{(1)}, \Delta_{(2)}) ( \Delta^{-1}(\delta_1(k)k^{-1}) \Delta^{1/2}(k^{-1}\delta_2(k))) \\
+ 2 \mathcal{D}_{2,2}(\Delta_{(1)}, \Delta_{(2)}) ( \Delta^{-1}(\delta_1(k)k^{-1}) (k^{-1}\delta_2(k))) \\
+ 2 \mathcal{D}_{2,2}(\Delta_{(1)}, \Delta_{(2)}) ( \Delta^{-3/2}(\delta_1(k)k^{-1}) \Delta^{1/2}(k^{-1}\delta_2(k))) \\
+ 2 \mathcal{D}_{2,2}(\Delta_{(1)}, \Delta_{(2)}) ( \Delta^{-3/2}(\delta_1(k)k^{-1}) (k^{-1}\delta_2(k))) \\
+ 4 \mathcal{D}_{3,1}(\Delta_{(1)}, \Delta_{(2)}) ( \Delta^{-2}(\delta_1(k)k^{-1}) \Delta^{1/2}(k^{-1}\delta_2(k))) \\
+ 4 \mathcal{D}_{3,1}(\Delta_{(1)}, \Delta_{(2)}) ( \Delta^{-2}(\delta_1(k)k^{-1}) (k^{-1}\delta_2(k))) \\
+ 4 \mathcal{D}_{3,1}(\Delta_{(1)}, \Delta_{(2)}) ( \Delta^{-5/2}(\delta_1(k)k^{-1}) \Delta^{1/2}(k^{-1}\delta_2(k))) \\
+ 4 \mathcal{D}_{3,1}(\Delta_{(1)}, \Delta_{(2)}) ( \Delta^{-5/2}(\delta_1(k)k^{-1}) (k^{-1}\delta_2(k))) \\
- 4 \mathcal{D}_{2,1}(\Delta_{(1)}, \Delta_{(2)}) ( \Delta^{-1}(\delta_1(k)k^{-1}) \Delta^{1/2} (k^{-1}\delta_2(k))) \\
- 4 \mathcal{D}_{2,1}(\Delta_{(1)}, \Delta_{(2)}) ( \Delta^{-1}(\delta_1(k)k^{-1})  (k^{-1}\delta_2(k))) \\
- 4 \mathcal{D}_{2,1}(\Delta_{(1)}, \Delta_{(2)}) ( \Delta^{-3/2}(\delta_1(k)k^{-1})  \Delta^{1/2}(k^{-1}\delta_2(k))) \\
- 4 \mathcal{D}_{2,1}(\Delta_{(1)}, \Delta_{(2)}) ( \Delta^{-3/2}(\delta_1(k)k^{-1})  (k^{-1}\delta_2(k))) \\
-  \mathcal{D}_{1,2}(\Delta_{(1)}, \Delta_{(2)}) ( (\delta_1(k)k^{-1}) \Delta^{1/2} (k^{-1}\delta_2(k))) \\
-  \mathcal{D}_{1,2}(\Delta_{(1)}, \Delta_{(2)}) (\Delta^{-1/2} (\delta_1(k)k^{-1}) \Delta^{1/2} (k^{-1}\delta_2(k))) \\
-  \mathcal{D}_{1,2}(\Delta_{(1)}, \Delta_{(2)}) ( \Delta^{-1/2}(\delta_1(k)k^{-1}) (k^{-1}\delta_2(k))) \\
-  \mathcal{D}_{1,2}(\Delta_{(1)}, \Delta_{(2)}) ( (\delta_1(k)k^{-1}) (k^{-1}\delta_2(k))) \\
-  \mathcal{D}_{2,1}(\Delta_{(1)}, \Delta_{(2)}) ( \Delta^{-3/2}(\delta_1(k)k^{-1}) \Delta^{1/2}(k^{-1}\delta_2(k))) \\
-  \mathcal{D}_{2,1}(\Delta_{(1)}, \Delta_{(2)}) ( \Delta^{-1}(\delta_1(k)k^{-1}) \Delta^{1/2}(k^{-1}\delta_2(k))) \\
-  \mathcal{D}_{2,1}(\Delta_{(1)}, \Delta_{(2)}) ( \Delta^{-3/2}(\delta_1(k)k^{-1}) (k^{-1}\delta_2(k))) \\
-  \mathcal{D}_{2,1}(\Delta_{(1)}, \Delta_{(2)}) ( \Delta^{-1}(\delta_1(k)k^{-1}) (k^{-1}\delta_2(k))) \\
-  \mathcal{D}_{2,1}(\Delta_{(1)}, \Delta_{(2)}) ( \Delta^{-1}(\delta_1(k)k^{-1}) (k^{-1}\delta_2(k))) \\
-  \mathcal{D}_{2,1}(\Delta_{(1)}, \Delta_{(2)}) ( \Delta^{-1}(\delta_1(k)k^{-1}) \Delta^{1/2}(k^{-1}\delta_2(k))) \\
-  \mathcal{D}_{2,1}(\Delta_{(1)}, \Delta_{(2)}) ( \Delta^{-3/2}(\delta_1(k)k^{-1}) (k^{-1}\delta_2(k))) \\
-  \mathcal{D}_{2,1}(\Delta_{(1)}, \Delta_{(2)}) ( \Delta^{-3/2}(\delta_1(k)k^{-1}) \Delta^{1/2}(k^{-1}\delta_2(k))) \\
+ \mathcal{D}_{1,1}(\Delta_{(1)}, \Delta_{(2)}) ( \Delta^{-1/2}(\delta_1(k)k^{-1}) \Delta^{1/2}(k^{-1}\delta_2(k))) \\
+ \mathcal{D}_{1,1}(\Delta_{(1)}, \Delta_{(2)}) ( (\delta_1(k)k^{-1}) \Delta^{1/2}(k^{-1}\delta_2(k))) \\
+ \mathcal{D}_{1,1}(\Delta_{(1)}, \Delta_{(2)}) ( \Delta^{-1/2}(\delta_1(k)k^{-1}) (k^{-1}\delta_2(k))) \\
+ \mathcal{D}_{1,1}(\Delta_{(1)}, \Delta_{(2)}) ( (\delta_1(k)k^{-1}) (k^{-1}\delta_2(k))) \\
+2 \mathcal{D}_{2,2}(\Delta_{(1)}, \Delta_{(2)}) ( \Delta^{-1}(\delta_2(k)k^{-1}) \Delta^{1/2}(k^{-1}\delta_1(k))) \\
+ 2 \mathcal{D}_{2,2}(\Delta_{(1)}, \Delta_{(2)}) ( \Delta^{-1}(\delta_2(k)k^{-1}) (k^{-1}\delta_1(k))) \\
+ 2 \mathcal{D}_{2,2}(\Delta_{(1)}, \Delta_{(2)}) ( \Delta^{-3/2}(\delta_2(k)k^{-1}) \Delta^{1/2}(k^{-1}\delta_1(k))) \\
+ 2 \mathcal{D}_{2,2}(\Delta_{(1)}, \Delta_{(2)}) ( \Delta^{-3/2}(\delta_2(k)k^{-1}) (k^{-1}\delta_1(k))) \\
+ 4 \mathcal{D}_{3,1}(\Delta_{(1)}, \Delta_{(2)}) ( \Delta^{-2}(\delta_2(k)k^{-1}) \Delta^{1/2}(k^{-1}\delta_1(k))) \\
+ 4 \mathcal{D}_{3,1}(\Delta_{(1)}, \Delta_{(2)}) ( \Delta^{-2}(\delta_2(k)k^{-1}) (k^{-1}\delta_1(k))) \\
+ 4 \mathcal{D}_{3,1}(\Delta_{(1)}, \Delta_{(2)}) ( \Delta^{-5/2}(\delta_2(k)k^{-1}) \Delta^{1/2}(k^{-1}\delta_1(k))) \\
+ 4 \mathcal{D}_{3,1}(\Delta_{(1)}, \Delta_{(2)}) ( \Delta^{-5/2}(\delta_2(k)k^{-1}) (k^{-1}\delta_1(k))) \\
- 4 \mathcal{D}_{2,1}(\Delta_{(1)}, \Delta_{(2)}) ( \Delta^{-1}(\delta_2(k)k^{-1}) \Delta^{1/2} (k^{-1}\delta_1(k))) \\
- 4 \mathcal{D}_{2,1}(\Delta_{(1)}, \Delta_{(2)}) ( \Delta^{-1}(\delta_2(k)k^{-1})  (k^{-1}\delta_1(k))) \\
- 4 \mathcal{D}_{2,1}(\Delta_{(1)}, \Delta_{(2)}) ( \Delta^{-3/2}(\delta_2(k)k^{-1})  \Delta^{1/2}(k^{-1}\delta_1(k))) \\
- 4 \mathcal{D}_{2,1}(\Delta_{(1)}, \Delta_{(2)}) ( \Delta^{-3/2}(\delta_2(k)k^{-1})  (k^{-1}\delta_1(k))) \\
-  \mathcal{D}_{1,2}(\Delta_{(1)}, \Delta_{(2)}) ( (\delta_2(k)k^{-1}) \Delta^{1/2} (k^{-1}\delta_1(k))) \\
-  \mathcal{D}_{1,2}(\Delta_{(1)}, \Delta_{(2)}) (\Delta^{-1/2} (\delta_2(k)k^{-1}) \Delta^{1/2} (k^{-1}\delta_1(k))) \\
-  \mathcal{D}_{1,2}(\Delta_{(1)}, \Delta_{(2)}) ( \Delta^{-1/2}(\delta_2(k)k^{-1}) (k^{-1}\delta_1(k))) \\
-  \mathcal{D}_{1,2}(\Delta_{(1)}, \Delta_{(2)}) ( (\delta_2(k)k^{-1}) (k^{-1}\delta_1(k))) \\
-  \mathcal{D}_{2,1}(\Delta_{(1)}, \Delta_{(2)}) ( \Delta^{-3/2}(\delta_2(k)k^{-1}) \Delta^{1/2}(k^{-1}\delta_1(k))) \\
-  \mathcal{D}_{2,1}(\Delta_{(1)}, \Delta_{(2)}) ( \Delta^{-1}(\delta_2(k)k^{-1}) \Delta^{1/2}(k^{-1}\delta_1(k))) \\
-  \mathcal{D}_{2,1}(\Delta_{(1)}, \Delta_{(2)}) ( \Delta^{-3/2}(\delta_2(k)k^{-1}) (k^{-1}\delta_1(k))) \\
-  \mathcal{D}_{2,1}(\Delta_{(1)}, \Delta_{(2)}) ( \Delta^{-1}(\delta_2(k)k^{-1}) (k^{-1}\delta_1(k))) \\
-  \mathcal{D}_{2,1}(\Delta_{(1)}, \Delta_{(2)}) ( \Delta^{-1}(\delta_2(k)k^{-1}) (k^{-1}\delta_1(k))) \\
-  \mathcal{D}_{2,1}(\Delta_{(1)}, \Delta_{(2)}) ( \Delta^{-1}(\delta_2(k)k^{-1}) \Delta^{1/2}(k^{-1}\delta_1(k))) \\
-  \mathcal{D}_{2,1}(\Delta_{(1)}, \Delta_{(2)}) ( \Delta^{-3/2}(\delta_2(k)k^{-1}) (k^{-1}\delta_1(k))) \\
-  \mathcal{D}_{2,1}(\Delta_{(1)}, \Delta_{(2)}) ( \Delta^{-3/2}(\delta_2(k)k^{-1}) \Delta^{1/2}(k^{-1}\delta_1(k))) \\
+ \mathcal{D}_{1,1}(\Delta_{(1)}, \Delta_{(2)}) ( \Delta^{-1/2}(\delta_2(k)k^{-1}) \Delta^{1/2}(k^{-1}\delta_1(k))) \\
+ \mathcal{D}_{1,1}(\Delta_{(1)}, \Delta_{(2)}) ( (\delta_2(k)k^{-1}) \Delta^{1/2}(k^{-1}\delta_1(k))) \\
+ \mathcal{D}_{1,1}(\Delta_{(1)}, \Delta_{(2)}) ( \Delta^{-1/2}(\delta_2(k)k^{-1}) (k^{-1}\delta_1(k))) \\
+ \mathcal{D}_{1,1}(\Delta_{(1)}, \Delta_{(2)}) ( (\delta_2(k)k^{-1}) (k^{-1}\delta_1(k)))  \Big ) \\
-i \tau_2 \Big ( 
-\mathcal{D}_{1,2}(\Delta_{(1)}, \Delta_{(2)})( \Delta^{-1/2}(\delta_1(k)k^{-1}) \Delta^{1/2}(k^{-1}\delta_2(k))) \\
- \mathcal{D}_{1,2}(\Delta_{(1)}, \Delta_{(2)})( (\delta_1(k)k^{-1}) \Delta^{1/2}(k^{-1}\delta_2(k))) \\
- \mathcal{D}_{1,2}(\Delta_{(1)}, \Delta_{(2)})( \Delta^{-1/2}(\delta_1(k)k^{-1}) (k^{-1}\delta_2(k))) \\
- \mathcal{D}_{1,2}(\Delta_{(1)}, \Delta_{(2)})( (\delta_1(k)k^{-1}) (k^{-1}\delta_2(k))) \\
- \mathcal{D}_{2,1}(\Delta_{(1)}, \Delta_{(2)})( \Delta^{-3/2}(\delta_1(k)k^{-1}) \Delta^{1/2}(k^{-1}\delta_2(k))) \\
- \mathcal{D}_{2,1}(\Delta_{(1)}, \Delta_{(2)})( \Delta^{-1}(\delta_1(k)k^{-1}) \Delta^{1/2}(k^{-1}\delta_2(k))) \\
- \mathcal{D}_{2,1}(\Delta_{(1)}, \Delta_{(2)})( \Delta^{-3/2}(\delta_1(k)k^{-1}) (k^{-1}\delta_2(k))) \\
- \mathcal{D}_{2,1}(\Delta_{(1)}, \Delta_{(2)})( \Delta^{-1}(\delta_1(k)k^{-1}) (k^{-1}\delta_2(k))) \\
+  \mathcal{D}_{2,1}(\Delta_{(1)}, \Delta_{(2)})( \Delta^{-1}(\delta_1(k)k^{-1}) (k^{-1}\delta_2(k))) \\
+ \mathcal{D}_{2,1}(\Delta_{(1)}, \Delta_{(2)})( \Delta^{-1}(\delta_1(k)k^{-1}) \Delta^{1/2}(k^{-1}\delta_2(k))) \\
+ \mathcal{D}_{2,1}(\Delta_{(1)}, \Delta_{(2)})( \Delta^{-3/2}(\delta_1(k)k^{-1}) (k^{-1}\delta_2(k))) \\
+ \mathcal{D}_{2,1}(\Delta_{(1)}, \Delta_{(2)})( \Delta^{-3/2}(\delta_1(k)k^{-1}) \Delta^{1/2}(k^{-1}\delta_2(k))) \\
+ \mathcal{D}_{1,1}(\Delta_{(1)}, \Delta_{(2)})( \Delta^{-1/2}(\delta_1(k)k^{-1/2}) \Delta^{1/2}(k^{-1}\delta_2(k))) \\
+ \mathcal{D}_{1,1}(\Delta_{(1)}, \Delta_{(2)})( (\delta_1(k)k^{-1}) \Delta^{1/2}(k^{-1}\delta_2(k))) \\
+ \mathcal{D}_{1,1}(\Delta_{(1)}, \Delta_{(2)})( \Delta^{-1/2}(\delta_1(k)k^{-1}) (k^{-1}\delta_2(k))) \\
+ \mathcal{D}_{1,1}(\Delta_{(1)}, \Delta_{(2)})( (\delta_1(k)k^{-1}) (k^{-1}\delta_2(k))) \\
+\mathcal{D}_{1,2}(\Delta_{(1)}, \Delta_{(2)})( \Delta^{-1/2}(\delta_2(k)k^{-1}) \Delta^{1/2}(k^{-1}\delta_1(k))) \\
+ \mathcal{D}_{1,2}(\Delta_{(1)}, \Delta_{(2)})( (\delta_2(k)k^{-1}) \Delta^{1/2}(k^{-1}\delta_1(k))) \\
+ \mathcal{D}_{1,2}(\Delta_{(1)}, \Delta_{(2)})( \Delta^{-1/2}(\delta_2(k)k^{-1}) (k^{-1}\delta_1(k))) \\
+ \mathcal{D}_{1,2}(\Delta_{(1)}, \Delta_{(2)})( (\delta_2(k)k^{-1}) (k^{-1}\delta_1(k))) \\
+ \mathcal{D}_{2,1}(\Delta_{(1)}, \Delta_{(2)})( \Delta^{-3/2}(\delta_2(k)k^{-1}) \Delta^{1/2}(k^{-1}\delta_1(k))) \\
+ \mathcal{D}_{2,1}(\Delta_{(1)}, \Delta_{(2)})( \Delta^{-1}(\delta_2(k)k^{-1}) \Delta^{1/2}(k^{-1}\delta_1(k))) \\
+ \mathcal{D}_{2,1}(\Delta_{(1)}, \Delta_{(2)})( \Delta^{-3/2}(\delta_2(k)k^{-1}) (k^{-1}\delta_1(k))) \\
+ \mathcal{D}_{2,1}(\Delta_{(1)}, \Delta_{(2)})( \Delta^{-1}(\delta_2(k)k^{-1}) (k^{-1}\delta_1(k))) \\
-  \mathcal{D}_{2,1}(\Delta_{(1)}, \Delta_{(2)})( \Delta^{-1}(\delta_2(k)k^{-1}) (k^{-1}\delta_1(k))) \\
- \mathcal{D}_{2,1}(\Delta_{(1)}, \Delta_{(2)})( \Delta^{-1}(\delta_2(k)k^{-1}) \Delta^{1/2}(k^{-1}\delta_1(k))) \\
- \mathcal{D}_{2,1}(\Delta_{(1)}, \Delta_{(2)})( \Delta^{-3/2}(\delta_2(k)k^{-1}) (k^{-1}\delta_1(k))) \\
- \mathcal{D}_{2,1}(\Delta_{(1)}, \Delta_{(2)})( \Delta^{-3/2}(\delta_2(k)k^{-1}) \Delta^{1/2}(k^{-1}\delta_1(k))) \\
- \mathcal{D}_{1,1}(\Delta_{(1)}, \Delta_{(2)})( \Delta^{-1/2}(\delta_2(k)k^{-1/2}) \Delta^{1/2}(k^{-1}\delta_1(k))) \\
- \mathcal{D}_{1,1}(\Delta_{(1)}, \Delta_{(2)})( (\delta_2(k)k^{-1}) \Delta^{1/2}(k^{-1}\delta_1(k))) \\
- \mathcal{D}_{1,1}(\Delta_{(1)}, \Delta_{(2)})( \Delta^{-1/2}(\delta_2(k)k^{-1}) (k^{-1}\delta_1(k))) \\
- \mathcal{D}_{1,1}(\Delta_{(1)}, \Delta_{(2)})( (\delta_2(k)k^{-1}) (k^{-1}\delta_1(k))) 
\Big ).
$ \\

Putting together the latter with \eqref{twob0inttwo}, up to an overall factor of $\pi$, we get\\

\begin{eqnarray} \label{threeb0inttwo}
&&g_1(\Delta)(k^{-1} \delta_1^2(k)) + g_2(\Delta)(k^{-2} \delta_1(k)^2) \nonumber \\
&+& G(\Delta_{(1)}, \Delta_{(2)}) (  (\delta_1(k) k^{-1}) (k^{-1}\delta_1(k))) \nonumber \\
&+& |\tau|^2 g_1(\Delta)(k^{-1} \delta_2^2(k)) +|\tau|^2 g_2(\Delta)(k^{-2} \delta_2(k)^2) \nonumber \\
&+& |\tau|^2 G(\Delta_{(1)}, \Delta_{(2)}) (  (\delta_2(k) k^{-1}) (k^{-1}\delta_2(k))) \nonumber \\
&+& \tau_1 g_1(\Delta)(k^{-1} \delta_1\delta_2(k)) + \tau_1 g_2(\Delta)(k^{-2} \delta_1(k)\delta_2(k)) \nonumber \\
&+& \tau_1 G(\Delta_{(1)}, \Delta_{(2)}) (  (\delta_1(k) k^{-1}) (k^{-1}\delta_2(k))  ) \nonumber \\
&+& \tau_1 g_1(\Delta)(k^{-1} \delta_2\delta_1(k)) + \tau_1 g_2(\Delta)(k^{-2} \delta_2(k)\delta_1(k)) \nonumber \\
&+& \tau_1 G(\Delta_{(1)}, \Delta_{(2)}) (  (\delta_2(k) k^{-1}) (k^{-1}\delta_1(k))  ) \nonumber \\
&-& i \tau_2 L(\Delta_{(1)}, \Delta_{(2)}) ((\delta_1(k)k^{-1}) (k^{-1}\delta_2(k))) \nonumber \\
&+& i \tau_2 L(\Delta_{(1)}, \Delta_{(2)}) ((\delta_2(k)k^{-1}) (k^{-1}\delta_1(k))),
\end{eqnarray}
where as given by formulas \eqref{gone} and \eqref{gtwo}:

\begin{eqnarray} 
g_1(u)=  \frac{-1+u^2-2u \log u}{(-1+u^{1/2})^3(1+u^{1/2})^2}, \nonumber
\end{eqnarray}
\begin{eqnarray} 
g_2(u)= 2\frac{-1+u^2-2u \log u}{(-1+u)^3}, \nonumber 
\end{eqnarray} 
the function $G$ is defined by
\begin{eqnarray} 
G(u,v)&:=& 2 \mathcal{D}_{2,2}(u,v) u^{-1} v^{1/2} +2\mathcal{D}_{2,2}(u,v) u^{-1} + 2\mathcal{D}_{2,2}(u,v) u^{-3/2} v^{1/2} \nonumber \\
&&+ 2 \mathcal{D}_{2,2}(u,v) u^{-3/2} + 4\mathcal{D}_{3,1}(u,v) u^{-2} v^{1/2} +4 \mathcal{D}_{3,1}(u,v) u^{-2} \nonumber \\
&&+4 \mathcal{D}_{3,1}(u,v) u^{-5/2} v^{1/2}+4\mathcal{D}_{3,1}(u,v) u^{-5/2}-4\mathcal{D}_{2,1}(u,v) u^{-1} v^{1/2} \nonumber \\
&& -4\mathcal{D}_{2,1}(u,v) u^{-1} -4 \mathcal{D}_{2,1}(u,v) u^{-3/2} v^{1/2} -4 \mathcal{D}_{2,1}(u,v) u^{-3/2} \nonumber \\
&& - \mathcal{D}_{1,2}(u,v) v^{1/2} - \mathcal{D}_{1,2}(u,v) u^{-1/2} v^{1/2} - \mathcal{D}_{1,2}(u,v) u^{-1/2} \nonumber \\
&& - \mathcal{D}_{1,2}(u,v) - \mathcal{D}_{2,1}(u,v) u^{-3/2} v^{1/2} - \mathcal{D}_{2,1}(u,v) u^{-1} v^{1/2} \nonumber \\
&& - \mathcal{D}_{2,1}(u,v) u^{-3/2} - \mathcal{D}_{2,1}(u,v) u^{-1} - \mathcal{D}_{2,1}(u,v) u^{-1} \nonumber \\
&&- \mathcal{D}_{2,1}(u,v) u^{-1} v^{1/2} - \mathcal{D}_{2,1}(u,v) u^{-3/2} - \mathcal{D}_{2,1}(u,v) u^{-3/2} v^{1/2} \nonumber \\
&& + \mathcal{D}_{1,1}(u,v) u^{-1/2} v^{1/2} +\mathcal{D}_{1,1}(u,v)  v^{1/2} + \mathcal{D}_{1,1}(u,v) u^{-1/2}  \nonumber \\
&&+\mathcal{D}_{1,1}(u,v) \nonumber
\end{eqnarray}
\begin{eqnarray} \label{G}
&=& -(\sqrt{u} (u (-1 + v)^2 (-1 + u v (-4 + u (4 + v))) \log(1/u) + (-1 + u)  \nonumber \\
&& ((1 + u (-2 + v)) (-1 + v) (-1 + u v) (1 + u v) + (-1 + u) v \nonumber \\
&& (-1 + u (-4 + v (4 + u v))) \log v)))/((-1 + \sqrt{u})^2 (1 + \sqrt{u}) (-1 + \sqrt{v})^2 \nonumber \\
&&(1 + \sqrt{v}) (-1 + u v)^3), \nonumber
\end{eqnarray}

and

\begin{eqnarray} 
L(u,v)&:=& - \mathcal{D}_{1,2}(u,v) u^{-1/2}v^{1/2} - \mathcal{D}_{1,2}(u,v) v^{1/2} - \mathcal{D}_{1,2}(u,v) u^{-1/2} \nonumber \\
&& - \mathcal{D}_{1,2}(u,v) -\mathcal{D}_{2,1}(u,v) u^{-3/2}v^{1/2} -\mathcal{D}_{2,1}(u,v) u^{-1}v^{1/2} \nonumber \\
&& - \mathcal{D}_{2,1}(u,v) u^{-3/2} - \mathcal{D}_{2,1}(u,v) u^{-1} + \mathcal{D}_{2,1}(u,v) u^{-1} \nonumber \\
&&+ \mathcal{D}_{2,1}(u,v) u^{-1}v^{1/2} + \mathcal{D}_{2,1}(u,v) u^{-3/2} + \mathcal{D}_{2,1}(u,v) u^{-3/2}v^{1/2} \nonumber \\
&& + \mathcal{D}_{1,1}(u,v) u^{-1/2}v^{1/2} + \mathcal{D}_{1,1}(u,v) v^{1/2} + \mathcal{D}_{1,1}(u,v) u^{-1/2} \nonumber \\
&& + \mathcal{D}_{1,1}(u,v) \nonumber
\end{eqnarray}
\begin{eqnarray} \label{L}
&=& (\sqrt{u} (u (-1 + v)^2 \log(1/u) + (-1 + u) ((-1 + v) (-1 + u v) + (v - u v) \nonumber \\
&& \log v)))/((-1 + \sqrt{u})^2 (1 + \sqrt{u}) (-1 + \sqrt{v})^2 (1 + \sqrt{v}) (-1 + u v)). \nonumber
\end{eqnarray}

\section{The Scalar Curvature in Terms of $\log (k)$} \label{scacurlogsec}
In order to express the scalar curvature of $(\mathbb{T}_\theta^2, \tau, k)$ in terms of $\log k$ 
we need to find some identities that relate $k^{-1} \delta_i\delta_j(k)$ and $k^{-2} \delta_i(k)^2$, for $i,j=1,2$, to 
$\log k$. This is carried out in the following lemma.

\newtheorem{derlog}[main1]{Lemma}
\begin{derlog} \label{derlog}
For $i,j=1,2,$ we have
\begin{eqnarray} \label{derlog2}
k^{-2} \delta_i(k) \delta_j(k) = 4 \frac{\Delta - \Delta^{1/2}}{\log \Delta}(\delta_i(\log k)) 
\frac{\Delta^{1/2}-1}{\log \Delta} (\delta_j(\log k)). 
\end{eqnarray}
Also we have
\begin{eqnarray} \label{derlog1}
k^{-1} \delta_i \delta_j (k) &=&   2 \frac{\Delta^{1/2}-1}{ \log \Delta} (\delta_i \delta_j (\log k)) +  g(\Delta_{(1)}, \Delta_{(2)})(\delta_j(\log k) \delta_i(\log k)) +\nonumber \\
&& g(\Delta_{(1)}, \Delta_{(2)})(\delta_i(\log k) \delta_j(\log k)), 
\end{eqnarray}
where
\begin{equation} \label{gfunction}
g(u, v):= 4 \frac{(\sqrt{uv}-1) \log u - (\sqrt{u}-1) \log (uv)}{\log v \log u \log (uv)}, \nonumber
\end{equation}
and $\Delta_{(i)}$ signifies the action of $\Delta$ on the $i$-th factor of the product.
\begin{proof}
We note that the following identity from \cite{contre} will be used in the proof of both identities:
\[ k^{-1} \delta_j(k) = 2\frac{\Delta^{1/2}-1}{\log \Delta} (\delta_j(\log k)).\]
First we prove \eqref{derlog2}:
\begin{eqnarray}
k^{-2} \delta_i(k) \delta_j(k) &=& k^{-1} 2\frac{\Delta^{1/2}-1}{\log \Delta}(\delta_i(\log k)) \, \delta_j(k) \nonumber \\
&=& 4 \frac{\Delta - \Delta^{1/2}}{\log \Delta}(\delta_i(\log k)) \frac{\Delta^{1/2}-1}{\log \Delta} (\delta_j(\log k)). \nonumber
\end{eqnarray} 
To prove \eqref{derlog1}, we write
\begin{eqnarray}
k^{-1} \delta_i \delta_j (k) &=&  \int_0^1 \Delta^{s/2} \delta_i \delta_j (\log k) \, ds + \nonumber \\
&& \int_0^1  \Delta^{s/2} (\delta_j(\log k)) 2 \frac{\Delta^{s/2}-1}{\log\Delta}(\delta_i(\log k)) \,ds + \nonumber \\
&& \int_0^1  \Delta^{s/2} (\delta_i(\log k)) 2\frac{\Delta^{s/2}-1}{\log\Delta}(\delta_j(\log k))\,ds \nonumber \\
&=& 2 \frac{\Delta^{1/2}-1}{ \log \Delta} (\delta_i \delta_j (\log k)) + \nonumber \\
&&  g(\Delta_{(1)}, \Delta_{(2)})(\delta_j(\log k) \delta_i(\log k)) + \nonumber \\
&& g(\Delta_{(1)}, \Delta_{(2)})(\delta_i(\log k) \delta_j(\log k)). \nonumber
\end{eqnarray}
\end{proof}
\end{derlog}

\subsection{The terms corresponding to $k \partial^* \partial k$.}
Applying  Lemma \ref{derlog} to the local expression \eqref{threeb0int}, we can write it in terms of $\log{k}$ as follows:
\begin{eqnarray}
\eqref{threeb0int}&=&f_1(\Delta)(2 \frac{\Delta^{1/2}-1}{ \log \Delta}(\delta_1^2(\log{k})))  \nonumber \\
&&+  f_1(\Delta)(2 g(\Delta_{(1)}, \Delta_{(2)}) (\delta_1(\log k) \delta_1(\log k)) )  \nonumber \\
&&+f_2(\Delta)(4\frac{\Delta - \Delta^{1/2}}{\log \Delta}(\delta_1(\log k)) \frac{\Delta^{1/2}-1}{\log \Delta}(\delta_1(\log k))) \nonumber \\
&&+F(\Delta_{(1)}, \Delta_{(2)})((-2\frac{\Delta^{-1/2}-1}{\log \Delta}(\delta_1(\log k)))(2\frac{\Delta^{1/2}-1}{\log \Delta}(\delta_1(\log k)))) \nonumber 
\end{eqnarray}
\begin{eqnarray}
&&+ |\tau|^2 f_1(\Delta)(2 \frac{\Delta^{1/2}-1}{ \log \Delta}(\delta_2^2(\log{k})))  \nonumber \\
&&+ |\tau|^2 f_1(\Delta)(2 g(\Delta_{(1)}, \Delta_{(2)}) (\delta_2(\log k) \delta_2(\log k)) )  \nonumber \\
&&+|\tau|^2 f_2(\Delta)(4\frac{\Delta - \Delta^{1/2}}{\log \Delta}(\delta_2(\log k)) \frac{\Delta^{1/2}-1}{\log \Delta}(\delta_2(\log k))) \nonumber \\
&&+|\tau|^2F(\Delta_{(1)}, \Delta_{(2)})((-2\frac{\Delta^{-1/2}-1}{\log \Delta}(\delta_2(\log k)))(2\frac{\Delta^{1/2}-1}{\log \Delta}(\delta_2(\log k))))  \nonumber 
\end{eqnarray}
\begin{eqnarray}
&& +\tau_1 f_1(\Delta)(2 \frac{\Delta^{1/2}-1}{ \log \Delta}(\delta_1\delta_2(\log{k})))  \nonumber \\
&&+\tau_1  f_1(\Delta)( g(\Delta_{(1)}, \Delta_{(2)}) (\delta_2(\log k) \delta_1(\log k)) )  \nonumber \\
&&+\tau_1  f_1(\Delta)( g(\Delta_{(1)}, \Delta_{(2)}) (\delta_1(\log k) \delta_2(\log k)) )  \nonumber \\
&&+\tau_1 f_2(\Delta)(4\frac{\Delta - \Delta^{1/2}}{\log \Delta}(\delta_1(\log k)) \frac{\Delta^{1/2}-1}{\log \Delta}(\delta_2(\log k))) \nonumber \\
&&+\tau_1F(\Delta_{(1)}, \Delta_{(2)})((-2\frac{\Delta^{-1/2}-1}{\log \Delta}(\delta_1(\log k)))(2\frac{\Delta^{1/2}-1}{\log \Delta}(\delta_2(\log k)))) \nonumber 
\end{eqnarray}
\begin{eqnarray} \label{locexp}
&& +\tau_1 f_1(\Delta)(2 \frac{\Delta^{1/2}-1}{ \log \Delta}(\delta_2\delta_1(\log{k})))  \nonumber \\
&&+\tau_1  f_1(\Delta)( g(\Delta_{(1)}, \Delta_{(2)}) (\delta_1(\log k) \delta_2(\log k)) )  \nonumber \\
&&+\tau_1  f_1(\Delta)( g(\Delta_{(1)}, \Delta_{(2)}) (\delta_2(\log k) \delta_1(\log k)) )  \nonumber \\
&&+\tau_1 f_2(\Delta)(4\frac{\Delta - \Delta^{1/2}}{\log \Delta}(\delta_2(\log k)) \frac{\Delta^{1/2}-1}{\log \Delta}(\delta_1(\log k))) \nonumber \\
&&+\tau_1F(\Delta_{(1)}, \Delta_{(2)})((-2\frac{\Delta^{-1/2}-1}{\log \Delta}(\delta_2(\log k)))(2\frac{\Delta^{1/2}-1}{\log \Delta}(\delta_1(\log k)))). \nonumber 
\end{eqnarray}
Now, writing the latter in terms of $\log \Delta$ and considering an overall factor of -1 (\emph{cf.} Subsection \ref{recipe}), up to an overall factor of 
$\frac{\pi}{\tau_2}$, we obtain the following expression for the first component of the scalar curvature of the perturbed spectral triple attached to 
$(\mathbb{T}_\theta^2, \tau, k)$:
\begin{eqnarray}\label{localexp}
&&K(\log \Delta) \big (\delta_1^2(\log k) + |\tau|^2 \delta_2^2(\log k) +2 \tau_1 \delta_1\delta_2(\log k) \big )+  \nonumber \\
&& H(\log \Delta_{(1)}, \log \Delta_{(2)}) \big ( \delta_1(\log k) \delta_1(\log k) + |\tau|^2 \delta_2(\log k) \delta_2(\log k) + \nonumber \\ 
&& \qquad \qquad \qquad \qquad \qquad \qquad \tau_1 \delta_1(\log k) \delta_2(\log k) + \tau_1 \delta_2(\log k) \delta_1(\log k) \big), \nonumber 
\end{eqnarray}
where
\[ K(x) := - 2 f_1(e^x) \frac{e^{x/2}-1}{x} = \frac{2 e^{x/2} (2 + e^x (-2 + x) + x)}{(-1 + e^x)^2 x} ,\]
and
\begin{eqnarray}
H(s, t) &:=& -   2 f_1(e^{s+t}) g(e^s,e^t)-4 f_2(e^{s+t}) \frac{e^s-e^{s/2}}{s} \frac{e^{t/2}-1}{t} + \nonumber \\
&& 4 F(e^s, e^t) \frac{e^{-s/2}-1}{s} \frac{e^{t/2}-1}{t} \nonumber \\
&=& \nonumber
\end{eqnarray}
\begin{eqnarray} \label{H}
-\frac{-t (s + t) \cosh{s} +  s (s + t) \cosh{t} - (s - t) (s + t + \sinh{s} + \sinh{t} - \sinh(s + t))}{s t (s + t)\sinh(s/2) \sinh(t/2) \sinh^2 ((s + t)/2) }.\nonumber 
\end{eqnarray}

\subsection{The terms corresponding to $ \partial^* k^2 \partial$.}

We also apply Lemma \ref{derlog} to the local expression \eqref{threeb0inttwo} and obtain the following:

 \begin{eqnarray}
\eqref{threeb0inttwo}  &=&g_1(\Delta)(2 \frac{\Delta^{1/2}-1}{ \log \Delta}(\delta_1^2(\log{k})))  \nonumber \\
&&+  g_1(\Delta)(2 g(\Delta_{(1)}, \Delta_{(2)}) (\delta_1(\log k) \delta_1(\log k)) )  \nonumber \\
&&+g_2(\Delta)(4\frac{\Delta - \Delta^{1/2}}{\log \Delta}(\delta_1(\log k)) \frac{\Delta^{1/2}-1}{\log \Delta}(\delta_1(\log k))) \nonumber \\
&&+G(\Delta_{(1)}, \Delta_{(2)})((-2\frac{\Delta^{-1/2}-1}{\log \Delta}(\delta_1(\log k)))(2\frac{\Delta^{1/2}-1}{\log \Delta}(\delta_1(\log k)))) \nonumber 
\end{eqnarray}
\begin{eqnarray}
&&+ |\tau|^2 g_1(\Delta)(2 \frac{\Delta^{1/2}-1}{ \log \Delta}(\delta_2^2(\log{k})))  \nonumber \\
&&+ |\tau|^2 g_1(\Delta)(2 g(\Delta_{(1)}, \Delta_{(2)}) (\delta_2(\log k) \delta_2(\log k)) )  \nonumber \\
&&+|\tau|^2 g_2(\Delta)(4\frac{\Delta - \Delta^{1/2}}{\log \Delta}(\delta_2(\log k)) \frac{\Delta^{1/2}-1}{\log \Delta}(\delta_2(\log k))) \nonumber \\
&&+|\tau|^2 G(\Delta_{(1)}, \Delta_{(2)})((-2\frac{\Delta^{-1/2}-1}{\log \Delta}(\delta_2(\log k)))(2\frac{\Delta^{1/2}-1}{\log \Delta}(\delta_2(\log k))))  \nonumber 
\end{eqnarray}
\begin{eqnarray}
&& +\tau_1 g_1(\Delta)(2 \frac{\Delta^{1/2}-1}{ \log \Delta}(\delta_1\delta_2(\log{k})))  \nonumber \\
&&+\tau_1  g_1(\Delta)( g(\Delta_{(1)}, \Delta_{(2)}) (\delta_2(\log k) \delta_1(\log k)) )  \nonumber \\
&&+\tau_1  g_1(\Delta)( g(\Delta_{(1)}, \Delta_{(2)}) (\delta_1(\log k) \delta_2(\log k)) )  \nonumber \\
&&+\tau_1 g_2(\Delta)(4\frac{\Delta - \Delta^{1/2}}{\log \Delta}(\delta_1(\log k)) \frac{\Delta^{1/2}-1}{\log \Delta}(\delta_2(\log k))) \nonumber \\
&&+\tau_1 G(\Delta_{(1)}, \Delta_{(2)})((-2\frac{\Delta^{-1/2}-1}{\log \Delta}(\delta_1(\log k)))(2\frac{\Delta^{1/2}-1}{\log \Delta}(\delta_2(\log k)))) \nonumber 
\end{eqnarray}
\begin{eqnarray} \label{locexp}
&& +\tau_1 g_1(\Delta)(2 \frac{\Delta^{1/2}-1}{ \log \Delta}(\delta_2\delta_1(\log{k})))  \nonumber \\
&&+\tau_1  g_1(\Delta)( g(\Delta_{(1)}, \Delta_{(2)}) (\delta_1(\log k) \delta_2(\log k)) )  \nonumber \\
&&+\tau_1  g_1(\Delta)( g(\Delta_{(1)}, \Delta_{(2)}) (\delta_2(\log k) \delta_1(\log k)) )  \nonumber \\
&&+\tau_1 g_2(\Delta)(4\frac{\Delta - \Delta^{1/2}}{\log \Delta}(\delta_2(\log k)) \frac{\Delta^{1/2}-1}{\log \Delta}(\delta_1(\log k))) \nonumber \\
&&+\tau_1G(\Delta_{(1)}, \Delta_{(2)})((-2\frac{\Delta^{-1/2}-1}{\log \Delta}(\delta_2(\log k)))(2\frac{\Delta^{1/2}-1}{\log \Delta}(\delta_1(\log k)))) \nonumber 
\end{eqnarray}
\begin{eqnarray}
&& - i \tau_2 L(\Delta_{(1)}, \Delta_{(2)}) ((-2\frac{\Delta^{-1/2}-1}{\log \Delta}(\delta_1(\log k)))
(2\frac{\Delta^{1/2}-1}{\log \Delta}(\delta_2(\log k)))) \nonumber \\
&& + i \tau_2 L(\Delta_{(1)}, \Delta_{(2)}) ((-2\frac{\Delta^{-1/2}-1}{\log \Delta}(\delta_2(\log k)))
(2\frac{\Delta^{1/2}-1}{\log \Delta}(\delta_1(\log k)))). \nonumber 
\end{eqnarray}
Now we write the latter in terms of $\log \Delta$, and after considering an overall factor of $-1$, up to an overall factor of $\frac{\pi}{\tau_2}$, 
we obtain the following expression for the second component of the scalar curvature of the perturbed spectral triple attached to 
$(\mathbb{T}_\theta^2, \tau, k)$:
\begin{eqnarray}\label{localexp}
&&S(\log \Delta) \big (\delta_1^2(\log k) + |\tau|^2 \delta_2^2(\log k) +2 \tau_1 \delta_1\delta_2(\log k) \big )+  \nonumber \\
&& T(\log \Delta_{(1)}, \log \Delta_{(2)}) \big ( \delta_1(\log k) \delta_1(\log k) + |\tau|^2 \delta_2(\log k) \delta_2(\log k) + \nonumber \\ 
&& \qquad \qquad \qquad \qquad \qquad \qquad \tau_1 \delta_1(\log k) \delta_2(\log k) + \tau_1 \delta_2(\log k) \delta_1(\log k) \big) -\nonumber \\
&& i \tau_2 W(\log \Delta_{(1)}, \log \Delta_{(2)}) \big (  \delta_1(\log k) \delta_2 (\log k) - \delta_2(\log k) \delta_1(\log k)   \big ), \nonumber
\end{eqnarray}
where
\begin{equation} \label{S}
S(x) := - 2 g_1(e^x) \frac{e^{x/2}-1}{x} = -\frac{4 e^x (-x + \sinh x)}{(-1 + e^{x/2})^2 (1 + e^{x/2})^2 x}, \nonumber
\end{equation}
\begin{eqnarray}
T(s, t) &:=&  -  2 g_1(e^{s+t}) g(e^s,e^t)-4 g_2(e^{s+t}) \frac{e^s-e^{s/2}}{s} \frac{e^{t/2}-1}{t} + \nonumber \\
&& 4 G(e^s, e^t) \frac{e^{-s/2}-1}{s} \frac{e^{t/2}-1}{t} \nonumber 
\end{eqnarray}
\begin{eqnarray} \label{T}
= - \cosh ((s+t)/2) \times \qquad \qquad \qquad \qquad \qquad \qquad \qquad \qquad \qquad \qquad \quad \quad \quad \quad \nonumber\\
 \frac{-t (s + t) \cosh{s} +  s (s + t) \cosh{t} - (s - t) (s + t + \sinh{s} + \sinh{t} - \sinh(s + t))}{s t (s + t)\sinh(s/2) \sinh(t/2) \sinh^2 ((s + t)/2) }, \nonumber 
\end{eqnarray}
and
\begin{eqnarray} \label{W}
W(s,t)&:=& +4 L(e^s, e^t)\frac{e^{-s/2} - 1}{s}\frac{e^{t/2} - 1}{t} \nonumber \\
&=& -4 \frac{ -s - t + t \cosh s + s \cosh t + \sinh s + \sinh t - \sinh (s + t)}{s t (\sinh s + \sinh t - \sinh (s + t))} \nonumber \\
&=& \frac{-s - t + t \cosh s + s \cosh t + \sinh s + \sinh t - \sinh (s + t)}{st\sinh (s/2) \sinh (t/2) \sinh ((s + t)/2)}. \nonumber 
\end{eqnarray}

\subsection{The scalar curvature.}

We collect the results of this paper in the following theorems. They are also independently proved  by Connes and Moscovici in \cite{conmos2}.  Note that in our final formulas we have considered an overall minus sign which comes from the 
change of sign initially considered in the Cauchy integral formula \eqref{cauchyint}.

\newtheorem{curfull}[main1]{Theorem}
\begin{curfull} \label{curfull}
Let $\theta$ be an irrational number, $\tau$ a complex number in the upper half plane representing the conformal class 
of a metric on $T_\theta^2$, and $k$ an invertible positive element in $A_\theta^\infty$ playing the role of the Weyl factor. Then the 
scalar curvature $R$ of the perturbed spectral triple attached to $(T_\theta^2, \tau, k)$, up to an overall factor of $-\frac{\pi}{\tau_2}$, 
is equal to

\begin{eqnarray} \label{finalexpungraded}
&& R_1(\log \Delta) \big ( \delta_1^2(\log k) + |\tau|^2 \delta_2^2(\log k) +2 \tau_1 \delta_1\delta_2(\log k)  \big ) \nonumber \\
&+&  R_2(\log \Delta_{(1)}, \log \Delta_{(2)}) \Big (  \delta_1(\log k) \delta_1(\log k) + |\tau|^2 \delta_2(\log k) \delta_2(\log k) + \nonumber \\
&& \qquad \qquad \qquad \qquad  \qquad  \qquad       \tau_1 \big (\delta_1(\log k) \delta_2(\log k) + \delta_2(\log k) \delta_1(\log k) \big )   \Big )
\nonumber \\
&-&  i  W(\log \Delta_{(1)}, \log \Delta_{(2)})  \Big ( \tau_2 \big ( \delta_1(\log k) \delta_2 (\log k) - \delta_2(\log k) \delta_1(\log k) \big ) \Big ), \nonumber
\end{eqnarray} 
where 
\begin{eqnarray}
R_1(x) :=K(x)+S(x)= -\frac{2 \coth (x/4)}{x} + \frac{1}{2 \sinh^2 (x/4)}=\frac{\frac{1}{2}- \frac{\sinh(x/2)}{x}}{\sinh^2(x/4)}, \nonumber
\end{eqnarray} 
$R_2(s, t):= H(s,t)+T(s,t) $
$= - (1+ \cosh ((s+t)/2) ) \times $
\begin{eqnarray}
\,\,\, \frac{-t (s + t) \cosh{s} +  s (s + t) \cosh{t} - (s - t) (s + t + \sinh{s} + \sinh{t} - \sinh(s + t))}{s t (s + t)\sinh(s/2) \sinh(t/2) \sinh^2 ((s + t)/2) }, \nonumber 
\end{eqnarray}
and 
\begin{eqnarray} 
W(s,t)=  \frac{-s - t + t \cosh s + s \cosh t + \sinh s + \sinh t - \sinh (s + t)}{st\sinh (s/2) \sinh (t/2) \sinh ((s + t)/2)}. \nonumber 
\end{eqnarray}
\end{curfull}

\newtheorem{curgraded}[main1]{Theorem}
\begin{curgraded} Assuming the hypotheses of Theorem \ref{curfull},  the chiral 
scalar curvature $R^\gamma$ of the perturbed graded spectral triple attached to $(T_\theta^2, \tau, k)$, up to an overall factor of $-\frac{\pi}{\tau_2}$, 
is given by

\begin{eqnarray} \label{finalexpgraded}
 && R^\gamma_1(\log \Delta) \big ( \delta_1^2(\log k) + |\tau|^2 \delta_2^2(\log k) +2 \tau_1 \delta_1\delta_2(\log k)  \big ) \nonumber \\
&+&  R^\gamma_2(\log \Delta_{(1)}, \log \Delta_{(2)}) \Big (  \delta_1(\log k) \delta_1(\log k) + |\tau|^2 \delta_2(\log k) \delta_2(\log k) + \nonumber \\
&& \qquad \qquad \qquad \qquad  \qquad  \qquad       \tau_1 \big (\delta_1(\log k) \delta_2(\log k) + \delta_2(\log k) \delta_1(\log k) \big )   \Big )
\nonumber \\
&+&  i  W(\log \Delta_{(1)}, \log \Delta_{(2)})  \Big ( \tau_2 \big ( \delta_1(\log k) \delta_2 (\log k) - \delta_2(\log k) \delta_1(\log k) \big ) \Big ), \nonumber 
\end{eqnarray} 
where 
\begin{eqnarray}
R^\gamma_1(x) :=K(x)-S(x)= \frac{x + 2 \sinh(x/2)}{x + x \cosh(x/2)}= \frac{\frac{1}{2}+\frac{\sinh(x/2)}{x}}{\cosh^2 (x/4)},  \nonumber
\end{eqnarray} 
$R^\gamma_2(s, t):= H(s,t)-T(s,t)   $
$=- (1- \cosh ((s+t)/2) ) \times$
\begin{eqnarray}
\,\,\, \frac{-t (s + t) \cosh{s} +  s (s + t) \cosh{t} - (s - t) (s + t + \sinh{s} + \sinh{t} - \sinh(s + t))}{s t (s + t)\sinh(s/2) \sinh(t/2) \sinh^2 ((s + t)/2) }, \nonumber 
\end{eqnarray}
and 
\begin{eqnarray} 
W(s,t)=  \frac{-s - t + t \cosh s + s \cosh t + \sinh s + \sinh t - \sinh (s + t)}{st\sinh (s/2) \sinh (t/2) \sinh ((s + t)/2)}. \nonumber 
\end{eqnarray}
\end{curgraded}

\newtheorem{commcase}[main1]{Remark}
\begin{commcase}
We note that the above local expressions $R$ and $R^\gamma$ for the scalar curvature of 
$(\mathbb{T}_\theta^2, \tau, k)$, reduce to the scalar curvature of the ordinary two torus  when $\theta = 0$. Namely, since 
\[\lim_{x \to 0} R_1(x) = -\frac{1}{3},\]
\[ \lim_{x \to 0} R_1^\gamma(x) = 1,\]
\[ \lim_{s,t \to 0} R_2(s,t) = \lim_{s,t \to 0} R^\gamma_2(s,t) = 0, \]
and
\[  \lim_{s,t \to 0} W(s,t) = - \frac{2}{3},\]
in the commutative case, the expressions for $R$ and $R^\gamma$ stated in the above theorems, reduce to 
constant multiples of 
\[  \frac{1}{\tau_2}\delta_1^2(\log k) + \frac{|\tau|^2}{\tau_2} \delta_2^2(\log k) + 2 \frac{\tau_1}{\tau_2} \delta_1\delta_2(\log k). \]
\end{commcase}

\section{Appendices}
As we mentioned before, after direct computations the $b_2$ terms of the two parts of the Laplacian $\triangle$ 
attached to $(\mathbb{T}_\theta^2, \tau, k)$, have quite lengthy formulas which, for the convenience of the reader, are 
recorded here.
   
\appendix
\section{ The $b_2$ term of $k \partial^* \partial k$} \label{firstb2}
The $b_2$ term of the first operator, namely $k \partial^* \partial k$, is equal to \\

\noindent 
$-b_0k\delta_1^2(k)b_0-
2\tau_1b_0k\delta_1\delta_2(k)b_0
-|\tau|^2b_0k\delta_2^2
(k)b_0+\\
6\xi_1^2b_0^2k^2\delta_1(k)^2b_0+\xi_1^2b_0^2k^
2\delta_1^2(k)b_0k+
5\xi_1^2b_0^2k^3\delta_1^2(k)b_0+\\
2\xi_1^2b_0k\delta_1(k)b_0\delta_1(k)b_0k+
6\xi_1^2b_0k\delta_1(k)b_0k\delta_1(k)b_0+\\
6\tau_1\xi_1^2b_0^2k^2\delta_1(k)\delta_2(k)b_0+
6\tau_1\xi_1^2b_0^2k^2\delta_2(k)\delta_1(k)b_0+\\
2\tau_1\xi_1^2b_0^2k^2\delta_1\delta_2(k)b_0k+
10\tau_1\xi_1^2b_0^2k^3\delta_1\delta_2(k)b_0+\\
2\tau_1\xi_1^2b_0k\delta_1(k)b_0\delta_2(k)b_0k+
6\tau_1\xi_1^2b_0k\delta_1(k)b_0k\delta_2(k)b_0+\\
2\tau_1\xi_1^2b_0k\delta_2(k)b_0\delta_1(k)b_0k+
6\tau_1\xi_1^2b_0k\delta_2(k)b_0k\delta_1(k)b_0+\\
12\tau_1\xi_1\xi_2b_0^2k^2\delta_1(k)^2b_0+
2\tau_1\xi_1\xi_2b_0^2k^2\delta_1^2(k)b_0k+\\
10\tau_1\xi_1\xi_2b_0^2k^3\delta_1^2(k)b_0+
4\tau_1\xi_1\xi_2b_0k\delta_1(k)b_0\delta_1(k)b_0k+\\
12\tau_1\xi_1\xi_2b_0k\delta_1(k)b_0k\delta_1(k)b_0+
4\tau_1^2\xi_1^2b_0^2k^2\delta_2(k)^2b_0+\\
4\tau_1^2\xi_1^2b_0^2k^3\delta_2^2(k)b_0+
4\tau_1^2\xi_1^2b_0k\delta_2(k)b_0k\delta_2(k)b_0+\\
8\tau_1^2\xi_1\xi_2b_0^2k^2\delta_1(k)\delta_2(k)b_0+
8\tau_1^2\xi_1\xi_2b_0^2k^2\delta_2(k)\delta_1(k)b_0+\\
4\tau_1^2\xi_1\xi_2b_0^2k^2\delta_1\delta_2(k)b_0k+
12\tau_1^2\xi_1\xi_2b_0^2k^3\delta_1\delta_2(k)b_0+\\
4\tau_1^2\xi_1\xi_2b_0k\delta_1(k)b_0\delta_2(k)b_0k+
8\tau_1^2\xi_1\xi_2b_0k\delta_1(k)b_0k\delta_2(k)b_0+\\
4\tau_1^2\xi_1\xi_2b_0k\delta_2(k)b_0\delta_1(k)b_0k+
8\tau_1^2\xi_1\xi_2b_0k\delta_2(k)b_0k\delta_1(k)b_0+\\
4\tau_1^2\xi_2^2b_0^2k^2\delta_1(k)^2b_0+
4\tau_1^2\xi_2^2b_0^2k^3\delta_1^2(k)b_0+\\
4\tau_1^2\xi_2^2b_0k\delta_1(k)b_0k\delta_1(k)b_0+
2|\tau|^2\xi_1^2b_0^2k^2\delta_2(k)^2b_0+\\
|\tau|^2\xi_1^2b_0^2k^2\delta_2^2(k)b_0k+
|\tau|^2\xi_1^2b_0^2k^3\delta_2^2(k)b_0+\\
2|\tau|^2\xi_1^2b_0k\delta_2(k)b_0\delta_2(k)b_0k+
2|\tau|^2\xi_1^2b_0k\delta_2(k)b_0k\delta_2(k)b_0+\\
4|\tau|^2\xi_1\xi_2b_0^2k^2\delta_1(k)\delta_2(k)b_0+
4|\tau|^2\xi_1\xi_2b_0^2k^2\delta_2(k)\delta_1(k)b_0+\\
8|\tau|^2\xi_1\xi_2b_0^2k^3\delta_1\delta_2(k)b_0+
4|\tau|^2\xi_1\xi_2b_0k\delta_1(k)b_0k\delta_2(k)b_0+\\
4|\tau|^2\xi_1\xi_2b_0k\delta_2(k)b_0k\delta_1(k)b_0+
2|\tau|^2\xi_2^2b_0^2k^2\delta_1(k)^2b_0+\\
|\tau|^2\xi_2^2b_0^2k^2\delta_1^2(k)b_0k+
|\tau|^2\xi_2^2b_0^2k^3\delta_1^2(k)b_0+\\
2|\tau|^2\xi_2^2b_0k\delta_1(k)b_0\delta_1(k)b_0k+
2|\tau|^2\xi_2^2b_0k\delta_1(k)b_0k\delta_1(k)b_0+\\
12\tau_1|\tau|^2\xi_1\xi_2b_0^2k^2\delta_2(k)^2b_0+
2\tau_1|\tau|^2\xi_1\xi_2b_0^2k^2\delta_2^2(k)b_0k+\\
10\tau_1|\tau|^2\xi_1\xi_2b_0^2k^3\delta_2^2(k)b_0+
4\tau_1|\tau|^2\xi_1\xi_2b_0k\delta_2(k)b_0\delta_2(k)b_0k+\\
12\tau_1|\tau|^2\xi_1\xi_2b_0k\delta_2(k)b_0k\delta_2(k)b_0+
6\tau_1|\tau|^2\xi_2^2b_0^2k^2\delta_1(k)\delta_2(k)b_0+\\
6\tau_1|\tau|^2\xi_2^2b_0^2k^2\delta_2(k)\delta_1(k)b_0+
2\tau_1|\tau|^2\xi_2^2b_0^2k^2\delta_1\delta_2(k)b_0k+\\
10\tau_1|\tau|^2\xi_2^2b_0^2k^3\delta_1\delta_2(k)b_0+
2\tau_1|\tau|^2\xi_2^2b_0k\delta_1(k)b_0\delta_2(k)b_0k+\\
6\tau_1|\tau|^2\xi_2^2b_0k\delta_1(k)b_0k\delta_2(k)b_0+
2\tau_1|\tau|^2\xi_2^2b_0k\delta_2(k)b_0\delta_1(k)b_0k+\\
6\tau_1|\tau|^2\xi_2^2b_0k\delta_2(k)b_0k\delta_1(k)b_0+
6|\tau|^4\xi_2^2b_0^2k^2\delta_2(k)^2b_0+\\
|\tau|^4\xi_2^2b_0^2k^2\delta_2^2(k)b_0k+
5|\tau|^4\xi_2^2b_0^2k^3\delta_2^2(k)b_0+\\
2|\tau|^4\xi_2^2b_0k\delta_2(k)b_0\delta_2(k)b_0k+
6|\tau|^4\xi_2^2b_0k\delta_2(k)b_0k\delta_2(k)b_0
-\\8\xi_1^4b_0^3k^4\delta_1(k)^2b_0
-4\xi_1^4b_0^3k^4\delta_1^2(k)b_0k
-\\4\xi_1^4b_0^3k^5\delta_1^2(k)b_0-
6\xi_1^4b_0^2k^2\delta_1(k)b_0k\delta_1(k)b_0k
-\\10\xi_1^4b_0^2k^2\delta_1(k)b_0k^2\delta_1(k)b_0
-10\xi_1^4b_0^2k^3\delta_1(k)b_0\delta_1(k)b_0k-\\
14\xi_1^4b_0^2k^3\delta_1(k)b_0k\delta_1(k)b_0
-4\xi_1^4b_0k\delta_1(k)b_0^2k^2\delta_1(k)b_0k
-\\4\xi_1^4b_0k\delta_1(k)b_0^2k^3\delta_1(k)b_0-
8\tau_1\xi_1^4b_0^3k^4\delta_1(k)\delta_2(k)b_0
-\\8\tau_1\xi_1^4b_0^3k^4\delta_2(k)\delta_1(k)b_0
-8\tau_1\xi_1^4b_0^3k^4\delta_1\delta_2(k)b_0k-\\
8\tau_1\xi_1^4b_0^3k^5\delta_1\delta_2(k)b_0-
6\tau_1\xi_1^4b_0^2k^2\delta_1(k)b_0k\delta_2(k)b_0k
-\\10\tau_1\xi_1^4b_0^2k^2\delta_1(k)b_0k^2\delta_2(k)b_0-
6\tau_1\xi_1^4b_0^2k^2\delta_2(k)b_0k\delta_1(k)b_0k
-\\10\tau_1\xi_1^4b_0^2k^2\delta_2(k)b_0k^2\delta_1(k)b_0
-10\tau_1\xi_1^4b_0^2k^3\delta_1(k)b_0\delta_2(k)b_0k-\\
14\tau_1\xi_1^4b_0^2k^3\delta_1(k)b_0k\delta_2(k)b_0
-10\tau_1\xi_1^4b_0^2k^3\delta_2(k)b_0\delta_1(k)b_0k
-\\14\tau_1\xi_1^4b_0^2k^3\delta_2(k)b_0k\delta_1(k)b_0-
4\tau_1\xi_1^4b_0k\delta_1(k)b_0^2k^2\delta_2(k)b_0k
-\\4\tau_1\xi_1^4b_0k\delta_1(k)b_0^2k^3\delta_2(k)b_0
-4\tau_1\xi_1^4b_0k\delta_2(k)b_0^2k^2\delta_1(k)b_0k-\\
4\tau_1\xi_1^4b_0k\delta_2(k)b_0^2k^3\delta_1(k)b_0-
32\tau_1\xi_1^3\xi_2b_0^3k^4\delta_1(k)^2b_0-\\
16\tau_1\xi_1^3\xi_2b_0^3k^4\delta_1^2(k)b_0k-
16\tau_1\xi_1^3\xi_2b_0^3k^5\delta_1^2(k)b_0
-\\24\tau_1\xi_1^3\xi_2b_0^2k^2\delta_1(k)b_0k\delta_1(k)b_0k-
40\tau_1\xi_1^3\xi_2b_0^2k^2\delta_1(k)b_0k^2\delta_1(k)b_0-\\
40\tau_1\xi_1^3\xi_2b_0^2k^3\delta_1(k)b_0\delta_1(k)b_0k-
56\tau_1\xi_1^3\xi_2b_0^2k^3\delta_1(k)b_0k\delta_1(k)b_0-\\
16\tau_1\xi_1^3\xi_2b_0k\delta_1(k)b_0^2k^2\delta_1(k)b_0k-
16\tau_1\xi_1^3\xi_2b_0k\delta_1(k)b_0^2k^3\delta_1(k)b_0
-\\8\tau_1^2\xi_1^4b_0^3k^4\delta_2(k)^2b_0-
4\tau_1^2\xi_1^4b_0^3k^4\delta_2^2(k)b_0k-\\
4\tau_1^2\xi_1^4b_0^3k^5\delta_2^2(k)b_0-
4\tau_1^2\xi_1^4b_0^2k^2\delta_2(k)b_0k\delta_2(k)b_0k-\\
8\tau_1^2\xi_1^4b_0^2k^2\delta_2(k)b_0k^2\delta_2(k)b_0-
8\tau_1^2\xi_1^4b_0^2k^3\delta_2(k)b_0\delta_2(k)b_0k-\\
12\tau_1^2\xi_1^4b_0^2k^3\delta_2(k)b_0k\delta_2(k)b_0
-4\tau_1^2\xi_1^4b_0k\delta_2(k)b_0^2k^2\delta_2(k)b_0k-\\
4\tau_1^2\xi_1^4b_0k\delta_2(k)b_0^2k^3\delta_2(k)b_0
-24\tau_1^2\xi_1^3\xi_2b_0^3k^4\delta_1(k)\delta_2(k)b_0
-\\24\tau_1^2\xi_1^3\xi_2b_0^3k^4\delta_2(k)\delta_1(k)b_0-
24\tau_1^2\xi_1^3\xi_2b_0^3k^4\delta_1\delta_2(k)b_0k-\\
24\tau_1^2\xi_1^3\xi_2b_0^3k^5\delta_1\delta_2(k)b_0
-20\tau_1^2\xi_1^3\xi_2b_0^2k^2\delta_1(k)b_0k\delta_2(k)b_0k-\\
32\tau_1^2\xi_1^3\xi_2b_0^2k^2\delta_1(k)b_0k^2\delta_2(k)b_0
-20\tau_1^2\xi_1^3\xi_2b_0^2k^2\delta_2(k)b_0k\delta_1(k)b_0k
-\\32\tau_1^2\xi_1^3\xi_2b_0^2k^2\delta_2(k)b_0k^2\delta_1(k)b_0-
32\tau_1^2\xi_1^3\xi_2b_0^2k^3\delta_1(k)b_0\delta_2(k)b_0k-\\
44\tau_1^2\xi_1^3\xi_2b_0^2k^3\delta_1(k)b_0k\delta_2(k)b_0
-32\tau_1^2\xi_1^3\xi_2b_0^2k^3\delta_2(k)b_0\delta_1(k)b_0k-\\
44\tau_1^2\xi_1^3\xi_2b_0^2k^3\delta_2(k)b_0k\delta_1(k)b_0
-12\tau_1^2\xi_1^3\xi_2b_0k\delta_1(k)b_0^2k^2\delta_2(k)b_0k
-\\12\tau_1^2\xi_1^3\xi_2b_0k\delta_1(k)b_0^2k^3\delta_2(k)b_0-
12\tau_1^2\xi_1^3\xi_2b_0k\delta_2(k)b_0^2k^2\delta_1(k)b_0k
-\\12\tau_1^2\xi_1^3\xi_2b_0k\delta_2(k)b_0^2k^3\delta_1(k)b_0-
40\tau_1^2\xi_1^2\xi_2^2b_0^3k^4\delta_1(k)^2b_0-\\
20\tau_1^2\xi_1^2\xi_2^2b_0^3k^4\delta_1^2(k)b_0k-
20\tau_1^2\xi_1^2\xi_2^2b_0^3k^5\delta_1^2(k)b_0-\\
28\tau_1^2\xi_1^2\xi_2^2b_0^2k^2\delta_1(k)b_0k\delta_1(k)b_0k-
48\tau_1^2\xi_1^2\xi_2^2b_0^2k^2\delta_1(k)b_0k^2\delta_1(k)b_0-\\
48\tau_1^2\xi_1^2\xi_2^2b_0^2k^3\delta_1(k)b_0\delta_1(k)b_0k-
68\tau_1^2\xi_1^2\xi_2^2b_0^2k^3\delta_1(k)b_0k\delta_1(k)b_0-\\
20\tau_1^2\xi_1^2\xi_2^2b_0k\delta_1(k)b_0^2k^2\delta_1(k)b_0k
-20\tau_1^2\xi_1^2\xi_2^2b_0k\delta_1(k)b_0^2k^3\delta_1(k)b_0
-\\2|\tau|^2\xi_1^4b_0^2k^2\delta_2(k)b_0k\delta_2(k)b_0k-
2|\tau|^2\xi_1^4b_0^2k^2\delta_2(k)b_0k^2\delta_2(k)b_0
-\\2|\tau|^2\xi_1^4b_0^2k^3\delta_2(k)b_0\delta_2(k)b_0k
-2|\tau|^2\xi_1^4b_0^2k^3\delta_2(k)b_0k\delta_2(k)b_0-\\
8|\tau|^2\xi_1^3\xi_2b_0^3k^4\delta_1(k)\delta_2(k)b_0-8
|\tau|^2\xi_1^3\xi_2b_0^3k^4\delta_2(k)\delta_1(k)b_0-\\
8|\tau|^2\xi_1^3\xi_2b_0^3k^4\delta_1\delta_2(k)b_0k-
8|\tau|^2\xi_1^3\xi_2b_0^3k^5\delta_1\delta_2(k)b_0-\\4
|\tau|^2\xi_1^3\xi_2b_0^2k^2\delta_1(k)b_0k\delta_2(k)b_0k-8
|\tau|^2\xi_1^3\xi_2b_0^2k^2\delta_1(k)b_0k^2\delta_2(k)b_0-\\
4|\tau|^2\xi_1^3\xi_2b_0^2k^2\delta_2(k)b_0k\delta_1(k)b_0k-
8|\tau|^2\xi_1^3\xi_2b_0^2k^2\delta_2(k)b_0k^2\delta_1(k)b_0-\\8
|\tau|^2\xi_1^3\xi_2b_0^2k^3\delta_1(k)b_0\delta_2(k)b_0k-
12|\tau|^2\xi_1^3\xi_2b_0^2k^3\delta_1(k)b_0k\delta_2(k)b_0
-\\8|\tau|^2\xi_1^3\xi_2b_0^2k^3\delta_2(k)b_0\delta_1(k)b_0k-
12|\tau|^2\xi_1^3\xi_2b_0^2k^3\delta_2(k)b_0k\delta_1(k)b_0-\\
4|\tau|^2\xi_1^3\xi_2b_0k\delta_1(k)b_0^2k^2\delta_2(k)b_0k-4
|\tau|^2\xi_1^3\xi_2b_0k\delta_1(k)b_0^2k^3\delta_2(k)b_0-\\
4|\tau|^2\xi_1^3\xi_2b_0k\delta_2(k)b_0^2k^2\delta_1(k)b_0k-
4|\tau|^2\xi_1^3\xi_2b_0k\delta_2(k)b_0^2k^3\delta_1(k)b_0-\\8
|\tau|^2\xi_1^2\xi_2^2b_0^3k^4\delta_1(k)^2b_0-
4|\tau|^2\xi_1^2\xi_2^2b_0^3k^4\delta_1^2(k)b_0k-\\
4|\tau|^2\xi_1^2\xi_2^2b_0^3k^5\delta_1^2(k)b_0
-8|\tau|^2\xi_1^2\xi_2^2b_0^2k^2\delta_1(k)b_0k\delta_1(k)b_0k-\\
12|\tau|^2\xi_1^2\xi_2^2b_0^2k^2\delta_1(k)b_0k^2\delta_1(k)b_0-12
|\tau|^2\xi_1^2\xi_2^2b_0^2k^3\delta_1(k)b_0\delta_1(k)b_0k-\\
16|\tau|^2\xi_1^2\xi_2^2b_0^2k^3\delta_1(k)b_0k\delta_1(k)b_0-
4|\tau|^2\xi_1^2\xi_2^2b_0k\delta_1(k)b_0^2k^2\delta_1(k)b_0k
-\\4|\tau|^2\xi_1^2\xi_2^2b_0k\delta_1(k)b_0^2k^3\delta_1(k)b_0-
16\tau_1^3\xi_1^3\xi_2b_0^3k^4\delta_2(k)^2b_0-\\
8\tau_1^3\xi_1^3\xi_2b_0^3k^4\delta_2^2(k)b_0k
-8\tau_1^3\xi_1^3\xi_2b_0^3k^5\delta_2^2(k)b_0-\\
8\tau_1^3\xi_1^3\xi_2b_0^2k^2\delta_2(k)b_0k\delta_2(k)b_0k-
16\tau_1^3\xi_1^3\xi_2b_0^2k^2\delta_2(k)b_0k^2\delta_2(k)b_0
-\\16\tau_1^3\xi_1^3\xi_2b_0^2k^3\delta_2(k)b_0\delta_2(k)b_0k-
24\tau_1^3\xi_1^3\xi_2b_0^2k^3\delta_2(k)b_0k\delta_2(k)b_0-\\
8\tau_1^3\xi_1^3\xi_2b_0k\delta_2(k)b_0^2k^2\delta_2(k)b_0k-8
\tau_1^3\xi_1^3\xi_2b_0k\delta_2(k)b_0^2k^3\delta_2(k)b_0-\\16
\tau_1^3\xi_1^2\xi_2^2b_0^3k^4\delta_1(k)\delta_2(k)b_0-
16\tau_1^3\xi_1^2\xi_2^2b_0^3k^4\delta_2(k)\delta_1(k)b_0
-\\16\tau_1^3\xi_1^2\xi_2^2b_0^3k^4\delta_1\delta_2(k)b_0k
-16\tau_1^3\xi_1^2\xi_2^2b_0^3k^5\delta_1\delta_2(k)b_0-\\
16\tau_1^3\xi_1^2\xi_2^2b_0^2k^2\delta_1(k)b_0k\delta_2(k)b_0k-
24\tau_1^3\xi_1^2\xi_2^2b_0^2k^2\delta_1(k)b_0k^2\delta_2(k)b_0-\\
16\tau_1^3\xi_1^2\xi_2^2b_0^2k^2\delta_2(k)b_0k\delta_1(k)b_0k-
24\tau_1^3\xi_1^2\xi_2^2b_0^2k^2\delta_2(k)b_0k^2\delta_1(k)b_0-\\
24\tau_1^3\xi_1^2\xi_2^2b_0^2k^3\delta_1(k)b_0\delta_2(k)b_0k-
32\tau_1^3\xi_1^2\xi_2^2b_0^2k^3\delta_1(k)b_0k\delta_2(k)b_0-\\
24\tau_1^3\xi_1^2\xi_2^2b_0^2k^3\delta_2(k)b_0\delta_1(k)b_0k-32
\tau_1^3\xi_1^2\xi_2^2b_0^2k^3\delta_2(k)b_0k\delta_1(k)b_0-\\8
\tau_1^3\xi_1^2\xi_2^2b_0k\delta_1(k)b_0^2k^2\delta_2(k)b_0k-
8\tau_1^3\xi_1^2\xi_2^2b_0k\delta_1(k)b_0^2k^3\delta_2(k)b_0
-\\8\tau_1^3\xi_1^2\xi_2^2b_0k\delta_2(k)b_0^2k^2\delta_1(k)b_0k
-8\tau_1^3\xi_1^2\xi_2^2b_0k\delta_2(k)b_0^2k^3\delta_1(k)b_0
-\\
16\tau_1^3\xi_1\xi_2^3b_0^3k^4\delta_1(k)^2b_0-
8\tau_1^3\xi_1\xi_2^3b_0^3k^4\delta_1^2(k)b_0k
-\\8\tau_1^3\xi_1\xi_2^3b_0^3k^5\delta_1^2(k)b_0-
8\tau_1^3\xi_1\xi_2^3b_0^2k^2\delta_1(k)b_0k\delta_1(k)b_0k
-\\16\tau_1^3\xi_1\xi_2^3b_0^2k^2\delta_1(k)b_0k^2\delta_1(k)b_0
-16\tau_1^3\xi_1\xi_2^3b_0^2k^3\delta_1(k)b_0\delta_1(k)b_0k-\\
24\tau_1^3\xi_1\xi_2^3b_0^2k^3\delta_1(k)b_0k\delta_1(k)b_0
-8\tau_1^3\xi_1\xi_2^3b_0k\delta_1(k)b_0^2k^2\delta_1(k)b_0k-\\
8\tau_1^3\xi_1\xi_2^3b_0k\delta_1(k)b_0^2k^3\delta_1(k)b_0-
16\tau_1|\tau|^2\xi_1^3\xi_2b_0^3k^4\delta_2(k)^2b_0-\\
8\tau_1|\tau|^2\xi_1^3\xi_2b_0^3k^4\delta_2^2(k)b_0k-8
\tau_1|\tau|^2\xi_1^3\xi_2b_0^3k^5\delta_2^2(k)b_0-\\
16\tau_1|\tau|^2\xi_1^3\xi_2b_0^2k^2\delta_2(k)b_0k\delta_2(k)b_0k
-24\tau_1|\tau|^2\xi_1^3\xi_2b_0^2k^2\delta_2(k)b_0k^2\delta_2(k)b_0
-\\24\tau_1|\tau|^2\xi_1^3\xi_2b_0^2k^3\delta_2(k)b_0\delta_2(k)b_0k-
32\tau_1|\tau|^2\xi_1^3\xi_2b_0^2k^3\delta_2(k)b_0k\delta_2(k)b_0-\\
8\tau_1|\tau|^2\xi_1^3\xi_2b_0k\delta_2(k)b_0^2k^2\delta_2(k)b_0k-
8\tau_1|\tau|^2\xi_1^3\xi_2b_0k\delta_2(k)b_0^2k^3\delta_2(k)b_0-\\
32\tau_1|\tau|^2\xi_1^2\xi_2^2b_0^3k^4\delta_1(k)\delta_2(k)b_0
-32\tau_1|\tau|^2\xi_1^2\xi_2^2b_0^3k^4\delta_2(k)\delta_1(k)b_0-\\
32\tau_1|\tau|^2\xi_1^2\xi_2^2b_0^3k^4\delta_1\delta_2(k)b_0k-
32\tau_1|\tau|^2\xi_1^2\xi_2^2b_0^3k^5\delta_1\delta_2(k)b_0-\\
20\tau_1|\tau|^2\xi_1^2\xi_2^2b_0^2k^2\delta_1(k)b_0k\delta_2(k)b_0k-
36\tau_1|\tau|^2\xi_1^2\xi_2^2b_0^2k^2\delta_1(k)b_0k^2\delta_2(k)b_0-\\
20\tau_1|\tau|^2\xi_1^2\xi_2^2b_0^2k^2\delta_2(k)b_0k\delta_1(k)b_0k-
36\tau_1|\tau|^2\xi_1^2\xi_2^2b_0^2k^2\delta_2(k)b_0k^2\delta_1(k)b_0-\\
36\tau_1|\tau|^2\xi_1^2\xi_2^2b_0^2k^3\delta_1(k)b_0\delta_2(k)b_0k-
52\tau_1|\tau|^2\xi_1^2\xi_2^2b_0^2k^3\delta_1(k)b_0k\delta_2(k)b_0-\\
36\tau_1|\tau|^2\xi_1^2\xi_2^2b_0^2k^3\delta_2(k)b_0\delta_1(k)b_0k-
52\tau_1|\tau|^2\xi_1^2\xi_2^2b_0^2k^3\delta_2(k)b_0k\delta_1(k)b_0-\\
16\tau_1|\tau|^2\xi_1^2\xi_2^2b_0k\delta_1(k)b_0^2k^2\delta_2(k)b_0k-
16\tau_1|\tau|^2\xi_1^2\xi_2^2b_0k\delta_1(k)b_0^2k^3\delta_2(k)b_0-\\
16\tau_1|\tau|^2\xi_1^2\xi_2^2b_0k\delta_2(k)b_0^2k^2\delta_1(k)b_0k-
16\tau_1|\tau|^2\xi_1^2\xi_2^2b_0k\delta_2(k)b_0^2k^3\delta_1(k)b_0
-\\16\tau_1|\tau|^2\xi_1\xi_2^3b_0^3k^4\delta_1(k)^2b_0-
8\tau_1|\tau|^2\xi_1\xi_2^3b_0^3k^4\delta_1^2(k)b_0k-\\
8\tau_1|\tau|^2\xi_1\xi_2^3b_0^3k^5\delta_1^2(k)b_0-
16\tau_1|\tau|^2\xi_1\xi_2^3b_0^2k^2\delta_1(k)b_0k\delta_1(k)b_0k
-\\24\tau_1|\tau|^2\xi_1\xi_2^3b_0^2k^2\delta_1(k)b_0k^2\delta_1(k)b_0-
24\tau_1|\tau|^2\xi_1\xi_2^3b_0^2k^3\delta_1(k)b_0\delta_1(k)b_0k-\\32
\tau_1|\tau|^2\xi_1\xi_2^3b_0^2k^3\delta_1(k)b_0k\delta_1(k)b_0-
8\tau_1|\tau|^2\xi_1\xi_2^3b_0k\delta_1(k)b_0^2k^2\delta_1(k)b_0k-\\
8\tau_1|\tau|^2\xi_1\xi_2^3b_0k\delta_1(k)b_0^2k^3\delta_1(k)b_0-
40\tau_1^2|\tau|^2\xi_1^2\xi_2^2b_0^3k^4\delta_2(k)^2b_0-\\
20\tau_1^2|\tau|^2\xi_1^2\xi_2^2b_0^3k^4\delta_2^2(k)b_0k-
20\tau_1^2|\tau|^2\xi_1^2\xi_2^2b_0^3k^5\delta_2^2(k)b_0-\\
28\tau_1^2|\tau|^2\xi_1^2\xi_2^2b_0^2k^2\delta_2(k)b_0k\delta_2(k)b_0k-
48\tau_1^2|\tau|^2\xi_1^2\xi_2^2b_0^2k^2\delta_2(k)b_0k^2\delta_2(k)b_0-\\
48\tau_1^2|\tau|^2\xi_1^2\xi_2^2b_0^2k^3\delta_2(k)b_0\delta_2(k)b_0k-
68\tau_1^2|\tau|^2\xi_1^2\xi_2^2b_0^2k^3\delta_2(k)b_0k\delta_2(k)b_0-\\
20\tau_1^2|\tau|^2\xi_1^2\xi_2^2b_0k\delta_2(k)b_0^2k^2\delta_2(k)b_0k-
20\tau_1^2|\tau|^2\xi_1^2\xi_2^2b_0k\delta_2(k)b_0^2k^3\delta_2(k)b_0-\\
24\tau_1^2|\tau|^2\xi_1\xi_2^3b_0^3k^4\delta_1(k)\delta_2(k)b_0-
24\tau_1^2|\tau|^2\xi_1\xi_2^3b_0^3k^4\delta_2(k)\delta_1(k)b_0
-\\24\tau_1^2|\tau|^2\xi_1\xi_2^3b_0^3k^4\delta_1\delta_2(k)b_0k-
24\tau_1^2|\tau|^2\xi_1\xi_2^3b_0^3k^5\delta_1\delta_2(k)b_0
-\\20\tau_1^2|\tau|^2\xi_1\xi_2^3b_0^2k^2\delta_1(k)b_0k\delta_2(k)b_0k-
32\tau_1^2|\tau|^2\xi_1\xi_2^3b_0^2k^2\delta_1(k)b_0k^2\delta_2(k)b_0
-\\20\tau_1^2|\tau|^2\xi_1\xi_2^3b_0^2k^2\delta_2(k)b_0k\delta_1(k)b_0k-
32\tau_1^2|\tau|^2\xi_1\xi_2^3b_0^2k^2\delta_2(k)b_0k^2\delta_1(k)b_0
-\\32\tau_1^2|\tau|^2\xi_1\xi_2^3b_0^2k^3\delta_1(k)b_0\delta_2(k)b_0k-
44\tau_1^2|\tau|^2\xi_1\xi_2^3b_0^2k^3\delta_1(k)b_0k\delta_2(k)b_0
-\\32\tau_1^2|\tau|^2\xi_1\xi_2^3b_0^2k^3\delta_2(k)b_0\delta_1(k)b_0k-
44\tau_1^2|\tau|^2\xi_1\xi_2^3b_0^2k^3\delta_2(k)b_0k\delta_1(k)b_0
-\\12\tau_1^2|\tau|^2\xi_1\xi_2^3b_0k\delta_1(k)b_0^2k^2\delta_2(k)b_0k-
12\tau_1^2|\tau|^2\xi_1\xi_2^3b_0k\delta_1(k)b_0^2k^3\delta_2(k)b_0-\\
12\tau_1^2|\tau|^2\xi_1\xi_2^3b_0k\delta_2(k)b_0^2k^2\delta_1(k)b_0k-
12\tau_1^2|\tau|^2\xi_1\xi_2^3b_0k\delta_2(k)b_0^2k^3\delta_1(k)b_0-\\
8\tau_1^2|\tau|^2\xi_2^4b_0^3k^4\delta_1(k)^2b_0
-4\tau_1^2|\tau|^2\xi_2^4b_0^3k^4\delta_1^2(k)b_0k-\\
4\tau_1^2|\tau|^2\xi_2^4b_0^3k^5\delta_1^2(k)b_0-
4\tau_1^2|\tau|^2\xi_2^4b_0^2k^2\delta_1(k)b_0k\delta_1(k)b_0k-\\
8\tau_1^2|\tau|^2\xi_2^4b_0^2k^2\delta_1(k)b_0k^2\delta_1(k)b_0
-8\tau_1^2|\tau|^2\xi_2^4b_0^2k^3\delta_1(k)b_0\delta_1(k)b_0k-\\
12\tau_1^2|\tau|^2\xi_2^4b_0^2k^3\delta_1(k)b_0k\delta_1(k)b_0-
4\tau_1^2|\tau|^2\xi_2^4b_0k\delta_1(k)b_0^2k^2\delta_1(k)b_0k-\\
4\tau_1^2|\tau|^2\xi_2^4b_0k\delta_1(k)b_0^2k^3\delta_1(k)b_0-
8|\tau|^4\xi_1^2\xi_2^2b_0^3k^4\delta_2(k)^2b_0-\\
4|\tau|^4\xi_1^2\xi_2^2b_0^3k^4\delta_2^2(k)b_0k-
4|\tau|^4\xi_1^2\xi_2^2b_0^3k^5\delta_2^2(k)b_0-\\
8|\tau|^4\xi_1^2\xi_2^2b_0^2k^2\delta_2(k)b_0k\delta_2(k)b_0k
-12|\tau|^4\xi_1^2\xi_2^2b_0^2k^2\delta_2(k)b_0k^2\delta_2(k)b_0-\\
12|\tau|^4\xi_1^2\xi_2^2b_0^2k^3\delta_2(k)b_0\delta_2(k)b_0k-
16|\tau|^4\xi_1^2\xi_2^2b_0^2k^3\delta_2(k)b_0k\delta_2(k)b_0-\\
4|\tau|^4\xi_1^2\xi_2^2b_0k\delta_2(k)b_0^2k^2\delta_2(k)b_0k-
4|\tau|^4\xi_1^2\xi_2^2b_0k\delta_2(k)b_0^2k^3\delta_2(k)b_0-\\
8|\tau|^4\xi_1\xi_2^3b_0^3k^4\delta_1(k)\delta_2(k)b_0-
8|\tau|^4\xi_1\xi_2^3b_0^3k^4\delta_2(k)\delta_1(k)b_0-\\
8|\tau|^4\xi_1\xi_2^3b_0^3k^4\delta_1\delta_2(k)b_0k-
8|\tau|^4\xi_1\xi_2^3b_0^3k^5\delta_1\delta_2(k)b_0-\\
4|\tau|^4\xi_1\xi_2^3b_0^2k^2\delta_1(k)b_0k\delta_2(k)b_0k-
8|\tau|^4\xi_1\xi_2^3b_0^2k^2\delta_1(k)b_0k^2\delta_2(k)b_0-\\
4|\tau|^4\xi_1\xi_2^3b_0^2k^2\delta_2(k)b_0k\delta_1(k)b_0k-
8|\tau|^4\xi_1\xi_2^3b_0^2k^2\delta_2(k)b_0k^2\delta_1(k)b_0-\\
8|\tau|^4\xi_1\xi_2^3b_0^2k^3\delta_1(k)b_0\delta_2(k)b_0k
-12|\tau|^4\xi_1\xi_2^3b_0^2k^3\delta_1(k)b_0k\delta_2(k)b_0-\\
8|\tau|^4\xi_1\xi_2^3b_0^2k^3\delta_2(k)b_0\delta_1(k)b_0k-
12|\tau|^4\xi_1\xi_2^3b_0^2k^3\delta_2(k)b_0k\delta_1(k)b_0-\\
4|\tau|^4\xi_1\xi_2^3b_0k\delta_1(k)b_0^2k^2\delta_2(k)b_0k
-4|\tau|^4\xi_1\xi_2^3b_0k\delta_1(k)b_0^2k^3\delta_2(k)b_0-\\
4|\tau|^4\xi_1\xi_2^3b_0k\delta_2(k)b_0^2k^2\delta_1(k)b_0k-
4|\tau|^4\xi_1\xi_2^3b_0k\delta_2(k)b_0^2k^3\delta_1(k)b_0-\\
2|\tau|^4\xi_2^4b_0^2k^2\delta_1(k)b_0k\delta_1(k)b_0k
-2|\tau|^4\xi_2^4b_0^2k^2\delta_1(k)b_0k^2\delta_1(k)b_0-\\
2|\tau|^4\xi_2^4b_0^2k^3\delta_1(k)b_0\delta_1(k)b_0k-
2|\tau|^4\xi_2^4b_0^2k^3\delta_1(k)b_0k\delta_1(k)b_0+\\
8\xi_1^6b_0^3k^4\delta_1(k)b_0k\delta_1(k)b_0k+
8\xi_1^6b_0^3k^4\delta_1(k)b_0k^2\delta_1(k)b_0+\\
8\xi_1^6b_0^3k^5\delta_1(k)b_0\delta_1(k)b_0k+
8\xi_1^6b_0^3k^5\delta_1(k)b_0k\delta_1(k)b_0+\\
4\xi_1^6b_0^2k^2\delta_1(k)b_0^2k^3\delta_1(k)b_0k+
4\xi_1^6b_0^2k^2\delta_1(k)b_0^2k^4\delta_1(k)b_0+\\
4\xi_1^6b_0^2k^3\delta_1(k)b_0^2k^2\delta_1(k)b_0k+
4\xi_1^6b_0^2k^3\delta_1(k)b_0^2k^3\delta_1(k)b_0-\\
32\tau_1|\tau|^4\xi_1\xi_2^3b_0^3k^4\delta_2(k)^2b_0
-16\tau_1|\tau|^4\xi_1\xi_2^3b_0^3k^4\delta_2^2(k)b_0k-\\
16\tau_1|\tau|^4\xi_1\xi_2^3b_0^3k^5\delta_2^2(k)b_0-24
\tau_1|\tau|^4\xi_1\xi_2^3b_0^2k^2\delta_2(k)b_0k\delta_2(k)b_0k-\\
40\tau_1|\tau|^4\xi_1\xi_2^3b_0^2k^2\delta_2(k)b_0k^2\delta_2(k)b_0-
40\tau_1|\tau|^4\xi_1\xi_2^3b_0^2k^3\delta_2(k)b_0\delta_2(k)b_0k-\\
56\tau_1|\tau|^4\xi_1\xi_2^3b_0^2k^3\delta_2(k)b_0k\delta_2(k)b_0
-16\tau_1|\tau|^4\xi_1\xi_2^3b_0k\delta_2(k)b_0^2k^2\delta_2(k)b_0k-\\
16\tau_1|\tau|^4\xi_1\xi_2^3b_0k\delta_2(k)b_0^2k^3\delta_2(k)b_0-
8\tau_1|\tau|^4\xi_2^4b_0^3k^4\delta_1(k)\delta_2(k)b_0-\\
8\tau_1|\tau|^4\xi_2^4b_0^3k^4\delta_2(k)\delta_1(k)b_0-
8\tau_1|\tau|^4\xi_2^4b_0^3k^4\delta_1\delta_2(k)b_0k-\\
8\tau_1|\tau|^4\xi_2^4b_0^3k^5\delta_1\delta_2(k)b_0-
6\tau_1|\tau|^4\xi_2^4b_0^2k^2\delta_1(k)b_0k\delta_2(k)b_0k-\\
10\tau_1|\tau|^4\xi_2^4b_0^2k^2\delta_1(k)b_0k^2\delta_2(k)b_0
-6\tau_1|\tau|^4\xi_2^4b_0^2k^2\delta_2(k)b_0k\delta_1(k)b_0k-\\
10\tau_1|\tau|^4\xi_2^4b_0^2k^2\delta_2(k)b_0k^2\delta_1(k)b_0-
10\tau_1|\tau|^4\xi_2^4b_0^2k^3\delta_1(k)b_0\delta_2(k)b_0k-\\
14\tau_1|\tau|^4\xi_2^4b_0^2k^3\delta_1(k)b_0k\delta_2(k)b_0-
10\tau_1|\tau|^4\xi_2^4b_0^2k^3\delta_2(k)b_0\delta_1(k)b_0k-\\
14\tau_1|\tau|^4\xi_2^4b_0^2k^3\delta_2(k)b_0k\delta_1(k)b_0-
4\tau_1|\tau|^4\xi_2^4b_0k\delta_1(k)b_0^2k^2\delta_2(k)b_0k-\\
4\tau_1|\tau|^4\xi_2^4b_0k\delta_1(k)b_0^2k^3\delta_2(k)b_0-
4\tau_1|\tau|^4\xi_2^4b_0k\delta_2(k)b_0^2k^2\delta_1(k)b_0k-\\
4\tau_1|\tau|^4\xi_2^4b_0k\delta_2(k)b_0^2k^3\delta_1(k)b_0+
8\tau_1\xi_1^6b_0^3k^4\delta_1(k)b_0k\delta_2(k)b_0k+\\
8\tau_1\xi_1^6b_0^3k^4\delta_1(k)b_0k^2\delta_2(k)b_0+
8\tau_1\xi_1^6b_0^3k^4\delta_2(k)b_0k\delta_1(k)b_0k+\\
8\tau_1\xi_1^6b_0^3k^4\delta_2(k)b_0k^2\delta_1(k)b_0+
8\tau_1\xi_1^6b_0^3k^5\delta_1(k)b_0\delta_2(k)b_0k+\\
8\tau_1\xi_1^6b_0^3k^5\delta_1(k)b_0k\delta_2(k)b_0+
8\tau_1\xi_1^6b_0^3k^5\delta_2(k)b_0\delta_1(k)b_0k+\\
8\tau_1\xi_1^6b_0^3k^5\delta_2(k)b_0k\delta_1(k)b_0+
4\tau_1\xi_1^6b_0^2k^2\delta_1(k)b_0^2k^3\delta_2(k)b_0k+\\
4\tau_1\xi_1^6b_0^2k^2\delta_1(k)b_0^2k^4\delta_2(k)b_0+
4\tau_1\xi_1^6b_0^2k^2\delta_2(k)b_0^2k^3\delta_1(k)b_0k+\\
4\tau_1\xi_1^6b_0^2k^2\delta_2(k)b_0^2k^4\delta_1(k)b_0+
4\tau_1\xi_1^6b_0^2k^3\delta_1(k)b_0^2k^2\delta_2(k)b_0k+\\
4\tau_1\xi_1^6b_0^2k^3\delta_1(k)b_0^2k^3\delta_2(k)b_0+
4\tau_1\xi_1^6b_0^2k^3\delta_2(k)b_0^2k^2\delta_1(k)b_0k+\\
4\tau_1\xi_1^6b_0^2k^3\delta_2(k)b_0^2k^3\delta_1(k)b_0+
48\tau_1\xi_1^5\xi_2b_0^3k^4\delta_1(k)b_0k\delta_1(k)b_0k+\\
48\tau_1\xi_1^5\xi_2b_0^3k^4\delta_1(k)b_0k^2\delta_1(k)b_0+
48\tau_1\xi_1^5\xi_2b_0^3k^5\delta_1(k)b_0\delta_1(k)b_0k+\\
48\tau_1\xi_1^5\xi_2b_0^3k^5\delta_1(k)b_0k\delta_1(k)b_0+
24\tau_1\xi_1^5\xi_2b_0^2k^2\delta_1(k)b_0^2k^3\delta_1(k)b_0k+\\
24\tau_1\xi_1^5\xi_2b_0^2k^2\delta_1(k)b_0^2k^4\delta_1(k)b_0+
24\tau_1\xi_1^5\xi_2b_0^2k^3\delta_1(k)b_0^2k^2\delta_1(k)b_0k+\\
24\tau_1\xi_1^5\xi_2b_0^2k^3\delta_1(k)b_0^2k^3\delta_1(k)b_0-
8|\tau|^6\xi_2^4b_0^3k^4\delta_2(k)^2b_0-
4|\tau|^6\xi_2^4b_0^3k^4\delta_2^2(k)b_0k-\\
4|\tau|^6\xi_2^4b_0^3k^5\delta_2^2(k)b_0-
6|\tau|^6\xi_2^4b_0^2k^2\delta_2(k)b_0k\delta_2(k)b_0k-\\
10|\tau|^6\xi_2^4b_0^2k^2\delta_2(k)b_0k^2\delta_2(k)b_0-
10|\tau|^6\xi_2^4b_0^2k^3\delta_2(k)b_0\delta_2(k)b_0k-\\14|\tau|^
6\xi_2^4b_0^2k^3\delta_2(k)b_0k\delta_2(k)b_0-4|\tau|^6\xi_2^4b_0k
\delta_2(k)b_0^2k^2\delta_2(k)b_0k-\\4|\tau|^6\xi_2^4b_0k\delta_2(k)b_0^2k^
3\delta_2(k)b_0+8\tau_1^2\xi_1^6b_0^3k^4\delta_2(k)b_0k\delta_2(k)b_0k+
\\8\tau_1^2\xi_1^6b_0^3k^4\delta_2(k)b_0k^2\delta_2(k)b_0+
8\tau_1^2\xi_1^
6b_0^3k^5\delta_2(k)b_0\delta_2(k)b_0k+\\8\tau_1^2\xi_1^6b_0^3k^5\delta_2(k)b_0k\delta_2(k)b_0+4\tau_1^2\xi_1^6b_0^2k^2\delta_2(k)b_0^2k^
3\delta_2(k)b_0k+\\
4\tau_1^2\xi_1^6b_0^2k^2\delta_2(k)b_0^2k^4\delta_2(k)b_0+4\tau_1^2\xi_1^6b_0^2k^3\delta_2(k)b_0^2k^2\delta_2(k)b_0k+\\4\tau_1^
2\xi_1^6b_0^2k^3\delta_2(k)b_0^2k^3\delta_2(k)b_0+
40\tau_1^2\xi_1^5\xi_2b_0^3k^4\delta_1(k)b_0k\delta_2(k)b_0k+\\40\tau_1^2\xi_1^5\xi_2b_0^3k^
4\delta_1(k)b_0k^2\delta_2(k)b_0+40\tau_1^2\xi_1^5\xi_2b_0^3k^4\delta_2(k)b_0k\delta_1(k)b_0k+\\40\tau_1^2\xi_1^5\xi_2b_0^3k^4\delta_2(k)b_0k^
2\delta_1(k)b_0+40\tau_1^2\xi_1^5\xi_2b_0^3k^5\delta_1(k)b_0\delta_2(k)b_0k+\\40\tau_1^2\xi_1^5\xi_2b_0^3k^5\delta_1(k)b_0k\delta_2(k)b_0+40\tau_1^
2\xi_1^5\xi_2b_0^3k^5\delta_2(k)b_0\delta_1(k)b_0k+\\40\tau_1^2\xi_1^5\xi_2b_0^3k^5\delta_2(k)b_0k\delta_1(k)b_0+20\tau_1^2\xi_1^5\xi_2b_0^2k^
2\delta_1(k)b_0^2k^3\delta_2(k)b_0k+\\20\tau_1^2\xi_1^5\xi_2b_0^2k^2\delta_1(k)b_0^2k^4\delta_2(k)b_0+20\tau_1^2\xi_1^5\xi_2b_0^2k^2\delta_2(k)b_0^
2k^3\delta_1(k)b_0k+\\20\tau_1^2\xi_1^5\xi_2b_0^2k^2\delta_2(k)b_0^2k^4\delta_1(k)b_0+20\tau_1^2\xi_1^5\xi_2b_0^2k^3\delta_1(k)b_0^2k^
2\delta_2(k)b_0k+\\20\tau_1^2\xi_1^5\xi_2b_0^2k^3\delta_1(k)b_0^2k^3\delta_2(k)b_0+20\tau_1^2\xi_1^5\xi_2b_0^2k^3\delta_2(k)b_0^2k^2\delta_1(k)b_0k+\\
20\tau_1^2\xi_1^5\xi_2b_0^2k^3\delta_2(k)b_0^2k^3\delta_1(k)b_0+104\tau_1^2\xi_1^4\xi_2^2b_0^3k^4\delta_1(k)b_0k\delta_1(k)b_0k+\\104\tau_1^2\xi_1^
4\xi_2^2b_0^3k^4\delta_1(k)b_0k^2\delta_1(k)b_0+104\tau_1^2\xi_1^4\xi_2^2b_0^3k^5\delta_1(k)b_0\delta_1(k)b_0k+\\104\tau_1^2\xi_1^4\xi_2^2b_0^3k^
5\delta_1(k)b_0k\delta_1(k)b_0+52\tau_1^2\xi_1^4\xi_2^2b_0^2k^2\delta_1(k)b_0^2k^3\delta_1(k)b_0k+\\52\tau_1^2\xi_1^4\xi_2^2b_0^2k^2\delta_1(k)b_0^
2k^4\delta_1(k)b_0+52\tau_1^2\xi_1^4\xi_2^2b_0^2k^3\delta_1(k)b_0^2k^2\delta_1(k)b_0k+\\52\tau_1^2\xi_1^4\xi_2^2b_0^2k^3\delta_1(k)b_0^2k^
3\delta_1(k)b_0+8|\tau|^2\xi_1^5\xi_2b_0^3k^4\delta_1(k)b_0k\delta_2(k)b_0k+\\8|\tau|^2\xi_1^5\xi_2b_0^3k^4\delta_1(k)b_0k^2\delta_2(k)b_0+8|\tau|^2\xi_1^
5\xi_2b_0^3k^4\delta_2(k)b_0k\delta_1(k)b_0k+\\8|\tau|^2\xi_1^5\xi_2b_0^3k^4\delta_2(k)b_0k^2\delta_1(k)b_0+8|\tau|^2\xi_1^5\xi_2b_0^3k^
5\delta_1(k)b_0\delta_2(k)b_0k+\\8|\tau|^2\xi_1^5\xi_2b_0^3k^5\delta_1(k)b_0k\delta_2(k)b_0+8|\tau|^2\xi_1^5\xi_2b_0^3k^5\delta_2(k)b_0\delta_1(k)b_0k+\\8|\tau|^
2\xi_1^5\xi_2b_0^3k^5\delta_2(k)b_0k\delta_1(k)b_0+4|\tau|^2\xi_1^5\xi_2b_0^2k^2\delta_1(k)b_0^2k^3\delta_2(k)b_0k+\\4|\tau|^2\xi_1^5\xi_2b_0^2k^
2\delta_1(k)b_0^2k^4\delta_2(k)b_0+4|\tau|^2\xi_1^5\xi_2b_0^2k^2\delta_2(k)b_0^2k^3\delta_1(k)b_0k+\\4|\tau|^2\xi_1^5\xi_2b_0^2k^2\delta_2(k)b_0^2k^
4\delta_1(k)b_0+4|\tau|^2\xi_1^5\xi_2b_0^2k^3\delta_1(k)b_0^2k^2\delta_2(k)b_0k+\\4|\tau|^2\xi_1^5\xi_2b_0^2k^3\delta_1(k)b_0^2k^3\delta_2(k)b_0+4|\tau|^2\xi_1^
5\xi_2b_0^2k^3\delta_2(k)b_0^2k^2\delta_1(k)b_0k+\\4|\tau|^2\xi_1^5\xi_2b_0^2k^3\delta_2(k)b_0^2k^3\delta_1(k)b_0+16|\tau|^2\xi_1^4\xi_2^2b_0^3k^
4\delta_1(k)b_0k\delta_1(k)b_0k+\\16|\tau|^2\xi_1^4\xi_2^2b_0^3k^4\delta_1(k)b_0k^2\delta_1(k)b_0+16|\tau|^2\xi_1^4\xi_2^2b_0^3k^
5\delta_1(k)b_0\delta_1(k)b_0k+\\16|\tau|^2\xi_1^4\xi_2^2b_0^3k^5\delta_1(k)b_0k\delta_1(k)b_0+8|\tau|^2\xi_1^4\xi_2^2b_0^2k^2\delta_1(k)b_0^2k^
3\delta_1(k)b_0k+\\8|\tau|^2\xi_1^4\xi_2^2b_0^2k^2\delta_1(k)b_0^2k^4\delta_1(k)b_0+8|\tau|^2\xi_1^4\xi_2^2b_0^2k^3\delta_1(k)b_0^2k^2\delta_1(k)b_0k+\\8|\tau|^
2\xi_1^4\xi_2^2b_0^2k^3\delta_1(k)b_0^2k^3\delta_1(k)b_0+32\tau_1^3\xi_1^5\xi_2b_0^3k^4\delta_2(k)b_0k\delta_2(k)b_0k+\\32\tau_1^3\xi_1^5\xi_2b_0^
3k^4\delta_2(k)b_0k^2\delta_2(k)b_0+32\tau_1^3\xi_1^5\xi_2b_0^3k^5\delta_2(k)b_0\delta_2(k)b_0k+\\32\tau_1^3\xi_1^5\xi_2b_0^3k^
5\delta_2(k)b_0k\delta_2(k)b_0+16\tau_1^3\xi_1^5\xi_2b_0^2k^2\delta_2(k)b_0^2k^3\delta_2(k)b_0k+\\16\tau_1^3\xi_1^5\xi_2b_0^2k^2\delta_2(k)b_0^2k^
4\delta_2(k)b_0+16\tau_1^3\xi_1^5\xi_2b_0^2k^3\delta_2(k)b_0^2k^2\delta_2(k)b_0k+\\16\tau_1^3\xi_1^5\xi_2b_0^2k^3\delta_2(k)b_0^2k^3\delta_2(k)b_0+
64\tau_1^3\xi_1^4\xi_2^2b_0^3k^4\delta_1(k)b_0k\delta_2(k)b_0k+\\64\tau_1^3\xi_1^4\xi_2^2b_0^3k^4\delta_1(k)b_0k^2\delta_2(k)b_0+64\tau_1^3\xi_1^
4\xi_2^2b_0^3k^4\delta_2(k)b_0k\delta_1(k)b_0k+\\64\tau_1^3\xi_1^4\xi_2^2b_0^3k^4\delta_2(k)b_0k^2\delta_1(k)b_0+64\tau_1^3\xi_1^4\xi_2^2b_0^3k^
5\delta_1(k)b_0\delta_2(k)b_0k+\\64\tau_1^3\xi_1^4\xi_2^2b_0^3k^5\delta_1(k)b_0k\delta_2(k)b_0+64\tau_1^3\xi_1^4\xi_2^2b_0^3k^
5\delta_2(k)b_0\delta_1(k)b_0k+\\64\tau_1^3\xi_1^4\xi_2^2b_0^3k^5\delta_2(k)b_0k\delta_1(k)b_0+32\tau_1^3\xi_1^4\xi_2^2b_0^2k^2\delta_1(k)b_0^2k^
3\delta_2(k)b_0k+\\32\tau_1^3\xi_1^4\xi_2^2b_0^2k^2\delta_1(k)b_0^2k^4\delta_2(k)b_0+32\tau_1^3\xi_1^4\xi_2^2b_0^2k^2\delta_2(k)b_0^2k^
3\delta_1(k)b_0k+\\32\tau_1^3\xi_1^4\xi_2^2b_0^2k^2\delta_2(k)b_0^2k^4\delta_1(k)b_0+32\tau_1^3\xi_1^4\xi_2^2b_0^2k^3\delta_1(k)b_0^2k^
2\delta_2(k)b_0k+\\32\tau_1^3\xi_1^4\xi_2^2b_0^2k^3\delta_1(k)b_0^2k^3\delta_2(k)b_0+32\tau_1^3\xi_1^4\xi_2^2b_0^2k^3\delta_2(k)b_0^2k^
2\delta_1(k)b_0k+\\32\tau_1^3\xi_1^4\xi_2^2b_0^2k^3\delta_2(k)b_0^2k^3\delta_1(k)b_0+96\tau_1^3\xi_1^3\xi_2^3b_0^3k^4\delta_1(k)b_0k\delta_1(k)b_0k+\\
96\tau_1^3\xi_1^3\xi_2^3b_0^3k^4\delta_1(k)b_0k^2\delta_1(k)b_0+96\tau_1^3\xi_1^3\xi_2^3b_0^3k^5\delta_1(k)b_0\delta_1(k)b_0k+\\96\tau_1^3\xi_1^
3\xi_2^3b_0^3k^5\delta_1(k)b_0k\delta_1(k)b_0+48\tau_1^3\xi_1^3\xi_2^3b_0^2k^2\delta_1(k)b_0^2k^3\delta_1(k)b_0k+\\48\tau_1^3\xi_1^3\xi_2^3b_0^2k^
2\delta_1(k)b_0^2k^4\delta_1(k)b_0+48\tau_1^3\xi_1^3\xi_2^3b_0^2k^3\delta_1(k)b_0^2k^2\delta_1(k)b_0k+\\48\tau_1^3\xi_1^3\xi_2^3b_0^2k^
3\delta_1(k)b_0^2k^3\delta_1(k)b_0+16\tau_1|\tau|^2\xi_1^5\xi_2b_0^3k^4\delta_2(k)b_0k\delta_2(k)b_0k+\\16\tau_1|\tau|^2\xi_1^5\xi_2b_0^3k^4\delta_2(k)b_0k^
2\delta_2(k)b_0+16\tau_1|\tau|^2\xi_1^5\xi_2b_0^3k^5\delta_2(k)b_0\delta_2(k)b_0k+\\16\tau_1|\tau|^2\xi_1^5\xi_2b_0^3k^5\delta_2(k)b_0k\delta_2(k)b_0+
8\tau_1|\tau|^2\xi_1^5\xi_2b_0^2k^2\delta_2(k)b_0^2k^3\delta_2(k)b_0k+\\8\tau_1|\tau|^2\xi_1^5\xi_2b_0^2k^2\delta_2(k)b_0^2k^4\delta_2(k)b_0+8\tau_1|\tau|^
2\xi_1^5\xi_2b_0^2k^3\delta_2(k)b_0^2k^2\delta_2(k)b_0k+\\8\tau_1|\tau|^2\xi_1^5\xi_2b_0^2k^3\delta_2(k)b_0^2k^3\delta_2(k)b_0+56\tau_1|\tau|^2\xi_1^4\xi_2^
2b_0^3k^4\delta_1(k)b_0k\delta_2(k)b_0k+\\56\tau_1|\tau|^2\xi_1^4\xi_2^2b_0^3k^4\delta_1(k)b_0k^2\delta_2(k)b_0+56\tau_1|\tau|^2\xi_1^4\xi_2^2b_0^3k^
4\delta_2(k)b_0k\delta_1(k)b_0k+\\56\tau_1|\tau|^2\xi_1^4\xi_2^2b_0^3k^4\delta_2(k)b_0k^2\delta_1(k)b_0+56\tau_1|\tau|^2\xi_1^4\xi_2^2b_0^3k^
5\delta_1(k)b_0\delta_2(k)b_0k+\\56\tau_1|\tau|^2\xi_1^4\xi_2^2b_0^3k^5\delta_1(k)b_0k\delta_2(k)b_0+56\tau_1|\tau|^2\xi_1^4\xi_2^2b_0^3k^
5\delta_2(k)b_0\delta_1(k)b_0k+\\56\tau_1|\tau|^2\xi_1^4\xi_2^2b_0^3k^5\delta_2(k)b_0k\delta_1(k)b_0+28\tau_1|\tau|^2\xi_1^4\xi_2^2b_0^2k^2\delta_1(k)b_0^2k^
3\delta_2(k)b_0k+\\28\tau_1|\tau|^2\xi_1^4\xi_2^2b_0^2k^2\delta_1(k)b_0^2k^4\delta_2(k)b_0+28\tau_1|\tau|^2\xi_1^4\xi_2^2b_0^2k^2\delta_2(k)b_0^2k^
3\delta_1(k)b_0k+\\28\tau_1|\tau|^2\xi_1^4\xi_2^2b_0^2k^2\delta_2(k)b_0^2k^4\delta_1(k)b_0+28\tau_1|\tau|^2\xi_1^4\xi_2^2b_0^2k^3\delta_1(k)b_0^2k^
2\delta_2(k)b_0k+\\28\tau_1|\tau|^2\xi_1^4\xi_2^2b_0^2k^3\delta_1(k)b_0^2k^3\delta_2(k)b_0+28\tau_1|\tau|^2\xi_1^4\xi_2^2b_0^2k^3\delta_2(k)b_0^2k^
2\delta_1(k)b_0k+\\28\tau_1|\tau|^2\xi_1^4\xi_2^2b_0^2k^3\delta_2(k)b_0^2k^3\delta_1(k)b_0+64\tau_1|\tau|^2\xi_1^3\xi_2^3b_0^3k^
4\delta_1(k)b_0k\delta_1(k)b_0k+\\64\tau_1|\tau|^2\xi_1^3\xi_2^3b_0^3k^4\delta_1(k)b_0k^2\delta_1(k)b_0+64\tau_1|\tau|^2\xi_1^3\xi_2^3b_0^3k^
5\delta_1(k)b_0\delta_1(k)b_0k+\\64\tau_1|\tau|^2\xi_1^3\xi_2^3b_0^3k^5\delta_1(k)b_0k\delta_1(k)b_0+32\tau_1|\tau|^2\xi_1^3\xi_2^3b_0^2k^2\delta_1(k)b_0^2k^
3\delta_1(k)b_0k+\\32\tau_1|\tau|^2\xi_1^3\xi_2^3b_0^2k^2\delta_1(k)b_0^2k^4\delta_1(k)b_0+32\tau_1|\tau|^2\xi_1^3\xi_2^3b_0^2k^3\delta_1(k)b_0^2k^
2\delta_1(k)b_0k+\\32\tau_1|\tau|^2\xi_1^3\xi_2^3b_0^2k^3\delta_1(k)b_0^2k^3\delta_1(k)b_0+32\tau_1^4\xi_1^4\xi_2^2b_0^3k^
4\delta_2(k)b_0k\delta_2(k)b_0k+\\32\tau_1^4\xi_1^4\xi_2^2b_0^3k^4\delta_2(k)b_0k^2\delta_2(k)b_0+32\tau_1^4\xi_1^4\xi_2^2b_0^3k^
5\delta_2(k)b_0\delta_2(k)b_0k+\\32\tau_1^4\xi_1^4\xi_2^2b_0^3k^5\delta_2(k)b_0k\delta_2(k)b_0+16\tau_1^4\xi_1^4\xi_2^2b_0^2k^2\delta_2(k)b_0^2k^
3\delta_2(k)b_0k+\\16\tau_1^4\xi_1^4\xi_2^2b_0^2k^2\delta_2(k)b_0^2k^4\delta_2(k)b_0+16\tau_1^4\xi_1^4\xi_2^2b_0^2k^3\delta_2(k)b_0^2k^
2\delta_2(k)b_0k+\\16\tau_1^4\xi_1^4\xi_2^2b_0^2k^3\delta_2(k)b_0^2k^3\delta_2(k)b_0+32\tau_1^4\xi_1^3\xi_2^3b_0^3k^4\delta_1(k)b_0k\delta_2(k)b_0k+\\
32\tau_1^4\xi_1^3\xi_2^3b_0^3k^4\delta_1(k)b_0k^2\delta_2(k)b_0+32\tau_1^4\xi_1^3\xi_2^3b_0^3k^4\delta_2(k)b_0k\delta_1(k)b_0k+\\32\tau_1^4\xi_1^
3\xi_2^3b_0^3k^4\delta_2(k)b_0k^2\delta_1(k)b_0+32\tau_1^4\xi_1^3\xi_2^3b_0^3k^5\delta_1(k)b_0\delta_2(k)b_0k+\\32\tau_1^4\xi_1^3\xi_2^3b_0^3k^
5\delta_1(k)b_0k\delta_2(k)b_0+32\tau_1^4\xi_1^3\xi_2^3b_0^3k^5\delta_2(k)b_0\delta_1(k)b_0k+\\32\tau_1^4\xi_1^3\xi_2^3b_0^3k^
5\delta_2(k)b_0k\delta_1(k)b_0+16\tau_1^4\xi_1^3\xi_2^3b_0^2k^2\delta_1(k)b_0^2k^3\delta_2(k)b_0k+\\16\tau_1^4\xi_1^3\xi_2^3b_0^2k^2\delta_1(k)b_0^
2k^4\delta_2(k)b_0+16\tau_1^4\xi_1^3\xi_2^3b_0^2k^2\delta_2(k)b_0^2k^3\delta_1(k)b_0k+\\16\tau_1^4\xi_1^3\xi_2^3b_0^2k^2\delta_2(k)b_0^2k^
4\delta_1(k)b_0+16\tau_1^4\xi_1^3\xi_2^3b_0^2k^3\delta_1(k)b_0^2k^2\delta_2(k)b_0k+\\16\tau_1^4\xi_1^3\xi_2^3b_0^2k^3\delta_1(k)b_0^2k^
3\delta_2(k)b_0+16\tau_1^4\xi_1^3\xi_2^3b_0^2k^3\delta_2(k)b_0^2k^2\delta_1(k)b_0k+\\16\tau_1^4\xi_1^3\xi_2^3b_0^2k^3\delta_2(k)b_0^2k^
3\delta_1(k)b_0+32\tau_1^4\xi_1^2\xi_2^4b_0^3k^4\delta_1(k)b_0k\delta_1(k)b_0k+\\32\tau_1^4\xi_1^2\xi_2^4b_0^3k^4\delta_1(k)b_0k^2\delta_1(k)b_0+
32\tau_1^4\xi_1^2\xi_2^4b_0^3k^5\delta_1(k)b_0\delta_1(k)b_0k+\\32\tau_1^4\xi_1^2\xi_2^4b_0^3k^5\delta_1(k)b_0k\delta_1(k)b_0+16\tau_1^4\xi_1^2\xi_2^
4b_0^2k^2\delta_1(k)b_0^2k^3\delta_1(k)b_0k+\\16\tau_1^4\xi_1^2\xi_2^4b_0^2k^2\delta_1(k)b_0^2k^4\delta_1(k)b_0+16\tau_1^4\xi_1^2\xi_2^4b_0^2k^
3\delta_1(k)b_0^2k^2\delta_1(k)b_0k+\\16\tau_1^4\xi_1^2\xi_2^4b_0^2k^3\delta_1(k)b_0^2k^3\delta_1(k)b_0+80\tau_1^2|\tau|^2\xi_1^4\xi_2^2b_0^3k^
4\delta_2(k)b_0k\delta_2(k)b_0k+\\80\tau_1^2|\tau|^2\xi_1^4\xi_2^2b_0^3k^4\delta_2(k)b_0k^2\delta_2(k)b_0+80\tau_1^2|\tau|^2\xi_1^4\xi_2^2b_0^3k^
5\delta_2(k)b_0\delta_2(k)b_0k+\\80\tau_1^2|\tau|^2\xi_1^4\xi_2^2b_0^3k^5\delta_2(k)b_0k\delta_2(k)b_0+40\tau_1^2|\tau|^2\xi_1^4\xi_2^2b_0^2k^2\delta_2(k)b_0^
2k^3\delta_2(k)b_0k+\\40\tau_1^2|\tau|^2\xi_1^4\xi_2^2b_0^2k^2\delta_2(k)b_0^2k^4\delta_2(k)b_0+40\tau_1^2|\tau|^2\xi_1^4\xi_2^2b_0^2k^3\delta_2(k)b_0^2k^
2\delta_2(k)b_0k+\\40\tau_1^2|\tau|^2\xi_1^4\xi_2^2b_0^2k^3\delta_2(k)b_0^2k^3\delta_2(k)b_0+112\tau_1^2|\tau|^2\xi_1^3\xi_2^3b_0^3k^
4\delta_1(k)b_0k\delta_2(k)b_0k+\\112\tau_1^2|\tau|^2\xi_1^3\xi_2^3b_0^3k^4\delta_1(k)b_0k^2\delta_2(k)b_0+112\tau_1^2|\tau|^2\xi_1^3\xi_2^3b_0^3k^
4\delta_2(k)b_0k\delta_1(k)b_0k+\\112\tau_1^2|\tau|^2\xi_1^3\xi_2^3b_0^3k^4\delta_2(k)b_0k^2\delta_1(k)b_0+112\tau_1^2|\tau|^2\xi_1^3\xi_2^3b_0^3k^
5\delta_1(k)b_0\delta_2(k)b_0k+\\112\tau_1^2|\tau|^2\xi_1^3\xi_2^3b_0^3k^5\delta_1(k)b_0k\delta_2(k)b_0+112\tau_1^2|\tau|^2\xi_1^3\xi_2^3b_0^3k^
5\delta_2(k)b_0\delta_1(k)b_0k+\\112\tau_1^2|\tau|^2\xi_1^3\xi_2^3b_0^3k^5\delta_2(k)b_0k\delta_1(k)b_0+56\tau_1^2|\tau|^2\xi_1^3\xi_2^3b_0^2k^2\delta_1(k)b_0^
2k^3\delta_2(k)b_0k+\\56\tau_1^2|\tau|^2\xi_1^3\xi_2^3b_0^2k^2\delta_1(k)b_0^2k^4\delta_2(k)b_0+56\tau_1^2|\tau|^2\xi_1^3\xi_2^3b_0^2k^2\delta_2(k)b_0^2k^
3\delta_1(k)b_0k+\\56\tau_1^2|\tau|^2\xi_1^3\xi_2^3b_0^2k^2\delta_2(k)b_0^2k^4\delta_1(k)b_0+56\tau_1^2|\tau|^2\xi_1^3\xi_2^3b_0^2k^3\delta_1(k)b_0^2k^
2\delta_2(k)b_0k+\\56\tau_1^2|\tau|^2\xi_1^3\xi_2^3b_0^2k^3\delta_1(k)b_0^2k^3\delta_2(k)b_0+56\tau_1^2|\tau|^2\xi_1^3\xi_2^3b_0^2k^3\delta_2(k)b_0^2k^
2\delta_1(k)b_0k+\\56\tau_1^2|\tau|^2\xi_1^3\xi_2^3b_0^2k^3\delta_2(k)b_0^2k^3\delta_1(k)b_0+80\tau_1^2|\tau|^2\xi_1^2\xi_2^4b_0^3k^
4\delta_1(k)b_0k\delta_1(k)b_0k+\\80\tau_1^2|\tau|^2\xi_1^2\xi_2^4b_0^3k^4\delta_1(k)b_0k^2\delta_1(k)b_0+80\tau_1^2|\tau|^2\xi_1^2\xi_2^4b_0^3k^
5\delta_1(k)b_0\delta_1(k)b_0k+\\80\tau_1^2|\tau|^2\xi_1^2\xi_2^4b_0^3k^5\delta_1(k)b_0k\delta_1(k)b_0+40\tau_1^2|\tau|^2\xi_1^2\xi_2^4b_0^2k^2\delta_1(k)b_0^
2k^3\delta_1(k)b_0k+\\40\tau_1^2|\tau|^2\xi_1^2\xi_2^4b_0^2k^2\delta_1(k)b_0^2k^4\delta_1(k)b_0+40\tau_1^2|\tau|^2\xi_1^2\xi_2^4b_0^2k^3\delta_1(k)b_0^2k^
2\delta_1(k)b_0k+\\40\tau_1^2|\tau|^2\xi_1^2\xi_2^4b_0^2k^3\delta_1(k)b_0^2k^3\delta_1(k)b_0+8|\tau|^4\xi_1^4\xi_2^2b_0^3k^4\delta_2(k)b_0k\delta_2(k)b_0k+\\
8|\tau|^4\xi_1^4\xi_2^2b_0^3k^4\delta_2(k)b_0k^2\delta_2(k)b_0+8|\tau|^4\xi_1^4\xi_2^2b_0^3k^5\delta_2(k)b_0\delta_2(k)b_0k+\\8|\tau|^4\xi_1^4\xi_2^2b_0^3k^
5\delta_2(k)b_0k\delta_2(k)b_0+4|\tau|^4\xi_1^4\xi_2^2b_0^2k^2\delta_2(k)b_0^2k^3\delta_2(k)b_0k+\\4|\tau|^4\xi_1^4\xi_2^2b_0^2k^2\delta_2(k)b_0^2k^
4\delta_2(k)b_0+4|\tau|^4\xi_1^4\xi_2^2b_0^2k^3\delta_2(k)b_0^2k^2\delta_2(k)b_0k+\\4|\tau|^4\xi_1^4\xi_2^2b_0^2k^3\delta_2(k)b_0^2k^3\delta_2(k)b_0+16|\tau|^
4\xi_1^3\xi_2^3b_0^3k^4\delta_1(k)b_0k\delta_2(k)b_0k+\\16|\tau|^4\xi_1^3\xi_2^3b_0^3k^4\delta_1(k)b_0k^2\delta_2(k)b_0+16|\tau|^4\xi_1^3\xi_2^3b_0^3k^
4\delta_2(k)b_0k\delta_1(k)b_0k+\\16|\tau|^4\xi_1^3\xi_2^3b_0^3k^4\delta_2(k)b_0k^2\delta_1(k)b_0+16|\tau|^4\xi_1^3\xi_2^3b_0^3k^
5\delta_1(k)b_0\delta_2(k)b_0k+\\16|\tau|^4\xi_1^3\xi_2^3b_0^3k^5\delta_1(k)b_0k\delta_2(k)b_0+16|\tau|^4\xi_1^3\xi_2^3b_0^3k^5\delta_2(k)b_0\delta_1(k)b_0k+\\
16|\tau|^4\xi_1^3\xi_2^3b_0^3k^5\delta_2(k)b_0k\delta_1(k)b_0+8|\tau|^4\xi_1^3\xi_2^3b_0^2k^2\delta_1(k)b_0^2k^3\delta_2(k)b_0k+\\8|\tau|^4\xi_1^3\xi_2^3b_0^
2k^2\delta_1(k)b_0^2k^4\delta_2(k)b_0+8|\tau|^4\xi_1^3\xi_2^3b_0^2k^2\delta_2(k)b_0^2k^3\delta_1(k)b_0k+\\8|\tau|^4\xi_1^3\xi_2^3b_0^2k^2\delta_2(k)b_0^2k^
4\delta_1(k)b_0+8|\tau|^4\xi_1^3\xi_2^3b_0^2k^3\delta_1(k)b_0^2k^2\delta_2(k)b_0k+\\8|\tau|^4\xi_1^3\xi_2^3b_0^2k^3\delta_1(k)b_0^2k^3\delta_2(k)b_0+8|\tau|^
4\xi_1^3\xi_2^3b_0^2k^3\delta_2(k)b_0^2k^2\delta_1(k)b_0k+\\8|\tau|^4\xi_1^3\xi_2^3b_0^2k^3\delta_2(k)b_0^2k^3\delta_1(k)b_0+8|\tau|^4\xi_1^2\xi_2^4b_0^3k^
4\delta_1(k)b_0k\delta_1(k)b_0k+\\8|\tau|^4\xi_1^2\xi_2^4b_0^3k^4\delta_1(k)b_0k^2\delta_1(k)b_0+8|\tau|^4\xi_1^2\xi_2^4b_0^3k^
5\delta_1(k)b_0\delta_1(k)b_0k+\\8|\tau|^4\xi_1^2\xi_2^4b_0^3k^5\delta_1(k)b_0k\delta_1(k)b_0+4|\tau|^4\xi_1^2\xi_2^4b_0^2k^2\delta_1(k)b_0^2k^
3\delta_1(k)b_0k+\\4|\tau|^4\xi_1^2\xi_2^4b_0^2k^2\delta_1(k)b_0^2k^4\delta_1(k)b_0+4|\tau|^4\xi_1^2\xi_2^4b_0^2k^3\delta_1(k)b_0^2k^2\delta_1(k)b_0k+\\4|\tau|^
4\xi_1^2\xi_2^4b_0^2k^3\delta_1(k)b_0^2k^3\delta_1(k)b_0+96\tau_1^3|\tau|^2\xi_1^3\xi_2^3b_0^3k^4\delta_2(k)b_0k\delta_2(k)b_0k+\\96\tau_1^3|\tau|^2\xi_1^
3\xi_2^3b_0^3k^4\delta_2(k)b_0k^2\delta_2(k)b_0+96\tau_1^3|\tau|^2\xi_1^3\xi_2^3b_0^3k^5\delta_2(k)b_0\delta_2(k)b_0k+\\96\tau_1^3|\tau|^2\xi_1^3\xi_2^3b_0^
3k^5\delta_2(k)b_0k\delta_2(k)b_0+48\tau_1^3|\tau|^2\xi_1^3\xi_2^3b_0^2k^2\delta_2(k)b_0^2k^3\delta_2(k)b_0k+\\48\tau_1^3|\tau|^2\xi_1^3\xi_2^3b_0^2k^
2\delta_2(k)b_0^2k^4\delta_2(k)b_0+48\tau_1^3|\tau|^2\xi_1^3\xi_2^3b_0^2k^3\delta_2(k)b_0^2k^2\delta_2(k)b_0k+\\48\tau_1^3|\tau|^2\xi_1^3\xi_2^3b_0^2k^
3\delta_2(k)b_0^2k^3\delta_2(k)b_0+64\tau_1^3|\tau|^2\xi_1^2\xi_2^4b_0^3k^4\delta_1(k)b_0k\delta_2(k)b_0k+\\64\tau_1^3|\tau|^2\xi_1^2\xi_2^4b_0^3k^
4\delta_1(k)b_0k^2\delta_2(k)b_0+64\tau_1^3|\tau|^2\xi_1^2\xi_2^4b_0^3k^4\delta_2(k)b_0k\delta_1(k)b_0k+\\64\tau_1^3|\tau|^2\xi_1^2\xi_2^4b_0^3k^
4\delta_2(k)b_0k^2\delta_1(k)b_0+64\tau_1^3|\tau|^2\xi_1^2\xi_2^4b_0^3k^5\delta_1(k)b_0\delta_2(k)b_0k+\\64\tau_1^3|\tau|^2\xi_1^2\xi_2^4b_0^3k^
5\delta_1(k)b_0k\delta_2(k)b_0+64\tau_1^3|\tau|^2\xi_1^2\xi_2^4b_0^3k^5\delta_2(k)b_0\delta_1(k)b_0k+\\64\tau_1^3|\tau|^2\xi_1^2\xi_2^4b_0^3k^
5\delta_2(k)b_0k\delta_1(k)b_0+32\tau_1^3|\tau|^2\xi_1^2\xi_2^4b_0^2k^2\delta_1(k)b_0^2k^3\delta_2(k)b_0k+\\32\tau_1^3|\tau|^2\xi_1^2\xi_2^4b_0^2k^
2\delta_1(k)b_0^2k^4\delta_2(k)b_0+32\tau_1^3|\tau|^2\xi_1^2\xi_2^4b_0^2k^2\delta_2(k)b_0^2k^3\delta_1(k)b_0k+\\32\tau_1^3|\tau|^2\xi_1^2\xi_2^4b_0^2k^
2\delta_2(k)b_0^2k^4\delta_1(k)b_0+32\tau_1^3|\tau|^2\xi_1^2\xi_2^4b_0^2k^3\delta_1(k)b_0^2k^2\delta_2(k)b_0k+\\32\tau_1^3|\tau|^2\xi_1^2\xi_2^4b_0^2k^
3\delta_1(k)b_0^2k^3\delta_2(k)b_0+32\tau_1^3|\tau|^2\xi_1^2\xi_2^4b_0^2k^3\delta_2(k)b_0^2k^2\delta_1(k)b_0k+\\32\tau_1^3|\tau|^2\xi_1^2\xi_2^4b_0^2k^
3\delta_2(k)b_0^2k^3\delta_1(k)b_0+32\tau_1^3|\tau|^2\xi_1\xi_2^5b_0^3k^4\delta_1(k)b_0k\delta_1(k)b_0k+\\32\tau_1^3|\tau|^2\xi_1\xi_2^5b_0^3k^
4\delta_1(k)b_0k^2\delta_1(k)b_0+32\tau_1^3|\tau|^2\xi_1\xi_2^5b_0^3k^5\delta_1(k)b_0\delta_1(k)b_0k+\\32\tau_1^3|\tau|^2\xi_1\xi_2^5b_0^3k^
5\delta_1(k)b_0k\delta_1(k)b_0+16\tau_1^3|\tau|^2\xi_1\xi_2^5b_0^2k^2\delta_1(k)b_0^2k^3\delta_1(k)b_0k+\\16\tau_1^3|\tau|^2\xi_1\xi_2^5b_0^2k^
2\delta_1(k)b_0^2k^4\delta_1(k)b_0+16\tau_1^3|\tau|^2\xi_1\xi_2^5b_0^2k^3\delta_1(k)b_0^2k^2\delta_1(k)b_0k+\\16\tau_1^3|\tau|^2\xi_1\xi_2^5b_0^2k^
3\delta_1(k)b_0^2k^3\delta_1(k)b_0+64\tau_1|\tau|^4\xi_1^3\xi_2^3b_0^3k^4\delta_2(k)b_0k\delta_2(k)b_0k+\\64\tau_1|\tau|^4\xi_1^3\xi_2^3b_0^3k^
4\delta_2(k)b_0k^2\delta_2(k)b_0+64\tau_1|\tau|^4\xi_1^3\xi_2^3b_0^3k^5\delta_2(k)b_0\delta_2(k)b_0k+\\64\tau_1|\tau|^4\xi_1^3\xi_2^3b_0^3k^
5\delta_2(k)b_0k\delta_2(k)b_0+32\tau_1|\tau|^4\xi_1^3\xi_2^3b_0^2k^2\delta_2(k)b_0^2k^3\delta_2(k)b_0k+\\32\tau_1|\tau|^4\xi_1^3\xi_2^3b_0^2k^
2\delta_2(k)b_0^2k^4\delta_2(k)b_0+32\tau_1|\tau|^4\xi_1^3\xi_2^3b_0^2k^3\delta_2(k)b_0^2k^2\delta_2(k)b_0k+\\32\tau_1|\tau|^4\xi_1^3\xi_2^3b_0^2k^
3\delta_2(k)b_0^2k^3\delta_2(k)b_0+56\tau_1|\tau|^4\xi_1^2\xi_2^4b_0^3k^4\delta_1(k)b_0k\delta_2(k)b_0k+\\56\tau_1|\tau|^4\xi_1^2\xi_2^4b_0^3k^
4\delta_1(k)b_0k^2\delta_2(k)b_0+56\tau_1|\tau|^4\xi_1^2\xi_2^4b_0^3k^4\delta_2(k)b_0k\delta_1(k)b_0k+\\56\tau_1|\tau|^4\xi_1^2\xi_2^4b_0^3k^
4\delta_2(k)b_0k^2\delta_1(k)b_0+56\tau_1|\tau|^4\xi_1^2\xi_2^4b_0^3k^5\delta_1(k)b_0\delta_2(k)b_0k+\\56\tau_1|\tau|^4\xi_1^2\xi_2^4b_0^3k^
5\delta_1(k)b_0k\delta_2(k)b_0+56\tau_1|\tau|^4\xi_1^2\xi_2^4b_0^3k^5\delta_2(k)b_0\delta_1(k)b_0k+\\56\tau_1|\tau|^4\xi_1^2\xi_2^4b_0^3k^
5\delta_2(k)b_0k\delta_1(k)b_0+28\tau_1|\tau|^4\xi_1^2\xi_2^4b_0^2k^2\delta_1(k)b_0^2k^3\delta_2(k)b_0k+\\28\tau_1|\tau|^4\xi_1^2\xi_2^4b_0^2k^
2\delta_1(k)b_0^2k^4\delta_2(k)b_0+28\tau_1|\tau|^4\xi_1^2\xi_2^4b_0^2k^2\delta_2(k)b_0^2k^3\delta_1(k)b_0k+\\28\tau_1|\tau|^4\xi_1^2\xi_2^4b_0^2k^
2\delta_2(k)b_0^2k^4\delta_1(k)b_0+28\tau_1|\tau|^4\xi_1^2\xi_2^4b_0^2k^3\delta_1(k)b_0^2k^2\delta_2(k)b_0k+\\28\tau_1|\tau|^4\xi_1^2\xi_2^4b_0^2k^
3\delta_1(k)b_0^2k^3\delta_2(k)b_0+28\tau_1|\tau|^4\xi_1^2\xi_2^4b_0^2k^3\delta_2(k)b_0^2k^2\delta_1(k)b_0k+\\28\tau_1|\tau|^4\xi_1^2\xi_2^4b_0^2k^
3\delta_2(k)b_0^2k^3\delta_1(k)b_0+16\tau_1|\tau|^4\xi_1\xi_2^5b_0^3k^4\delta_1(k)b_0k\delta_1(k)b_0k+\\16\tau_1|\tau|^4\xi_1\xi_2^5b_0^3k^4\delta_1(k)b_0k^
2\delta_1(k)b_0+16\tau_1|\tau|^4\xi_1\xi_2^5b_0^3k^5\delta_1(k)b_0\delta_1(k)b_0k+\\16\tau_1|\tau|^4\xi_1\xi_2^5b_0^3k^5\delta_1(k)b_0k\delta_1(k)b_0+
8\tau_1|\tau|^4\xi_1\xi_2^5b_0^2k^2\delta_1(k)b_0^2k^3\delta_1(k)b_0k+\\8\tau_1|\tau|^4\xi_1\xi_2^5b_0^2k^2\delta_1(k)b_0^2k^4\delta_1(k)b_0+8\tau_1|\tau|^
4\xi_1\xi_2^5b_0^2k^3\delta_1(k)b_0^2k^2\delta_1(k)b_0k+\\8\tau_1|\tau|^4\xi_1\xi_2^5b_0^2k^3\delta_1(k)b_0^2k^3\delta_1(k)b_0+104\tau_1^2|\tau|^4\xi_1^
2\xi_2^4b_0^3k^4\delta_2(k)b_0k\delta_2(k)b_0k+\\104\tau_1^2|\tau|^4\xi_1^2\xi_2^4b_0^3k^4\delta_2(k)b_0k^2\delta_2(k)b_0+104\tau_1^2|\tau|^4\xi_1^2\xi_2^
4b_0^3k^5\delta_2(k)b_0\delta_2(k)b_0k+\\104\tau_1^2|\tau|^4\xi_1^2\xi_2^4b_0^3k^5\delta_2(k)b_0k\delta_2(k)b_0+52\tau_1^2|\tau|^4\xi_1^2\xi_2^4b_0^2k^
2\delta_2(k)b_0^2k^3\delta_2(k)b_0k+\\52\tau_1^2|\tau|^4\xi_1^2\xi_2^4b_0^2k^2\delta_2(k)b_0^2k^4\delta_2(k)b_0+52\tau_1^2|\tau|^4\xi_1^2\xi_2^4b_0^2k^
3\delta_2(k)b_0^2k^2\delta_2(k)b_0k+\\52\tau_1^2|\tau|^4\xi_1^2\xi_2^4b_0^2k^3\delta_2(k)b_0^2k^3\delta_2(k)b_0+40\tau_1^2|\tau|^4\xi_1\xi_2^5b_0^3k^
4\delta_1(k)b_0k\delta_2(k)b_0k+\\40\tau_1^2|\tau|^4\xi_1\xi_2^5b_0^3k^4\delta_1(k)b_0k^2\delta_2(k)b_0+40\tau_1^2|\tau|^4\xi_1\xi_2^5b_0^3k^
4\delta_2(k)b_0k\delta_1(k)b_0k+\\40\tau_1^2|\tau|^4\xi_1\xi_2^5b_0^3k^4\delta_2(k)b_0k^2\delta_1(k)b_0+40\tau_1^2|\tau|^4\xi_1\xi_2^5b_0^3k^
5\delta_1(k)b_0\delta_2(k)b_0k+\\40\tau_1^2|\tau|^4\xi_1\xi_2^5b_0^3k^5\delta_1(k)b_0k\delta_2(k)b_0+40\tau_1^2|\tau|^4\xi_1\xi_2^5b_0^3k^
5\delta_2(k)b_0\delta_1(k)b_0k+\\40\tau_1^2|\tau|^4\xi_1\xi_2^5b_0^3k^5\delta_2(k)b_0k\delta_1(k)b_0+20\tau_1^2|\tau|^4\xi_1\xi_2^5b_0^2k^2\delta_1(k)b_0^2k^
3\delta_2(k)b_0k+\\20\tau_1^2|\tau|^4\xi_1\xi_2^5b_0^2k^2\delta_1(k)b_0^2k^4\delta_2(k)b_0+20\tau_1^2|\tau|^4\xi_1\xi_2^5b_0^2k^2\delta_2(k)b_0^2k^
3\delta_1(k)b_0k+\\20\tau_1^2|\tau|^4\xi_1\xi_2^5b_0^2k^2\delta_2(k)b_0^2k^4\delta_1(k)b_0+20\tau_1^2|\tau|^4\xi_1\xi_2^5b_0^2k^3\delta_1(k)b_0^2k^
2\delta_2(k)b_0k+\\20\tau_1^2|\tau|^4\xi_1\xi_2^5b_0^2k^3\delta_1(k)b_0^2k^3\delta_2(k)b_0+20\tau_1^2|\tau|^4\xi_1\xi_2^5b_0^2k^3\delta_2(k)b_0^2k^
2\delta_1(k)b_0k+\\20\tau_1^2|\tau|^4\xi_1\xi_2^5b_0^2k^3\delta_2(k)b_0^2k^3\delta_1(k)b_0+8\tau_1^2|\tau|^4\xi_2^6b_0^3k^4\delta_1(k)b_0k\delta_1(k)b_0k+\\
8\tau_1^2|\tau|^4\xi_2^6b_0^3k^4\delta_1(k)b_0k^2\delta_1(k)b_0+8\tau_1^2|\tau|^4\xi_2^6b_0^3k^5\delta_1(k)b_0\delta_1(k)b_0k+\\8\tau_1^2|\tau|^4\xi_2^6b_0^3k^
5\delta_1(k)b_0k\delta_1(k)b_0+4\tau_1^2|\tau|^4\xi_2^6b_0^2k^2\delta_1(k)b_0^2k^3\delta_1(k)b_0k+\\4\tau_1^2|\tau|^4\xi_2^6b_0^2k^2\delta_1(k)b_0^2k^
4\delta_1(k)b_0+4\tau_1^2|\tau|^4\xi_2^6b_0^2k^3\delta_1(k)b_0^2k^2\delta_1(k)b_0k+\\4\tau_1^2|\tau|^4\xi_2^6b_0^2k^3\delta_1(k)b_0^2k^3\delta_1(k)b_0+16|\tau|^
6\xi_1^2\xi_2^4b_0^3k^4\delta_2(k)b_0k\delta_2(k)b_0k+\\16|\tau|^6\xi_1^2\xi_2^4b_0^3k^4\delta_2(k)b_0k^2\delta_2(k)b_0+16|\tau|^6\xi_1^2\xi_2^4b_0^3k^
5\delta_2(k)b_0\delta_2(k)b_0k+\\16|\tau|^6\xi_1^2\xi_2^4b_0^3k^5\delta_2(k)b_0k\delta_2(k)b_0+8|\tau|^6\xi_1^2\xi_2^4b_0^2k^2\delta_2(k)b_0^2k^
3\delta_2(k)b_0k+\\8|\tau|^6\xi_1^2\xi_2^4b_0^2k^2\delta_2(k)b_0^2k^4\delta_2(k)b_0+8|\tau|^6\xi_1^2\xi_2^4b_0^2k^3\delta_2(k)b_0^2k^2\delta_2(k)b_0k+\\8|\tau|^
6\xi_1^2\xi_2^4b_0^2k^3\delta_2(k)b_0^2k^3\delta_2(k)b_0+8|\tau|^6\xi_1\xi_2^5b_0^3k^4\delta_1(k)b_0k\delta_2(k)b_0k+\\8|\tau|^6\xi_1\xi_2^5b_0^3k^
4\delta_1(k)b_0k^2\delta_2(k)b_0+8|\tau|^6\xi_1\xi_2^5b_0^3k^4\delta_2(k)b_0k\delta_1(k)b_0k+\\8|\tau|^6\xi_1\xi_2^5b_0^3k^4\delta_2(k)b_0k^2\delta_1(k)b_0+
8|\tau|^6\xi_1\xi_2^5b_0^3k^5\delta_1(k)b_0\delta_2(k)b_0k+\\8|\tau|^6\xi_1\xi_2^5b_0^3k^5\delta_1(k)b_0k\delta_2(k)b_0+8|\tau|^6\xi_1\xi_2^5b_0^3k^
5\delta_2(k)b_0\delta_1(k)b_0k+\\8|\tau|^6\xi_1\xi_2^5b_0^3k^5\delta_2(k)b_0k\delta_1(k)b_0+4|\tau|^6\xi_1\xi_2^5b_0^2k^2\delta_1(k)b_0^2k^3\delta_2(k)b_0k+\\
4|\tau|^6\xi_1\xi_2^5b_0^2k^2\delta_1(k)b_0^2k^4\delta_2(k)b_0+4|\tau|^6\xi_1\xi_2^5b_0^2k^2\delta_2(k)b_0^2k^3\delta_1(k)b_0k+\\4|\tau|^6\xi_1\xi_2^5b_0^2k^
2\delta_2(k)b_0^2k^4\delta_1(k)b_0+4|\tau|^6\xi_1\xi_2^5b_0^2k^3\delta_1(k)b_0^2k^2\delta_2(k)b_0k+\\4|\tau|^6\xi_1\xi_2^5b_0^2k^3\delta_1(k)b_0^2k^
3\delta_2(k)b_0+4|\tau|^6\xi_1\xi_2^5b_0^2k^3\delta_2(k)b_0^2k^2\delta_1(k)b_0k+\\4|\tau|^6\xi_1\xi_2^5b_0^2k^3\delta_2(k)b_0^2k^3\delta_1(k)b_0+48\tau_1|\tau|^
6\xi_1\xi_2^5b_0^3k^4\delta_2(k)b_0k\delta_2(k)b_0k+\\48\tau_1|\tau|^6\xi_1\xi_2^5b_0^3k^4\delta_2(k)b_0k^2\delta_2(k)b_0+48\tau_1|\tau|^6\xi_1\xi_2^5b_0^
3k^5\delta_2(k)b_0\delta_2(k)b_0k+\\48\tau_1|\tau|^6\xi_1\xi_2^5b_0^3k^5\delta_2(k)b_0k\delta_2(k)b_0+24\tau_1|\tau|^6\xi_1\xi_2^5b_0^2k^2\delta_2(k)b_0^2k^
3\delta_2(k)b_0k+\\24\tau_1|\tau|^6\xi_1\xi_2^5b_0^2k^2\delta_2(k)b_0^2k^4\delta_2(k)b_0+24\tau_1|\tau|^6\xi_1\xi_2^5b_0^2k^3\delta_2(k)b_0^2k^
2\delta_2(k)b_0k+\\24\tau_1|\tau|^6\xi_1\xi_2^5b_0^2k^3\delta_2(k)b_0^2k^3\delta_2(k)b_0+8\tau_1|\tau|^6\xi_2^6b_0^3k^4\delta_1(k)b_0k\delta_2(k)b_0k+\\
8\tau_1|\tau|^6\xi_2^6b_0^3k^4\delta_1(k)b_0k^2\delta_2(k)b_0+8\tau_1|\tau|^6\xi_2^6b_0^3k^4\delta_2(k)b_0k\delta_1(k)b_0k+\\8\tau_1|\tau|^6\xi_2^6b_0^3k^
4\delta_2(k)b_0k^2\delta_1(k)b_0+8\tau_1|\tau|^6\xi_2^6b_0^3k^5\delta_1(k)b_0\delta_2(k)b_0k+\\8\tau_1|\tau|^6\xi_2^6b_0^3k^5\delta_1(k)b_0k\delta_2(k)b_0+
8\tau_1|\tau|^6\xi_2^6b_0^3k^5\delta_2(k)b_0\delta_1(k)b_0k+\\8\tau_1|\tau|^6\xi_2^6b_0^3k^5\delta_2(k)b_0k\delta_1(k)b_0+4\tau_1|\tau|^6\xi_2^6b_0^2k^
2\delta_1(k)b_0^2k^3\delta_2(k)b_0k+\\4\tau_1|\tau|^6\xi_2^6b_0^2k^2\delta_1(k)b_0^2k^4\delta_2(k)b_0+4\tau_1|\tau|^6\xi_2^6b_0^2k^2\delta_2(k)b_0^2k^
3\delta_1(k)b_0k+\\4\tau_1|\tau|^6\xi_2^6b_0^2k^2\delta_2(k)b_0^2k^4\delta_1(k)b_0+4\tau_1|\tau|^6\xi_2^6b_0^2k^3\delta_1(k)b_0^2k^2\delta_2(k)b_0k+\\
4\tau_1|\tau|^6\xi_2^6b_0^2k^3\delta_1(k)b_0^2k^3\delta_2(k)b_0+4\tau_1|\tau|^6\xi_2^6b_0^2k^3\delta_2(k)b_0^2k^2\delta_1(k)b_0k+\\4\tau_1|\tau|^6\xi_2^6b_0^
2k^3\delta_2(k)b_0^2k^3\delta_1(k)b_0+8|\tau|^8\xi_2^6b_0^3k^4\delta_2(k)b_0k\delta_2(k)b_0k+\\8|\tau|^8\xi_2^6b_0^3k^4\delta_2(k)b_0k^2\delta_2(k)b_0+8|\tau|^
8\xi_2^6b_0^3k^5\delta_2(k)b_0\delta_2(k)b_0k+\\8|\tau|^8\xi_2^6b_0^3k^5\delta_2(k)b_0k\delta_2(k)b_0+4|\tau|^8\xi_2^6b_0^2k^2\delta_2(k)b_0^2k^
3\delta_2(k)b_0k+\\
4|\tau|^8\xi_2^6b_0^2k^2\delta_2(k)b_0^2k^4\delta_2(k)b_0+4|\tau|^8\xi_2^6b_0^2k^3\delta_2(k)b_0^2k^2\delta_2(k)b_0k+\\4|\tau|^8\xi_2^6b_0^
2k^3\delta_2(k)b_0^2k^3\delta_2(k)b_0.$\\

\section{The $b_2$ term of  $\partial^* k^2 \partial$} \label{secondb2}
The $b_2$ term of the second operator, namely $\partial^* k^2 \partial$, is equal to \\

\noindent
$
+\xi_1^2b_0\delta_1(k^2)b_0\delta_1(k^2)b_0                        
+\tau_1\xi_1^2b_0\delta_1(k^2)b_0\delta_2(k^2)b_0\\                        
+\tau_1\xi_1^2b_0\delta_2(k^2)b_0\delta_1(k^2)b_0                       
+2\tau_1\xi_1\xi_2b_0\delta_1(k^2)b_0\delta_1(k^2)b_0\\                        
+i\tau_2\xi_1^2b_0\delta_1(k^2)b_0\delta_2(k^2)b_0                        
+i\tau_2\xi_1^2b_0\delta_2(k^2)b_0\delta_1(k^2)b_0\\                        
-2i\tau_2\xi_1\xi_2b_0\delta_1(k^2)b_0\delta_1(k^2)b_0                        
+2\xi_1^2b_0^2k^2\delta_1^2(k^2)b_0\\                        
+\xi_1^2b_0\delta_1(k^2)b_0k\delta_1(k)b_0                        
+\xi_1^2b_0\delta_1(k^2)b_0\delta_1(k)b_0k\\                        
+\tau_1^2\xi_1^2b_0\delta_2(k^2)b_0\delta_2(k^2)b_0                        
+\tau_1^2\xi_1\xi_2b_0\delta_1(k^2)b_0\delta_2(k^2)b_0\\                        
+\tau_1^2\xi_1\xi_2b_0\delta_2(k^2)b_0\delta_1(k^2)b_0                       
+\tau_1^2\xi_2^2b_0\delta_1(k^2)b_0\delta_1(k^2)b_0\\                        
+2i\tau_1\tau_2\xi_1^2b_0\delta_2(k^2)b_0\delta_2(k^2)b_0                       
-2i\tau_1\tau_2\xi_2^2b_0\delta_1(k^2)b_0\delta_1(k^2)b_0\\                        
+4\tau_1\xi_1^2b_0^2k^2\delta_1\delta_2(k^2)b_0                       
+\tau_1\xi_1^2b_0\delta_1(k^2)b_0k\delta_2(k)b_0\\                       
 +\tau_1\xi_1^2b_0\delta_1(k^2)b_0\delta_2(k)b_0k                        
 +\tau_1\xi_1^2b_0\delta_2(k^2)b_0k\delta_1(k)b_0\\                        
 +\tau_1\xi_1^2b_0\delta_2(k^2)b_0\delta_1(k)b_0k                       
 +4\tau_1\xi_1\xi_2b_0^2k^2\delta_1^2(k^2)b_0\\                       
 +2\tau_1\xi_1\xi_2b_0\delta_1(k^2)b_0k\delta_1(k)b_0                      
 +2\tau_1\xi_1\xi_2b_0\delta_1(k^2)b_0\delta_1(k)b_0k\\                        
 -\tau_2^2\xi_1^2b_0\delta_2(k^2)b_0\delta_2(k^2)b_0                       
 +\tau_2^2\xi_1\xi_2b_0\delta_1(k^2)b_0\delta_2(k^2)b_0\\                        
 +\tau_2^2\xi_1\xi_2b_0\delta_2(k^2)b_0\delta_1(k^2)b_0                       
 -\tau_2^2\xi_2^2b_0\delta_1(k^2)b_0\delta_1(k^2)b_0\\                        
+2i\tau_2\xi_1^2b_0^2k^2\delta_1\delta_2(k^2)b_0                       
-i\tau_2\xi_1^2b_0\delta_1(k^2)b_0k\delta_2(k)b_0\\                        
-i\tau_2\xi_1^2b_0\delta_1(k^2)b_0\delta_2(k)b_0k                      
+i\tau_2\xi_1^2b_0\delta_2(k^2)b_0k\delta_1(k)b_0\\                        
+i\tau_2\xi_1^2b_0\delta_2(k^2)b_0\delta_1(k)b_0k                       
-2i\tau_2\xi_1\xi_2b_0^2k^2\delta_1^2(k^2)b_0\\                        
+|\tau|^2\xi_1\xi_2b_0\delta_1(k^2)b_0\delta_2(k^2)b_0                       
+|\tau|^2\xi_1\xi_2b_0\delta_2(k^2)b_0\delta_1(k^2)b_0\\                        
+\xi_1^2b_0^2k^3\delta_1^2(k)b_0                        
+2\xi_1^2b_0^2k^2\delta_1(k)^2b_0\\                        
+\xi_1^2b_0^2k^2\delta_1^2(k)b_0k                        
+2\tau_1^2\xi_1^2b_0^2k^2\delta_2^2(k^2)b_0\\                        
+4\tau_1^2\xi_1\xi_2b_0^2k^2\delta_1\delta_2(k^2)b_0                        
+2\tau_1^2\xi_1\xi_2b_0\delta_1(k^2)b_0k\delta_2(k)b_0\\                        
+2\tau_1^2\xi_1\xi_2b_0\delta_1(k^2)b_0\delta_2(k)b_0k                        
+2\tau_1^2\xi_1\xi_2b_0\delta_2(k^2)b_0k\delta_1(k)b_0\\                        
+2\tau_1^2\xi_1\xi_2b_0\delta_2(k^2)b_0\delta_1(k)b_0k                        
+2\tau_1^2\xi_2^2b_0^2k^2\delta_1^2(k^2)b_0\\                        
+2i\tau_1\tau_2\xi_1^2b_0^2k^2\delta_2^2(k^2)b_0                       
-2i\tau_1\tau_2\xi_1\xi_2b_0\delta_1(k^2)b_0k\delta_2(k)b_0\\                        
-2i\tau_1\tau_2\xi_1\xi_2b_0\delta_1(k^2)b_0\delta_2(k)b_0k                        
+2i\tau_1\tau_2\xi_1\xi_2b_0\delta_2(k^2)b_0k\delta_1(k)b_0\\                        
+2i\tau_1\tau_2\xi_1\xi_2b_0\delta_2(k^2)b_0\delta_1(k)b_0k                        
-2i\tau_1\tau_2\xi_2^2b_0^2k^2\delta_1^2(k^2)b_0\\                       
+2\tau_1|\tau|^2\xi_1\xi_2b_0\delta_2(k^2)b_0\delta_2(k^2)b_0                        
+\tau_1|\tau|^2\xi_2^2b_0\delta_1(k^2)b_0\delta_2(k^2)b_0\\                        
+\tau_1|\tau|^2\xi_2^2b_0\delta_2(k^2)b_0\delta_1(k^2)b_0                       
 +2\tau_1\xi_1^2b_0^2k^3\delta_1\delta_2(k)b_0\\                        
 +2\tau_1\xi_1^2b_0^2k^2\delta_1(k)\delta_2(k)b_0                       
 +2\tau_1\xi_1^2b_0^2k^2\delta_2(k)\delta_1(k)b_0\\                        
 +2\tau_1\xi_1^2b_0^2k^2\delta_1\delta_2(k)b_0k                      
 +2\tau_1\xi_1\xi_2b_0^2k^3\delta_1^2(k)b_0\\                        
 +4\tau_1\xi_1\xi_2b_0^2k^2\delta_1(k)^2b_0                        
 +2\tau_1\xi_1\xi_2b_0^2k^2\delta_1^2(k)b_0k\\                        
 +2i\tau_2|\tau|^2\xi_1\xi_2b_0\delta_2(k^2)b_0\delta_2(k^2)b_0                        
 -i\tau_2|\tau|^2\xi_2^2b_0\delta_1(k^2)b_0\delta_2(k^2)b_0\\                        
 -i \tau_2|\tau|^2\xi_2^2b_0\delta_2(k^2)b_0\delta_1(k^2)b_0                       
 +|\tau|^2\xi_1^2b_0\delta_2(k^2)b_0k\delta_2(k)b_0\\                        
 +|\tau|^2\xi_1^2b_0\delta_2(k^2)b_0\delta_2(k)b_0k                        
 +4|\tau|^2\xi_1\xi_2b_0^2k^2\delta_1\delta_2(k^2)b_0\\                        
 +|\tau|^2\xi_2^2b_0\delta_1(k^2)b_0k\delta_1(k)b_0                        
 +|\tau|^2\xi_2^2b_0\delta_1(k^2)b_0\delta_1(k)b_0k\\                        
 +4\tau_1^2\xi_1\xi_2b_0^2k^3\delta_1\delta_2(k)b_0                        
 +4\tau_1^2\xi_1\xi_2b_0^2k^2\delta_1(k)\delta_2(k)b_0\\                        
 +4\tau_1^2\xi_1\xi_2b_0^2k^2\delta_2(k)\delta_1(k)b_0                        
 +4\tau_1^2\xi_1\xi_2b_0^2k^2\delta_1\delta_2(k)b_0k\\                        
 +4\tau_1|\tau|^2\xi_1\xi_2b_0^2k^2\delta_2^2(k^2)b_0                        
 +2\tau_1|\tau|^2\xi_1\xi_2b_0\delta_2(k^2)b_0k\delta_2(k)b_0\\                        
 +2\tau_1|\tau|^2\xi_1\xi_2b_0\delta_2(k^2)b_0\delta_2(k)b_0k                        
 +4\tau_1|\tau|^2\xi_2^2b_0^2k^2\delta_1\delta_2(k^2)b_0\\                        
 +\tau_1|\tau|^2\xi_2^2b_0\delta_1(k^2)b_0k\delta_2(k)b_0                        
 +\tau_1|\tau|^2\xi_2^2b_0\delta_1(k^2)b_0\delta_2(k)b_0k\\                        
 +\tau_1|\tau|^2\xi_2^2b_0\delta_2(k^2)b_0k\delta_1(k)b_0                        
 +\tau_1|\tau|^2\xi_2^2b_0\delta_2(k^2)b_0\delta_1(k)b_0k\\                        
 +2i\tau_2|\tau|^2\xi_1\xi_2b_0^2k^2\delta_2^2(k^2)b_0                        
 -2i\tau_2|\tau|^2\xi_2^2b_0^2k^2\delta_1\delta_2(k^2)b_0\\                        
 -i\tau_2|\tau|^2\xi_2^2b_0\delta_1(k^2)b_0k\delta_2(k)b_0                        
 -i \tau_2|\tau|^2\xi_2^2b_0\delta_1(k^2)b_0\delta_2(k)b_0k\\                        
 +i\tau_2|\tau|^2\xi_2^2b_0\delta_2(k^2)b_0k\delta_1(k)b_0                        
 +i\tau_2|\tau|^2\xi_2^2b_0\delta_2(k^2)b_0\delta_1(k)b_0k\\                        
 +|\tau|^4\xi_2^2b_0\delta_2(k^2)b_0\delta_2(k^2)b_0                        
 +|\tau|^2\xi_1^2b_0^2k^3\delta_2^2(k)b_0\\                        
 +2|\tau|^2\xi_1^2b_0^2k^2\delta_2(k)^2b_0                       
 +|\tau|^2\xi_1^2b_0^2k^2\delta_2^2(k)b_0k\\                        
 +|\tau|^2\xi_2^2b_0^2k^3\delta_1^2(k)b_0                        
 +2|\tau|^2\xi_2^2b_0^2k^2\delta_1(k)^2b_0\\                        
 +|\tau|^2\xi_2^2b_0^2k^2\delta_1^2(k)b_0k                        
 +2\tau_1|\tau|^2\xi_1\xi_2b_0^2k^3\delta_2^2(k)b_0\\                       
 +4\tau_1|\tau|^2\xi_1\xi_2b_0^2k^2\delta_2(k)^2b_0                        
 +2\tau_1|\tau|^2\xi_1\xi_2b_0^2k^2\delta_2^2(k)b_0k\\                        
 +2\tau_1|\tau|^2\xi_2^2b_0^2k^3\delta_1\delta_2(k)b_0                        
 +2\tau_1|\tau|^2\xi_2^2b_0^2k^2\delta_1(k)\delta_2(k)b_0\\                        
 +2\tau_1|\tau|^2\xi_2^2b_0^2k^2\delta_2(k)\delta_1(k)b_0                        
 +2\tau_1|\tau|^2\xi_2^2b_0^2k^2\delta_1\delta_2(k)b_0k\\                        
 +2|\tau|^4\xi_2^2b_0^2k^2\delta_2^2(k^2)b_0                        
 +|\tau|^4\xi_2^2b_0\delta_2(k^2)b_0k\delta_2(k)b_0\\                        
 +|\tau|^4\xi_2^2b_0\delta_2(k^2)b_0\delta_2(k)b_0k                       
 +|\tau|^4\xi_2^2b_0^2k^3\delta_2^2(k)b_0\\                        
 +2|\tau|^4\xi_2^2b_0^2k^2\delta_2(k)^2b_0                        
 +|\tau|^4\xi_2^2b_0^2k^2\delta_2^2(k)b_0k\\                        
 -2\xi_1^4b_0^2k^3\delta_1(k)b_0\delta_1(k^2)b_0                        
 -2\xi_1^4b_0^2k^2\delta_1(k)b_0k\delta_1(k^2)b_0\\                        
 -2\xi_1^4b_0^2k^2\delta_1(k^2)b_0k\delta_1(k)b_0                        
 -2\xi_1^4b_0^2k^2\delta_1(k^2)b_0\delta_1(k)b_0k\\                        
 -2\xi_1^4b_0\delta_1(k^2)b_0^2k^3\delta_1(k)b_0                       
 -2\xi_1^4b_0\delta_1(k^2)b_0^2k^2\delta_1(k)b_0k\\                        
 -2\tau_1\xi_1^4b_0^2k^3\delta_1(k)b_0\delta_2(k^2)b_0                        
 -2\tau_1\xi_1^4b_0^2k^3\delta_2(k)b_0\delta_1(k^2)b_0\\                        
 -2\tau_1\xi_1^4b_0^2k^2\delta_1(k)b_0k\delta_2(k^2)b_0                       
 -2\tau_1\xi_1^4b_0^2k^2\delta_2(k)b_0k\delta_1(k^2)b_0\\                        
 -2\tau_1\xi_1^4b_0^2k^2\delta_1(k^2)b_0k\delta_2(k)b_0                        
 -2\tau_1\xi_1^4b_0^2k^2\delta_1(k^2)b_0\delta_2(k)b_0k\\                        
 -2\tau_1\xi_1^4b_0^2k^2\delta_2(k^2)b_0k\delta_1(k)b_0                       
 -2\tau_1\xi_1^4b_0^2k^2\delta_2(k^2)b_0\delta_1(k)b_0k\\                        
 -2\tau_1\xi_1^4b_0\delta_1(k^2)b_0^2k^3\delta_2(k)b_0                       
 -2\tau_1\xi_1^4b_0\delta_1(k^2)b_0^2k^2\delta_2(k)b_0k\\                       
 -2\tau_1\xi_1^4b_0\delta_2(k^2)b_0^2k^3\delta_1(k)b_0                        
 -2\tau_1\xi_1^4b_0\delta_2(k^2)b_0^2k^2\delta_1(k)b_0k\\                        
 -8\tau_1\xi_1^3\xi_2b_0^2k^3\delta_1(k)b_0\delta_1(k^2)b_0                       
  -8\tau_1\xi_1^3\xi_2b_0^2k^2\delta_1(k)b_0k\delta_1(k^2)b_0\\                        
-8\tau_1\xi_1^3\xi_2b_0^2k^2\delta_1(k^2)b_0k\delta_1(k)b_0                      
-8\tau_1\xi_1^3\xi_2b_0^2k^2\delta_1(k^2)b_0\delta_1(k)b_0k\\                        
-8\tau_1\xi_1^3\xi_2b_0\delta_1(k^2)b_0^2k^3\delta_1(k)b_0                       
-8\tau_1\xi_1^3\xi_2b_0\delta_1(k^2)b_0^2k^2\delta_1(k)b_0k\\                        
-2i\tau_2\xi_1^4b_0^2k^3\delta_1(k)b_0\delta_2(k^2)b_0                       
-2i\tau_2\xi_1^4b_0^2k^2\delta_1(k)b_0k\delta_2(k^2)b_0\\                        
-2i\tau_2\xi_1^4b_0^2k^2\delta_2(k^2)b_0k\delta_1(k)b_0                    
-2i\tau_2\xi_1^4b_0^2k^2\delta_2(k^2)b_0\delta_1(k)b_0k\\                        
-2i\tau_2\xi_1^4b_0\delta_2(k^2)b_0^2k^3\delta_1(k)b_0                       
-2i\tau_2\xi_1^4b_0\delta_2(k^2)b_0^2k^2\delta_1(k)b_0k\\                        
+2i\tau_2\xi_1^3\xi_2b_0^2k^3\delta_1(k)b_0\delta_1(k^2)b_0                       
+2i\tau_2\xi_1^3\xi_2b_0^2k^2\delta_1(k)b_0k\delta_1(k^2)b_0\\                       
+2i\tau_2\xi_1^3\xi_2b_0^2k^2\delta_1(k^2)b_0k\delta_1(k)b_0                        
+2i\tau_2\xi_1^3\xi_2b_0^2k^2\delta_1(k^2)b_0\delta_1(k)b_0k\\                        
+2i\tau_2\xi_1^3\xi_2b_0\delta_1(k^2)b_0^2k^3\delta_1(k)b_0                        
-2i\tau_2\xi_1^3\xi_2b_0\delta_1(k^2)b_0^2k^2\delta_1(k)b_0k\\                        
-4\xi_1^4b_0^3k^5\delta_1^2(k)b_0                        
-8\xi_1^4b_0^3k^4\delta_1(k)^2b_0\\                        
-4\xi_1^4b_0^3k^4\delta_1^2(k)b_0k                        
-6\xi_1^4b_0^2k^3\delta_1(k)b_0k\delta_1(k)b_0\\                        
-6\xi_1^4b_0^2k^3\delta_1(k)b_0\delta_1(k)b_0k                        
-6\xi_1^4b_0^2k^2\delta_1(k)b_0k^2\delta_1(k)b_0\\                        
-6\xi_1^4b_0^2k^2\delta_1(k)b_0k\delta_1(k)b_0k                        
-2\tau_1^2\xi_1^4b_0^2k^3\delta_2(k)b_0\delta_2(k^2)b_0\\                        
-2\tau_1^2\xi_1^4b_0^2k^2\delta_2(k)b_0k\delta_2(k^2)b_0                        
-2\tau_1^2\xi_1^4b_0^2k^2\delta_2(k^2)b_0k\delta_2(k)b_0\\                        
-2\tau_1^2\xi_1^4b_0^2k^2\delta_2(k^2)b_0\delta_2(k)b_0k                        
-2\tau_1^2\xi_1^4b_0\delta_2(k^2)b_0^2k^3\delta_2(k)b_0\\                        
-2\tau_1^2\xi_1^4b_0\delta_2(k^2)b_0^2k^2\delta_2(k)b_0k                        
-6\tau_1^2\xi_1^3\xi_2b_0^2k^3\delta_1(k)b_0\delta_2(k^2)b_0\\                        
-6\tau_1^2\xi_1^3\xi_2b_0^2k^3\delta_2(k)b_0\delta_1(k^2)b_0                        
-6\tau_1^2\xi_1^3\xi_2b_0^2k^2\delta_1(k)b_0k\delta_2(k^2)b_0\\                       
 -6\tau_1^2\xi_1^3\xi_2b_0^2k^2\delta_2(k)b_0k\delta_1(k^2)b_0                        
 -6\tau_1^2\xi_1^3\xi_2b_0^2k^2\delta_1(k^2)b_0k\delta_2(k)b_0\\                        
 -6\tau_1^2\xi_1^3\xi_2b_0^2k^2\delta_1(k^2)b_0\delta_2(k)b_0k                       
  -6\tau_1^2\xi_1^3\xi_2b_0^2k^2\delta_2(k^2)b_0k\delta_1(k)b_0\\                        
  -6\tau_1^2\xi_1^3\xi_2b_0^2k^2\delta_2(k^2)b_0\delta_1(k)b_0k                        
  -6\tau_1^2\xi_1^3\xi_2b_0\delta_1(k^2)b_0^2k^3\delta_2(k)b_0\\                        
  -6\tau_1^2\xi_1^3\xi_2b_0\delta_1(k^2)b_0^2k^2\delta_2(k)b_0k                       
   -6\tau_1^2\xi_1^3\xi_2b_0\delta_2(k^2)b_0^2k^3\delta_1(k)b_0\\                        
   -6\tau_1^2\xi_1^3\xi_2b_0\delta_2(k^2)b_0^2k^2\delta_1(k)b_0k                        
   -10\tau_1^2\xi_1^2\xi_2^2b_0^2k^3\delta_1(k)b_0\delta_1(k^2)b_0\\                        
   -10\tau_1^2\xi_1^2\xi_2^2b_0^2k^2\delta_1(k)b_0k\delta_1(k^2)b_0                        
   -10\tau_1^2\xi_1^2\xi_2^2b_0^2k^2\delta_1(k^2)b_0k\delta_1(k)b_0\\                        
   -10\tau_1^2\xi_1^2\xi_2^2b_0^2k^2\delta_1(k^2)b_0\delta_1(k)b_0k                        
   -10\tau_1^2\xi_1^2\xi_2^2b_0\delta_1(k^2)b_0^2k^3\delta_1(k)b_0\\                        
   -10\tau_1^2\xi_1^2\xi_2^2b_0\delta_1(k^2)b_0^2k^2\delta_1(k)b_0k                        
   -2i\tau_1\tau_2\xi_1^4b_0^2k^3\delta_2(k)b_0\delta_2(k^2)b_0\\                        
   -2i\tau_1\tau_2\xi_1^4b_0^2k^2\delta_2(k)b_0k\delta_2(k^2)b_0                        
   -2i\tau_1\tau_2\xi_1^4b_0^2k^2\delta_2(k^2)b_0k\delta_2(k)b_0\\                        
   -2i\tau_1\tau_2\xi_1^4b_0^2k^2\delta_2(k^2)b_0\delta_2(k)b_0k                        
   -2i\tau_1\tau_2\xi_1^4b_0\delta_2(k^2)b_0^2k^3\delta_2(k)b_0\\                        
   -2i\tau_1\tau_2\xi_1^4b_0\delta_2(k^2)b_0^2k^2\delta_2(k)b_0k                        
   -6i\tau_1\tau_2\xi_1^3\xi_2b_0^2k^3\delta_1(k)b_0\delta_2(k^2)b_0\\                        
   +2i\tau_1\tau_2\xi_1^3\xi_2b_0^2k^3\delta_2(k)b_0\delta_1(k^2)b_0                        
   -6i\tau_1\tau_2\xi_1^3\xi_2b_0^2k^2\delta_1(k)b_0k\delta_2(k^2)b_0\\                        
   +2i\tau_1\tau_2\xi_1^3\xi_2b_0^2k^2\delta_2(k)b_0k\delta_1(k^2)b_0                        
   -2i\tau_1\tau_2\xi_1^3\xi_2b_0^2k^2\delta_1(k^2)b_0k\delta_2(k)b_0\\                        
   +2i\tau_1\tau_2\xi_1^3\xi_2b_0^2k^2\delta_1(k^2)b_0\delta_2(k)b_0k                        
   -6i\tau_1\tau_2\xi_1^3\xi_2b_0^2k^2\delta_2(k^2)b_0k\delta_1(k)b_0\\                        
   -6i\tau_1\tau_2\xi_1^3\xi_2b_0^2k^2\delta_2(k^2)b_0\delta_1(k)b_0k                        
   +2i\tau_1\tau_2\xi_1^3\xi_2b_0\delta_1(k^2)b_0^2k^3\delta_2(k)b_0\\                        
   +2i\tau_1\tau_2\xi_1^3\xi_2b_0\delta_1(k^2)b_0^2k^2\delta_2(k)b_0k                        
   -6i\tau_1\tau_2\xi_1^3\xi_2b_0\delta_2(k^2)b_0^2k^3\delta_1(k)b_0\\                        
   -6i\tau_1\tau_2\xi_1^3\xi_2b_0\delta_2(k^2)b_0^2k^2\delta_1(k)b_0k                        
   +6i\tau_1\tau_2\xi_1^2\xi_2^2b_0^2k^3\delta_1(k)b_0\delta_1(k^2)b_0\\                        
   +6i\tau_1\tau_2\xi_1^2\xi_2^2b_0^2k^2\delta_1(k)b_0k\delta_1(k^2)b_0                        
   +6i\tau_1\tau_2\xi_1^2\xi_2^2b_0^2k^2\delta_1(k^2)b_0k\delta_1(k)b_0\\                        
   +6i\tau_1\tau_2\xi_1^2\xi_2^2b_0^2k^2\delta_1(k^2)b_0\delta_1(k)b_0k                        
   +6i\tau_1\tau_2\xi_1^2\xi_2^2b_0\delta_1(k^2)b_0^2k^3\delta_1(k)b_0\\                        
   +6i\tau_1\tau_2\xi_1^2\xi_2^2b_0\delta_1(k^2)b_0^2k^2\delta_1(k)b_0k                        
   -8\tau_1\xi_1^4b_0^3k^5\delta_1\delta_2(k)b_0\\                        
   -8\tau_1\xi_1^4b_0^3k^4\delta_1(k)\delta_2(k)b_0                        
   -8\tau_1\xi_1^4b_0^3k^4\delta_2(k)\delta_1(k)b_0\\                        
   -8\tau_1\xi_1^4b_0^3k^4\delta_1\delta_2(k)b_0k                        
   -6\tau_1\xi_1^4b_0^2k^3\delta_1(k)b_0k\delta_2(k)b_0\\                        
   -6\tau_1\xi_1^4b_0^2k^3\delta_1(k)b_0\delta_2(k)b_0k                        
   -6\tau_1\xi_1^4b_0^2k^3\delta_2(k)b_0k\delta_1(k)b_0\\                        
-6\tau_1\xi_1^4b_0^2k^3\delta_2(k)b_0\delta_1(k)b_0k                        
-6\tau_1\xi_1^4b_0^2k^2\delta_1(k)b_0k^2\delta_2(k)b_0\\                        
-6\tau_1\xi_1^4b_0^2k^2\delta_1(k)b_0k\delta_2(k)b_0k                        
-6\tau_1\xi_1^4b_0^2k^2\delta_2(k)b_0k^2\delta_1(k)b_0\\                        
-6\tau_1\xi_1^4b_0^2k^2\delta_2(k)b_0k\delta_1(k)b_0k                        
-16\tau_1\xi_1^3\xi_2b_0^3k^5\delta_1^2(k)b_0\\                        
-32\tau_1\xi_1^3\xi_2b_0^3k^4\delta_1(k)^2b_0                        
-16\tau_1\xi_1^3\xi_2b_0^3k^4\delta_1^2(k)b_0k\\                        
-24\tau_1\xi_1^3\xi_2b_0^2k^3\delta_1(k)b_0k\delta_1(k)b_0                        
-24\tau_1\xi_1^3\xi_2b_0^2k^3\delta_1(k)b_0\delta_1(k)b_0k\\                        
-24\tau_1\xi_1^3\xi_2b_0^2k^2\delta_1(k)b_0k^2\delta_1(k)b_0                        
-24\tau_1\xi_1^3\xi_2b_0^2k^2\delta_1(k)b_0k\delta_1(k)b_0k\\                        
-2|\tau|^2\xi_1^3\xi_2b_0^2k^3\delta_1(k)b_0\delta_2(k^2)b_0                        
-2|\tau|^2\xi_1^3\xi_2b_0^2k^3\delta_2(k)b_0\delta_1(k^2)b_0\\                        
-2|\tau|^2\xi_1^3\xi_2b_0^2k^2\delta_1(k)b_0k\delta_2(k^2)b_0                        
-2|\tau|^2\xi_1^3\xi_2b_0^2k^2\delta_2(k)b_0k\delta_1(k^2)b_0\\                        
-2|\tau|^2\xi_1^3\xi_2b_0^2k^2\delta_1(k^2)b_0k\delta_2(k)b_0                        
-2|\tau|^2\xi_1^3\xi_2b_0^2k^2\delta_1(k^2)b_0\delta_2(k)b_0k\\                        
-2|\tau|^2\xi_1^3\xi_2b_0^2k^2\delta_2(k^2)b_0k\delta_1(k)b_0                        
-2|\tau|^2\xi_1^3\xi_2b_0^2k^2\delta_2(k^2)b_0\delta_1(k)b_0k\\                        
-2|\tau|^2\xi_1^3\xi_2b_0\delta_1(k^2)b_0^2k^3\delta_2(k)b_0                        
-2|\tau|^2\xi_1^3\xi_2b_0\delta_1(k^2)b_0^2k^2\delta_2(k)b_0k\\                        
-2|\tau|^2\xi_1^3\xi_2b_0\delta_2(k^2)b_0^2k^3\delta_1(k)b_0                        
-2|\tau|^2\xi_1^3\xi_2b_0\delta_2(k^2)b_0^2k^2\delta_1(k)b_0k\\                        
-2|\tau|^2\xi_1^2\xi_2^2b_0^2k^3\delta_1(k)b_0\delta_1(k^2)b_0                        
-2|\tau|^2\xi_1^2\xi_2^2b_0^2k^2\delta_1(k)b_0k\delta_1(k^2)b_0\\                        
-2|\tau|^2\xi_1^2\xi_2^2b_0^2k^2\delta_1(k^2)b_0k\delta_1(k)b_0                        
-2|\tau|^2\xi_1^2\xi_2^2b_0^2k^2\delta_1(k^2)b_0\delta_1(k)b_0k\\                        
-2|\tau|^2\xi_1^2\xi_2^2b_0\delta_1(k^2)b_0^2k^3\delta_1(k)b_0                        
-2|\tau|^2\xi_1^2\xi_2^2b_0\delta_1(k^2)b_0^2k^2\delta_1(k)b_0k\\                        
-4\tau_1^3\xi_1^3\xi_2b_0^2k^3\delta_2(k)b_0\delta_2(k^2)b_0                        
-4\tau_1^3\xi_1^3\xi_2b_0^2k^2\delta_2(k)b_0k\delta_2(k^2)b_0\\                        
-4\tau_1^3\xi_1^3\xi_2b_0^2k^2\delta_2(k^2)b_0k\delta_2(k)b_0                        
-4\tau_1^3\xi_1^3\xi_2b_0^2k^2\delta_2(k^2)b_0\delta_2(k)b_0k\\                        
-4\tau_1^3\xi_1^3\xi_2b_0\delta_2(k^2)b_0^2k^3\delta_2(k)b_0                        
-4\tau_1^3\xi_1^3\xi_2b_0\delta_2(k^2)b_0^2k^2\delta_2(k)b_0k\\                        
-4\tau_1^3\xi_1^2\xi_2^2b_0^2k^3\delta_1(k)b_0\delta_2(k^2)b_0                        
-4\tau_1^3\xi_1^2\xi_2^2b_0^2k^3\delta_2(k)b_0\delta_1(k^2)b_0\\                        
-4\tau_1^3\xi_1^2\xi_2^2b_0^2k^2\delta_1(k)b_0k\delta_2(k^2)b_0                        
-4\tau_1^3\xi_1^2\xi_2^2b_0^2k^2\delta_2(k)b_0k\delta_1(k^2)b_0\\                        
-4\tau_1^3\xi_1^2\xi_2^2b_0^2k^2\delta_1(k^2)b_0k\delta_2(k)b_0                        
-4\tau_1^3\xi_1^2\xi_2^2b_0^2k^2\delta_1(k^2)b_0\delta_2(k)b_0k\\                        
-4\tau_1^3\xi_1^2\xi_2^2b_0^2k^2\delta_2(k^2)b_0k\delta_1(k)b_0                        
-4\tau_1^3\xi_1^2\xi_2^2b_0^2k^2\delta_2(k^2)b_0\delta_1(k)b_0k\\                        
-4\tau_1^3\xi_1^2\xi_2^2b_0\delta_1(k^2)b_0^2k^3\delta_2(k)b_0                        
-4\tau_1^3\xi_1^2\xi_2^2b_0\delta_1(k^2)b_0^2k^2\delta_2(k)b_0k\\                        
-4\tau_1^3\xi_1^2\xi_2^2b_0\delta_2(k^2)b_0^2k^3\delta_1(k)b_0                        
-4\tau_1^3\xi_1^2\xi_2^2b_0\delta_2(k^2)b_0^2k^2\delta_1(k)b_0k\\                        
-4\tau_1^3\xi_1\xi_2^3b_0^2k^3\delta_1(k)b_0\delta_1(k^2)b_0                        
-4\tau_1^3\xi_1\xi_2^3b_0^2k^2\delta_1(k)b_0k\delta_1(k^2)b_0\\                        
-4\tau_1^3\xi_1\xi_2^3b_0^2k^2\delta_1(k^2)b_0k\delta_1(k)b_0                        
-4\tau_1^3\xi_1\xi_2^3b_0^2k^2\delta_1(k^2)b_0\delta_1(k)b_0k\\                       
-4\tau_1^3\xi_1\xi_2^3b_0\delta_1(k^2)b_0^2k^3\delta_1(k)b_0                        
-4\tau_1^3\xi_1\xi_2^3b_0\delta_1(k^2)b_0^2k^2\delta_1(k)b_0k\\                        
-4i\tau_1^2\tau_2\xi_1^3\xi_2b_0^2k^3\delta_2(k)b_0\delta_2(k^2)b_0                        
-4i\tau_1^2\tau_2\xi_1^3\xi_2b_0^2k^2\delta_2(k)b_0k\delta_2(k^2)b_0\\                        
-4i\tau_1^2\tau_2\xi_1^3\xi_2b_0^2k^2\delta_2(k^2)b_0k\delta_2(k)b_0                        
-4i\tau_1^2\tau_2\xi_1^3\xi_2b_0^2k^2\delta_2(k^2)b_0\delta_2(k)b_0k\\                        
-4i\tau_1^2\tau_2\xi_1^3\xi_2b_0\delta_2(k^2)b_0^2k^3\delta_2(k)b_0                        
-4i\tau_1^2\tau_2\xi_1^3\xi_2b_0\delta_2(k^2)b_0^2k^2\delta_2(k)b_0k\\                        
-4i\tau_1^2\tau_2\xi_1^2\xi_2^2b_0^2k^3\delta_1(k)b_0\delta_2(k^2)b_0                        
+4i\tau_1^2\tau_2\xi_1^2\xi_2^2b_0^2k^3\delta_2(k)b_0\delta_1(k^2)b_0\\                        
-4i\tau_1^2\tau_2\xi_1^2\xi_2^2b_0^2k^2\delta_1(k)b_0k\delta_2(k^2)b_0                        
+4i\tau_1^2\tau_2\xi_1^2\xi_2^2b_0^2k^2\delta_2(k)b_0k\delta_1(k^2)b_0\\                        
+4i\tau_1^2\tau_2\xi_1^2\xi_2^2b_0^2k^
2\delta_1(k^2)b_0k\delta_2(k)b_0                        
+4i\tau_1^2\tau_2\xi_1^2\xi_2^2b_0^2k^2\delta_1(k^2)b_0\delta_2(k)b_0k\\                        
-4i\tau_1^2\tau_2\xi_1^2\xi_2^2b_0^2k^2\delta_2(k^2)b_0k\delta_1(k)b_0                       
-4i\tau_1^2\tau_2\xi_1^2\xi_2^2b_0^2k^2\delta_2(k^2)b_0\delta_1(k)b_0k\\                        
+4i\tau_1^2\tau_2\xi_1^2\xi_2^2b_0\delta_1(k^2)b_0^2k^3\delta_2(k)b_0                        
+4i\tau_1^2\tau_2\xi_1^2\xi_2^2b_0\delta_1(k^2)b_0^2k^2\delta_2(k)b_0k\\                        
-4i\tau_1^2\tau_2\xi_1^2\xi_2^2b_0\delta_2(k^2)b_0^2k^3\delta_1(k)b_0                       
-4i\tau_1^2\tau_2\xi_1^2\xi_2^2b_0\delta_2(k^2)b_0^2k^2\delta_1(k)b_0k\\                        
+4i\tau_1^2\tau_2\xi_1\xi_2^3b_0^2k^3\delta_1(k)b_0\delta_1(k^2)b_0                        
+4i\tau_1^2\tau_2\xi_1\xi_2^3b_0^2k^2\delta_1(k)b_0k\delta_1(k^2)b_0\\                        
+4i\tau_1^2\tau_2\xi_1\xi_2^3b_0^2k^2\delta_1(k^2)b_0k\delta_1(k)b_0                       
+4i\tau_1^2\tau_2\xi_1\xi_2^3b_0^2k^2\delta_1(k^2)b_0\delta_1(k)b_0k\\                        
+4i\tau_1^2\tau_2\xi_1\xi_2^3b_0\delta_1(k^2)b_0^2k^3\delta_1(k)b_0                        
+4i\tau_1^2\tau_2\xi_1\xi_2^3b_0\delta_1(k^2)b_0^2k^2\delta_1(k)b_0k\\                        
-4\tau_1^2\xi_1^4b_0^3k^5\delta_2^2(k)b_0                        
-8\tau_1^2\xi_1^4b_0^3k^4\delta_2(k)^2b_0\\                        
-4\tau_1^2\xi_1^4b_0^3k^4\delta_2^2(k)b_0k                        
-4\tau_1^2\xi_1^4b_0^2k^3\delta_2(k)b_0k\delta_2(k)b_0\\                        
-4\tau_1^2\xi_1^4b_0^2k^3\delta_2(k)b_0\delta_2(k)b_0k                        
-4\tau_1^2\xi_1^4b_0^2k^2\delta_2(k)b_0k^2\delta_2(k)b_0\\                        
-4\tau_1^2\xi_1^4b_0^2k^2\delta_2(k)b_0k\delta_2(k)b_0k                        
-24\tau_1^2\xi_1^3\xi_2b_0^3k^5\delta_1\delta_2(k)b_0\\                        
-24\tau_1^2\xi_1^3\xi_2b_0^3k^4\delta_1(k)\delta_2(k)b_0                        
-24\tau_1^2\xi_1^3\xi_2b_0^3k^4\delta_2(k)\delta_1(k)b_0\\                        
-24\tau_1^2\xi_1^3\xi_2b_0^3k^4\delta_1\delta_2(k)b_0k                        
-20\tau_1^2\xi_1^3\xi_2b_0^2k^3\delta_1(k)b_0k\delta_2(k)b_0\\                        
-20\tau_1^2\xi_1^3\xi_2b_0^2k^3\delta_1(k)b_0\delta_2(k)b_0k                        
-20\tau_1^2\xi_1^3\xi_2b_0^2k^3\delta_2(k)b_0k\delta_1(k)b_0\\                        
-20\tau_1^2\xi_1^3\xi_2b_0^2k^3\delta_2(k)b_0\delta_1(k)b_0k                        
-20\tau_1^2\xi_1^3\xi_2b_0^2k^2\delta_1(k)b_0k^2\delta_2(k)b_0\\                        
-20\tau_1^2\xi_1^3\xi_2b_0^2k^2\delta_1(k)b_0k\delta_2(k)b_0k                        
-20\tau_1^2\xi_1^3\xi_2b_0^2k^2\delta_2(k)b_0k^2\delta_1(k)b_0\\                        
-20\tau_1^2\xi_1^3\xi_2b_0^2k^2\delta_2(k)b_0k\delta_1(k)b_0k                        
-20\tau_1^2\xi_1^2\xi_2^2b_0^3k^5\delta_1^2(k)b_0\\                        
-40\tau_1^2\xi_1^2\xi_2^2b_0^3k^4\delta_1(k)^2b_0                        
-20\tau_1^2\xi_1^2\xi_2^2b_0^3k^4\delta_1^2(k)b_0k\\                        
-28\tau_1^2\xi_1^2\xi_2^2b_0^2k^3\delta_1(k)b_0k\delta_1(k)b_0                        
-28\tau_1^2\xi_1^2\xi_2^2b_0^2k^3\delta_1(k)b_0\delta_1(k)b_0k\\                        
-28\tau_1^2\xi_1^2\xi_2^2b_0^2k^2\delta_1(k)b_0k^2\delta_1(k)b_0                        
-28\tau_1^2\xi_1^2\xi_2^2b_0^2k^2\delta_1(k)b_0k\delta_1(k)b_0k\\                        
-4\tau_1|\tau|^2\xi_1^3\xi_2b_0^2k^3\delta_2(k)b_0\delta_2(k^2)b_0                        
-4\tau_1|\tau|^2\xi_1^3\xi_2b_0^2k^2\delta_2(k)b_0k\delta_2(k^2)b_0\\                        
-4\tau_1|\tau|^2\xi_1^3\xi_2b_0^2k^2\delta_2(k^2)b_0k\delta_2(k)b_0                        
-4\tau_1|\tau|^2\xi_1^3\xi_2b_0^2k^2\delta_2(k^2)b_0\delta_2(k)b_0k\\                        
-4\tau_1|\tau|^2\xi_1^3\xi_2b_0\delta_2(k^2)b_0^2k^3\delta_2(k)b_0                      
-4\tau_1|\tau|^2\xi_1^3\xi_2b_0\delta_2(k^2)b_0^2k^2\delta_2(k)b_0k\\                        
-8\tau_1|\tau|^2\xi_1^2\xi_2^2b_0^2k^3\delta_1(k)b_0\delta_2(k^2)b_0                        
-8\tau_1|\tau|^2\xi_1^2\xi_2^2b_0^2k^3\delta_2(k)b_0\delta_1(k^2)b_0\\                        
-8\tau_1|\tau|^2\xi_1^2\xi_2^2b_0^2k^2\delta_1(k)b_0k\delta_2(k^2)b_0                        
-8\tau_1|\tau|^2\xi_1^2\xi_2^2b_0^2k^2\delta_2(k)b_0k\delta_1(k^2)b_0\\                        
-8\tau_1|\tau|^2\xi_1^2\xi_2^2b_0^2k^2\delta_1(k^2)b_0k\delta_2(k)b_0                        
-8\tau_1|\tau|^2\xi_1^2\xi_2^2b_0^2k^2\delta_1(k^2)b_0\delta_2(k)b_0k\\                        
-8\tau_1|\tau|^2\xi_1^2\xi_2^2b_0^2k^2\delta_2(k^2)b_0k\delta_1(k)b_0                        
-8\tau_1|\tau|^2\xi_1^2\xi_2^2b_0^2k^2\delta_2(k^2)b_0\delta_1(k)b_0k\\                        
-8\tau_1|\tau|^2\xi_1^2\xi_2^2b_0\delta_1(k^2)b_0^2k^3\delta_2(k)b_0                        
-8\tau_1|\tau|^2\xi_1^2\xi_2^2b_0\delta_1(k^2)b_0^2k^2\delta_2(k)b_0k\\                        
-8\tau_1|\tau|^2\xi_1^2\xi_2^2b_0\delta_2(k^2)b_0^2k^3\delta_1(k)b_0                        
-8\tau_1|\tau|^2\xi_1^2\xi_2^2b_0\delta_2(k^2)b_0^2k^2\delta_1(k)b_0k\\                        
-4\tau_1|\tau|^2\xi_1\xi_2^3b_0^2k^3\delta_1(k)b_0\delta_1(k^2)b_0                        
-4\tau_1|\tau|^2\xi_1\xi_2^3b_0^2k^2\delta_1(k)b_0k\delta_1(k^2)b_0\\                        
-4\tau_1|\tau|^2\xi_1\xi_2^3b_0^2k^2\delta_1(k^2)b_0k\delta_1(k)b_0                        
-4\tau_1|\tau|^2\xi_1\xi_2^3b_0^2k^2\delta_1(k^2)b_0\delta_1(k)b_0k\\                        
-4\tau_1|\tau|^2\xi_1\xi_2^3b_0\delta_1(k^2)b_0^2k^3\delta_1(k)b_0                        
-4\tau_1|\tau|^2\xi_1\xi_2^3b_0\delta_1(k^2)b_0^2k^2\delta_1(k)b_0k\\                        
-2i\tau_2|\tau|^2\xi_1^3\xi_2b_0^2k^3\delta_2(k)b_0\delta_2(k^2)b_0                        
-2i\tau_2|\tau|^2\xi_1^3\xi_2b_0^2k^2\delta_2(k)b_0k\delta_2(k^2)b_0\\                        
-2i\tau_2|\tau|^2\xi_1^3\xi_2b_0^2k^2\delta_2(k^2)b_0k\delta_2(k)b_0                        
-2i\tau_2|\tau|^2\xi_1^3\xi_2b_0^2k^2\delta_2(k^2)b_0\delta_2(k)b_0k\\                        
-2i\tau_2|\tau|^2\xi_1^3\xi_2b_0\delta_2(k^2)b_0^2k^3\delta_2(k)b_0                        
-2i\tau_2|\tau|^2\xi_1^3\xi_2b_0\delta_2(k^2)b_0^2k^2\delta_2(k)b_0k\\                        
-2i\tau_2|\tau|^2\xi_1^2\xi_2^2b_0^2k^3\delta_1(k)b_0\delta_2(k^2)b_0                        
+2i\tau_2|\tau|^2\xi_1^2\xi_2^2b_0^2k^3\delta_2(k)b_0\delta_1(k^2)b_0\\                        
-2i\tau_2|\tau|^2\xi_1^2\xi_2^2b_0^2k^2\delta_1(k)b_0k\delta_2(k^2)b_0                        
+2i\tau_2|\tau|^2\xi_1^2\xi_2^2b_0^2k^2\delta_2(k)b_0k\delta_1(k^2)b_0\\                        
+2i\tau_2|\tau|^2\xi_1^2\xi_2^2b_0^2k^2\delta_1(k^2)b_0k\delta_2(k)b_0                        
+2i\tau_2|\tau|^2\xi_1^2\xi_2^2b_0^2k^2\delta_1(k^2)b_0\delta_2(k)b_0k\\                        
+2i\tau_2|\tau|^2\xi_1^2\xi_2^2b_0^2k^2\delta_2(k^2)b_0k\delta_1(k)b_0                        
-2i\tau_2|\tau|^2\xi_1^2\xi_2^2b_0^2k^2\delta_2(k^2)b_0\delta_1(k)b_0k\\                        
+2i\tau_2|\tau|^2\xi_1^2\xi_2^2b_0\delta_1(k^2)b_0^2k^3\delta_2(k)b_0                        
+2i\tau_2|\tau|^2\xi_1^2\xi_2^2b_0\delta_1(k^2)b_0^2k^2\delta_2(k)b_0k\\                        
-2i\tau_2|\tau|^2\xi_1^2\xi_2^2b_0\delta_2(k^2)b_0^2k^3\delta_1(k)b_0                        
-2i\tau_2|\tau|^2\xi_1^2\xi_2^2b_0\delta_2(k^2)b_0^2k^2\delta_1(k)b_0k\\                        
+2i\tau_2|\tau|^2\xi_1\xi_2^3b_0^2k^3\delta_1(k)b_0\delta_1(k^2)b_0                        
+2i\tau_2|\tau|^2\xi_1\xi_2^3b_0^2k^2\delta_1(k)b_0k\delta_1(k^2)b_0\\                        
+2i\tau_2|\tau|^2\xi_1\xi_2^3b_0^2k^2\delta_1(k^2)b_0k\delta_1(k)b_0                        
+2i\tau_2|\tau|^2\xi_1\xi_2^3b_0^2k^2\delta_1(k^2)b_0\delta_1(k)b_0k\\                        
+2i\tau_2|\tau|^2\xi_1\xi_2^3b_0\delta_1(k^2)b_0^2k^3\delta_1(k)b_0                        
+2i\tau_2|\tau|^2\xi_1\xi_2^3b_0\delta_1(k^2)b_0^2k^2\delta_1(k)b_0k\\                        
-2|\tau|^2\xi_1^4b_0^2k^3\delta_2(k)b_0k\delta_2(k)b_0                        
-2|\tau|^2\xi_1^4b_0^2k^3\delta_2(k)b_0\delta_2(k)b_0k\\                        
-2|\tau|^2\xi_1^4b_0^2k^2\delta_2(k)b_0k^2\delta_2(k)b_0                       
-2|\tau|^2\xi_1^4b_0^2k^2\delta_2(k)b_0k\delta_2(k)b_0k\\                        
-8|\tau|^2\xi_1^3\xi_2b_0^3k^5\delta_1\delta_2(k)b_0                        
-8|\tau|^2\xi_1^3\xi_2b_0^3k^4\delta_1(k)\delta_2(k)b_0\\                        
-8|\tau|^2\xi_1^3\xi_2b_0^3k^4\delta_2(k)\delta_1(k)b_0                        
-8|\tau|^2\xi_1^3\xi_2b_0^3k^4\delta_1\delta_2(k)b_0k\\                        
-4|\tau|^2\xi_1^3\xi_2b_0^2k^3\delta_1(k)b_0k\delta_2(k)b_0                        
-4|\tau|^2\xi_1^3\xi_2b_0^2k^3\delta_1(k)b_0\delta_2(k)b_0k\\                        
-4|\tau|^2\xi_1^3\xi_2b_0^2k^3\delta_2(k)b_0k\delta_1(k)b_0                        
-4|\tau|^2\xi_1^3\xi_2b_0^2k^3\delta_2(k)b_0\delta_1(k)b_0k\\                        
-4|\tau|^2\xi_1^3\xi_2b_0^2k^2\delta_1(k)b_0k^2\delta_2(k)b_0                        
-4|\tau|^2\xi_1^3\xi_2b_0^2k^2\delta_1(k)b_0k\delta_2(k)b_0k\\                        
-4|\tau|^2\xi_1^3\xi_2b_0^2k^2\delta_2(k)b_0k^2\delta_1(k)b_0                        
-4|\tau|^2\xi_1^3\xi_2b_0^2k^2\delta_2(k)b_0k\delta_1(k)b_0k\\                        
-4|\tau|^2\xi_1^2\xi_2^2b_0^3k^5\delta_1^2(k)b_0                        
-8|\tau|^2\xi_1^2\xi_2^2b_0^3k^4\delta_1(k)^2b_0\\                                                
-8|\tau|^2\xi_1^2\xi_2^2b_0^2k^3\delta_1(k)b_0k\delta_1(k)b_0                        
-8|\tau|^2\xi_1^2\xi_2^2b_0^2k^3\delta_1(k)b_0\delta_1(k)b_0k\\                        
-8|\tau|^2\xi_1^2\xi_2^2b_0^2k^2\delta_1(k)b_0k^2\delta_1(k)b_0                        
-8|\tau|^2\xi_1^2\xi_2^2b_0^2k^2\delta_1(k)b_0k\delta_1(k)b_0k\\                        
-8\tau_1^3\xi_1^3\xi_2b_0^3k^5\delta_2^2(k)b_0                        
-16\tau_1^3\xi_1^3\xi_2b_0^3k^4\delta_2(k)^2b_0\\                        
-8\tau_1^3\xi_1^3\xi_2b_0^3k^4\delta_2^2(k)b_0k                        
-8\tau_1^3\xi_1^3\xi_2b_0^2k^3\delta_2(k)b_0k\delta_2(k)b_0\\                       
 -8\tau_1^3\xi_1^3\xi_2b_0^2k^3\delta_2(k)b_0\delta_2(k)b_0k                       
 -8\tau_1^3\xi_1^3\xi_2b_0^2k^2\delta_2(k)b_0k^2\delta_2(k)b_0\\                        
 -8\tau_1^3\xi_1^3\xi_2b_0^2k^2\delta_2(k)b_0k\delta_2(k)b_0k                       
 -16\tau_1^3\xi_1^2\xi_2^2b_0^3k^5\delta_1\delta_2(k)b_0\\                        
 -16\tau_1^3\xi_1^2\xi_2^2b_0^3k^4\delta_1(k)\delta_2(k)b_0                        
 -16\tau_1^3\xi_1^2\xi_2^2b_0^3k^4\delta_2(k)\delta_1(k)b_0\\                        
 -16\tau_1^3\xi_1^2\xi_2^2b_0^3k^4\delta_1\delta_2(k)b_0k                        
 -16\tau_1^3\xi_1^2\xi_2^2b_0^2k^3\delta_1(k)b_0k\delta_2(k)b_0\\                        
 -16\tau_1^3\xi_1^2\xi_2^2b_0^2k^3\delta_1(k)b_0\delta_2(k)b_0k                       
 -16\tau_1^3\xi_1^2\xi_2^2b_0^2k^3\delta_2(k)b_0k\delta_1(k)b_0\\                        
 -16\tau_1^3\xi_1^2\xi_2^2b_0^2k^3\delta_2(k)b_0\delta_1(k)b_0k                       
 -16\tau_1^3\xi_1^2\xi_2^2b_0^2k^2\delta_1(k)b_0k^2\delta_2(k)b_0\\                        
 -16\tau_1^3\xi_1^2\xi_2^2b_0^2k^2\delta_1(k)b_0k\delta_2(k)b_0k                     
  -16\tau_1^3\xi_1^2\xi_2^2b_0^2k^2\delta_2(k)b_0k^2\delta_1(k)b_0\\                        
  -16\tau_1^3\xi_1^2\xi_2^2b_0^2k^2\delta_2(k)b_0k\delta_1(k)b_0k                        
  -8\tau_1^3\xi_1\xi_2^3b_0^3k^5\delta_1^2(k)b_0\\                        
  -16\tau_1^3\xi_1\xi_2^3b_0^3k^4\delta_1(k)^2b_0                     
   -8\tau_1^3\xi_1\xi_2^3b_0^3k^4\delta_1^2(k)b_0k\\                        
   -8\tau_1^3\xi_1\xi_2^3b_0^2k^3\delta_1(k)b_0k\delta_1(k)b_0                        
   -8\tau_1^3\xi_1\xi_2^3b_0^2k^3\delta_1(k)b_0\delta_1(k)b_0k\\                        
   -8\tau_1^3\xi_1\xi_2^3b_0^2k^2\delta_1(k)b_0k^2\delta_1(k)b_0
-4|\tau|^2\xi_1^2\xi_2^2b_0^3k^4\delta_1^2(k)b_0k\\                        
-8\tau_1^3\xi_1\xi_2^3b_0^2k^2\delta_1(k)b_0k\delta_1(k)b_0k                        
-10\tau_1^2|\tau|^2\xi_1^2\xi_2^2b_0^2k^3\delta_2(k)b_0\delta_2(k^2)b_0\\                        
-10\tau_1^2|\tau|^2\xi_1^2\xi_2^2b_0^2k^2\delta_2(k)b_0k\delta_2(k^2)b_0                       
 -10\tau_1^2|\tau|^2\xi_1^2\xi_2^2b_0^2k^2\delta_2(k^2)b_0k\delta_2(k)b_0\\                        
 -10\tau_1^2|\tau|^2\xi_1^2\xi_2^2b_0^2k^2\delta_2(k^2)b_0\delta_2(k)b_0k                      
 -10\tau_1^2|\tau|^2\xi_1^2\xi_2^2b_0\delta_2(k^2)b_0^2k^3\delta_2(k)b_0\\                        
 -10\tau_1^2|\tau|^2\xi_1^2\xi_2^2b_0\delta_2(k^2)b_0^2k^2\delta_2(k)b_0k                       
 -6\tau_1^2|\tau|^2\xi_1\xi_2^3b_0^2k^3\delta_1(k)b_0\delta_2(k^2)b_0\\                        
 -6\tau_1^2|\tau|^2\xi_1\xi_2^3b_0^2k^3\delta_2(k)b_0\delta_1(k^2)b_0                        
 -6\tau_1^2|\tau|^2\xi_1\xi_2^3b_0^2k^2\delta_1(k)b_0k\delta_2(k^2)b_0\\                        
 -6\tau_1^2|\tau|^2\xi_1\xi_2^3b_0^2k^2\delta_2(k)b_0k\delta_1(k^2)b_0                        
 -6\tau_1^2|\tau|^2\xi_1\xi_2^3b_0^2k^2\delta_1(k^2)b_0k\delta_2(k)b_0\\                        
 -6\tau_1^2|\tau|^2\xi_1\xi_2^3b_0^2k^2\delta_1(k^2)b_0\delta_2(k)b_0k                        
 -6\tau_1^2|\tau|^2\xi_1\xi_2^3b_0^2k^2\delta_2(k^2)b_0k\delta_1(k)b_0\\                        
 -6\tau_1^2|\tau|^2\xi_1\xi_2^3b_0^2k^2\delta_2(k^2)b_0\delta_1(k)b_0k                        
 -6\tau_1^2|\tau|^2\xi_1\xi_2^3b_0\delta_1(k^2)b_0^2k^3\delta_2(k)b_0\\                        
 -6\tau_1^2|\tau|^2\xi_1\xi_2^3b_0\delta_1(k^2)b_0^2k^2\delta_2(k)b_0k                        
 -6\tau_1^2|\tau|^2\xi_1\xi_2^3b_0\delta_2(k^2)b_0^2k^3\delta_1(k)b_0\\                        
 -6\tau_1^2|\tau|^2\xi_1\xi_2^3b_0\delta_2(k^2)b_0^2k^2\delta_1(k)b_0k                        
 -2\tau_1^2|\tau|^2\xi_2^4b_0^2k^3\delta_1(k)b_0\delta_1(k^2)b_0\\                       
  -2\tau_1^2|\tau|^2\xi_2^4b_0^2k^2\delta_1(k)b_0k\delta_1(k^2)b_0                       
   -2\tau_1^2|\tau|^2\xi_2^4b_0^2k^2\delta_1(k^2)b_0k\delta_1(k)b_0\\                        
   -2\tau_1^2|\tau|^2\xi_2^4b_0^2k^2\delta_1(k^2)b_0\delta_1(k)b_0k                        
   -2\tau_1^2|\tau|^2\xi_2^4b_0\delta_1(k^2)b_0^2k^3\delta_1(k)b_0\\                        
   -2\tau_1^2|\tau|^2\xi_2^4b_0\delta_1(k^2)b_0^2k^2\delta_1(k)b_0k                       
    -6i\tau_1\tau_2|\tau|^2\xi_1^2\xi_2^2b_0^2k^3\delta_2(k)b_0\delta_2(k^2)b_0\\                        
    -6i\tau_1\tau_2|\tau|^2\xi_1^2\xi_2^2b_0^2k^2\delta_2(k)b_0k\delta_2(k^2)b_0                        
    -6i\tau_1\tau_2|\tau|^2\xi_1^2\xi_2^2b_0^2k^2\delta_2(k^2)b_0k\delta_2(k)b_0\\                        
    -6i\tau_1\tau_2|\tau|^2\xi_1^2\xi_2^2b_0^2k^2\delta_2(k^2)b_0\delta_2(k)b_0k                        
    -6i\tau_1\tau_2|\tau|^2\xi_1^2\xi_2^2b_0\delta_2(k^2)b_0^2k^3\delta_2(k)b_0\\                        
-6i\tau_1\tau_2|\tau|^2\xi_1^2\xi_2^2b_0\delta_2(k^2)b_0^2k^2\delta_2(k)b_0k                        
-2i\tau_1\tau_2|\tau|^2\xi_1\xi_2^3b_0^2k^3\delta_1(k)b_0\delta_2(k^2)b_0\\                        
+6i\tau_1\tau_2|\tau|^2\xi_1\xi_2^3b_0^2k^3\delta_2(k)b_0\delta_1(k^2)b_0                        
-2i\tau_1\tau_2|\tau|^2\xi_1\xi_2^3b_0^2k^2\delta_1(k)b_0k\delta_2(k^2)b_0\\                        
+6i\tau_1\tau_2|\tau|^2\xi_1\xi_2^3b_0^2k^
2\delta_2(k)b_0k\delta_1(k^2)b_0                        
+6i\tau_1\tau_2|\tau|^2\xi_1\xi_2^3b_0^2k^2\delta_1(k^2)b_0k\delta_2(k)b_0\\                        
+6i\tau_1\tau_2|\tau|^2\xi_1\xi_2^3b_0^2k^2\delta_1(k^2)b_0\delta_2(k)b_0k                        
-2i\tau_1\tau_2|\tau|^2\xi_1\xi_2^3b_0^2k^2\delta_2(k^2)b_0k\delta_1(k)b_0\\                        
-2i\tau_1\tau_2|\tau|^2\xi_1\xi_2^3b_0^2k^2\delta_2(k^2)b_0\delta_1(k)b_0k                        
+6i\tau_1\tau_2|\tau|^2\xi_1\xi_2^3b_0\delta_1(k^2)b_0^2k^3\delta_2(k)b_0\\                        
+6i\tau_1\tau_2|\tau|^2\xi_1\xi_2^3b_0\delta_1(k^2)b_0^2k^2\delta_2(k)b_0k                        
-2i\tau_1\tau_2|\tau|^2\xi_1\xi_2^3b_0\delta_2(k^2)b_0^2k^3\delta_1(k)b_0\\                        
-2i\tau_1\tau_2|\tau|^2\xi_1\xi_2^3b_0\delta_2(k^2)b_0^2k^2\delta_1(k)b_0k                        
+2i\tau_1\tau_2|\tau|^2\xi_2^4b_0^2k^3\delta_1(k)b_0\delta_1(k^2)b_0\\                        
+2i\tau_1\tau_2|\tau|^2\xi_2^4b_0^2k^2\delta_1(k)b_0k\delta_1(k^2)b_0                       
+2i\tau_1\tau_2|\tau|^2\xi_2^4b_0^2k^2\delta_1(k^2)b_0k\delta_1(k)b_0\\                        
+2i\tau_1\tau_2|\tau|^2\xi_2^4b_0^2k^2\delta_1(k^2)b_0\delta_1(k)b_0k                        
+2i\tau_1\tau_2|\tau|^2\xi_2^4b_0\delta_1(k^2)b_0^2k^3\delta_1(k)b_0\\                        
+2i\tau_1\tau_2|\tau|^2\xi_2^4b_0\delta_1(k^2)b_0^2k^2\delta_1(k)b_0k                        
-8\tau_1|\tau|^2\xi_1^3\xi_2b_0^3k^5\delta_2^2(k)b_0\\                        
-16\tau_1|\tau|^2\xi_1^3\xi_2b_0^3k^4\delta_2(k)^2b_0                       
-8\tau_1|\tau|^2\xi_1^3\xi_2b_0^3k^4\delta_2^2(k)b_0k\\                        
-16\tau_1|\tau|^2\xi_1^3\xi_2b_0^2k^3\delta_2(k)b_0k\delta_2(k)b_0                        
-16\tau_1|\tau|^2\xi_1^3\xi_2b_0^2k^3\delta_2(k)b_0\delta_2(k)b_0k\\                        
-16\tau_1|\tau|^2\xi_1^3\xi_2b_0^2k^2\delta_2(k)b_0k^2\delta_2(k)b_0                       
-16\tau_1|\tau|^2\xi_1^3\xi_2b_0^2k^2\delta_2(k)b_0k\delta_2(k)b_0k\\                        
-32\tau_1|\tau|^2\xi_1^2\xi_2^2b_0^3k^5\delta_1\delta_2(k)b_0                        
-32\tau_1|\tau|^2\xi_1^2\xi_2^2b_0^3k^4\delta_1(k)\delta_2(k)b_0\\                        
-32\tau_1|\tau|^2\xi_1^2\xi_2^2b_0^3k^4\delta_2(k)\delta_1(k)b_0                        
-32\tau_1|\tau|^2\xi_1^2\xi_2^2b_0^3k^4\delta_1\delta_2(k)b_0k\\                        
-20\tau_1|\tau|^2\xi_1^2\xi_2^2b_0^2k^3\delta_1(k)b_0k\delta_2(k)b_0                        
-20\tau_1|\tau|^2\xi_1^2\xi_2^2b_0^2k^3\delta_1(k)b_0\delta_2(k)b_0k\\                        
-20\tau_1|\tau|^2\xi_1^2\xi_2^2b_0^2k^3\delta_2(k)b_0k\delta_1(k)b_0                        
-20\tau_1|\tau|^2\xi_1^2\xi_2^2b_0^2k^3\delta_2(k)b_0\delta_1(k)b_0k\\                        
-20\tau_1|\tau|^2\xi_1^2\xi_2^2b_0^2k^2\delta_1(k)b_0k^2\delta_2(k)b_0                        
-20\tau_1|\tau|^2\xi_1^2\xi_2^2b_0^2k^2\delta_1(k)b_0k\delta_2(k)b_0k\\                        
-20\tau_1|\tau|^2\xi_1^2\xi_2^2b_0^2k^2\delta_2(k)b_0k^2\delta_1(k)b_0                        
-20\tau_1|\tau|^2\xi_1^2\xi_2^2b_0^2k^2\delta_2(k)b_0k\delta_1(k)b_0k\\                        
-8\tau_1|\tau|^2\xi_1\xi_2^3b_0^3k^5\delta_1^2(k)b_0                        
-16\tau_1|\tau|^2\xi_1\xi_2^3b_0^3k^4\delta_1(k)^2b_0\\                        
-8\tau_1|\tau|^2\xi_1\xi_2^3b_0^3k^4\delta_1^2(k)b_0k                        
-16\tau_1|\tau|^2\xi_1\xi_2^3b_0^2k^3\delta_1(k)b_0k\delta_1(k)b_0\\                        
-16\tau_1|\tau|^2\xi_1\xi_2^3b_0^2k^3\delta_1(k)b_0\delta_1(k)b_0k                        
-16\tau_1|\tau|^2\xi_1\xi_2^3b_0^2k^2\delta_1(k)b_0k^2\delta_1(k)b_0\\                        
-16\tau_1|\tau|^2\xi_1\xi_2^3b_0^2k^2\delta_1(k)b_0k\delta_1(k)b_0k                        
-2|\tau|^4\xi_1^2\xi_2^2b_0^2k^3\delta_2(k)b_0\delta_2(k^2)b_0\\                        
-2|\tau|^4\xi_1^2\xi_2^2b_0^2k^2\delta_2(k)b_0k\delta_2(k^2)b_0                        
-2|\tau|^4\xi_1^2\xi_2^2b_0^2k^2\delta_2(k^2)b_0k\delta_2(k)b_0\\                        
-2|\tau|^4\xi_1^2\xi_2^2b_0^2k^2\delta_2(k^2)b_0\delta_2(k)b_0k                        
-2|\tau|^4\xi_1^2\xi_2^2b_0\delta_2(k^2)b_0^2k^3\delta_2(k)b_0\\                        
-2|\tau|^4\xi_1^2\xi_2^2b_0\delta_2(k^2)b_0^2k^2\delta_2(k)b_0k                        
-2|\tau|^4\xi_1\xi_2^3b_0^2k^3\delta_1(k)b_0\delta_2(k^2)b_0\\                        
-2|\tau|^4\xi_1\xi_2^3b_0^2k^3\delta_2(k)b_0\delta_1(k^2)b_0                       
 -2|\tau|^4\xi_1\xi_2^3b_0^2k^2\delta_1(k)b_0k\delta_2(k^2)b_0\\                        
 -2|\tau|^4\xi_1\xi_2^3b_0^2k^2\delta_2(k)b_0k\delta_1(k^2)b_0                        
 -2|\tau|^4\xi_1\xi_2^3b_0^2k^2\delta_1(k^2)b_0k\delta_2(k)b_0\\                        
 -2|\tau|^4\xi_1\xi_2^3b_0^2k^2\delta_1(k^2)b_0\delta_2(k)b_0k                        
 -2|\tau|^4\xi_1\xi_2^3b_0^2k^2\delta_2(k^2)b_0k\delta_1(k)b_0\\                        
 -2|\tau|^4\xi_1\xi_2^3b_0^2k^2\delta_2(k^2)b_0\delta_1(k)b_0k                        
 -2|\tau|^4\xi_1\xi_2^3b_0\delta_1(k^2)b_0^2k^3\delta_2(k)b_0\\                        
 -2|\tau|^4\xi_1\xi_2^3b_0\delta_1(k^2)b_0^2k^2\delta_2(k)b_0k                        
 -2|\tau|^4\xi_1\xi_2^3b_0\delta_2(k^2)b_0^2k^3\delta_1(k)b_0\\                        
 -2|\tau|^4\xi_1\xi_2^3b_0\delta_2(k^2)b_0^2k^2\delta_1(k)b_0k                        
 -20\tau_1^2|\tau|^2\xi_1^2\xi_2^2b_0^3k^5\delta_2^2(k)b_0\\                        
 -40\tau_1^2|\tau|^2\xi_1^2\xi_2^2b_0^3k^4\delta_2(k)^2b_0                        
 -20\tau_1^2|\tau|^2\xi_1^2\xi_2^2b_0^3k^4\delta_2^2(k)b_0k\\                        
 -28\tau_1^2|\tau|^2\xi_1^2\xi_2^2b_0^2k^3\delta_2(k)b_0k\delta_2(k)b_0                        
 -28\tau_1^2|\tau|^2\xi_1^2\xi_2^2b_0^2k^3\delta_2(k)b_0\delta_2(k)b_0k\\                        
 -28\tau_1^2|\tau|^2\xi_1^2\xi_2^2b_0^2k^2\delta_2(k)b_0k^2\delta_2(k)b_0                        
 -28\tau_1^2|\tau|^2\xi_1^2\xi_2^2b_0^2k^2\delta_2(k)b_0k\delta_2(k)b_0k\\                        
 -24\tau_1^2|\tau|^2\xi_1\xi_2^3b_0^3k^5\delta_1\delta_2(k)b_0                        
 -24\tau_1^2|\tau|^2\xi_1\xi_2^3b_0^3k^4\delta_1(k)\delta_2(k)b_0\\                        
 -24\tau_1^2|\tau|^2\xi_1\xi_2^3b_0^3k^4\delta_2(k)\delta_1(k)b_0                        
 -24\tau_1^2|\tau|^2\xi_1\xi_2^3b_0^3k^4\delta_1\delta_2(k)b_0k\\                        
 -20\tau_1^2|\tau|^2\xi_1\xi_2^3b_0^2k^3\delta_1(k)b_0k\delta_2(k)b_0                        
 -20\tau_1^2|\tau|^2\xi_1\xi_2^3b_0^2k^3\delta_1(k)b_0\delta_2(k)b_0k\\                        
 -20\tau_1^2|\tau|^2\xi_1\xi_2^3b_0^2k^3\delta_2(k)b_0k\delta_1(k)b_0                        
 -20\tau_1^2|\tau|^2\xi_1\xi_2^3b_0^2k^3\delta_2(k)b_0\delta_1(k)b_0k\\                        
 -20\tau_1^2|\tau|^2\xi_1\xi_2^3b_0^2k^2\delta_1(k)b_0k^2\delta_2(k)b_0                        
 -20\tau_1^2|\tau|^2\xi_1\xi_2^3b_0^2k^2\delta_1(k)b_0k\delta_2(k)b_0k\\                        
 -20\tau_1^2|\tau|^2\xi_1\xi_2^3b_0^2k^2\delta_2(k)b_0k^2\delta_1(k)b_0                        
 -20\tau_1^2|\tau|^2\xi_1\xi_2^3b_0^2k^2\delta_2(k)b_0k\delta_1(k)b_0k\\                        
 -4\tau_1^2|\tau|^2\xi_2^4b_0^3k^5\delta_1^2(k)b_0                       
 -8\tau_1^2|\tau|^2\xi_2^4b_0^3k^4\delta_1(k)^2b_0\\                        
 -4\tau_1^2|\tau|^2\xi_2^4b_0^3k^4\delta_1^2(k)b_0k                        
 -4\tau_1^2|\tau|^2\xi_2^4b_0^2k^3\delta_1(k)b_0k\delta_1(k)b_0\\                        
 -4\tau_1^2|\tau|^2\xi_2^4b_0^2k^3\delta_1(k)b_0\delta_1(k)b_0k                        
 -4\tau_1^2|\tau|^2\xi_2^4b_0^2k^2\delta_1(k)b_0k^2\delta_1(k)b_0\\                        
 -4\tau_1^2|\tau|^2\xi_2^4b_0^2k^2\delta_1(k)b_0k\delta_1(k)b_0k                        
 -8\tau_1|\tau|^4\xi_1\xi_2^3b_0^2k^3\delta_2(k)b_0\delta_2(k^2)b_0\\                        
 -8\tau_1|\tau|^4\xi_1\xi_2^3b_0^2k^2\delta_2(k)b_0k\delta_2(k^2)b_0                        
 -8\tau_1|\tau|^4\xi_1\xi_2^3b_0^2k^2\delta_2(k^2)b_0k\delta_2(k)b_0\\                        
 -8\tau_1|\tau|^4\xi_1\xi_2^3b_0^2k^2\delta_2(k^2)b_0\delta_2(k)b_0k                        
 -8\tau_1|\tau|^4\xi_1\xi_2^3b_0\delta_2(k^2)b_0^2k^3\delta_2(k)b_0\\                        
 -8\tau_1|\tau|^4\xi_1\xi_2^3b_0\delta_2(k^2)b_0^2k^2\delta_2(k)b_0k                        
 -2\tau_1|\tau|^4\xi_2^4b_0^2k^3\delta_1(k)b_0\delta_2(k^2)b_0\\                        
 -2\tau_1|\tau|^4\xi_2^4b_0^2k^3\delta_2(k)b_0\delta_1(k^2)b_0                        
 -2\tau_1|\tau|^4\xi_2^4b_0^2k^2\delta_1(k)b_0k\delta_2(k^2)b_0\\                        
 -2\tau_1|\tau|^4\xi_2^4b_0^2k^2\delta_2(k)b_0k\delta_1(k^2)b_0                        
 -2\tau_1|\tau|^4\xi_2^4b_0^2k^2\delta_1(k^2)b_0k\delta_2(k)b_0\\                        
 -2\tau_1|\tau|^4\xi_2^4b_0^2k^2\delta_1(k^2)b_0\delta_2(k)b_0k                        
 -2\tau_1|\tau|^4\xi_2^4b_0^2k^2\delta_2(k^2)b_0k\delta_1(k)b_0\\                        
 -2\tau_1|\tau|^4\xi_2^4b_0^2k^2\delta_2(k^2)b_0\delta_1(k)b_0k                        
 -2\tau_1|\tau|^4\xi_2^4b_0\delta_1(k^2)b_0^2k^3\delta_2(k)b_0\\                        
 -2\tau_1|\tau|^4\xi_2^4b_0\delta_1(k^2)b_0^2k^2\delta_2(k)b_0k                        
-2\tau_1|\tau|^4\xi_2^4b_0\delta_2(k^2)b_0^2k^3\delta_1(k)b_0\\                        
-2\tau_1|\tau|^4\xi_2^4b_0\delta_2(k^2)b_0^2k^2\delta_1(k)b_0k                        
-2i\tau_2|\tau|^4\xi_1\xi_2^3b_0^2k^3\delta_2(k)b_0\delta_2(k^2)b_0\\                        
-2i\tau_2|\tau|^4\xi_1\xi_2^3b_0^2k^2\delta_2(k)b_0k\delta_2(k^2)b_0                        
-2i\tau_2|\tau|^4\xi_1\xi_2^3b_0^2k^2\delta_2(k^2)b_0k\delta_2(k)b_0\\                        
-2i\tau_2|\tau|^4\xi_1\xi_2^3b_0^2k^2\delta_2(k^2)b_0\delta_2(k)b_0k                        
-2i\tau_2|\tau|^4\xi_1\xi_2^3b_0\delta_2(k^2)b_0^2k^3\delta_2(k)b_0\\                        
-2i\tau_2|\tau|^4\xi_1\xi_2^3b_0\delta_2(k^2)b_0^2k^2\delta_2(k)b_0k                        
+2i\tau_2|\tau|^4\xi_2^4b_0^2k^3\delta_2(k)b_0\delta_1(k^2)b_0\\                        
+2i\tau_2|\tau|^4\xi_2^4b_0^2k^2\delta_2(k)b_0k\delta_1(k^2)b_0                        
+2i\tau_2|\tau|^4\xi_2^4b_0^2k^2\delta_1(k^2)b_0k\delta_2(k)b_0\\                        
+2i\tau_2|\tau|^4\xi_2^4b_0^2k^2\delta_1(k^2)b_0\delta_2(k)b_0k                       
+2i\tau_2|\tau|^4\xi_2^4b_0\delta_1(k^2)b_0^2k^3\delta_2(k)b_0\\                        
+2i\tau_2|\tau|^4\xi_2^4b_0\delta_1(k^2)b_0^2k^2\delta_2(k)b_0k                        
-4|\tau|^4\xi_1^2\xi_2^2b_0^3k^5\delta_2^2(k)b_0\\                        
-8|\tau|^4\xi_1^2\xi_2^2b_0^3k^4\delta_2(k)^2b_0                        
-4|\tau|^4\xi_1^2\xi_2^2b_0^3k^4\delta_2^2(k)b_0k\\                        
-8|\tau|^4\xi_1^2\xi_2^2b_0^2k^3\delta_2(k)b_0k\delta_2(k)b_0                        
-8|\tau|^4\xi_1^2\xi_2^2b_0^2k^3\delta_2(k)b_0\delta_2(k)b_0k\\                       
 -8|\tau|^4\xi_1^2\xi_2^2b_0^2k^2\delta_2(k)b_0k^2\delta_2(k)b_0                        
 -8|\tau|^4\xi_1^2\xi_2^2b_0^2k^2\delta_2(k)b_0k\delta_2(k)b_0k\\                        
 -8|\tau|^4\xi_1\xi_2^3b_0^3k^5\delta_1\delta_2(k)b_0                        
 -8|\tau|^4\xi_1\xi_2^3b_0^3k^4\delta_1(k)\delta_2(k)b_0\\                        
 -8|\tau|^4\xi_1\xi_2^3b_0^3k^4\delta_2(k)\delta_1(k)b_0                        
 -8|\tau|^4\xi_1\xi_2^3b_0^3k^4\delta_1\delta_2(k)b_0k\\                        
 -4|\tau|^4\xi_1\xi_2^3b_0^2k^3\delta_1(k)b_0k\delta_2(k)b_0                        
 -4|\tau|^4\xi_1\xi_2^3b_0^2k^3\delta_1(k)b_0\delta_2(k)b_0k\\                        
 -4|\tau|^4\xi_1\xi_2^3b_0^2k^3\delta_2(k)b_0k\delta_1(k)b_0                        
 -4|\tau|^4\xi_1\xi_2^3b_0^2k^3\delta_2(k)b_0\delta_1(k)b_0k\\                        
 -4|\tau|^4\xi_1\xi_2^3b_0^2k^2\delta_1(k)b_0k^2\delta_2(k)b_0                        
 -4|\tau|^4\xi_1\xi_2^3b_0^2k^2\delta_1(k)b_0k\delta_2(k)b_0k\\                        
 -4|\tau|^4\xi_1\xi_2^3b_0^2k^2\delta_2(k)b_0k^2\delta_1(k)b_0                        
 -4|\tau|^4\xi_1\xi_2^3b_0^2k^2\delta_2(k)b_0k\delta_1(k)b_0k\\                        
 -2|\tau|^4\xi_2^4b_0^2k^3\delta_1(k)b_0k\delta_1(k)b_0                        
 -2|\tau|^4\xi_2^4b_0^2k^3\delta_1(k)b_0\delta_1(k)b_0k\\                        
 -2|\tau|^4\xi_2^4b_0^2k^2\delta_1(k)b_0k^2\delta_1(k)b_0                        
 -2|\tau|^4\xi_2^4b_0^2k^2\delta_1(k)b_0k\delta_1(k)b_0k\\                        
 -16\tau_1|\tau|^4\xi_1\xi_2^3b_0^3k^5\delta_2^2(k)b_0                        
 -32\tau_1|\tau|^4\xi_1\xi_2^3b_0^3k^4\delta_2(k)^2b_0\\                        
 -16\tau_1|\tau|^4\xi_1\xi_2^3b_0^3k^4\delta_2^2(k)b_0k                        
 -24\tau_1|\tau|^4\xi_1\xi_2^3b_0^2k^3\delta_2(k)b_0k\delta_2(k)b_0\\                        
 -24\tau_1|\tau|^4\xi_1\xi_2^3b_0^2k^3\delta_2(k)b_0\delta_2(k)b_0k                        
 -24\tau_1|\tau|^4\xi_1\xi_2^3b_0^2k^2\delta_2(k)b_0k^2\delta_2(k)b_0\\                        
 -24\tau_1|\tau|^4\xi_1\xi_2^3b_0^2k^2\delta_2(k)b_0k\delta_2(k)b_0k                        
 -8\tau_1|\tau|^4\xi_2^4b_0^3k^5\delta_1\delta_2(k)b_0\\                        
 -8\tau_1|\tau|^4\xi_2^4b_0^3k^4\delta_1(k)\delta_2(k)b_0                        
 -8\tau_1|\tau|^4\xi_2^4b_0^3k^4\delta_2(k)\delta_1(k)b_0\\                        
 -8\tau_1|\tau|^4\xi_2^4b_0^3k^4\delta_1\delta_2(k)b_0k                       
 -6\tau_1|\tau|^4\xi_2^4b_0^2k^3\delta_1(k)b_0k\delta_2(k)b_0\\                        
 -6\tau_1|\tau|^4\xi_2^4b_0^2k^3\delta_1(k)b_0\delta_2(k)b_0k                       
 -6\tau_1|\tau|^4\xi_2^4b_0^2k^3\delta_2(k)b_0k\delta_1(k)b_0\\                        
 -6\tau_1|\tau|^4\xi_2^4b_0^2k^3\delta_2(k)b_0\delta_1(k)b_0k                       
 -6\tau_1|\tau|^4\xi_2^4b_0^2k^2\delta_1(k)b_0k^2\delta_2(k)b_0\\                        
 -6\tau_1|\tau|^4\xi_2^4b_0^2k^2\delta_1(k)b_0k\delta_2(k)b_0k                       
 -6\tau_1|\tau|^4\xi_2^4b_0^2k^2\delta_2(k)b_0k^2\delta_1(k)b_0\\                        
 -6\tau_1|\tau|^4\xi_2^4b_0^2k^2\delta_2(k)b_0k\delta_1(k)b_0k                       
 -2|\tau|^6\xi_2^4b_0^2k^3\delta_2(k)b_0\delta_2(k^2)b_0\\                        
 -2|\tau|^6\xi_2^4b_0^2k^2\delta_2(k)b_0k\delta_2(k^2)b_0                      
 -2|\tau|^6\xi_2^4b_0^2k^2\delta_2(k^2)b_0k\delta_2(k)b_0\\                        
 -2|\tau|^6\xi_2^4b_0^2k^2\delta_2(k^2)b_0\delta_2(k)b_0k                        
 -2|\tau|^6\xi_2^4b_0\delta_2(k^2)b_0^2k^3\delta_2(k)b_0\\                        
 -2|\tau|^6\xi_2^4b_0\delta_2(k^2)b_0^2k^2\delta_2(k)b_0k                       
 +8\xi_1^6b_0^3k^5\delta_1(k)b_0k\delta_1(k)b_0\\                        
 +8\xi_1^6b_0^3k^5\delta_1(k)b_0\delta_1(k)b_0k                        
 +8\xi_1^6b_0^3k^4\delta_1(k)b_0k^2\delta_1(k)b_0\\                       
 +8\xi_1^6b_0^3k^4\delta_1(k)b_0k\delta_1(k)b_0k                        
 +4\xi_1^6b_0^2k^3\delta_1(k)b_0^2k^3\delta_1(k)b_0\\                        
 +4\xi_1^6b_0^2k^3\delta_1(k)b_0^2k^2\delta_1(k)b_0k                        
 +4\xi_1^6b_0^2k^2\delta_1(k)b_0^2k^4\delta_1(k)b_0\\                        
 +4\xi_1^6b_0^2k^2\delta_1(k)b_0^2k^3\delta_1(k)b_0k                        
 +8\tau_1\xi_1^6b_0^3k^5\delta_1(k)b_0k\delta_2(k)b_0\\                        
 +8\tau_1\xi_1^6b_0^3k^5\delta_1(k)b_0\delta_2(k)b_0k                      
  +8\tau_1\xi_1^6b_0^3k^5\delta_2(k)b_0k\delta_1(k)b_0\\                        
  +8\tau_1\xi_1^6b_0^3k^5\delta_2(k)b_0\delta_1(k)b_0k                        
  +8\tau_1\xi_1^6b_0^3k^4\delta_1(k)b_0k^2\delta_2(k)b_0\\                        
  +8\tau_1\xi_1^6b_0^3k^4\delta_1(k)b_0k\delta_2(k)b_0k                        
  +8\tau_1\xi_1^6b_0^3k^4\delta_2(k)b_0k^2\delta_1(k)b_0\\                        
  +8\tau_1\xi_1^6b_0^3k^4\delta_2(k)b_0k\delta_1(k)b_0k                        
  +4\tau_1\xi_1^6b_0^2k^3\delta_1(k)b_0^2k^3\delta_2(k)b_0\\                        
  +4\tau_1\xi_1^6b_0^2k^3\delta_1(k)b_0^2k^2\delta_2(k)b_0k                        
  +4\tau_1\xi_1^6b_0^2k^3\delta_2(k)b_0^2k^3\delta_1(k)b_0\\                        
  +4\tau_1\xi_1^6b_0^2k^3\delta_2(k)b_0^2k^2\delta_1(k)b_0k                        
  +4\tau_1\xi_1^6b_0^2k^2\delta_1(k)b_0^2k^4\delta_2(k)b_0\\                        
  +4\tau_1\xi_1^6b_0^2k^2\delta_1(k)b_0^2k^3\delta_2(k)b_0k                        
  +4\tau_1\xi_1^6b_0^2k^2\delta_2(k)b_0^2k^4\delta_1(k)b_0\\                        
  +4\tau_1\xi_1^6b_0^2k^2\delta_2(k)b_0^2k^3\delta_1(k)b_0k                        
  +48\tau_1\xi_1^5\xi_2b_0^3k^5\delta_1(k)b_0k\delta_1(k)b_0\\                        
  +48\tau_1\xi_1^5\xi_2b_0^3k^5\delta_1(k)b_0\delta_1(k)b_0k                       
   +48\tau_1\xi_1^5\xi_2b_0^3k^4\delta_1(k)b_0k^2\delta_1(k)b_0\\                        
   +48\tau_1\xi_1^5\xi_2b_0^3k^4\delta_1(k)b_0k\delta_1(k)b_0k                        
   +24\tau_1\xi_1^5\xi_2b_0^2k^3\delta_1(k)b_0^2k^3\delta_1(k)b_0\\                        
   +24\tau_1\xi_1^5\xi_2b_0^2k^3\delta_1(k)b_0^2k^2\delta_1(k)b_0k                        
   +24\tau_1\xi_1^5\xi_2b_0^2k^2\delta_1(k)b_0^2k^4\delta_1(k)b_0\\                        
   +24\tau_1\xi_1^5\xi_2b_0^2k^2\delta_1(k)b_0^2k^3\delta_1(k)b_0k                       
   -4|\tau|^6\xi_2^4b_0^3k^5\delta_2^2(k)b_0\\                        
   -8|\tau|^6\xi_2^4b_0^3k^4\delta_2(k)^2b_0                        
   -4|\tau|^6\xi_2^4b_0^3k^4\delta_2^2(k)b_0k\\                        
-6|\tau|^6\xi_2^4b_0^2k^3\delta_2(k)b_0k\delta_2(k)b_0                        
-6|\tau|^6\xi_2^4b_0^2k^3\delta_2(k)b_0\delta_2(k)b_0k\\                        
-6|\tau|^6\xi_2^4b_0^2k^2\delta_2(k)b_0k^2\delta_2(k)b_0                        
-6|\tau|^6\xi_2^4b_0^2k^2\delta_2(k)b_0k\delta_2(k)b_0k\\                        
+8\tau_1^2\xi_1^6b_0^3k^5\delta_2(k)b_0k\delta_2(k)b_0                        
+8\tau_1^2\xi_1^6b_0^3k^5\delta_2(k)b_0\delta_2(k)b_0k\\                        
+8\tau_1^2\xi_1^6b_0^3k^4\delta_2(k)b_0k^2\delta_2(k)b_0                        
+8\tau_1^2\xi_1^6b_0^3k^4\delta_2(k)b_0k\delta_2(k)b_0k\\                        
+4\tau_1^2\xi_1^6b_0^2k^3\delta_2(k)b_0^2k^3\delta_2(k)b_0                        
+4\tau_1^2\xi_1^6b_0^2k^3\delta_2(k)b_0^2k^2\delta_2(k)b_0k\\                        
+4\tau_1^2\xi_1^6b_0^2k^2\delta_2(k)b_0^2k^4\delta_2(k)b_0                        
+4\tau_1^2\xi_1^6b_0^2k^2\delta_2(k)b_0^2k^3\delta_2(k)b_0k\\                        
+40\tau_1^2\xi_1^5\xi_2b_0^3k^5\delta_1(k)b_0k\delta_2(k)b_0                        
+40\tau_1^2\xi_1^5\xi_2b_0^3k^5\delta_1(k)b_0\delta_2(k)b_0k\\                        
+40\tau_1^2\xi_1^5\xi_2b_0^3k^5\delta_2(k)b_0k\delta_1(k)b_0                        
+40\tau_1^2\xi_1^5\xi_2b_0^3k^5\delta_2(k)b_0\delta_1(k)b_0k\\                        
+40\tau_1^2\xi_1^5\xi_2b_0^3k^4\delta_1(k)b_0k^2\delta_2(k)b_0                        
+40\tau_1^2\xi_1^5\xi_2b_0^3k^4\delta_1(k)b_0k\delta_2(k)b_0k\\                        
+40\tau_1^2\xi_1^5\xi_2b_0^3k^4\delta_2(k)b_0k^2\delta_1(k)b_0                        
+40\tau_1^2\xi_1^5\xi_2b_0^3k^4\delta_2(k)b_0k\delta_1(k)b_0k\\                        
+20\tau_1^2\xi_1^5\xi_2b_0^2k^3\delta_1(k)b_0^2k^3\delta_2(k)b_0                        
+20\tau_1^2\xi_1^5\xi_2b_0^2k^3\delta_1(k)b_0^2k^2\delta_2(k)b_0k\\                        
+20\tau_1^2\xi_1^5\xi_2b_0^2k^3\delta_2(k)b_0^2k^3\delta_1(k)b_0                       
+20\tau_1^2\xi_1^5\xi_2b_0^2k^3\delta_2(k)b_0^2k^2\delta_1(k)b_0k\\                        
+20\tau_1^2\xi_1^5\xi_2b_0^2k^2\delta_1(k)b_0^2k^4\delta_2(k)b_0                       
+20\tau_1^2\xi_1^5\xi_2b_0^2k^2\delta_1(k)b_0^2k^3\delta_2(k)b_0k\\                        
+20\tau_1^2\xi_1^5\xi_2b_0^2k^2\delta_2(k)b_0^2k^4\delta_1(k)b_0                        
+20\tau_1^2\xi_1^5\xi_2b_0^2k^2\delta_2(k)b_0^2k^3\delta_1(k)b_0k\\                        
+104\tau_1^2\xi_1^4\xi_2^2b_0^3k^5\delta_1(k)b_0k\delta_1(k)b_0                       
+104\tau_1^2\xi_1^4\xi_2^2b_0^3k^5\delta_1(k)b_0\delta_1(k)b_0k\\                        
+104\tau_1^2\xi_1^4\xi_2^2b_0^3k^4\delta_1(k)b_0k^2\delta_1(k)b_0                       
+104\tau_1^2\xi_1^4\xi_2^2b_0^3k^4\delta_1(k)b_0k\delta_1(k)b_0k\\                        
+52\tau_1^2\xi_1^4\xi_2^2b_0^2k^3\delta_1(k)b_0^2k^3\delta_1(k)b_0                        
+52\tau_1^2\xi_1^4\xi_2^2b_0^2k^3\delta_1(k)b_0^2k^2\delta_1(k)b_0k\\                        
+52\tau_1^2\xi_1^4\xi_2^2b_0^2k^2\delta_1(k)b_0^2k^4\delta_1(k)b_0                       
+52\tau_1^2\xi_1^4\xi_2^2b_0^2k^2\delta_1(k)b_0^2k^3\delta_1(k)b_0k\\                        
+8|\tau|^2\xi_1^5\xi_2b_0^3k^5\delta_1(k)b_0k\delta_2(k)b_0                       
+8|\tau|^2\xi_1^5\xi_2b_0^3k^5\delta_1(k)b_0\delta_2(k)b_0k\\                        
+8|\tau|^2\xi_1^5\xi_2b_0^3k^5\delta_2(k)b_0k\delta_1(k)b_0                        
+8|\tau|^2\xi_1^5\xi_2b_0^3k^5\delta_2(k)b_0\delta_1(k)b_0k\\                        
+8|\tau|^2\xi_1^5\xi_2b_0^3k^4\delta_1(k)b_0k^2\delta_2(k)b_0                        
+8|\tau|^2\xi_1^5\xi_2b_0^3k^4\delta_1(k)b_0k\delta_2(k)b_0k\\                        
+8|\tau|^2\xi_1^5\xi_2b_0^3k^4\delta_2(k)b_0k^2\delta_1(k)b_0                        
+8|\tau|^2\xi_1^5\xi_2b_0^3k^4\delta_2(k)b_0k\delta_1(k)b_0k\\                        
+4|\tau|^2\xi_1^5\xi_2b_0^2k^3\delta_1(k)b_0^2k^3\delta_2(k)b_0                        
+4|\tau|^2\xi_1^5\xi_2b_0^2k^3\delta_1(k)b_0^2k^2\delta_2(k)b_0k\\                        
+4|\tau|^2\xi_1^5\xi_2b_0^2k^3\delta_2(k)b_0^2k^3\delta_1(k)b_0                       
+4|\tau|^2\xi_1^5\xi_2b_0^2k^3\delta_2(k)b_0^2k^2\delta_1(k)b_0k\\                        
+4|\tau|^2\xi_1^5\xi_2b_0^2k^2\delta_1(k)b_0^2k^4\delta_2(k)b_0                       
+4|\tau|^2\xi_1^5\xi_2b_0^2k^2\delta_1(k)b_0^2k^3\delta_2(k)b_0k\\                        
+4|\tau|^2\xi_1^5\xi_2b_0^2k^2\delta_2(k)b_0^2k^4\delta_1(k)b_0                        
+4|\tau|^2\xi_1^5\xi_2b_0^2k^2\delta_2(k)b_0^2k^3\delta_1(k)b_0k\\                        
+16|\tau|^2\xi_1^4\xi_2^2b_0^3k^5\delta_1(k)b_0k\delta_1(k)b_0                        
+16|\tau|^2\xi_1^4\xi_2^2b_0^3k^5\delta_1(k)b_0\delta_1(k)b_0k\\                        
+16|\tau|^2\xi_1^4\xi_2^2b_0^3k^4\delta_1(k)b_0k^2\delta_1(k)b_0                        
+16|\tau|^2\xi_1^4\xi_2^2b_0^3k^4\delta_1(k)b_0k\delta_1(k)b_0k\\                        
+8|\tau|^2\xi_1^4\xi_2^2b_0^2k^3\delta_1(k)b_0^2k^3\delta_1(k)b_0                        
+8|\tau|^2\xi_1^4\xi_2^2b_0^2k^3\delta_1(k)b_0^2k^2\delta_1(k)b_0k\\                        
+8|\tau|^2\xi_1^4\xi_2^2b_0^2k^2\delta_1(k)b_0^2k^4\delta_1(k)b_0                        
+8|\tau|^2\xi_1^4\xi_2^2b_0^2k^2\delta_1(k)b_0^2k^3\delta_1(k)b_0k\\                        
+32\tau_1^3\xi_1^5\xi_2b_0^3k^5\delta_2(k)b_0k\delta_2(k)b_0                        
+32\tau_1^3\xi_1^5\xi_2b_0^3k^5\delta_2(k)b_0\delta_2(k)b_0k\\                        
+32\tau_1^3\xi_1^5\xi_2b_0^3k^4\delta_2(k)b_0k^2\delta_2(k)b_0                        
+32\tau_1^3\xi_1^5\xi_2b_0^3k^4\delta_2(k)b_0k\delta_2(k)b_0k\\                        
+16\tau_1^3\xi_1^5\xi_2b_0^2k^3\delta_2(k)b_0^2k^3\delta_2(k)b_0                       
 +16\tau_1^3\xi_1^5\xi_2b_0^2k^3\delta_2(k)b_0^2k^2\delta_2(k)b_0k\\                        
 +16\tau_1^3\xi_1^5\xi_2b_0^2k^2\delta_2(k)b_0^2k^4\delta_2(k)b_0                       
  +16\tau_1^3\xi_1^5\xi_2b_0^2k^2\delta_2(k)b_0^2k^3\delta_2(k)b_0k\\                        
  +64\tau_1^3\xi_1^4\xi_2^2b_0^3k^5\delta_1(k)b_0k\delta_2(k)b_0                        
  +64\tau_1^3\xi_1^4\xi_2^2b_0^3k^5\delta_1(k)b_0\delta_2(k)b_0k\\                        
  +64\tau_1^3\xi_1^4\xi_2^2b_0^3k^5\delta_2(k)b_0k\delta_1(k)b_0                        
  +64\tau_1^3\xi_1^4\xi_2^2b_0^3k^5\delta_2(k)b_0\delta_1(k)b_0k\\                        
  +64\tau_1^3\xi_1^4\xi_2^2b_0^3k^4\delta_1(k)b_0k^2\delta_2(k)b_0                        
  +64\tau_1^3\xi_1^4\xi_2^2b_0^3k^4\delta_1(k)b_0k\delta_2(k)b_0k\\                        
  +64\tau_1^3\xi_1^4\xi_2^2b_0^3k^4\delta_2(k)b_0k^2\delta_1(k)b_0                        
  +64\tau_1^3\xi_1^4\xi_2^2b_0^3k^4\delta_2(k)b_0k\delta_1(k)b_0k\\                        
  +32\tau_1^3\xi_1^4\xi_2^2b_0^2k^3\delta_1(k)b_0^2k^3\delta_2(k)b_0                        
  +32\tau_1^3\xi_1^4\xi_2^2b_0^2k^3\delta_1(k)b_0^2k^2\delta_2(k)b_0k\\                        
  +32\tau_1^3\xi_1^4\xi_2^2b_0^2k^3\delta_2(k)b_0^2k^3\delta_1(k)b_0                        
  +32\tau_1^3\xi_1^4\xi_2^2b_0^2k^3\delta_2(k)b_0^2k^2\delta_1(k)b_0k\\                        
  +32\tau_1^3\xi_1^4\xi_2^2b_0^2k^2\delta_1(k)b_0^2k^4\delta_2(k)b_0                        
  +32\tau_1^3\xi_1^4\xi_2^2b_0^2k^2\delta_1(k)b_0^2k^3\delta_2(k)b_0k\\                        
  +32\tau_1^3\xi_1^4\xi_2^2b_0^2k^2\delta_2(k)b_0^2k^4\delta_1(k)b_0                        
  +32\tau_1^3\xi_1^4\xi_2^2b_0^2k^2\delta_2(k)b_0^2k^3\delta_1(k)b_0k\\                        
  +96\tau_1^3\xi_1^3\xi_2^3b_0^3k^5\delta_1(k)b_0k\delta_1(k)b_0                        
  +96\tau_1^3\xi_1^3\xi_2^3b_0^3k^5\delta_1(k)b_0\delta_1(k)b_0k\\                       
   +96\tau_1^3\xi_1^3\xi_2^3b_0^3k^4\delta_1(k)b_0k^2\delta_1(k)b_0                        
   +96\tau_1^3\xi_1^3\xi_2^3b_0^3k^4\delta_1(k)b_0k\delta_1(k)b_0k\\                        
   +48\tau_1^3\xi_1^3\xi_2^3b_0^2k^3\delta_1(k)b_0^2k^3\delta_1(k)b_0                        
   +48\tau_1^3\xi_1^3\xi_2^3b_0^2k^3\delta_1(k)b_0^2k^2\delta_1(k)b_0k\\                        
   +48\tau_1^3\xi_1^3\xi_2^3b_0^2k^2\delta_1(k)b_0^2k^4\delta_1(k)b_0                        
   +48\tau_1^3\xi_1^3\xi_2^3b_0^2k^2\delta_1(k)b_0^2k^3\delta_1(k)b_0k\\                       
+16\tau_1|\tau|^2\xi_1^5\xi_2b_0^3k^5\delta_2(k)b_0k\delta_2(k)b_0                        
+16\tau_1|\tau|^2\xi_1^5\xi_2b_0^3k^
5\delta_2(k)b_0\delta_2(k)b_0k\\                        
+16\tau_1|\tau|^2\xi_1^5\xi_2b_0^3k^
4\delta_2(k)b_0k^2\delta_2(k)b_0                        
+16\tau_1|\tau|^2\xi_1^5\xi_2b_0^3k^
4\delta_2(k)b_0k\delta_2(k)b_0k\\                        
+8\tau_1|\tau|^2\xi_1^5\xi_2b_0^2k^
3\delta_2(k)b_0^2k^3\delta_2(k)b_0                        
+8\tau_1|\tau|^2\xi_1^5\xi_2b_0^2k^
3\delta_2(k)b_0^2k^2\delta_2(k)b_0k\\                        
+8\tau_1|\tau|^2\xi_1^5\xi_2b_0^2k^
2\delta_2(k)b_0^2k^4\delta_2(k)b_0                        
+8\tau_1|\tau|^2\xi_1^5\xi_2b_0^2k^
2\delta_2(k)b_0^2k^3\delta_2(k)b_0k\\                        
+56\tau_1|\tau|^2\xi_1^4\xi_2^2b_0^3k^
5\delta_1(k)b_0k\delta_2(k)b_0                        
+56\tau_1|\tau|^2\xi_1^4\xi_2^2b_0^3k^
5\delta_1(k)b_0\delta_2(k)b_0k\\                        
+56\tau_1|\tau|^2\xi_1^4\xi_2^2b_0^3k^
5\delta_2(k)b_0k\delta_1(k)b_0                        
+56\tau_1|\tau|^2\xi_1^4\xi_2^2b_0^3k^
5\delta_2(k)b_0\delta_1(k)b_0k\\                        
+56\tau_1|\tau|^2\xi_1^4\xi_2^2b_0^3k^
4\delta_1(k)b_0k^2\delta_2(k)b_0                        
+56\tau_1|\tau|^2\xi_1^4\xi_2^2b_0^3k^
4\delta_1(k)b_0k\delta_2(k)b_0k\\                        
+56\tau_1|\tau|^2\xi_1^4\xi_2^2b_0^3k^
4\delta_2(k)b_0k^2\delta_1(k)b_0                        
+56\tau_1|\tau|^2\xi_1^4\xi_2^2b_0^3k^
4\delta_2(k)b_0k\delta_1(k)b_0k\\                        
+28\tau_1|\tau|^2\xi_1^4\xi_2^2b_0^2k^
3\delta_1(k)b_0^2k^3\delta_2(k)b_0                        
+28\tau_1|\tau|^2\xi_1^4\xi_2^2b_0^2k^
3\delta_1(k)b_0^2k^2\delta_2(k)b_0k\\                        
+28\tau_1|\tau|^2\xi_1^4\xi_2^2b_0^2k^
3\delta_2(k)b_0^2k^3\delta_1(k)b_0                        
+28\tau_1|\tau|^2\xi_1^4\xi_2^2b_0^2k^
3\delta_2(k)b_0^2k^2\delta_1(k)b_0k\\                        
+28\tau_1|\tau|^2\xi_1^4\xi_2^2b_0^2k^
2\delta_1(k)b_0^2k^4\delta_2(k)b_0                        
+28\tau_1|\tau|^2\xi_1^4\xi_2^2b_0^2k^
2\delta_1(k)b_0^2k^3\delta_2(k)b_0k\\                        
+28\tau_1|\tau|^2\xi_1^4\xi_2^2b_0^2k^
2\delta_2(k)b_0^2k^4\delta_1(k)b_0                        
+28\tau_1|\tau|^2\xi_1^4\xi_2^2b_0^2k^
2\delta_2(k)b_0^2k^3\delta_1(k)b_0k\\                        
+64\tau_1|\tau|^2\xi_1^3\xi_2^3b_0^3k^
5\delta_1(k)b_0k\delta_1(k)b_0                        
+64\tau_1|\tau|^2\xi_1^3\xi_2^3b_0^3k^
5\delta_1(k)b_0\delta_1(k)b_0k\\                        
+64\tau_1|\tau|^2\xi_1^3\xi_2^3b_0^3k^
4\delta_1(k)b_0k^2\delta_1(k)b_0                        
+64\tau_1|\tau|^2\xi_1^3\xi_2^3b_0^3k^
4\delta_1(k)b_0k\delta_1(k)b_0k\\                        
+32\tau_1|\tau|^2\xi_1^3\xi_2^3b_0^2k^
3\delta_1(k)b_0^2k^3\delta_1(k)b_0                        
+32\tau_1|\tau|^2\xi_1^3\xi_2^3b_0^2k^
3\delta_1(k)b_0^2k^2\delta_1(k)b_0k\\                        
+32\tau_1|\tau|^2\xi_1^3\xi_2^3b_0^2k^
2\delta_1(k)b_0^2k^4\delta_1(k)b_0                        
+32\tau_1|\tau|^2\xi_1^3\xi_2^3b_0^2k^
2\delta_1(k)b_0^2k^3\delta_1(k)b_0k\\                        
+32\tau_1^4\xi_1^4\xi_2^2b_0^3k^
5\delta_2(k)b_0k\delta_2(k)b_0                        
+32\tau_1^4\xi_1^4\xi_2^2b_0^3k^
5\delta_2(k)b_0\delta_2(k)b_0k\\                        
+32\tau_1^4\xi_1^4\xi_2^2b_0^3k^
4\delta_2(k)b_0k^2\delta_2(k)b_0                        
+32\tau_1^4\xi_1^4\xi_2^2b_0^3k^
4\delta_2(k)b_0k\delta_2(k)b_0k\\                        
+16\tau_1^4\xi_1^4\xi_2^2b_0^2k^
3\delta_2(k)b_0^2k^3\delta_2(k)b_0                        
+16\tau_1^4\xi_1^4\xi_2^2b_0^2k^
3\delta_2(k)b_0^2k^2\delta_2(k)b_0k\\                        
+16\tau_1^4\xi_1^4\xi_2^2b_0^2k^
2\delta_2(k)b_0^2k^4\delta_2(k)b_0                        
+16\tau_1^4\xi_1^4\xi_2^2b_0^2k^
2\delta_2(k)b_0^2k^3\delta_2(k)b_0k\\                        
+32\tau_1^4\xi_1^3\xi_2^3b_0^3k^
5\delta_1(k)b_0k\delta_2(k)b_0                        
+32\tau_1^4\xi_1^3\xi_2^3b_0^3k^
5\delta_1(k)b_0\delta_2(k)b_0k\\                        
+32\tau_1^4\xi_1^3\xi_2^3b_0^3k^
5\delta_2(k)b_0k\delta_1(k)b_0                        
+32\tau_1^4\xi_1^3\xi_2^3b_0^3k^
5\delta_2(k)b_0\delta_1(k)b_0k\\                        
+32\tau_1^4\xi_1^3\xi_2^3b_0^3k^
4\delta_1(k)b_0k^2\delta_2(k)b_0                        
+32\tau_1^4\xi_1^3\xi_2^3b_0^3k^
4\delta_1(k)b_0k\delta_2(k)b_0k\\                        
+32\tau_1^4\xi_1^3\xi_2^3b_0^3k^
4\delta_2(k)b_0k^2\delta_1(k)b_0                        
+32\tau_1^4\xi_1^3\xi_2^3b_0^3k^
4\delta_2(k)b_0k\delta_1(k)b_0k\\                        
+16\tau_1^4\xi_1^3\xi_2^3b_0^2k^
3\delta_1(k)b_0^2k^3\delta_2(k)b_0                       
+16\tau_1^4\xi_1^3\xi_2^3b_0^2k^
3\delta_1(k)b_0^2k^2\delta_2(k)b_0k\\                        
+16\tau_1^4\xi_1^3\xi_2^3b_0^2k^
3\delta_2(k)b_0^2k^3\delta_1(k)b_0                        
+16\tau_1^4\xi_1^3\xi_2^3b_0^2k^
3\delta_2(k)b_0^2k^2\delta_1(k)b_0k\\                        
+16\tau_1^4\xi_1^3\xi_2^3b_0^2k^
2\delta_1(k)b_0^2k^4\delta_2(k)b_0                       
+16\tau_1^4\xi_1^3\xi_2^3b_0^2k^
2\delta_1(k)b_0^2k^3\delta_2(k)b_0k\\                        
+16\tau_1^4\xi_1^3\xi_2^3b_0^2k^
2\delta_2(k)b_0^2k^4\delta_1(k)b_0                        
+16\tau_1^4\xi_1^3\xi_2^3b_0^2k^
2\delta_2(k)b_0^2k^3\delta_1(k)b_0k\\                        
+32\tau_1^4\xi_1^2\xi_2^4b_0^3k^
5\delta_1(k)b_0k\delta_1(k)b_0                        
+32\tau_1^4\xi_1^2\xi_2^4b_0^3k^
5\delta_1(k)b_0\delta_1(k)b_0k\\                        
+32\tau_1^4\xi_1^2\xi_2^4b_0^3k^
4\delta_1(k)b_0k^2\delta_1(k)b_0                        
+32\tau_1^4\xi_1^2\xi_2^4b_0^3k^
4\delta_1(k)b_0k\delta_1(k)b_0k\\                        
+16\tau_1^4\xi_1^2\xi_2^4b_0^2k^
3\delta_1(k)b_0^2k^3\delta_1(k)b_0                        
+16\tau_1^4\xi_1^2\xi_2^4b_0^2k^
3\delta_1(k)b_0^2k^2\delta_1(k)b_0k\\                        
+16\tau_1^4\xi_1^2\xi_2^4b_0^2k^
2\delta_1(k)b_0^2k^4\delta_1(k)b_0                        
+16\tau_1^4\xi_1^2\xi_2^4b_0^2k^
2\delta_1(k)b_0^2k^3\delta_1(k)b_0k\\                        
+80\tau_1^2|\tau|^2\xi_1^4\xi_2^2b_0^3k^
5\delta_2(k)b_0k\delta_2(k)b_0                        
+80\tau_1^2|\tau|^2\xi_1^4\xi_2^2b_0^3k^
5\delta_2(k)b_0\delta_2(k)b_0k\\                        
+80\tau_1^2|\tau|^2\xi_1^4\xi_2^2b_0^3k^
4\delta_2(k)b_0k^2\delta_2(k)b_0                        
+80\tau_1^2|\tau|^2\xi_1^4\xi_2^2b_0^3k^
4\delta_2(k)b_0k\delta_2(k)b_0k\\                        
+40\tau_1^2|\tau|^2\xi_1^4\xi_2^2b_0^2k^
3\delta_2(k)b_0^2k^3\delta_2(k)b_0                       
 +40\tau_1^2|\tau|^2\xi_1^4\xi_2^2b_0^2k^
3\delta_2(k)b_0^2k^2\delta_2(k)b_0k\\                        
+40\tau_1^2|\tau|^2\xi_1^4\xi_2^2b_0^2k^
2\delta_2(k)b_0^2k^4\delta_2(k)b_0                       
 +40\tau_1^2|\tau|^2\xi_1^4\xi_2^2b_0^2k^
2\delta_2(k)b_0^2k^3\delta_2(k)b_0k\\                        
+112\tau_1^2|\tau|^2\xi_1^3\xi_2^3b_0^3k^
5\delta_1(k)b_0k\delta_2(k)b_0                       
+112\tau_1^2|\tau|^2\xi_1^3\xi_2^3b_0^3k^
5\delta_1(k)b_0\delta_2(k)b_0k\\                        
+112\tau_1^2|\tau|^2\xi_1^3\xi_2^3b_0^3k^
5\delta_2(k)b_0k\delta_1(k)b_0                        
+112\tau_1^2|\tau|^2\xi_1^3\xi_2^3b_0^3k^
5\delta_2(k)b_0\delta_1(k)b_0k\\                        
+112\tau_1^2|\tau|^2\xi_1^3\xi_2^3b_0^3k^
4\delta_1(k)b_0k^2\delta_2(k)b_0                        
+112\tau_1^2|\tau|^2\xi_1^3\xi_2^3b_0^3k^
4\delta_1(k)b_0k\delta_2(k)b_0k\\                        
+112\tau_1^2|\tau|^2\xi_1^3\xi_2^3b_0^3k^
4\delta_2(k)b_0k^2\delta_1(k)b_0                        
+112\tau_1^2|\tau|^2\xi_1^3\xi_2^3b_0^3k^
4\delta_2(k)b_0k\delta_1(k)b_0k\\                        
+56\tau_1^2|\tau|^2\xi_1^3\xi_2^3b_0^2k^
3\delta_1(k)b_0^2k^3\delta_2(k)b_0                        
+56\tau_1^2|\tau|^2\xi_1^3\xi_2^3b_0^2k^
3\delta_1(k)b_0^2k^2\delta_2(k)b_0k\\                        
+56\tau_1^2|\tau|^2\xi_1^3\xi_2^3b_0^2k^
3\delta_2(k)b_0^2k^3\delta_1(k)b_0                        
+56\tau_1^2|\tau|^2\xi_1^3\xi_2^3b_0^2k^
3\delta_2(k)b_0^2k^2\delta_1(k)b_0k\\                       
 +56\tau_1^2|\tau|^2\xi_1^3\xi_2^3b_0^2k^
2\delta_1(k)b_0^2k^4\delta_2(k)b_0                        
+56\tau_1^2|\tau|^2\xi_1^3\xi_2^3b_0^2k^
2\delta_1(k)b_0^2k^3\delta_2(k)b_0k\\                        
+56\tau_1^2|\tau|^2\xi_1^3\xi_2^3b_0^2k^
2\delta_2(k)b_0^2k^4\delta_1(k)b_0                        
+56\tau_1^2|\tau|^2\xi_1^3\xi_2^3b_0^2k^
2\delta_2(k)b_0^2k^3\delta_1(k)b_0k\\                        
+80\tau_1^2|\tau|^2\xi_1^2\xi_2^4b_0^3k^
5\delta_1(k)b_0k\delta_1(k)b_0                        
+80\tau_1^2|\tau|^2\xi_1^2\xi_2^4b_0^3k^
5\delta_1(k)b_0\delta_1(k)b_0k\\                        
+80\tau_1^2|\tau|^2\xi_1^2\xi_2^4b_0^3k^
4\delta_1(k)b_0k^2\delta_1(k)b_0                        
+80\tau_1^2|\tau|^2\xi_1^2\xi_2^4b_0^3k^
4\delta_1(k)b_0k\delta_1(k)b_0k\\                        
+40\tau_1^2|\tau|^2\xi_1^2\xi_2^4b_0^2k^
3\delta_1(k)b_0^2k^3\delta_1(k)b_0                        
+40\tau_1^2|\tau|^2\xi_1^2\xi_2^4b_0^2k^
3\delta_1(k)b_0^2k^2\delta_1(k)b_0k\\                        
+40\tau_1^2|\tau|^2\xi_1^2\xi_2^4b_0^2k^
2\delta_1(k)b_0^2k^4\delta_1(k)b_0                       
 +40\tau_1^2|\tau|^2\xi_1^2\xi_2^4b_0^2k^
2\delta_1(k)b_0^2k^3\delta_1(k)b_0k\\                        
+8|\tau|^4\xi_1^4\xi_2^2b_0^3k^
5\delta_2(k)b_0k\delta_2(k)b_0                        
+8|\tau|^4\xi_1^4\xi_2^2b_0^3k^
5\delta_2(k)b_0\delta_2(k)b_0k\\                        
+8|\tau|^4\xi_1^4\xi_2^2b_0^3k^4\delta_2(k)b_0k^
2\delta_2(k)b_0                       
+8|\tau|^4\xi_1^4\xi_2^2b_0^3k^4\delta_2(k)b_0k\delta_2(k)b_0k\\                        
+4|\tau|^4\xi_1^4\xi_2^2b_0^2k^3\delta_2(k)b_0^2k^3\delta_2(k)b_0                        
+4|\tau|^4\xi_1^4\xi_2^2b_0^2k^3\delta_2(k)b_0^2k^2\delta_2(k)b_0k\\                        
+4|\tau|^4\xi_1^4\xi_2^2b_0^2k^2\delta_2(k)b_0^2k^4\delta_2(k)b_0                        
+4|\tau|^4\xi_1^4\xi_2^2b_0^2k^2\delta_2(k)b_0^2k^3\delta_2(k)b_0k\\                        
+16|\tau|^4\xi_1^3\xi_2^3b_0^3k^5\delta_1(k)b_0k\delta_2(k)b_0                        
+16|\tau|^4\xi_1^3\xi_2^3b_0^3k^5\delta_1(k)b_0\delta_2(k)b_0k\\                        
+16|\tau|^4\xi_1^3\xi_2^3b_0^3k^5\delta_2(k)b_0k\delta_1(k)b_0                        
+16|\tau|^4\xi_1^3\xi_2^3b_0^3k^5\delta_2(k)b_0\delta_1(k)b_0k\\                        
+16|\tau|^4\xi_1^3\xi_2^3b_0^3k^4\delta_1(k)b_0k^2\delta_2(k)b_0                        
+16|\tau|^4\xi_1^3\xi_2^3b_0^3k^4\delta_1(k)b_0k\delta_2(k)b_0k\\                        
+16|\tau|^4\xi_1^3\xi_2^3b_0^3k^4\delta_2(k)b_0k^2\delta_1(k)b_0                        
+16|\tau|^4\xi_1^3\xi_2^3b_0^3k^4\delta_2(k)b_0k\delta_1(k)b_0k\\                        
+8|\tau|^4\xi_1^3\xi_2^3b_0^2k^3\delta_1(k)b_0^2k^3\delta_2(k)b_0                        
+8|\tau|^4\xi_1^3\xi_2^3b_0^2k^3\delta_1(k)b_0^2k^2\delta_2(k)b_0k\\                       
 +8|\tau|^4\xi_1^3\xi_2^3b_0^2k^3\delta_2(k)b_0^2k^3\delta_1(k)b_0                       
 +8|\tau|^4\xi_1^3\xi_2^3b_0^2k^3\delta_2(k)b_0^2k^2\delta_1(k)b_0k\\                        
 +8|\tau|^4\xi_1^3\xi_2^3b_0^2k^2\delta_1(k)b_0^2k^4\delta_2(k)b_0                       
  +8|\tau|^4\xi_1^3\xi_2^3b_0^2k^2\delta_1(k)b_0^2k^3\delta_2(k)b_0k\\                        
  +8|\tau|^4\xi_1^3\xi_2^3b_0^2k^2\delta_2(k)b_0^2k^4\delta_1(k)b_0                        
  +8|\tau|^4\xi_1^3\xi_2^3b_0^2k^2\delta_2(k)b_0^2k^3\delta_1(k)b_0k\\                        
+8|\tau|^4\xi_1^2\xi_2^4b_0^3k^5\delta_1(k)b_0k\delta_1(k)b_0                        
+8|\tau|^4\xi_1^2\xi_2^4b_0^3k^5\delta_1(k)b_0\delta_1(k)b_0k\\                        
+8|\tau|^4\xi_1^2\xi_2^4b_0^3k^4\delta_1(k)b_0k^2\delta_1(k)b_0                        
+8|\tau|^4\xi_1^2\xi_2^4b_0^3k^4\delta_1(k)b_0k\delta_1(k)b_0k\\                        
+4|\tau|^4\xi_1^2\xi_2^4b_0^2k^3\delta_1(k)b_0^2k^3\delta_1(k)b_0                        
+4|\tau|^4\xi_1^2\xi_2^4b_0^2k^3\delta_1(k)b_0^2k^2\delta_1(k)b_0k\\                        
+4|\tau|^4\xi_1^2\xi_2^4b_0^2k^2\delta_1(k)b_0^2k^4\delta_1(k)b_0                        
+4|\tau|^4\xi_1^2\xi_2^4b_0^2k^2\delta_1(k)b_0^2k^3\delta_1(k)b_0k\\                        
+96\tau_1^3|\tau|^2\xi_1^3\xi_2^3b_0^3k^5\delta_2(k)b_0k\delta_2(k)b_0                        
+96\tau_1^3|\tau|^2\xi_1^3\xi_2^3b_0^3k^5\delta_2(k)b_0\delta_2(k)b_0k\\                        
+96\tau_1^3|\tau|^2\xi_1^3\xi_2^3b_0^3k^4\delta_2(k)b_0k^2\delta_2(k)b_0                        
+96\tau_1^3|\tau|^2\xi_1^3\xi_2^3b_0^3k^4\delta_2(k)b_0k\delta_2(k)b_0k\\                        
+48\tau_1^3|\tau|^2\xi_1^3\xi_2^3b_0^2k^3\delta_2(k)b_0^2k^3\delta_2(k)b_0                        
+48\tau_1^3|\tau|^2\xi_1^3\xi_2^3b_0^2k^3\delta_2(k)b_0^2k^2\delta_2(k)b_0k\\                        
+48\tau_1^3|\tau|^2\xi_1^3\xi_2^3b_0^2k^2\delta_2(k)b_0^2k^4\delta_2(k)b_0                        
+48\tau_1^3|\tau|^2\xi_1^3\xi_2^3b_0^2k^2\delta_2(k)b_0^2k^3\delta_2(k)b_0k\\                        
+64\tau_1^3|\tau|^2\xi_1^2\xi_2^4b_0^3k^5\delta_1(k)b_0k\delta_2(k)b_0                       
+64\tau_1^3|\tau|^2\xi_1^2\xi_2^4b_0^3k^
5\delta_1(k)b_0\delta_2(k)b_0k\\                        
+64\tau_1^3|\tau|^2\xi_1^2\xi_2^4b_0^3k^
5\delta_2(k)b_0k\delta_1(k)b_0                        
+64\tau_1^3|\tau|^2\xi_1^2\xi_2^4b_0^3k^
5\delta_2(k)b_0\delta_1(k)b_0k\\                        
+64\tau_1^3|\tau|^2\xi_1^2\xi_2^4b_0^3k^
4\delta_1(k)b_0k^2\delta_2(k)b_0                        
+64\tau_1^3|\tau|^2\xi_1^2\xi_2^4b_0^3k^
4\delta_1(k)b_0k\delta_2(k)b_0k\\                        
+64\tau_1^3|\tau|^2\xi_1^2\xi_2^4b_0^3k^
4\delta_2(k)b_0k^2\delta_1(k)b_0                        
+64\tau_1^3|\tau|^2\xi_1^2\xi_2^4b_0^3k^
4\delta_2(k)b_0k\delta_1(k)b_0k\\                        
+32\tau_1^3|\tau|^2\xi_1^2\xi_2^4b_0^2k^
3\delta_1(k)b_0^2k^3\delta_2(k)b_0                        
+32\tau_1^3|\tau|^2\xi_1^2\xi_2^4b_0^2k^
3\delta_1(k)b_0^2k^2\delta_2(k)b_0k\\                        
+32\tau_1^3|\tau|^2\xi_1^2\xi_2^4b_0^2k^
3\delta_2(k)b_0^2k^3\delta_1(k)b_0                        
+32\tau_1^3|\tau|^2\xi_1^2\xi_2^4b_0^2k^
3\delta_2(k)b_0^2k^2\delta_1(k)b_0k\\                        
+32\tau_1^3|\tau|^2\xi_1^2\xi_2^4b_0^2k^
2\delta_1(k)b_0^2k^4\delta_2(k)b_0                        
+32\tau_1^3|\tau|^2\xi_1^2\xi_2^4b_0^2k^
2\delta_1(k)b_0^2k^3\delta_2(k)b_0k\\                        
+32\tau_1^3|\tau|^2\xi_1^2\xi_2^4b_0^2k^
2\delta_2(k)b_0^2k^4\delta_1(k)b_0                        
+32\tau_1^3|\tau|^2\xi_1^2\xi_2^4b_0^2k^
2\delta_2(k)b_0^2k^3\delta_1(k)b_0k\\                        
+32\tau_1^3|\tau|^2\xi_1\xi_2^5b_0^3k^
5\delta_1(k)b_0k\delta_1(k)b_0                        
+32\tau_1^3|\tau|^2\xi_1\xi_2^5b_0^3k^
5\delta_1(k)b_0\delta_1(k)b_0k\\                        
+32\tau_1^3|\tau|^2\xi_1\xi_2^5b_0^3k^
4\delta_1(k)b_0k^2\delta_1(k)b_0                        
+32\tau_1^3|\tau|^2\xi_1\xi_2^5b_0^3k^
4\delta_1(k)b_0k\delta_1(k)b_0k\\                        
+16\tau_1^3|\tau|^2\xi_1\xi_2^5b_0^2k^
3\delta_1(k)b_0^2k^3\delta_1(k)b_0                        
+16\tau_1^3|\tau|^2\xi_1\xi_2^5b_0^2k^
3\delta_1(k)b_0^2k^2\delta_1(k)b_0k\\                        
+16\tau_1^3|\tau|^2\xi_1\xi_2^5b_0^2k^
2\delta_1(k)b_0^2k^4\delta_1(k)b_0                        
+16\tau_1^3|\tau|^2\xi_1\xi_2^5b_0^2k^
2\delta_1(k)b_0^2k^3\delta_1(k)b_0k\\                        
+64\tau_1|\tau|^4\xi_1^3\xi_2^3b_0^3k^
5\delta_2(k)b_0k\delta_2(k)b_0                        
+64\tau_1|\tau|^4\xi_1^3\xi_2^3b_0^3k^
5\delta_2(k)b_0\delta_2(k)b_0k\\                       
 +64\tau_1|\tau|^4\xi_1^3\xi_2^3b_0^3k^
4\delta_2(k)b_0k^2\delta_2(k)b_0                        
+64\tau_1|\tau|^4\xi_1^3\xi_2^3b_0^3k^
4\delta_2(k)b_0k\delta_2(k)b_0k\\                        
+32\tau_1|\tau|^4\xi_1^3\xi_2^3b_0^2k^
3\delta_2(k)b_0^2k^3\delta_2(k)b_0                        
+32\tau_1|\tau|^4\xi_1^3\xi_2^3b_0^2k^
3\delta_2(k)b_0^2k^2\delta_2(k)b_0k\\                        
+32\tau_1|\tau|^4\xi_1^3\xi_2^3b_0^2k^
2\delta_2(k)b_0^2k^4\delta_2(k)b_0                        
+32\tau_1|\tau|^4\xi_1^3\xi_2^3b_0^2k^
2\delta_2(k)b_0^2k^3\delta_2(k)b_0k\\                        
+56\tau_1|\tau|^4\xi_1^2\xi_2^4b_0^3k^
5\delta_1(k)b_0k\delta_2(k)b_0                       
+56\tau_1|\tau|^4\xi_1^2\xi_2^4b_0^3k^
5\delta_1(k)b_0\delta_2(k)b_0k\\                        
+56\tau_1|\tau|^4\xi_1^2\xi_2^4b_0^3k^
5\delta_2(k)b_0k\delta_1(k)b_0                       
 +56\tau_1|\tau|^4\xi_1^2\xi_2^4b_0^3k^
5\delta_2(k)b_0\delta_1(k)b_0k\\                        
+56\tau_1|\tau|^4\xi_1^2\xi_2^4b_0^3k^
4\delta_1(k)b_0k^2\delta_2(k)b_0                        
+56\tau_1|\tau|^4\xi_1^2\xi_2^4b_0^3k^
4\delta_1(k)b_0k\delta_2(k)b_0k\\                        
+56\tau_1|\tau|^4\xi_1^2\xi_2^4b_0^3k^
4\delta_2(k)b_0k^2\delta_1(k)b_0                        
+56\tau_1|\tau|^4\xi_1^2\xi_2^4b_0^3k^
4\delta_2(k)b_0k\delta_1(k)b_0k\\                       
 +28\tau_1|\tau|^4\xi_1^2\xi_2^4b_0^2k^
3\delta_1(k)b_0^2k^3\delta_2(k)b_0                        
+28\tau_1|\tau|^4\xi_1^2\xi_2^4b_0^2k^
3\delta_1(k)b_0^2k^2\delta_2(k)b_0k\\                        
+28\tau_1|\tau|^4\xi_1^2\xi_2^4b_0^2k^
3\delta_2(k)b_0^2k^3\delta_1(k)b_0                        
+28\tau_1|\tau|^4\xi_1^2\xi_2^4b_0^2k^
3\delta_2(k)b_0^2k^2\delta_1(k)b_0k\\                        
+28\tau_1|\tau|^4\xi_1^2\xi_2^4b_0^2k^
2\delta_1(k)b_0^2k^4\delta_2(k)b_0                        
+28\tau_1|\tau|^4\xi_1^2\xi_2^4b_0^2k^
2\delta_1(k)b_0^2k^3\delta_2(k)b_0k\\                        
+28\tau_1|\tau|^4\xi_1^2\xi_2^4b_0^2k^
2\delta_2(k)b_0^2k^4\delta_1(k)b_0                        
+28\tau_1|\tau|^4\xi_1^2\xi_2^4b_0^2k^
2\delta_2(k)b_0^2k^3\delta_1(k)b_0k\\                        
+16\tau_1|\tau|^4\xi_1\xi_2^5b_0^3k^
5\delta_1(k)b_0k\delta_1(k)b_0                       
 +16\tau_1|\tau|^4\xi_1\xi_2^5b_0^3k^
5\delta_1(k)b_0\delta_1(k)b_0k\\                        
+16\tau_1|\tau|^4\xi_1\xi_2^5b_0^3k^
4\delta_1(k)b_0k^2\delta_1(k)b_0                       
 +16\tau_1|\tau|^4\xi_1\xi_2^5b_0^3k^
4\delta_1(k)b_0k\delta_1(k)b_0k\\                        
+8\tau_1|\tau|^4\xi_1\xi_2^5b_0^2k^
3\delta_1(k)b_0^2k^3\delta_1(k)b_0                        
+8\tau_1|\tau|^4\xi_1\xi_2^5b_0^2k^
3\delta_1(k)b_0^2k^2\delta_1(k)b_0k\\                       
+8\tau_1|\tau|^4\xi_1\xi_2^5b_0^2k^
2\delta_1(k)b_0^2k^4\delta_1(k)b_0                        
+8\tau_1|\tau|^4\xi_1\xi_2^5b_0^2k^
2\delta_1(k)b_0^2k^3\delta_1(k)b_0k\\                        
+104\tau_1^2|\tau|^4\xi_1^2\xi_2^4b_0^3k^
5\delta_2(k)b_0k\delta_2(k)b_0                       
+104\tau_1^2|\tau|^4\xi_1^2\xi_2^4b_0^3k^
5\delta_2(k)b_0\delta_2(k)b_0k\\                        
+104\tau_1^2|\tau|^4\xi_1^2\xi_2^4b_0^3k^
4\delta_2(k)b_0k^2\delta_2(k)b_0                        
+104\tau_1^2|\tau|^4\xi_1^2\xi_2^4b_0^3k^
4\delta_2(k)b_0k\delta_2(k)b_0k\\                        
+52\tau_1^2|\tau|^4\xi_1^2\xi_2^4b_0^2k^
3\delta_2(k)b_0^2k^3\delta_2(k)b_0                        
+52\tau_1^2|\tau|^4\xi_1^2\xi_2^4b_0^2k^
3\delta_2(k)b_0^2k^2\delta_2(k)b_0k\\                        
+52\tau_1^2|\tau|^4\xi_1^2\xi_2^4b_0^2k^
2\delta_2(k)b_0^2k^4\delta_2(k)b_0                       
+52\tau_1^2|\tau|^4\xi_1^2\xi_2^4b_0^2k^
2\delta_2(k)b_0^2k^3\delta_2(k)b_0k\\                        
+40\tau_1^2|\tau|^4\xi_1\xi_2^5b_0^3k^
5\delta_1(k)b_0k\delta_2(k)b_0                        
+40\tau_1^2|\tau|^4\xi_1\xi_2^5b_0^3k^
5\delta_1(k)b_0\delta_2(k)b_0k\\                        
+40\tau_1^2|\tau|^4\xi_1\xi_2^5b_0^3k^
5\delta_2(k)b_0k\delta_1(k)b_0                        
+40\tau_1^2|\tau|^4\xi_1\xi_2^5b_0^3k^
5\delta_2(k)b_0\delta_1(k)b_0k\\                       
 +40\tau_1^2|\tau|^4\xi_1\xi_2^5b_0^3k^
4\delta_1(k)b_0k^2\delta_2(k)b_0                        
+40\tau_1^2|\tau|^4\xi_1\xi_2^5b_0^3k^
4\delta_1(k)b_0k\delta_2(k)b_0k\\                        
+40\tau_1^2|\tau|^4\xi_1\xi_2^5b_0^3k^
4\delta_2(k)b_0k^2\delta_1(k)b_0                       
 +40\tau_1^2|\tau|^4\xi_1\xi_2^5b_0^3k^
4\delta_2(k)b_0k\delta_1(k)b_0k\\                        
+20\tau_1^2|\tau|^4\xi_1\xi_2^5b_0^2k^
3\delta_1(k)b_0^2k^3\delta_2(k)b_0                        
+20\tau_1^2|\tau|^4\xi_1\xi_2^5b_0^2k^
3\delta_1(k)b_0^2k^2\delta_2(k)b_0k\\                        
+20\tau_1^2|\tau|^4\xi_1\xi_2^5b_0^2k^
3\delta_2(k)b_0^2k^3\delta_1(k)b_0                        
+20\tau_1^2|\tau|^4\xi_1\xi_2^5b_0^2k^
3\delta_2(k)b_0^2k^2\delta_1(k)b_0k\\                        
+20\tau_1^2|\tau|^4\xi_1\xi_2^5b_0^2k^
2\delta_1(k)b_0^2k^4\delta_2(k)b_0                       
+20\tau_1^2|\tau|^4\xi_1\xi_2^5b_0^2k^
2\delta_1(k)b_0^2k^3\delta_2(k)b_0k\\                        
+20\tau_1^2|\tau|^4\xi_1\xi_2^5b_0^2k^
2\delta_2(k)b_0^2k^4\delta_1(k)b_0                        
+20\tau_1^2|\tau|^4\xi_1\xi_2^5b_0^2k^
2\delta_2(k)b_0^2k^3\delta_1(k)b_0k\\                        
+8\tau_1^2|\tau|^4\xi_2^6b_0^3k^
5\delta_1(k)b_0k\delta_1(k)b_0                       
+8\tau_1^2|\tau|^4\xi_2^6b_0^3k^
5\delta_1(k)b_0\delta_1(k)b_0k\\                        
+8\tau_1^2|\tau|^4\xi_2^6b_0^3k^4\delta_1(k)b_0k^
2\delta_1(k)b_0                        
+8\tau_1^2|\tau|^4\xi_2^6b_0^3k^
4\delta_1(k)b_0k\delta_1(k)b_0k\\                        
+4\tau_1^2|\tau|^4\xi_2^6b_0^2k^3\delta_1(k)b_0^
2k^3\delta_1(k)b_0                        
+4\tau_1^2|\tau|^4\xi_2^6b_0^2k^3\delta_1(k)b_0^2k^
2\delta_1(k)b_0k\\                        
+4\tau_1^2|\tau|^4\xi_2^6b_0^2k^2\delta_1(k)b_0^2k^
4\delta_1(k)b_0                        
+4\tau_1^2|\tau|^4\xi_2^6b_0^2k^2\delta_1(k)b_0^2k^
3\delta_1(k)b_0k\\                        
+16|\tau|^6\xi_1^2\xi_2^4b_0^3k^
5\delta_2(k)b_0k\delta_2(k)b_0                       
 +16|\tau|^6\xi_1^2\xi_2^4b_0^3k^5\delta_2(k)b_0\delta_2(k)b_0k\\                        
 +16|\tau|^6\xi_1^2\xi_2^4b_0^3k^4\delta_2(k)b_0k^
2\delta_2(k)b_0                        
+16|\tau|^6\xi_1^2\xi_2^4b_0^3k^
4\delta_2(k)b_0k\delta_2(k)b_0k\\                        
+8|\tau|^6\xi_1^2\xi_2^4b_0^2k^3\delta_2(k)b_0^
2k^3\delta_2(k)b_0                        
+8|\tau|^6\xi_1^2\xi_2^4b_0^2k^3\delta_2(k)b_0^2k^
2\delta_2(k)b_0k\\                        
+8|\tau|^6\xi_1^2\xi_2^4b_0^2k^2\delta_2(k)b_0^2k^
4\delta_2(k)b_0                        
+8|\tau|^6\xi_1^2\xi_2^4b_0^2k^2\delta_2(k)b_0^2k^
3\delta_2(k)b_0k\\                        
+8|\tau|^6\xi_1\xi_2^5b_0^3k^5\delta_1(k)b_0k\delta_2(k)b_0                        
+8|\tau|^6\xi_1\xi_2^5b_0^3k^5\delta_1(k)b_0\delta_2(k)b_0k\\                        
+8|\tau|^6\xi_1\xi_2^5b_0^3k^5\delta_2(k)b_0k\delta_1(k)b_0                        
+8|\tau|^6\xi_1\xi_2^5b_0^3k^5\delta_2(k)b_0\delta_1(k)b_0k\\                        
+8|\tau|^6\xi_1\xi_2^5b_0^3k^4\delta_1(k)b_0k^2\delta_2(k)b_0                        
+8|\tau|^6\xi_1\xi_2^5b_0^3k^4\delta_1(k)b_0k\delta_2(k)b_0k\\                        
+8|\tau|^6\xi_1\xi_2^5b_0^3k^4\delta_2(k)b_0k^2\delta_1(k)b_0                        
+8|\tau|^6\xi_1\xi_2^5b_0^3k^4\delta_2(k)b_0k\delta_1(k)b_0k\\                        
+4|\tau|^6\xi_1\xi_2^5b_0^2k^3\delta_1(k)b_0^2k^3\delta_2(k)b_0                        
+4|\tau|^6\xi_1\xi_2^5b_0^2k^3\delta_1(k)b_0^2k^2\delta_2(k)b_0k\\                        
+4|\tau|^6\xi_1\xi_2^5b_0^2k^3\delta_2(k)b_0^2k^3\delta_1(k)b_0                        
+4|\tau|^6\xi_1\xi_2^5b_0^2k^3\delta_2(k)b_0^2k^2\delta_1(k)b_0k\\                        
+4|\tau|^6\xi_1\xi_2^5b_0^2k^2\delta_1(k)b_0^2k^4\delta_2(k)b_0                        
+4|\tau|^6\xi_1\xi_2^5b_0^2k^2\delta_1(k)b_0^2k^3\delta_2(k)b_0k\\                        
+4|\tau|^6\xi_1\xi_2^5b_0^2k^2\delta_2(k)b_0^2k^4\delta_1(k)b_0                        
+4|\tau|^6\xi_1\xi_2^5b_0^2k^2\delta_2(k)b_0^2k^3\delta_1(k)b_0k\\                        
+48\tau_1|\tau|^6\xi_1\xi_2^5b_0^3k^5\delta_2(k)b_0k\delta_2(k)b_0                        
+48\tau_1|\tau|^6\xi_1\xi_2^5b_0^3k^5\delta_2(k)b_0\delta_2(k)b_0k\\                       
 +48\tau_1|\tau|^6\xi_1\xi_2^5b_0^3k^4\delta_2(k)b_0k^2\delta_2(k)b_0                        
 +48\tau_1|\tau|^6\xi_1\xi_2^5b_0^3k^4\delta_2(k)b_0k\delta_2(k)b_0k\\                        
 +24\tau_1|\tau|^6\xi_1\xi_2^5b_0^2k^3\delta_2(k)b_0^2k^3\delta_2(k)b_0                       
  +24\tau_1|\tau|^6\xi_1\xi_2^5b_0^2k^3\delta_2(k)b_0^2k^2\delta_2(k)b_0k\\                        
  +24\tau_1|\tau|^6\xi_1\xi_2^5b_0^2k^2\delta_2(k)b_0^2k^4\delta_2(k)b_0                        
  +24\tau_1|\tau|^6\xi_1\xi_2^5b_0^2k^2\delta_2(k)b_0^2k^3\delta_2(k)b_0k\\                        
  +8\tau_1|\tau|^6\xi_2^6b_0^3k^5\delta_1(k)b_0k\delta_2(k)b_0                        
  +8\tau_1|\tau|^6\xi_2^6b_0^3k^5\delta_1(k)b_0\delta_2(k)b_0k\\                       
   +8\tau_1|\tau|^6\xi_2^6b_0^3k^5\delta_2(k)b_0k\delta_1(k)b_0                        
   +8\tau_1|\tau|^6\xi_2^6b_0^3k^5\delta_2(k)b_0\delta_1(k)b_0k\\                        
   +8\tau_1|\tau|^6\xi_2^6b_0^3k^4\delta_1(k)b_0k^2\delta_2(k)b_0                        
   +8\tau_1|\tau|^6\xi_2^6b_0^3k^4\delta_1(k)b_0k\delta_2(k)b_0k\\                        
   +8\tau_1|\tau|^6\xi_2^6b_0^3k^4\delta_2(k)b_0k^2\delta_1(k)b_0                        
   +8\tau_1|\tau|^6\xi_2^6b_0^3k^4\delta_2(k)b_0k\delta_1(k)b_0k\\                        
   +4\tau_1|\tau|^6\xi_2^6b_0^2k^3\delta_1(k)b_0^2k^3\delta_2(k)b_0                        
   +4\tau_1|\tau|^6\xi_2^6b_0^2k^3\delta_1(k)b_0^2k^2\delta_2(k)b_0k\\                        
   +4\tau_1|\tau|^6\xi_2^6b_0^2k^3\delta_2(k)b_0^2k^3\delta_1(k)b_0                       
    +4\tau_1|\tau|^6\xi_2^6b_0^2k^3\delta_2(k)b_0^2k^2\delta_1(k)b_0k\\                        
+4\tau_1|\tau|^6\xi_2^6b_0^2k^2\delta_1(k)b_0^2k^4\delta_2(k)b_0                       
+4\tau_1|\tau|^6\xi_2^6b_0^2k^2\delta_1(k)b_0^2k^3\delta_2(k)b_0k\\                        
+4\tau_1|\tau|^6\xi_2^6b_0^2k^2\delta_2(k)b_0^2k^4\delta_1(k)b_0                        
+4\tau_1|\tau|^6\xi_2^6b_0^2k^2\delta_2(k)b_0^2k^3\delta_1(k)b_0k\\                        
+8|\tau|^8\xi_2^6b_0^3k^5\delta_2(k)b_0k\delta_2(k)b_0                        
+8|\tau|^8\xi_2^6b_0^3k^5\delta_2(k)b_0\delta_2(k)b_0k\\                        
+8|\tau|^8\xi_2^6b_0^3k^4\delta_2(k)b_0k^2\delta_2(k)b_0                        
+8|\tau|^8\xi_2^6b_0^3k^4\delta_2(k)b_0k\delta_2(k)b_0k\\                        
+4|\tau|^8\xi_2^6b_0^2k^3\delta_2(k)b_0^2k^3\delta_2(k)b_0                        
+4|\tau|^8\xi_2^6b_0^2k^3\delta_2(k)b_0^2k^2\delta_2(k)b_0k\\                        
+4|\tau|^8\xi_2^6b_0^2k^2\delta_2(k)b_0^2k^4\delta_2(k)b_0                        
+4|\tau|^8\xi_2^6b_0^2k^2\delta_2(k)b_0^2k^3\delta_2(k)b_0k.$\\

\end{document}